

\input amstex
\expandafter\ifx\csname mathdefs.tex\endcsname\relax
  \expandafter\gdef\csname mathdefs.tex\endcsname{}
\else \message{Hey!  Apparently you were trying to
  \string twice.   This does not make sense.} 
\errmessage{Please edit your file (probably \jobname.tex) and remove
any duplicate ``\string\input'' lines} \fi




\catcode`\X=12\catcode`\@=11

\def\n@wcount{\alloc@0\count\countdef\insc@unt}
\def\n@wwrite{\alloc@7\write\chardef\sixt@@n}
\def\n@wread{\alloc@6\read\chardef\sixt@@n}
\def\r@s@t{\relax}\def\v@idline{\par}\def\@mputate#1/{#1}
\def\l@c@l#1X{\firstpart.#1}\def\gl@b@l#1X{#1}\def\t@d@l#1X{{}}

\def\crossrefs#1{\ifx\all#1\let\tr@ce=\all\else\def\tr@ce{#1,}\fi
   \n@wwrite\cit@tionsout\openout\cit@tionsout=\jobname.cit 
   \write\cit@tionsout{\tr@ce}\expandafter\setfl@gs\tr@ce,}
\def\setfl@gs#1,{\def\@{#1}\ifx\@\empty\let\next=\relax
   \else\let\next=\setfl@gs\expandafter\xdef
   \csname#1tr@cetrue\endcsname{}\fi\next}
\def\m@ketag#1#2{\expandafter\n@wcount\csname#2tagno\endcsname
     \csname#2tagno\endcsname=0\let\tail=\all\xdef\all{\tail#2,}
   \ifx#1\l@c@l\let\tail=\r@s@t\xdef\r@s@t{\csname#2tagno\endcsname=0\tail}\fi
   \expandafter\gdef\csname#2cite\endcsname##1{\expandafter
     \ifx\csname#2tag##1\endcsname\relax?\else\csname#2tag##1\endcsname\fi
     \expandafter\ifx\csname#2tr@cetrue\endcsname\relax\else
     \write\cit@tionsout{#2tag ##1 cited on page \folio.}\fi}
   \expandafter\gdef\csname#2page\endcsname##1{\expandafter
     \ifx\csname#2page##1\endcsname\relax?\else\csname#2page##1\endcsname\fi
     \expandafter\ifx\csname#2tr@cetrue\endcsname\relax\else
     \write\cit@tionsout{#2tag ##1 cited on page \folio.}\fi}
   \expandafter\gdef\csname#2tag\endcsname##1{\expandafter
      \ifx\csname#2check##1\endcsname\relax
      \expandafter\xdef\csname#2check##1\endcsname{}%
      \else\immediate\write16{Warning: #2tag ##1 used more than once.}\fi
      \multit@g{#1}{#2}##1/X%
      \write\t@gsout{#2tag ##1 assigned number \csname#2tag##1\endcsname\space
      on page \number\count0.}%
   \csname#2tag##1\endcsname}}
\def\multit@g#1#2#3/#4X{\def\t@mp{#4}\ifx\t@mp\empty%
      \global\advance\csname#2tagno\endcsname by 1 
      \expandafter\xdef\csname#2tag#3\endcsname
      {#1\number\csname#2tagno\endcsnameX}%
   \else\expandafter\ifx\csname#2last#3\endcsname\relax
      \expandafter\n@wcount\csname#2last#3\endcsname
      \global\advance\csname#2tagno\endcsname by 1 
      \expandafter\xdef\csname#2tag#3\endcsname
      {#1\number\csname#2tagno\endcsnameX}
      \write\t@gsout{#2tag #3 assigned number \csname#2tag#3\endcsname\space
      on page \number\count0.}\fi
   \global\advance\csname#2last#3\endcsname by 1
   \def\t@mp{\expandafter\xdef\csname#2tag#3/}%
   \expandafter\t@mp\@mputate#4\endcsname
   {\csname#2tag#3\endcsname\lastpart{\csname#2last#3\endcsname}}\fi}
\def\t@gs#1{\def\all{}\m@ketag#1e\m@ketag#1s\m@ketag\t@d@l p
   \m@ketag\gl@b@l r \n@wread\t@gsin
   \openin\t@gsin=\jobname.tgs \re@der \closein\t@gsin
   \n@wwrite\t@gsout\openout\t@gsout=\jobname.tgs }
\outer\def\localtags{\t@gs\l@c@l}
\outer\def\globaltags{\t@gs\gl@b@l}
\outer\def\newlocaltag#1{\m@ketag\l@c@l{#1}}
\outer\def\newglobaltag#1{\m@ketag\gl@b@l{#1}}

\newif\ifpr@ 
\def\m@kecs #1tag #2 assigned number #3 on page #4.%
   {\expandafter\gdef\csname#1tag#2\endcsname{#3}
   \expandafter\gdef\csname#1page#2\endcsname{#4}
   \ifpr@\expandafter\xdef\csname#1check#2\endcsname{}\fi}
\def\re@der{\ifeof\t@gsin\let\next=\relax\else
   \read\t@gsin to\t@gline\ifx\t@gline\v@idline\else
   \expandafter\m@kecs \t@gline\fi\let \next=\re@der\fi\next}
\def\pretags#1{\pr@true\pret@gs#1,,}
\def\pret@gs#1,{\def\@{#1}\ifx\@\empty\let\n@xtfile=\relax
   \else\let\n@xtfile=\pret@gs \openin\t@gsin=#1.tgs \message{#1} \re@der 
   \closein\t@gsin\fi \n@xtfile}

\newcount\sectno\sectno=0\newcount\subsectno\subsectno=0
\newif\ifultr@local \def\ultralocal{\ultr@localtrue}
\def\firstpart{\number\sectno}
\def\lastpart#1{\ifcase#1 \or a\or b\or c\or d\or e\or f\or g\or h\or 
   i\or k\or l\or m\or n\or o\or p\or q\or r\or s\or t\or u\or v\or w\or 
   x\or y\or z \fi}

\def\resetall{\global\advance\sectno by 1\subsectno=0
   \gdef\firstpart{\number\sectno}\r@s@t}
\def\resetsub{\global\advance\subsectno by 1
   \gdef\firstpart{\number\sectno.\number\subsectno}\r@s@t}
\def\newsection#1\par{\resetall\vskip0pt plus.3\vsize\penalty-250
   \vskip0pt plus-.3\vsize\bigskip\bigskip
   \message{#1}\leftline{\bf#1}\nobreak\bigskip}
\def\subsection#1\par{\ifultr@local\resetsub\fi
   \vskip0pt plus.2\vsize\penalty-250\vskip0pt plus-.2\vsize
   \bigskip\smallskip\message{#1}\leftline{\bf#1}\nobreak\medskip}

\def\t@gsoff#1,{\def\@{#1}\ifx\@\empty\let\next=\relax\else\let\next=\t@gsoff
   \def\@@{p}\ifx\@\@@\else
   \expandafter\gdef\csname#1cite\endcsname##1{\zeigen{##1}}
   \expandafter\gdef\csname#1page\endcsname##1{?}
   \expandafter\gdef\csname#1tag\endcsname##1{\zeigen{##1}}\fi\fi\next}
\def\verbatimtags{\ifx\all\relax\else\expandafter\t@gsoff\all,\fi}
\def\zeigen#1{\hbox{$\langle$}#1\hbox{$\rangle$}}

\def\(#1){\edef\dot@g{\ifmmode\ifinner(\hbox{\noexpand\etag{#1}})
   \else\noexpand\eqno(\hbox{\noexpand\etag{#1}})\fi
   \else(\noexpand\ecite{#1})\fi}\dot@g}

\newif\ifbr@ck
\def\eat#1{}
\def\[#1]{\br@cktrue[\br@cket#1'X]}
\def\br@cket#1'#2X{\def\temp{#2}\ifx\temp\empty\let\next\eat
   \else\let\next\br@cket\fi
   \ifbr@ck\br@ckfalse\br@ck@t#1,X\else\br@cktrue#1\fi\next#2X}
\def\br@ck@t#1,#2X{\def\temp{#2}\ifx\temp\empty\let\neext\eat
   \else\let\neext\br@ck@t\def\temp{,}\fi
   \def\teemp{#1}\ifx\teemp\empty\else\rcite{#1}\fi\temp\neext#2X}
\def\resetbr@cket{\gdef\[##1]{[\rtag{##1}]}}
\def\references{\resetbr@cket\newsection References\par}

\newtoks\symb@ls\newtoks\s@mb@ls\newtoks\p@gelist\n@wcount\ftn@mber
    \ftn@mber=1\newif\ifftn@mbers\ftn@mbersfalse\newif\ifbyp@ge\byp@gefalse
\def\defm@rk{\ifftn@mbers\n@mberm@rk\else\symb@lm@rk\fi}
\def\n@mberm@rk{\xdef\m@rk{{\the\ftn@mber}}%
    \global\advance\ftn@mber by 1 }
\def\rot@te#1{\let\temp=#1\global#1=\expandafter\r@t@te\the\temp,X}
\def\r@t@te#1,#2X{{#2#1}\xdef\m@rk{{#1}}}
\def\b@@st#1{{$^{#1}$}}\def\str@p#1{#1}
\def\symb@lm@rk{\ifbyp@ge\rot@te\p@gelist\ifnum\expandafter\str@p\m@rk=1 
    \s@mb@ls=\symb@ls\fi\write\f@nsout{\number\count0}\fi \rot@te\s@mb@ls}
\def\byp@ge{\byp@getrue\n@wwrite\f@nsin\openin\f@nsin=\jobname.fns 
    \n@wcount\currentp@ge\currentp@ge=0\p@gelist={0}
    \re@dfns\closein\f@nsin\rot@te\p@gelist
    \n@wread\f@nsout\openout\f@nsout=\jobname.fns }
\def\m@kelist#1X#2{{#1,#2}}
\def\re@dfns{\ifeof\f@nsin\let\next=\relax\else\read\f@nsin to \f@nline
    \ifx\f@nline\v@idline\else\let\t@mplist=\p@gelist
    \ifnum\currentp@ge=\f@nline
    \global\p@gelist=\expandafter\m@kelist\the\t@mplistX0
    \else\currentp@ge=\f@nline
    \global\p@gelist=\expandafter\m@kelist\the\t@mplistX1\fi\fi
    \let\next=\re@dfns\fi\next}
\def\symbols#1{\symb@ls={#1}\s@mb@ls=\symb@ls} 
\def\bigsymbol{\textstyle}
\symbols{\bigsymbol\ast,\dagger,\ddagger,\sharp,\flat,\natural,\star}
\def\ftnumbers{\ftn@mberstrue} \def\ftsymbols{\ftn@mbersfalse}
\def\paginal{\byp@ge} \def\resetftnumbers{\ftn@mber=1}
\def\ftnote#1{\defm@rk\expandafter\expandafter\expandafter\footnote
    \expandafter\b@@st\m@rk{#1}}

\long\def\jump#1\endjump{}
\def\ssum{\mathop{\lower .1em\hbox{$\textstyle\Sigma$}}\nolimits}

\def\qed{\nobreak\kern 1em \vrule height .5em width .5em depth 0em}
\def\newneq{\hbox{\rlap{\hbox to 1\wd9{\hss$=$\hss}}\raise .1em 
   \hbox to 1\wd9{\hss$\scriptscriptstyle/$\hss}}}
\def\subsetne{\setbox9 = \hbox{$\subset$}\mathrel{\hbox{\rlap
   {\lower .4em \newneq}\raise .13em \hbox{$\subset$}}}}
\def\supsetne{\setbox9 = \hbox{$\subset$}\mathrel{\hbox{\rlap
   {\lower .4em \newneq}\raise .13em \hbox{$\supset$}}}}

\def\vbar{\mathchoice{\vrule height6.3ptdepth-.5ptwidth.8pt\kern-.8pt}
   {\vrule height6.3ptdepth-.5ptwidth.8pt\kern-.8pt}
   {\vrule height4.1ptdepth-.35ptwidth.6pt\kern-.6pt}
   {\vrule height3.1ptdepth-.25ptwidth.5pt\kern-.5pt}}
\def\f@dge{\mathchoice{}{}{\mkern.5mu}{\mkern.8mu}}
\def\b@c#1#2{{\rm \mkern#2mu\vbar\mkern-#2mu#1}}
\def\b@b#1{{\rm I\mkern-3.5mu #1}}
\def\b@a#1#2{{\rm #1\mkern-#2mu\f@dge #1}}
\def\bb#1{{\count4=`#1 \advance\count4by-64 \ifcase\count4\or\b@a A{11.5}\or
   \b@b B\or\b@c C{5}\or\b@b D\or\b@b E\or\b@b F \or\b@c G{5}\or\b@b H\or
   \b@b I\or\b@c J{3}\or\b@b K\or\b@b L \or\b@b M\or\b@b N\or\b@c O{5} \or
   \b@b P\or\b@c Q{5}\or\b@b R\or\b@a S{8}\or\b@a T{10.5}\or\b@c U{5}\or
   \b@a V{12}\or\b@a W{16.5}\or\b@a X{11}\or\b@a Y{11.7}\or\b@a Z{7.5}\fi}}

\catcode`\X=11 \catcode`\@=12

\expandafter\ifx\csname citeadd.tex\endcsname\relax
\expandafter\gdef\csname citeadd.tex\endcsname{}
\else \message{Hey!  Apparently you were trying to
\string twice.   This does not make sense.} 
\errmessage{Please edit your file (probably \jobname.tex) and remove
any duplicate ``\string\input'' lines} \fi

\ifx\shlhetal\undefinedcontrolsequence\let\shlhetal\relax\fi
\sectno=-1   
\documentstyle {amsppt}
\localtags
\NoBlackBoxes
\define\mr{\medskip\roster}
\define\sn{\smallskip\noindent}
\define\mn{\medskip\noindent}
\define\bn{\bigskip\noindent}
\define\ub{\underbar}
\define\wilog{\text{without loss of generality}}
\define\ermn{\endroster\medskip\noindent}

\define\dbcu{\dsize\bigcup}
\define \nl{\newline}
\newbox\noforkbox \newdimen\forklinewidth
\forklinewidth=0.3pt   
\setbox0\hbox{$\textstyle\bigcup$}
\setbox1\hbox to \wd0{\hfil\vrule width \forklinewidth depth \dp0
                        height \ht0 \hfil}
\wd1=0 cm
\setbox\noforkbox\hbox{\box1\box0\relax}
\def\unionstick{\mathop{\copy\noforkbox}\limits}
\def\nonfork#1#2_#3{#1\unionstick_{\textstyle #3}#2}
\def\nonforkin#1#2_#3^#4{#1\unionstick_{\textstyle #3}^{\textstyle #4}#2}     
%
\setbox0\hbox{$\textstyle\bigcup$}
\setbox1\hbox to \wd0{\hfil{\sl /\/}\hfil}
\setbox2\hbox to \wd0{\hfil\vrule height \ht0 depth \dp0 width
                                \forklinewidth\hfil}
\wd1=0cm
\wd2=0cm
\newbox\doesforkbox
\setbox\doesforkbox\hbox{\box1\box0\relax}
\def\nunionstick{\mathop{\copy\doesforkbox}\limits}

\def\fork#1#2_#3{#1\nunionstick_{\textstyle #3}#2}
\def\forkin#1#2_#3^#4{#1\nunionstick_{\textstyle #3}^{\textstyle #4}#2}     
\topmatter
\title {Categoricity of an abstract \\
elementary class in \\
two successive cardinals \\
Sh576} \endtitle
\rightheadtext{Categoricity in two successive cardinals}
\author {Saharon Shelah \thanks{\null\newline
Partially supported by the United States-Israel Binational Science
Foundation. \null\newline
I thank Alice Leonhardt for the excellent typing. \newline
Done 10/94, (set theoretic, pseudo descriptive set theory)
22-29/11/94 (essentially \S1-\S2); more 9-14/12/95 \null\newline
(fulfilling \S1,\S3 eliminating use of the
maximal model in $\lambda^{+2}$, (new \S4) of categoricity \null\newline
in $\lambda^{+2}$ (see \S7,\S8), more 18-2 of 4/95 (redo \S3,\S10). 
\null\newline
\S10 written 1/2/95 \null\newline
First typed at Rutgers - 95/Feb/10 \null\newline
Latest Revision 98/May/15} \endthanks} \endauthor
\affil {Institute of Mathematics \\
The Hebrew University \\
Jerusalem, Israel
\medskip
Rutgers University \\
Department of Mathematics \\
New Brunswick, NJ  USA} \endaffil
\abstract {We investigate categoricity of abstract elementary classes
without any remnants of compactness (like non-definability of well ordering,
existence of E.M. models, or existence of large cardinals).  We prove 
(assuming a weak version of GCH around $\lambda$) that
if ${\frak K}$ is categorical in $\lambda,\lambda^+$,
$LS({\frak K}) \le \lambda$ and 
$1 \le I(\lambda^{++},{\frak K}) < 2^{\lambda^{++}}$ 
\underbar{then} ${\frak K}$ has a model in $\lambda^{+++}$.} \endabstract
\endtopmatter
\document  

\expandafter\ifx\csname alice2jlem.tex\endcsname\relax
  \expandafter\gdef\csname alice2jlem.tex\endcsname{}
\else \message{Hey!  Apparently you were trying to
\string  twice.   This does not make sense.}
\errmessage{Please edit your file (probably \jobname.tex) and remove
any duplicate ``\string\input'' lines} \fi

\expandafter\ifx\csname bib4plain.tex\endcsname\relax
  \expandafter\gdef\csname bib4plain.tex\endcsname{}
\else \message{Hey!  Apparently you were trying to \string twice.   This does not make sense.}
\errmessage{Please edit your file (probably \jobname.tex) and remove
any duplicate ``\string\input'' lines} \fi

\def\renewcommand{\newcommand}	       
\edef\cite{\the\catcode`@}%
\catcode`@ = 11
\let\@oldatcatcode = \cite
\chardef\@letter = 11
\chardef\@other = 12
%
%
%
%
\def\@innerdef#1#2{\edef#1{\expandafter\noexpand\csname #2\endcsname}}%
%
%
\@innerdef\@innernewcount{newcount}%
\@innerdef\@innernewdimen{newdimen}%
\@innerdef\@innernewif{newif}%
\@innerdef\@innernewwrite{newwrite}%
%
%
%
\def\@gobble#1{}%
%
%
%
\ifx\inputlineno\@undefined
   \let\@linenumber = \empty 
\else
   \def\@linenumber{\the\inputlineno:\space}%
\fi
%
%
%
\def\@futurenonspacelet#1{\def\cs{#1}%
   \afterassignment\@stepone\let\@nexttoken=
}%
\begingroup 
\def\\{\global\let\@stoken= }%
\\ 
\endgroup
\def\@stepone{\expandafter\futurelet\cs\@steptwo}%
\def\@steptwo{\expandafter\ifx\cs\@stoken\let\@@next=\@stepthree
   \else\let\@@next=\@nexttoken\fi \@@next}%
\def\@stepthree{\afterassignment\@stepone\let\@@next= }%
%
%
%
\def\@getoptionalarg#1{%
   \let\@optionaltemp = #1%
   \let\@optionalnext = \relax
   \@futurenonspacelet\@optionalnext\@bracketcheck
}%
%
%
\def\@bracketcheck{%
   \ifx [\@optionalnext
      \expandafter\@@getoptionalarg
   \else
      \let\@optionalarg = \empty
      \expandafter\@optionaltemp
   \fi
}%
\def\@@getoptionalarg[#1]{%
   \def\@optionalarg{#1}%
   \@optionaltemp
}%
%
%
%
\def\@nnil{\@nil}%
\def\@fornoop#1\@@#2#3{}%
\def\@for#1:=#2\do#3{%
   \edef\@fortmp{#2}%
   \ifx\@fortmp\empty \else
      \expandafter\@forloop#2,\@nil,\@nil\@@#1{#3}%
   \fi
}%
\def\@forloop#1,#2,#3\@@#4#5{\def#4{#1}\ifx #4\@nnil \else
       #5\def#4{#2}\ifx #4\@nnil \else#5\@iforloop #3\@@#4{#5}\fi\fi
}%
\def\@iforloop#1,#2\@@#3#4{\def#3{#1}\ifx #3\@nnil
       \let\@nextwhile=\@fornoop \else
      #4\relax\let\@nextwhile=\@iforloop\fi\@nextwhile#2\@@#3{#4}%
}%
%
%
%
\@innernewif\if@fileexists
\def\@testfileexistence{\@getoptionalarg\@finishtestfileexistence}%
\def\@finishtestfileexistence#1{%
   \begingroup
      \def\extension{#1}%
      \immediate\openin0 =
         \ifx\@optionalarg\empty\jobname\else\@optionalarg\fi
         \ifx\extension\empty \else .#1\fi
         \space
      \ifeof 0
         \global\@fileexistsfalse
      \else
         \global\@fileexiststrue
      \fi
      \immediate\closein0
   \endgroup
}%
%
%
%
%
\def\bibliographystyle#1{%
   \@readauxfile
   \@writeaux{\string\bibstyle{#1}}%
}%
\let\bibstyle = \@gobble
%
%
\let\bblfilebasename = \jobname
\def\bibliography#1{%
   \@readauxfile
   \@writeaux{\string\bibdata{#1}}%
   \@testfileexistence[\bblfilebasename]{bbl}%
   \if@fileexists
      \nobreak
      \@readbblfile
   \fi
}%
\let\bibdata = \@gobble
%
%
\def\nocite#1{%
   \@readauxfile
   \@writeaux{\string\citation{#1}}%
}%
\@innernewif\if@notfirstcitation
%
%
\def\cite{\@getoptionalarg\@cite}%
%
%
\def\@cite#1{%
   \let\@citenotetext = \@optionalarg
   \printcitestart
   \nocite{#1}%
   \@notfirstcitationfalse
   \@for \@citation :=#1\do
   {%
      \expandafter\@onecitation\@citation\@@
   }%
   \ifx\empty\@citenotetext\else
      \printcitenote{\@citenotetext}%
   \fi
   \printcitefinish
}%
\def\@onecitation#1\@@{%
   \if@notfirstcitation
      \printbetweencitations
   \fi
   \expandafter \ifx \csname\@citelabel{#1}\endcsname \relax
      \if@citewarning
         \message{\@linenumber Undefined citation `#1'.}%
      \fi
      \expandafter\gdef\csname\@citelabel{#1}\endcsname{%
\strut
\vadjust{\vskip-\dp\strutbox
\vbox to 0pt{\vss\parindent0cm \leftskip=\hsize 
\advance\leftskip3mm
\advance\hsize 4cm\strut\openup-4pt 
\rightskip 0cm plus 1cm minus 0.5cm ?  #1 ?\strut}}
         {\tt
            \escapechar = -1
            \nobreak\hskip0pt
            \expandafter\string\csname#1\endcsname
            \nobreak\hskip0pt
         }%
      }%
   \fi
   \csname\@citelabel{#1}\endcsname
   \@notfirstcitationtrue
}%
%
%
\def\@citelabel#1{b@#1}%
%
%
\def\@citedef#1#2{\expandafter\gdef\csname\@citelabel{#1}\endcsname{#2}}%
%
%
%
\def\@readbblfile{%
   \ifx\@itemnum\@undefined
      \@innernewcount\@itemnum
   \fi
   \begingroup
      \def\begin##1##2{%
         \setbox0 = \hbox{\biblabelcontents{##2}}%
         \biblabelwidth = \wd0
      }%
      \def\end##1{}
      %
      %
      \@itemnum = 0
      \def\bibitem{\@getoptionalarg\@bibitem}%
      \def\@bibitem{%
         \ifx\@optionalarg\empty
            \expandafter\@numberedbibitem
         \else
            \expandafter\@alphabibitem
         \fi
      }%
      \def\@alphabibitem##1{%
         \expandafter \xdef\csname\@citelabel{##1}\endcsname {\@optionalarg}%
         \ifx\biblabelprecontents\@undefined
            \let\biblabelprecontents = \relax
         \fi
         \ifx\biblabelpostcontents\@undefined
            \let\biblabelpostcontents = \hss
         \fi
         \@finishbibitem{##1}%
      }%
      \def\@numberedbibitem##1{%
         \advance\@itemnum by 1
         \expandafter \xdef\csname\@citelabel{##1}\endcsname{\number\@itemnum}%
         \ifx\biblabelprecontents\@undefined
            \let\biblabelprecontents = \hss
         \fi
         \ifx\biblabelpostcontents\@undefined
            \let\biblabelpostcontents = \relax
         \fi
         \@finishbibitem{##1}%
      }%
      \def\@finishbibitem##1{%
         \biblabelprint{\csname\@citelabel{##1}\endcsname}%
         \@writeaux{\string\@citedef{##1}{\csname\@citelabel{##1}\endcsname}}%
         \ignorespaces
      }%
      %
      %
      \let\em = \bblem
      \let\newblock = \bblnewblock
      \let\sc = \bblsc
      \frenchspacing
      \clubpenalty = 4000 \widowpenalty = 4000
      \tolerance = 10000 \hfuzz = .5pt
      \everypar = {\hangindent = \biblabelwidth
                      \advance\hangindent by \biblabelextraspace}%
      \bblrm
      \parskip = 1.5ex plus .5ex minus .5ex
      \biblabelextraspace = .5em
      \bblhook
      \input \bblfilebasename.bbl
   \endgroup
}%
%
%
\@innernewdimen\biblabelwidth
\@innernewdimen\biblabelextraspace
%
%
%
\def\biblabelprint#1{%
   \noindent
   \hbox to \biblabelwidth{%
      \biblabelprecontents
      \biblabelcontents{#1}%
      \biblabelpostcontents
   }%
   \kern\biblabelextraspace
}%
%
%
%
\def\biblabelcontents#1{{\bblrm [#1]}}%
%
%
\def\bblrm{\rm}%
%
%
\def\bblem{\it}%
%
%
\def\bblsc{\ifx\@scfont\@undefined
              \font\@scfont = cmcsc10
           \fi
           \@scfont
}%
%
%
\def\bblnewblock{\hskip .11em plus .33em minus .07em }%
%
%
\let\bblhook = \empty
%
%
%
\def\printcitestart{[}
\def\printcitefinish{]}
\def\printbetweencitations{, }
\def\printcitenote#1{, #1}
%
%
%
\let\citation = \@gobble
%
%
%
\@innernewcount\@numparams
%
%
\def\newcommand#1{%
   \def\@commandname{#1}%
   \@getoptionalarg\@continuenewcommand
}%
%
%
\def\@continuenewcommand{%
   \@numparams = \ifx\@optionalarg\empty 0\else\@optionalarg \fi \relax
   \@newcommand
}%
%
%
\def\@newcommand#1{%
   \def\@startdef{\expandafter\edef\@commandname}%
   \ifnum\@numparams=0
      \let\@paramdef = \empty
   \else
      \ifnum\@numparams>9
         \errmessage{\the\@numparams\space is too many parameters}%
      \else
         \ifnum\@numparams<0
            \errmessage{\the\@numparams\space is too few parameters}%
         \else
            \edef\@paramdef{%
               \ifcase\@numparams
                  \empty  No arguments.
               \or ####1%
               \or ####1####2%
               \or ####1####2####3%
               \or ####1####2####3####4%
               \or ####1####2####3####4####5%
               \or ####1####2####3####4####5####6%
               \or ####1####2####3####4####5####6####7%
               \or ####1####2####3####4####5####6####7####8%
               \or ####1####2####3####4####5####6####7####8####9%
               \fi
            }%
         \fi
      \fi
   \fi
   \expandafter\@startdef\@paramdef{#1}%
}%
%
%
%
%
\def\@readauxfile{%
   \if@auxfiledone \else 
      \global\@auxfiledonetrue
      \@testfileexistence{aux}%
      \if@fileexists
         \begingroup
            \endlinechar = -1
            \catcode`@ = 11
            \input \jobname.aux
         \endgroup
      \else
         \message{\@undefinedmessage}%
         \global\@citewarningfalse
      \fi
      \immediate\openout\@auxfile = \jobname.aux
   \fi
}%
%
%
\newif\if@auxfiledone
\ifx\noauxfile\@undefined \else \@auxfiledonetrue\fi
%
%
%
%
\@innernewwrite\@auxfile
\def\@writeaux#1{\ifx\noauxfile\@undefined \write\@auxfile{#1}\fi}%
%
%
%
\ifx\@undefinedmessage\@undefined
   \def\@undefinedmessage{No .aux file; I won't give you warnings about
                          undefined citations.}%
\fi
%
%
\@innernewif\if@citewarning
\ifx\noauxfile\@undefined \@citewarningtrue\fi
%
%
%
\catcode`@ = \@oldatcatcode


\def\widestnumber#1#2{}

\def\rm{\fam0 \tenrm}

\def\fakesubhead#1\endsubhead{\bigskip\noindent{\bf#1}\par}


%
%
%

%

\font\textrsfs=rsfs10
\font\scriptrsfs=rsfs7
\font\scriptscriptrsfs=rsfs5

\newfam\rsfsfam
\textfont\rsfsfam=\textrsfs
\scriptfont\rsfsfam=\scriptrsfs
\scriptscriptfont\rsfsfam=\scriptscriptrsfs

\edef\oldcatcodeofat{\the\catcode`\@}
\catcode`\@11

\def\Cal@@#1{\noaccents@ \fam \rsfsfam #1}

\catcode`\@\oldcatcodeofat

 \catcode`\@\active  
\newpage

\head {Annotated Content} \endhead \bigskip

\noindent
\S0 $\quad$ \underbar{Introduction}
\roster
\item "{${}$}"  [We give three versions of the main theorem in 
\scite{0.B} - \scite{0.D}.  In \scite{0.1} - \scite{0.28} we review the 
relevant knowledge of abstract elementary classes to help make this paper 
self-contained.  This includes the representation by PC-classes defined by 
omission of quantifier free types (\scite{0.8}, \scite{0.9}); types and
stability (based on $\le_{\frak K}$); and the equivalence of saturation to 
model homogeneity (\scite{0.19}).]
\endroster
\medskip

\noindent
\S1 $\quad$ \underbar{Weak Diamond}
\roster
\item "{${}$}"  [We mainly present necessary material on the weak diamond,
a combinatorial principle whose main variant holds for $\lambda^+$
if $2^\lambda < 2^{\lambda^+}$. \newline
We state cases provable from ZFC together with suitable cardinal arithmetic 
assumptions (\scite{1.2}).  We present applications of weak diamond to the
number of models of fixed cardinality (notably in \scite{1.6}).  We deal also
with the definable weak diamond, and introduce a smaller ideal of ``small" 
sets: UDmId$^{\Cal F}(\lambda)$ (\scite{1.8} - \scite{1.9}).]
\endroster
\medskip

\noindent
\S2 $\quad$ \underbar{First Attempts}
\roster
\item "{${}$}"  [We define the class $K^3_\lambda$ of triples $(M,N,a)$
representing types in ${\Cal S}(M)$ for $M \in K_\lambda$, and start to 
investigate it, dealing with the weak extension property, the extension 
property, minimality, reduced types (except for minimality, in the first 
order case, these hold trivially).  Our aims are to have the extension 
property or at least the weak extension property for all triples in 
$K^3_\lambda$, and the density of minimal triples.  The first property makes 
the model theory more like the first order case, and the
second is connected with categoricity.  We start by proving the weak extension
property under reasonable assumptions.  We prove the density of minimal 
triples under the strong assumption $K_{\lambda^{+3}} = \emptyset$ and an 
extra cardinal arithmetic assumption $(2^{\lambda^+} > \lambda^{++})$.  
In the end, under the additional assumption $K_{\lambda^{+3}} = 0$ we 
prove that all triples have the extension property and that we have 
disjoint amalgamation in $K_\lambda$.
Now the assumption $K_{\lambda^{+3}} = 0$ does no harm if we just want to
prove Theorem \scite{0.B}.  The reader willing to accept these assumptions 
may skip some proofs later.  The proof of the extension property makes 
essential use of categoricity in $\lambda^+$.]
\endroster
\medskip

\noindent
\S3 $\quad$ \underbar{Non-structure}
\roster
\item "{${}$}"  [We try to present clearly and in some generality the 
construction of many models in $\lambda^{++}$ based on knowledge of models 
of size $\lambda^+$, using weak diamond on $\lambda^+$ and on $\lambda^{++}$.
This is done by forming a tree
$\langle \bar M^\eta:\eta \in {}^{\lambda^{++}>}2 \rangle$ with $\bar M^\eta$
an $\le_{\frak K}$-increasing continuous sequence of members of $K_\lambda$
with limit $\dsize \bigcup_{i < \lambda^+} M^\eta_i$ increasing with $\eta$
(and an additional restriction).  Actually $\lambda^+$ can be replaced by a
regular uncountable $\lambda'$ (so $\|M^\eta_i\| = \lambda$ is replaced by
$\|M^\eta_i\| < \lambda'$)].
\endroster
\medskip

\noindent
\S4 $\quad$ \underbar{Minimal types}
\roster
\item "{${}$}"  [We prove that every member of $K^3_\lambda$ has the extension
property, by proving it for minimal triples.  We use: if $M_\ell \in 
K_\lambda$ and for every minimal \nl
$p \in {\Cal S}(M_0)$ the set 
${\Cal S}_{\ge p}(M_1)$ has cardinality $\le \lambda^+$, then the 
$M \in K_{\lambda^+}$ is saturated for minimal types and hence the number 
of minimal types in ${\Cal S}(M_1)$ is $\le \lambda^+$ (for 
$M_1 \in K_\lambda$), which is a step toward stability in $\lambda$.]
\endroster
\medskip

\noindent
\S5 $\quad$ \underbar{Inevitable types and stability in $\lambda$}
\roster
\item "{${}$}"  [We continue to ``climb the ladder", using the amount of
structure we already have (and sometimes categoricity) to get more.  We
start by assuming there are minimal types, and show that some minimal types 
are inevitable, construct $p_i \in {\Cal S}(N_i)$ minimal $(i \le \lambda^+)$ 
both strictly increasing continuous and with $p_0,p_\delta$ inevitable, and 
then as in the proof of the equivalence of saturativity and model 
homogeneity, we show $N_\delta$ 
is universal over $N_0$.  We can then deduce stability in $\lambda$, 
so the model in $\lambda^+$ is saturated.  Then we note that we 
have disjoint amalgamation in $K_\lambda$.]
\endroster
\medskip

\noindent
\S6 $\quad$ \underbar{A proof for ${\frak K}$ categorical in $\lambda^{+2}$}
\roster
\item "{${}$}"  [We give a shortcut to proving the main theorem by using
stronger assumptions.  If $I(\lambda^{+2},K) = 1$ and $I(\lambda^{+3},K) = 0$ 
then in some triple $(M,N,a) \in K_{\lambda^+}$, $a$ is essentially algebraic 
over $M$, i.e. this is a maximal triple.  Now first assuming for some pair
$M_0 \le_{\frak K} M_2$ in $K_\lambda$ we have unique amalgamation for every
possible $M_1$ with $M_0 \le_{\frak K} M_1 \in K_\lambda$ 
(and using stability), we get models in $\lambda^{++}$ thus contradicting 
the existence of maximal triples.  We
then use the methods of \S3 to prove there are enough cases of unique
amalgamation.]
\endroster
\medskip

\noindent
\S7 $\quad$ \underbar{Extensions and Conjugacy}
\roster
\item "{${}$}"  [We investigate types.  We prove that in ${\Cal S}(N),N \in 
K_\lambda$ that reduced implies inevitable, and that non-algebraic
extensions preserve the conjugacy classes (so solving the realize/
materialize problem).  Hence if $\langle N_i:i < \alpha \rangle$ is
$<_{\frak K}$-increasing in $K_\lambda$ and $\lambda$ divides $\alpha$ then
$N_\alpha$ is $(\lambda,\text{cf}(\alpha))$-saturated over $N_0$.]
\endroster
\medskip

\noindent
\S8 $\quad$ \underbar{Uniqueness of Amalgamation in ${\frak K}_\lambda$}
\roster
\item "{${}$}"  [We have by \S6 only $M^*_0 <_{\frak K} M^*_2$ in
${\frak K}_\lambda$ such that $M^*_0 \le_{\frak K} M_1 \Rightarrow M^*_0,
M_1,M^*_2$ has unique (disjoint) amalgamation.  Now if we have a
$\le_{\frak K}$-increasing continuous sequence $\langle N_i:i \le
\alpha \rangle$ such that $(N_i,N_{i+1}) \cong (M^*_0,M^*_2)$, 
we can amalgamate
$N_0,M_1,N_\alpha$ whenever $N_0 \le_{\frak K} M_1$, step by step.  So some
uniqueness is preserved and $N_\alpha$ can be any $(\lambda,\text{cf}
(\alpha))$-saturated model over $N_0$.  When we require also saturativity
of $M_i$ and of the resulting model, we get a nonforking relation denoted
NF$_{\lambda,\bar \delta}$.  We define the general
nonforking relation NF$_\lambda$ by closing NF$_{\lambda,\bar \delta}$
downward.  So we succeed to define a relation which should behave as a 
nice nonforking relation.
But we have to work to prove that this relation satisfies the expected 
properties, first for the ``saturated" version and then in the general 
case by a diagram chase.]
\endroster
\medskip

\noindent
\S9 $\quad$ \underbar{Nice extensions in $K_{\lambda^+}$}
\roster
\item "{${}$}"  [As we have a notion of ``nonforking" amalgamation in 
$K_\lambda$, we can use it to build $\le_{\frak K}$-extensions $M_1 \in 
K_{\lambda^+}$ for any given $M_0 \in K_{\lambda^+}$.  This defines naturally a two-place
relation $\le^*_{\lambda^+}$ on $K_{\lambda^+}$; ``being a nice
$\le_{\frak K}$-submodel".  We investigate it and variants.  In particular, we
prove the existence of disjoint amalgamation.]
\endroster
\medskip

\noindent
\S10 $\quad$ \underbar{Non-structure for $\le^*_{\lambda^+}$}
\roster
\item "{${}$}"  [Instead proving that all disjoint amalgamations in 
$K_\lambda$ are nonforking ones, we prove that on $K_{\lambda^+}$ the
relation $<^*_{\lambda^+}$
is the same as $\le^*_\lambda$, which is just as good for our purpose.  Toward
this we assume a failure and get many pairwise non-isomorphic models in
$K_{\lambda^{+2}}$, contradicting an assumption of \scite{0.B}(2).  
But once we have that $\le_{\frak K}$ agrees on $K_{\lambda^+}$ with 
$\le^*_{\lambda^+}$ we have disjoint amalgamation, which suffices for 
building a model in $K_{\lambda^{+3}}$.]
\endroster
\newpage

\head {\S0 Introduction} \endhead  \resetall
\bigskip

\noindent
Makowski \cite{Mw85} is a readable and good exposition concerning 
categoricity in abstract elementary classes around $\aleph_1$. 
\newline
Our primary concern is:
\demo{Problem \stag{0.A}} Can we have some 
(not necessarily much) classification theory for reasonable non-first order
classes ${\frak K}$ of models, with no uses of even traces of compactness 
and only mild set theoretic assumptions? 
\enddemo
\bn
Let me try to clarify the meaning of Problem \scite{0.A}. \newline
What is the meaning of ``\underbar{mild set theoretic assumptions}?"  We are
allowing requirements on cardinal arithmetic like GCH and weaker relatives.
Preferably, assumptions like diamonds and squares and even mild large 
cardinals will not be used (apart from cases provable in ZFC, or in ZFC plus
allowable assumptions).  

In fact we try to continue \cite{Sh:88}, where results about the 
number of non-isomorphic models in $\aleph_1$ and $\aleph_2$ of a 
sentence $\psi \in L_{\omega_1,\omega}$
are obtained.  Now in \cite{Sh:88} the theorem parallel to the present one
is proved assuming $2^{\aleph_0} < 2^{\aleph_1}$, so it is quite natural 
to use such assumptions here.

What is the meaning of ``some classification theory?"  While the 
dream is to have a classification theory as ``full" as the one obtained 
in \cite{Sh:c}, we will be glad to have theorems speaking just on having few
models in some cardinals or even categoricity and at least one in others.  
E.g. by \cite{Sh:88} if $\psi \in
L_{\omega_1,\omega}$ satisfies $1 \le I(\aleph_1,\psi) < 2^{\aleph_1}$ (and
$2^{\aleph_0} < 2^{\aleph_1}$) then $I(\aleph_2,\psi) > 0$.
Here $I(\mu,{\frak K}$) is the number of models in ${\frak K}$ of 
cardinality $\mu$, up to isomorphism.

What are ``reasonable non first order classes?"  This means we allow
classes of ``locally finite" or ``atomic" structures, or structures ``omitting
a type", or more generally the class
of models of a sentence in $L_{\kappa,\omega}$, (i.e. allowing conjunction 
$< \kappa$ but quantification only over a finite string) but not one 
restricting ourselves to e.g. well orderings.  In fact, we use ``abstract 
elementary classes" from \cite{Sh:88} (reviewed below).

What is the meaning of ``uses traces of compactness?"  For non
first order classes we cannot use the powerful compactness theorem, but 
there are many ways to get weak forms of it: one way is using large 
cardinals (compact cardinals in Makkai Shelah \cite{MaSh:285}, or 
just measurable cardinals as in Kolman Shelah \cite{KlSh:362}, or in 
\cite{Sh:472}).  Another way is to use
``non-definability of well ordering" which follows from the existence of
Ehrenfeucht-Mostowski models, and also from $\psi \in L_{\omega_1,\omega}$
having uncountable models (used extensively in \cite{Sh:88}).  
Our aim is to use \underbar{none of this} and 
we would like to see if any theory is left.

Above all, we hope the proofs will initiate classification theory in 
this case, so we hope the flavour will be one of
introducing and investigating notions of a model theoretic character.
Proofs of, say, a descriptive set theory character, will not satisfy this 
hope.

It seems to us that this goal is met here.  We prove:
\bigskip

\proclaim{Theorem \stag{0.B}}  
$(2^\lambda < 2^{\lambda^+} < 2^{\lambda^{++}})$.
Let ${\frak K}$ be an abstract elementary class. \newline
If ${\frak K}$ categorical in $\lambda,\lambda^+$ and $\lambda^{++}$ \ub{then}
$I(\lambda^{+3},{\frak K}) > 0$. 
\endproclaim
\mn
Here $\lambda^{+3} = \lambda^{+++}$ and in general $\aleph^{+ \beta}_\alpha$
means $\aleph_{\alpha + \beta}$. \nl
Of course, the categoricity in three successive cardinals is a strong
assumption.  Now note that in \cite{Sh:88}, the categoricity in $\aleph_0$
is gained ``freely", so the gap is smaller than seems at first glance.
Still it is better to have
\proclaim{Theorem \stag{0.C}}  $(2^\lambda < 2^{\lambda^+} < 
2^{\lambda^{++}})$.  Let ${\frak K}$ be an abstract elementary class. \newline
If ${\frak K}$ is categorical in $\lambda$ and $\lambda^+$, and $1 \le 
I(\lambda^{++},{\frak K}) < 2^{\lambda^{++}}$ \underbar{then}
$I(\lambda^{+3},{\frak K}) > 0$. \nl
Here we do not fully prove \scite{0.C} but we prove a slightly weaker 
version:
\endproclaim
\bigskip

\proclaim{Theorem \stag{0.D}}  If ${\frak K}$ is categorical in $\lambda$
and $\lambda^+$ with $\lambda > \aleph_0$ and $1 \le I(\lambda^{++},
{\frak K}) < \mu_{\text{wd}}(\lambda^{++})$, \underbar{then} $I(\lambda^{+3},
{\frak K}) > 0$.
\endproclaim
\mn
Note however:
\medskip
\roster
\item "{$(\alpha)$}"  A silly point; at exactly one point in the proof of
\scite{0.D} we assume 
$\lambda > \aleph_0$ (in the proof of \scite{4.2B}).  This
is silly as our intent is to prove for general $\lambda$ what we know for 
$\lambda = \aleph_0$ by \cite{Sh:88}; however, there we assume ${\frak K}$ is
PC$_{\aleph_1,\aleph_0}$, a reasonable assumption, but one which is not 
assumed here.  We shall complete this in \cite[\S2]{Sh:600},
so we do not mention the assumption $\lambda > \aleph_0$ in theorems relying 
on \scite{4.2B}. 
\sn
\item "{$(\beta)$}"  More seriously, at some point we assume toward a 
contradiction that $K_{\lambda^{+3}} = \emptyset$.  This is fine for proving 
theorem \scite{0.B}, but is not desirable if we want to develop a 
classification theory.  This will also be dealt with in \cite[\S2]{Sh:600}.
\sn
\item "{$(\gamma)$}"  Concerning $\mu_{\text{wd}}(\lambda^{+2})$:
in \scite{6.6} we get $I(\lambda^{++},K) \ge \mu_{\text{wd}}(\lambda^{++})$
instead of $I(\lambda^{++},K) = 2^{\lambda^{++}}$
(eliminated in \cite{Sh:600}, too).  The cardinal $\mu_{\text{wd}}(\chi)$ is 
defined in \scite{1.1}(6).  It is almost $2^\chi$; i.e. always 
$2^\chi < \mu_{\text{wd}}(\chi^+) \le 2^{\chi^+} \le \mu_{\text{wd}}
(\chi^+)^{\aleph_0}$.   Already
in \cite{Sh:87b} sometimes instead of getting $2^\chi$ we get only 
$\mu_{\text{wd}}(\chi)$ (again see \cite{Sh:600}).  If you 
are satisfied with $I(\chi,K) \ge \mu_{\text{wd}}(\chi)$ 
in those cases, you can ignore some proofs.
\endroster
\medskip

\noindent
We present below, as background, the following open questions which 
appeared in \cite{Sh:88}, for ${\frak K}$ an abstract elementary class, 
of course, e.g. the class of models
of $\psi \in L_{\lambda^+,\omega}$ with the relation $M \le_{\frak K} N$ 
being $M \prec_{\Cal L} N$ for ${\Cal L}$ a fragment of 
$L_{\kappa^+,\omega}$ to which $\psi$ belongs.  
In \cite{Sh:87a}, \cite{Sh:87b}, \cite{Sh:88} we prove:
\medskip
\roster
\item "{$(*)_3$}"  categoricity (of $\psi \in L_{\omega_1,\omega}(Q)$) in
$\aleph_1$ implies the existence of a model of $\psi$ of cardinality 
$\aleph_2$;
\sn
\item "{$(*)_4$}"  if $n > 0,2^{\aleph_0} < 2^{\aleph_1} < \cdots <
2^{\aleph_n},\psi \in L_{\omega_1,\omega}$ and $1 \le I(\aleph_\ell,\psi) 
< \mu_{wd}(\aleph_\ell)$ for $1 \le \ell \le n$, \ub{then} $\psi$ has a 
model of cardinality $\aleph_{n+1}$.
\endroster
\medskip

\noindent
Now the problems were:
\smallskip

\noindent
\underbar{Problem} 1) Prove $(*)_3,(*)_4$ in the context of an abstract 
elementary class ${\frak K}$ which is $PC_{\aleph_0}$.
\medskip

\noindent
\underbar{Problem} 2) Parallel results in ZFC; e.g. 
prove $(*)_3$ when $n = 1,2^{\aleph_0} = 2^{\aleph_1}$.
\medskip

\noindent
\underbar{Problem} 3) Construct examples; e.g. ${\frak K}$ 
(or $\psi \in L_{\omega_1,\omega}$),
categorical in $\aleph_0,\aleph_1,\dotsc,\aleph_n$ but not in
$\aleph_{n + 1}$.
\medskip

\noindent
\underbar{Problem} 4) If ${\frak K}$ is $PC_\lambda$ (abstract elementary
class), and is categorical in $\lambda$ and $\lambda^+$, does it 
necessarily have a model in $\lambda^{++}$? \newline

Concerning Problem 3, by Hart Shelah \cite{HaSh:323},2.10(2) + 3.8 
there is $\psi_n \in {\Cal L}_{\omega_1,\omega}$ categorical 
in $\aleph_0,\aleph_1,\dotsc,\aleph_{k-1}$, but not categorical in 
$\lambda$ if $2^\lambda > 2^{\aleph_{k-1}}$. \newline
The direct motivation for the present work, is that Grossberg asked me 
(Oct. 94) some questions in this neighborhood, in particular: 
\medskip

\noindent
\underbar{Problem} 5) Assume $K = \text{ Mod}(T)$, (i.e. $K$ is the class
of models of $T$), \newline
$T \subseteq L_{\omega_1,\omega},|T| = \lambda$,
$I(\lambda,K) = 1$ and $1 \le I(\lambda^+,K) < 2^{\lambda^+}$.  Does it
follow that $I(\lambda^{++},K) > 0$? \newline
We think of this as a test problem and much prefer a model theoretic to a
set theoretic solution.  This is closely related to Problem 4 above and to
\cite[Theorem 3.7]{Sh:88} (where we assume categoricity in $\lambda^+$, do
not require $2^\lambda < 2^{\lambda^+}$ \underbar{but} take 
$\lambda = \aleph_0$
or some similar cases) and \cite[Theorem 5.17(4)]{Sh:88} (and see
\cite[5.1,4.5]{Sh:88} on the assumptions) (there we require $2^\lambda
< 2^{\lambda^+},1 \le I(\lambda^+,K) < 2^{\lambda^+}$ \underbar{and}
$\lambda = \aleph_0$). \newline 
As said above, we are dealing with a closely related
problem.  Problem \scite{0.A} was stated {\it a posteriori} but is, I think, 
the real problem.
\medskip

In a first try we used more set theory, i.e. we used the definitional weak 
diamond on both $\lambda^+$ and $\lambda^{++}$
(see Definition \scite{1.7}) and things like ``a nice equivalence relation on 
${\Cal P}(\lambda)$ has either few or many classes" (see \S1).
\newline
Here we take a model theoretic approach.

We feel that this paper provides a reasonable positive solution to problem 
\scite{0.A}, with a classification theory flavor.  We shall continue in 
\cite{Sh:600} toward a parallel of \cite{Sh:87a}, \cite{Sh:87b}.
Grossberg and Shelah, in the mid-eighties, started to write a paper
\cite{GrSh:266} to prove that: if $\psi \in L_{\lambda^+,\omega}$ has models 
of arbitrarily large cardinality, and is categorical in $\mu^{+n}$ for each
$n$ and if $\mu \ge \lambda$ and $2^{\mu^{+n}} < 2^{\mu^{+n+1}}$ for 
$n < \omega$, \underbar{then} $\psi$ is categorical in every $\mu' \ge \mu$; 
this is a weak form of the upward part of Los' conjecture.  See Makkai Shelah 
\cite{MaSh:285} on $T \subseteq {\Cal L}_{\kappa,\omega}$ where $\kappa$ is a
compact cardinal; where we get downward and upward theorems for 
successor cardinals which are sufficiently bigger than $\kappa + |T|$.  
On the downward part, see \cite{KlSh:362},\cite{Sh:472} which
deals with a downward theorem for successor
cardinals which are sufficiently larger than $\kappa + |T|$ when the 
theory $T$ is in the logic ${\Cal L}_{\kappa,\omega}$ and $\kappa$ is 
measurable.  See also \cite{Sh:394} which deals with
abstract elementary classes with amalgamation, getting similar results with
no large cardinals. \newline
Part of \S1 and \S3 have a combinatorial character.  Most of the paper
forms the content of a course given in Fall '94 
(essentially without \S3, \S9, \S10).  The paper is written with an eye to
developing the model theory, rather than just proving Theorem \scite{0.B}. 
\bigskip

\noindent
\centerline {$* \qquad * \qquad *$}
\bigskip

\noindent
On abstract elementary classes see \cite{Sh:88} and \cite[II,\S3]{Sh:300}.  
To make the paper self-contained, we will review some relevant definitions 
and results.  We thank Gregory Cherlin for improving the writing of \S1-\S3
and for breaking \scite{3.16A}, \scite{3.16B} and \scite{3.16C} to 
three proofs. 
\newline
This work is continued in \cite{Sh:600}.
\newpage

\head{Abstract elementary classes} \endhead
\bigskip

\demo{\stag{0.1} Conventions}  ${\frak K} = (K,\le_{\frak K})$, where
$K$ is a class of $\tau$-models for a fixed vocabulary $\tau = \tau_K =
\tau_{\frak K}$ and $\le_{\frak K}$ is a two-place relation on the models in 
$K$.  We do not always 
strictly distinguish between $K$ and $(K,\le_{\frak K})$.  We shall 
assume that $K,\le_{\frak K}$ are fixed, and $M \le_{\frak K} N \Rightarrow 
M,N \in K$; and we assume that
the following axioms hold.  When we use $<$ in the sense of elementary
submodel for first order logic, we write $<_{{\Cal L}_{\omega,\omega}}$.
\enddemo
\bigskip

\definition{\stag{0.2} Definition}  \nl

$Ax\,0$: The validity of $M \in K$ or of $N \le_{\frak K} M$ depends on $N$
and $M$ only up to isomorphism - in the second case, isomorphism of the
pair. 
\medskip

$Ax\,I$: If $M \le_{\frak K} N$ then $M \subseteq N$ (i.e. $M$ is a submodel
of $N$).
\medskip

$Ax\,II$: $\le_{\frak K}$ is transitive and reflexive on $K$.
\medskip

$Ax\,III$: If $\lambda$ is a regular cardinal, $M_i \, (i < \lambda)$ is
$\le_{\frak K}$-increasing (i.e. $i < j < \lambda$ implies $M_i \le_{\frak K}
M_j$) and continuous (i.e. for limit ordinals $\delta < \lambda$ we have
\newline
$M_\delta = \dsize \bigcup_{i < \delta} M_i$) \underbar{then}
$M_0 \le_{\frak K} \dsize \bigcup_{i < \lambda} M_i$.
\medskip

$Ax\,IV$: If $\lambda$ is a regular cardinal, $M_i \, (i < \lambda)$ is
$\le_{\frak K}$-increasing continuous, \newline
$M_i \le_{\frak K} N$ \underbar{then}
$\dsize \bigcup_{i < \lambda} M_i \le_{\frak K} N$.
\medskip

$Ax\,V$:  If $M_0 \subseteq M_1$ and $M_\ell \le_{\frak K} N$ for $\ell =
0,1$, then $M_0 \le_{\frak K} M_1$.
\medskip

$Ax\,VI$:  There is a cardinal $\lambda$ such that: if $A \subseteq N$ and 
$|A| \le \lambda$ \underbar{then} for some $M \le_{\frak K} N$ we have 
$A \subseteq |M| \le \lambda$.  We define the L\"owenheim-Skolem number
$LS({\frak K})$ as the least such $\lambda$ with $\lambda \ge |\tau|$.
\enddefinition
\bigskip

\noindent
\underbar{Notation}:  $K_\lambda = \{ M \in K:\|M\| = \lambda\}$ and
$K_{< \lambda} = \dsize \bigcup_{\mu < \lambda} K_\mu$. \newline
${\Cal L}_{\omega,\omega}$ is first order logic. \newline
A theory in ${\Cal L}(\tau)$ is a set of sentences from ${\Cal L}(\tau)$.
\bigskip

\definition{\stag{0.3} Definition}  The embedding $f:N \rightarrow M$ 
is a ${\frak K}$-embedding or a $\le_{\frak K}$-embedding if its range is 
the universe of a model $N' \le_{\frak K} M$, 
(so $f:N \rightarrow N'$ is an isomorphism (onto)).
\enddefinition
\bigskip

\noindent
Very central in \cite{Sh:88}, but peripheral here is:
\definition{\stag{0.4} Definition}  1) For a logic ${\Cal L}$ and 
vocabulary $\tau,{\Cal L}(\tau)$ is the set of ${\Cal L}$-formulas in 
this vocabulary. \newline
2)  Let $T_1$ be a theory in 
${\Cal L}_{\omega,\omega}(\tau_1), \tau \subseteq \tau_1$ vocabularies,
$\Gamma$ a set of types in ${\Cal L}_{\omega,\omega}(\tau_1)$; (i.e.
for some $m$, a set of formulas $\varphi(x_0,\dotsc,x_{m-1}) \in 
{\Cal L}_{\omega,\omega}(\tau_1)$). \newline
3) $EC(T_1,\Gamma) = \{M:M \text{ a } \tau_1 \text{-model of } T_1
\text{ which omits every } p \in \Gamma\}$. \newline
4) $PC_\tau(T_1,\Gamma) = PC(T_1,\Gamma,\tau) = \{M:M \text{ is a } 
\tau \text{-reduct of some } M_1 \in EC(T_1,\Gamma)\}$. \nl
5)  We say that ${\frak K}$ is $PC^\mu_\lambda$ if for some 
$T_1,T_2,\Gamma_1,\Gamma_2$ and $\tau_1$ and $\tau_2$ we have: 
$T_\ell$ is a first order theory in the vocabulary $\tau_\ell,
\Gamma_\ell$ is a set of types and $K = PC(T_1,\Gamma_1,\tau_{\frak K})$ and
$\{(M,N):M \le_{\frak K} N,M,N \in K\} = PC(T_2,\Gamma_2,\tau')$ where
\newline
$\tau' = \tau_1 \cup \{P\}$ ($P$ a new one place predicate), $|T_\ell| 
\le \lambda,
|\Gamma_\ell| \le \mu$ for $\ell = 1,2$.  If $\mu = \lambda$, we may omit it.
\enddefinition
\bigskip

\demo{\stag{0.5} Example}  If $\tau_1 = \tau,T_1,\Gamma$ are as above, and
$(K,\le_{\frak K})$ are defined by \newline
$K =: EC(T_1,\Gamma)$, $\le_{\frak K}\, =: \,
\prec_{{\Cal L}_{\omega,\omega}}$, 
then the pair satisfies the Axioms from \scite{1.2} and 
$LS({\frak K}) \le |T_1| + \aleph_0$.
\enddemo
\bigskip

\demo{\stag{0.5X} Example}  $(V = L)$. 

Let $\text{cf}(\lambda) \ge \aleph_1,n < \omega$ \underbar{then} for some
$\psi \in L_{\lambda^+,\omega}$ we have: $\psi$ has no model of cardinality
$\lambda^{+(n+1)}$, and is categorical in $\lambda^{+n}$ (i.e. has one and 
only one model up to isomorphism). \newline
Let $M^* = (L_{\lambda^{+n}},\in,i)_{i \le \lambda}$ and let $\psi$ be

$$
\align
\bigwedge &\biggl\{ \varphi:\varphi \text{ is a first order sentence which }
M^* \text{ satisfies} \biggr\} \\
  &\wedge (\forall x) \left( x \in \lambda \equiv 
\dsize \bigvee_{i < \lambda} x = i \right).
\endalign
$$
\enddemo
\bigskip

\proclaim{\stag{0.6} Lemma}  Let $I$ be a 
directed set (i.e. partially ordered
by $\le$, such that any two elements have a common upper bound). \newline
1) If $M_t$ is defined for $t \in I$, and $t \le s \in I$ implies $M_t
\le_{\frak K} M_s$ \underbar{then} for every $t \in I$ we have
$M_t \le_{\frak K} \dsize \bigcup_{s \in I} M_s$. \newline
2) If in addition $t \in I$ implies $M_t \le_{\frak K} N$ \underbar{then}
$\dsize \bigcup_{s \in I} M_s \le_{\frak K} N$.
\endproclaim
\bigskip

\demo{Proof}  By induction on $|I|$ (simultaneously for 1) and 2)), or see
\cite[1.6]{Sh:88}.  \hfill$\square_{\scite{0.6}}$
\enddemo
\bigskip

\proclaim{\stag{0.7} Lemma}  \nl
1) Let $\tau_1 = \tau \cup \{F^n_i:i < LS({\frak K}) \text{ and }
n < \omega\},F^n_i$ an $n$-place function symbol (assuming, of course,
$F^n_i \notin \tau$).  If $M_1$ is an expansion of $M$ to a $\tau_1$-model
and $\bar a \in {}^n|M|$ for some $n$, let $M_{\bar a}$ be the submodel of
$M$ with universe $\{F^n_i(\bar a):i < LS(k)\}$.   \newline
Every model $M$ from $K$ can be expanded to a $\tau_1$-model $M_1$ such that:
\medskip
\roster
\item "{$(A)$}"  $M_{\bar a} \le_{\frak K} M$ for any $\bar a \in {}^n|M|$;
\sn
\item "{$(B)$}"  if $\bar a \in {}^n|M|$ then $\|M_{\bar a}\|
\le LS({\frak K})$;
\sn
\item "{$(C)$}"  if $\bar b$ is a subsequence of $\bar a$ (even up to
rearrangement), \underbar{then} $M_{\bar b} \le_{\frak K} M_{\bar a}$;
\sn
\item "{$(D)$}"  for every $N_1 \subseteq M_1$ we have $N_1 \restriction \tau
\le_{\frak K} M$.
\endroster
\medskip

\noindent
2) If $N \le_{\frak K} M$ is given, then we can choose the expansion $M_1$
so that clauses (A) - (D) hold and
\medskip
\roster
\item "{$(E)$}"  $N = N_1 \restriction \tau$ for some $N_1 \subseteq M_1$.
\endroster
\endproclaim
\bigskip

\demo{Proof}  We define by induction on $n$, the values of $F^n_i(\bar a)$
for every $i < LS({\frak K})$, \newline
$\bar a \in {}^n|M|$ such that $\bar a \subseteq N \Rightarrow M_{\bar a} 
\subseteq N$ when we are proving (2).  By $Ax\,VI$ there is an
$M_{\bar a} \le_{\frak K} M,\|M_{\bar a}\| \le LS({\frak K})$ such that
$|M_{\bar a}|$ includes

$$
\bigcup \{ M_{\bar b}:\bar b \text{ a subsequence of } \bar a
\text{ of length } < n\} \cup \bar a
$$
\mn
and $M_{\bar a}$ does not depend on
the order of $\bar a$.  Let $|M_{\bar a}| = \{c_i:i < i_0 \le 
LS({\frak K})\}$, and define $F^n_i(\bar a) = c_i$ for $i < i_0$ and $F^n_i
(\bar a) = c_0$ for $i_0 \le i < LS({\frak K})$ (so we can demand ``$F^n_i$ is
symmetric").

Clearly our conditions are satisfied: if $\bar b$ is a subsequence of
$\bar a$ then $M_{\bar b} \le_{\frak K} M_{\bar a}$ by $Ax\,V$.
\hfill$\square_{\scite{0.7}}$
\enddemo
\bigskip

\proclaim{\stag{0.8} Lemma}  1) There is a set $\Gamma$ of types in 
${\Cal L}_{\omega,\omega}(\tau_1)$ (where $\tau_1$ is as in Lemma \scite{0.7})
such that $K=PC(\emptyset,\Gamma,\tau)$, so it is a
$PC^{2^{LS({\frak K})}}_{LS({\frak K})}$-class. \newline
The types above consist of quantifier-free formulas. \nl
2) Moreover, if $N_1 \subseteq M_1 \in EC(\emptyset,\Gamma)$, and $N,M$ are
the $\tau$-reducts of $N_1,M_1$ respectively \underbar{then} $N \le_{\frak K} 
M$. Also, if
$N \le_{\frak K} M$ then there is an $\tau_1$-expansion $M_1$ of $M$ and 
a submodel $N_1$ of $M_1$ such that $M_1 \in EC(\emptyset,\Gamma)$ and
$N_1 \restriction \tau = N$. \newline
3) Also $\{(M,N):N \le_{\frak K} M\}$ is a 
$PC^{(2^{LS({\frak K})})}_{LS({\frak K})}$-class, hence ${\frak K}$ is as
well.
\endproclaim
\bigskip

\demo{Proof}  1) Let $\Gamma_n$ be the set of complete quantifier free
$n$-types $p(x_0,\dotsc,x_{n-1})$ in ${\Cal L}_{\omega,\omega}(\tau_1)$ 
such that:  if $M_1$ is an
$L_1$-model, $\bar a$ realizes $p$ in $M_1$ and $M$ is the $L$-reduct of
$M_1$, then $M_{\bar b} \le_{\frak K} M_{\bar a}$ for any subsequence of
$\bar b$ of $\bar a$.  Where 
$M_{\bar c}\,\,(\text{for }\bar c \in {}^m|M_1|)$ is the
submodel of $M$ whose universe is $\{F^m_i(\bar c):i < LS({\frak K})\}$
(and there are such submodels).
\medskip

Let $\Gamma$ be the set of $p$ which are complete quantifier free 
$n$-types for some $n < \omega$ in ${\Cal L}_{\omega,\omega}(\tau_1)$ and
which 
do not belong to $\Gamma_n$.  So if $M^1$ is in $PC(\emptyset,\Gamma,\tau_1)$ 
then by \scite{0.6} we have $M_1 \restriction \tau \in K$ hence
$PC(\emptyset,\Gamma,L) \subseteq K$ and by \scite{0.7}(1) we have
$K \subseteq PC(\emptyset,\Gamma,L)$. \newline
2) The first phrase is proven as in (1).  For the second phrase, use
\scite{0.7}(2). \newline
3) Follows from (2). \hfill$\square_{\scite{0.8}}$
\enddemo
\bigskip

\definition{\stag{0.9} Conclusion}  There is $\tau_1$ with 
$\tau \subseteq \tau_1$ and $|\tau_1| \le LS({\frak K})$ such that: 
for any $M \in K$ and 
any $\tau_1$-expansion $M_1$ of $M$ which is in $EC(\emptyset,\Gamma)$,

$$
N_1 \subseteq M_1 \Rightarrow N_1 \restriction \tau \le_{\frak K} M
$$

$$
N_1 \subseteq N_2 \subseteq M_1 \Rightarrow N_1 \restriction \tau 
\le_{\frak K} N_2 \restriction \tau.
$$
\enddefinition
\bigskip

\definition{\stag{0.10} Conclusion}  If for every 
$\alpha < (2^{LS({\frak K}})^+,K$ has a
model of cardinality $\ge \beth_\alpha$ \underbar{then} $K$ has a model in
every cardinality $\ge LS({\frak K})$.
\enddefinition
\bigskip

\demo{Proof}  Use \scite{0.8} and the value of the Hanf number for: 
models of a first order theory omitting a given set of types, for 
languages of cardinality $LS(K)$ (see \cite[VII,\S5]{Sh:c}). \newline
${}$ \hfill$\square_{\scite{0.10}}$
\enddemo
\bigskip

\definition{\stag{0.11} Definition}  For $\lambda$ regular 
$> LS({\frak K})$ and $N \in {\frak K}_\lambda$ we say 
$\bar N = \langle N_\alpha: \alpha < \lambda \rangle$ is a representation 
of $N$ if $\bar N$ is $\le_{\frak K}$-increasing 
continuous, $\| N_\alpha\| < \lambda$ and $N = \dsize 
\bigcup_{\alpha < \lambda} N_\alpha$.  If $\lambda = \mu^+$ then if not
said otherwise, we require $\|N_\alpha\| = \mu$.

How will we define types, and in particular, the set ${\Cal S}(M)$ of complete
types over $M$, when no formulas are present?  If we have a ``monster model"
${\frak C}$ we can use automorphisms; but any such ``monster" is far down 
the road.  So we will ``chase diagrams" in $K_\lambda$ (being careful, not to
use excessivley large models).  This gives us a relation of ``having the same
type" we call $E^{at}_\mu$, {\it but} this relation in general is
not transitive (if we do not have amalgamation in $K_\lambda$).  So
$E_\mu$ will be defined as the transitive closure of $E^{at}_\mu$.
\enddefinition
\bigskip

\definition{\stag{0.12} Definition}  1)  The two-place relation 
$E^{\text{at}}_M$ is defined on triples $(M,N,a)$ with $M$ fixed, $M
\le_{\frak K} N \in K_{\|M\|}$, and $a \in N$ by:

$$
\align
(M,N_1,a_1)E^{\text{at}}_M &(M,N_2,a_2) \text{ \underbar{iff} there is }
N \in K_\mu \text{ and } \le_{\frak K} \text{-embeddings} \\
  &f_\ell:N_\ell \rightarrow N \text{ for } \ell = 1,2 \text{ such that:} \\
  &f_1 \restriction M = \text{ id}_M = f_2 \restriction M \text{ and }
f_1(a_1) = f_2(a_2).
\endalign
$$
\mn
Let $E_M$ be the transitive closure of $E^{\text{at}}_M$. \nl
2) For $\mu \ge LS({\frak K})$ and $M \in K_\mu$ we define ${\Cal S}(M)$ as 
$\{\text{tp}(a,M,N):M \le_{\frak K} N \in K_\mu \text{ and } a \in N\}$ where
$\text{tp}(a,M,N) = (M,N,a)/E_M$. \nl
3)  We say ``$a$ realizes $p$ in $N$" if $a \in N,p \in {\Cal S}(M)$ and 
for some $N' \in K_\mu$ we have $M \le_{\frak K} N' \le_{\frak K} N,a \in N'$,
and $p = \text{ tp}(a,M,N')$. \newline
4) We say ``$a_2$ strongly realizes $(M,N^1,a^1)/E^{\text{at}}_M$ in $N$" if 
for some $N^2,a^2$ we have $M \le_{\frak K} N^2 \le_{\frak K} N,
a_2 \in N^2$, and $(M,N^1,a^1)\,E^{at}_M\,(M,N^2,a^2)$. \newline
5) We say $M_0 \in {\frak K}$ is an amalgamation base if letting
$\lambda = \|M_0\|$ we have: for every
$M_1,M_2 \in {\frak K}_\lambda$ and $\le_{\frak K}$-embeddings
$f_\ell:M_0 \rightarrow M_\ell$ (for $\ell = 1,2$) there is $M_3 \in
{\frak K}_\lambda$ and $\le_{\frak K}$-embeddings $g_\ell:M_\ell \rightarrow
M_3$ (for $\ell=1,2$) such that $g_1 \circ f_1 = g_2 \circ f_2$.
\newline
6) We say ${\frak K}$ is stable in $\lambda$ \underbar{if} $LS({\frak K})
\le \lambda$ and for all $M \in K_\lambda$ we have 
$|{\Cal S}(M)| \le \lambda$.
\enddefinition
\bigskip

\demo{\stag{0.13} Observation}  If $M$ is an amalgamation base then
$E_M = E^{\text{at}}_M$, and we have: 

$$
``a \text{ strongly realizes } (M,N,b)/E_M \text{ in } N" \text{ iff }
``a \text{ realizes } (M,N,b)/E_M.
$$
\enddemo
\bigskip

\definition{\stag{0.16} Definition}  \nl
1) ${\frak K}$ has the 
$\lambda$-amalgamation property or has amalgamation in $\lambda$, if every 
$M_0 \in {\frak K}_\lambda$ is an amalgamation base (see \scite{0.12}(5)
above). \newline
2) $N$ is universal over $M$ if: $M \le_{\frak K} N$ and if
$M \le_{\frak K} N' \in K_{\le\|N\|}$ then $N'$ can be 
$\le_{\frak K}$-embedded into $N$ over $M$. \newline
3) ${\frak K}$ has universal extensions in $\lambda$ \underbar{if} for every
$M \in K_\lambda$ there is $N$ such that:
\medskip
\roster
\item "{$(a)$}"  $M \le_{\frak K} N \in K_\lambda$
\sn
\item "{$(b)$}"  $N$ is universal over $M$.
\endroster
\medskip

\noindent
4) $N_1,N_2$ have a joint embedding if for some $N \in K$ there are
$\le_{\frak K}$-embeddings $h_\ell$ of $N_\ell$ into $N$ for $\ell =1,2$.
Let JEP$_\mu$ (JEP) means this holds for $N_1,N_2$ in $K_\mu$ (in $K$).
\newline
5) ${\frak K}_\lambda$ has unique amalgamation (or ${\frak K}$ has unique
amalgamation in $\lambda$) when: if $M^\ell_i \in K_\lambda$ for $\ell < 2,
i < 4$, and for $i=1,2$ we have $h^\ell_{i,0}$ is a 
$\le_{\frak K}$-embedding of $M^\ell_0$ into 
$M^\ell_i,h^\ell_{3,i}$ is a $\le_{\frak K}$-embedding of $M^\ell_i$ into
$M^\ell_3,h^\ell_{3,1} \circ h^\ell_{1,0} = h^\ell_{3,2} \circ h^\ell_{2,0}$
and Rang$(h^\ell_{3,1}) \cap \text{ Rang}(h^\ell_{3,2}) = \text{ Rang}(h^\ell
_{3,1} \circ h^\ell_{1,0})$ 
and for $i < 3,f_i$ is an isomorphism from $M^0_i$ onto $M^1_i,f_0 \subseteq
f_1,f_0 \subseteq f_2$ then for some $N \in K_\lambda$ there are
${\frak K}$-embedding $h_\ell:M^\ell_3 \rightarrow N$ such that
$h_0 \circ h^0_{3,i} = h_1 \circ h^1_{3,i}$, for $i=1,2$. \nl
6)  Let $p_\ell \in {\Cal S}(M_\ell)$ for 
$\ell=0,1$.  We say $p_0 \le p_1$ \underbar{if} $M_0 \le_{\frak K} M_1$ 
and for some $N$ and $a$ we have $M_1 \le_{\frak K} N \in 
K_{\|M_1\| + LS({\frak K})},a \in N$ and 
$\text{tp}(a,M_\ell,N) = p_\ell$ for $\ell = 0,1$.  We also write 
$p_0 = p_1 \restriction M_0$ ($p_0$ is unique knowing $M_0,M_1,p_1$
hence $p_1 \restriction M_0$ is well defined). 
\enddefinition
\bigskip

\proclaim{\stag{0.17} Claim}  \nl
1) If ${\frak K}$ is categorical in 
$\lambda$ and $LS({\frak K}) \le \lambda$ then: there is a model in 
$K_{\lambda^+}$ \underbar{iff} for some (equivalently, every) model 
$M \in K_\lambda$ there is $N$ such that $M <_{\frak K} N \in K_\lambda$ 
\underbar{iff} for some (equivalently every) $M \in K_\lambda$ there is
$N$ such that $M <_{\frak K} N \in K_{\lambda^+}$. 
\newline
2) If ${\frak K}$ has amalgamation in $\lambda,LS({\frak K}) \le \lambda$, 
and $M_0 \le_{\frak K} M_1$ are in $K_\lambda$ with $M_0 \le_{\frak K} N_0 
\in K_{\lambda^+}$ \underbar{then} we can find $h$ and $N_1$ such that 
$N_0 \le_{\frak K} N_1 \in K_{\lambda^+}$ and $h$ is a 
$\le_{\frak K}$-embedding of $M_1$ into 
$N_1$ extending $\text{id}_{M_0}$.   We can allow $N_0 \in K_\mu$ with 
$\mu > \lambda$ if ${\frak K}$ has amalgamation in every $\lambda' \in 
[\lambda,\mu)$.  \newline
3) Assume ${\frak K}$ has amalgamation in $\lambda$ and $LS({\frak K}) \le
\lambda$.  If $M_0 \le M_1$ are from
$K_\lambda$ and $p_0 \in {\Cal S}(M_0)$ \ub{then} we can find an 
extension $p_1 \in {\Cal S}(M_1)$ of $p_0$.
\endproclaim
\bigskip

\demo{Proof}  1)  As then we can choose by induction on $i < \lambda^+$ 
models \newline
$M_i \in K_\lambda,\le_{\frak K}$-increasing continuous, 
$M_i \ne M_{i+1}$, for $i = 0$ use $K_{\lambda^+} \ne \emptyset$, for \newline
$i$ limit take union, for
$i = j+1$ use the previous sentence; so \newline
$M_{\lambda^+} = \cup \{ M_i:i < \lambda^+\} \in K_{\lambda^+}$ as required.
\newline
2), 3)  Left to the reader.  \hfill$\square_{\scite{0.17}}$
\enddemo
\bigskip

\remark{\stag{0.17A} Remark}  We can here add the content of \scite{6.4}, 
\scite{6.5}, \scite{6.6}.
\endremark
\bigskip

\definition{\stag{0.18} Definition}  1) For $\lambda > LS({\frak K})$ we say
``$N \in {\frak K}$ is $\lambda$-saturated" \underbar{if} for every
$M \le_{\frak K} N$ of cardinality $< \lambda$, if $M \le_{\frak K} N' 
\in K_{< \lambda}$ and $a' \in N'$ then some $a \in N$ strongly realizes 
$(M,N',a')/E^{\text{at}}_M$ (in the interesting cases this suffices). \nl
2)  We say ``$N \in {\frak K}$ is $\lambda$-saturated above $\mu$ (or is
$\lambda$-saturated $\ge \mu$)" 
\underbar{if} above we restrict ourselves to $M$ of cardinality $\ge \mu$.
\newline
(Cf. $(\lambda,\kappa)$-saturated in Definition \scite{0.21}, \scite{0.22}
below).
\enddefinition
\bigskip

\demo{\stag{0.18A} Fact}  1) If $\mu + LS({\frak K}) < \lambda,\lambda$ is 
regular, ${\frak K}$ has the amalgamation property in every 
$\mu' \in [\mu,\lambda)$, and for all $M \in K_{[\mu,\lambda)}$ we have
$|{\Cal S}(M)| \le \lambda$ and $\lambda = \text{ cf}(\lambda)$, \ub{then} 
there is some $M \in K_\lambda$ saturated above $\mu$. \nl
2) Assume $\lambda > \mu \ge LS({\frak K})$ and $N \in K$ and ${\frak K}$ has
amalgamation in $K_\mu$, for every $\mu_1 \in [\mu,\lambda)$.  Then: $N$
is $\lambda$-saturated above $\mu$ \ub{iff} for every $M \le_{\frak K} N$ of
cardinality $< \lambda$, every $p \in {\Cal S}(M)$ is realized in $N$ 
(i.e. for some $a \in N$ we have tp$(a,M,N) = p$). \nl
3) If $LS({\frak K}) \le \mu_0 \le \mu'_0 \le \mu' \le \mu$ and $M$ is
$\mu$-saturated above $\mu_0$, \ub{then} 
it is $\mu'$-saturated above $\mu'_0$.
If $LS({\frak K}) \le \mu_0 < \mu$ \ub{then}: $M$
is $\mu$-saturated above $\mu_0$ \underbar{iff} for every $\lambda \in [\mu_0,
\mu),M$ is $\lambda^+$-saturated above $\lambda$. 
\enddemo
\bigskip

\definition{\stag{0.14} Definition}   
The type $p \in {\Cal S}(M)$ is local when: if for any directed partial
order $I$ and models $M_t \le_{\frak K} M$ for $t \in I$ with 
$I \models t \le s \Rightarrow M_t \le_{\frak K} M_s$ and $M = 
\dsize \bigcup_{t \in I} M_t$, and any $p' \in {\Cal S}(M)$ with 
$(p \restriction M_t = p' \restriction M_t)$ for all $t \in I,p = p'$.  
We say $M$ is local if every $p \in {\Cal S}(M)$ is, and ${\Cal K}$ is local 
if every $M \in {\Cal K}$ is.  We can add ``above $\mu$" as in 
Definition \scite{0.18}(2).
\enddefinition
\bigskip

\definition{\stag{0.15} Definition}  1) $I(\lambda,K) = I(\lambda,{\frak K})$
is the number of $M \in K_\lambda$ up to isomorphism.
2) ${\frak K}$ (or $K$) is categorical in $\lambda$ if $I(\lambda,K) = 1$.
\newline
3) $IE(\lambda,{\frak K}) = \text{ sup}\{|K'|:K' \subseteq K_\lambda
\text{ and for } M \ne N$ in $K'_\lambda,M$ is not 
$\le_{\frak K}$-embeddable into $N\}$.  If we write $IE(\lambda,K) \ge \mu$,
we mean that for some $K' \subseteq K_\lambda$ as above, $|K'| \ge \mu$, and
similarly for $= \mu$.  If there is a problem with attainment of the supremum
we shall say explicitly.
\enddefinition
\bigskip

\proclaim{\stag{0.19}  The Model-homogeneity = Saturativity Lemma}  Let 
$\lambda > \mu \ge LS({\frak K})$ and $M \in K$. \newline
1) $M$ is $\lambda$-saturated above $\mu$ \underbar{iff} $M$ is 
$({\Cal D}_{\frak K},\lambda)$-homogeneous above $\mu$, which means: for every
$N_1 \le_{\frak K} N_2 \in K$ such that $\mu \le \|N_1\| \le \|N_2\| < 
\lambda$ and $N_1 \le_{\frak K} M$ there is a $\le_{\frak K}$-embedding $f$ 
of $N_2$ into $M$ over $N_1$. \newline
2)  If $M_1,M_2 \in K_\lambda$ are $\lambda$-saturated above $\mu < \lambda$ 
and for some $N_1 \le_{\frak K} M_1,N_2 \le_{\frak K} M_2$, both 
of cardinality $\in [\mu,\lambda)$, we have $N_1 \cong N_2$ \ub{then} 
$M_1 \cong M_2$; in fact, any isomorphism $f$ from $N_1$ onto $N_2$ 
can be extended to an isomorphism from $M_1$ onto $M_2$. \newline
3) If in (2) we demand only ``$M_2$ is $\lambda$-saturated" and $M_1 \in
K_{\le \lambda}$ \underbar{then} $f$ can be extended to a 
$\le_{\frak K}$-embedding from $M_1$ into $M_2$. \newline
4) In part (2) instead of $N_1 \cong N_2$ it suffices to assume that
$N_1$ and $N_2$ can be $\le_{\frak K}$-embedded into some $N \in K$, which 
holds if ${\frak K}$ has the JEP.
\endproclaim
\bigskip

\demo{Proof}  1)  The ``if" direction is easy as $\lambda > LS({\frak K})$.  
Let us prove the other direction, so without loss of generality
$\|N_1\| = \mu$.

By \scite{0.16}(6) without loss of generality $\lambda$ is regular, moreover
$N_2$ has cardinality $\|N_1\|$. 

Let $|N_2| = \{ a_i:i < \kappa \}$, and we know $\kappa = \|N_1\| =
\|N_2\| < \lambda$.  We define by induction on $i \le \kappa,N^i_1,N^i_2,f_i$ 
such that:
\medskip
\roster
\item "{(a)}"  $N^i_1 \le_{\frak K} N^i_2$ and $\| N^i_2 \| < \lambda$
\sn
\item "{(b)}" $N^i_1$ is $\le_{\frak K}$-increasing continuous in $i$
\sn
\item "{(c)}" $N^i_2$ is $\le_{\frak K}$-increasing continuous in $i$
\sn
\item "{(d)}" $f_i$ is a $\le_{\frak K}$-embedding of $N^i_1$ into $M$
\sn
\item "{(e)}" $f_i$ is increasing continuous in $i$
\sn
\item "{(f)}" $a_i \in N^{i+1}_1$
\sn
\item "{(g)}" $N^0_1 = N_1,N^0_2 = N_2,f_0 = id_{N_1}$
\sn
\item "{(h)}" $N^i_1$ and $N^i_2$ has cardinality $\kappa$.  
\endroster
\medskip

\noindent
For $i = 0$, clause $(g)$ gives the definition.  For $i$ limit let: \newline
\smallskip
\noindent
$N^i_1 = \dsize \bigcup_{j < i} N^j_1$ and 
\sn
$N^i_2 = \dsize \bigcup_{j < i} N^j_2$ and 
\sn
$f_i = \dsize \bigcup_{j < i} f_j$. \newline
\smallskip
\noindent
Now (a)-(f) continues to hold by continuity (and $\|N^i_2\| < \lambda$ as
$\lambda$ is regular).
\medskip

For $i$ successor we use our assumption; more elaborately, let
$M^{i-1}_1 \le_{\frak K} M$ be $f_{i-1}(N^{i-1}_1)$ and $M^{i-1}_2,
g_{i-1}$ be such that $g_{i-1}$ is an isomorphism from $N^{i-1}_2$ \newline
onto $M^{i-1}_2$ extending
$f_{i-1}$, so $M^{i-1}_1 \le_{\frak K} M^{i-1}_2$ (but \wilog \nl
$M^{i-1}_2 \cap M = M^{i-1}_1$).  Now apply the saturation
assumption with $M,M^{i-1}_1$, \nl
$\text{tp}(g_{i-1}(a_{i_1}),M^{i-1}_1,M^{i-1}_2)$ 
here standing for $N,M,p$ there
(note: $a_{i-1} \in N_2 = N^0_2 \subseteq N^{i-1}_2$ and 
$\mu > \|N^{i-1}_2\| = \|M^{i-1}_2\|$ and $\|M^{i-1}_1\| =
\|N^{i-1}_1\| \ge \|N^0_1\| = \|N_1\| \ge \mu_0$ so the requirements and
also the cardinality in Definition \scite{0.18} holds).
So there is $b \in M$ such that tp$(b,M^{i-1}_1,M) = \text{ tp}(g_{i-1}
(a_{i-1}),M^{i-1}_1,M^{i-1}_2)$.  Moreover (if ${\frak K}$ has amalgamation
in $\mu$ the proof is slightly shorter) remembering the second sentence in
\scite{0.18}(1) which speaks about ``strongly realizes" there is $b \in M$
such that $b$ strongly realizes $(M^{i-1}_1,M^{i-1}_2,g_{i-1}(a_{i-1}))/
E^{\text{at}}_{M^{i-1}_1}$ in $M$.  This means (see Definition \scite{0.12}
(3)) that for some $M^{i,*}_1$ we have $M^{i-1}_1 \le_{\frak K} M^{i,*}_1
\le_{\frak K} M$ and $(M^{i-1}_1,M^{i-1}_2,g_{i-1}(a_{i-1}))
E^{\text{at}}_{M^{i-1}_1}(M^{i-1}_1,M^{i,*}_1 b)$.  This means (see
Definition \scite{0.12}(1) that $M'$ too has cardinality $\kappa$ and there
is $M^{2,*}_i \in K_\mu$ such that $M^{i-1}_1 \le_{\frak K} M^{i,*}_2$ and
there are $\le_{\frak K}$-embeddings $h^i_2,h^i_1$ of $M^{i-1}_2,M^{i,*}_1$
into $M^{i,*}_2$ over $M^{i-1}_1$ respectively, such that $h^i_2(g_{i-1}
(a_{i-1})) = h^i_1(b)$. \nl
Now changing names, \wilog \, $h^i_2$ is the identity.
\nl
Let $N^i_2,h_i$ be such that $N^{i-1}_2 <_{\frak K} N^i_2,h_i$ an 
isomorphism from $N^i_2$ onto $M^{i,*}_2$ extending $g_{i-1}$.  Let
$N^i_1 = h^{-1}_i \circ h^i_1(M^{i,*}_1),f_i = (h^i_1)^{-1} \circ 
(h_i \restriction N^i_1)$.
\medskip

We have carried the induction.  Now $f_\kappa$ is a $\le_{\frak K}$-
embedding of $N^\kappa_1$ into $M$ over $N_1$, but $|N_2| = \{ a_i:i <
\kappa \} \subseteq N^\kappa_1$, so $f_\kappa \restriction N_2:N_2
\rightarrow M$ is as required. \newline
2), 3)  By the hence and forth argument (or see \cite[II,\S3]{Sh:300} =
\cite[II,\S3]{Sh:e}). \newline
4) Easy, too.  \hfill$\square_{\scite{0.19}}$
\enddemo
\bigskip

\remark{\stag{0.20} Remark}  Note that by \scite{0.19}(2) if $M$ 
is $\mu$-saturated above $\mu_0$ and ${\frak K}$ has the JEP$_{\mu_0}$
then ${\frak K}$ has $\lambda$-amalgamation 
for each $\lambda \in [\mu_0,\mu)$.
\endremark
\bigskip

\definition{\stag{0.21} Definition}  Fix $\lambda \ge \kappa$ with $\kappa$
regular. \nl
1) We say that $N_1 \in K_\lambda$ is $(\lambda,\kappa)$-saturated over 
$N_0$ or that $(N_1,c)_{c \in N_0}$ is $(\lambda,\kappa)$-saturated (and
$\lambda$-saturated of cofinality $\kappa$ means $(\lambda,\kappa)$-saturated)
if: \newline
there is a sequence $\langle M_i:i < \kappa \rangle$ which
is $\le_{\frak K}$-increasing continuous with
$M_0 = N,M_\kappa = N_1$ and $M_{i+1} \in K_\lambda$
universal over $M_i$ (see \scite{0.16}(2)). \nl
2) If we omit $\kappa$, we mean $\kappa = \text{ cf}(\lambda);
(\lambda,\alpha)$-saturated means $(\lambda,\text{cf}(\alpha))$-saturated;
and $N_1$ is $(\lambda,1)$-saturated over $N_0$ means just
$N_0 \le_{\frak K} N_1$ are in $K_\lambda$. \newline
\enddefinition
\bigskip

\proclaim{\stag{0.22} Claim}  Fix $\lambda \ge \kappa$ with $\kappa$
regular. \nl
1) If $N_1$ is $(\lambda,\kappa)$-saturated over $N$, \ub{then} $N_1$ is 
unique over $N$. \newline
2) If $K_\lambda \ne \emptyset,\lambda \ge \kappa = \text{ cf}(\kappa)$ and
over every $M \in K_\lambda$ there is $N$ with 
$M \le_{\frak K} N \in K_\lambda$
universal over $M$, \underbar{then} for every $N \in K_\lambda$ there is
$N_1 \in K_\lambda$ which is $(\lambda,\kappa)$-saturated over $N$. \newline
3) If ${\frak K}_\lambda$ has amalgamation and ${\frak K}$ is stable in
$\lambda$ (i.e. $M \in K_\lambda \Rightarrow |{\Cal S}(M)| \le \lambda$)
\underbar{then} every $M \in K_\lambda$ has a universal extension (so part
(2)'s conclusion holds).
\endproclaim
\bigskip

\demo{Proof}  See \cite[Ch.II]{Sh:300} or check. \hfill$\square_{\scite{0.22}}$
\enddemo
\bigskip

\noindent
We do not need at present but recall from \cite{Sh:88}:
\proclaim{\stag{0.23} Claim}  There is $\tau' \supseteq \tau \cup 
\{P_0,P_1,P_2,c\}$ of cardinality $\le LS({\frak K})$ with $c$ an
individual constant, with $P_\ell$ unary predicates, and a set $\Gamma$ of 
quantifier free types such that:
\medskip
\roster
\item "{$(a)$}"  if $M' \in PC_{\tau'}(\emptyset,\Gamma)$ and $M_\ell =
(M' \restriction \tau) \restriction P^{M'}_\ell$ for $\ell =0,1,2$,
then \newline
$M_\ell \in K,M_0 \le_{\frak K} M_1,M_0 \le_{\frak K} M_2,c^{M'} \in
M_2$, and \newline
$N' \subseteq M' \Rightarrow (N' \restriction \tau) \restriction
P^{N'}_\ell \le_{\frak K} M_\ell$, and there is no $b \in M_1$ satisfying:
{\roster
\itemitem{ $\otimes$ }  for every $\bar a \in {}^{\omega >}(P^{M'}_0)$,
letting $N_{\bar a}$ be the $\tau'$-submodel of $M'$ generated by
$\bar a$ and $M^\ell_{\bar a} = M_\ell \restriction (M_\ell \cap N_{\bar a})$,
we have $M^\ell_{\bar a} \le_{\frak K} M_\ell$, and $b$ strongly realizes
tp$(c^{M'},M^0_{\bar a},M^2_{\bar a})$ in $M^1_{\bar a}$
\endroster}
\item "{$(b)$}"  if $M^*_\ell,c$ are as in $(a)$ and $M_0 = M_1 \cap M_2$.
then for some $M'$ we have clause (a).
\endroster
\endproclaim
\bigskip

\demo{Proof}  Should be clear; see \cite{Sh:88}.
\hfill$\square_{\scite{0.23}}$
\enddemo 
\bigskip

\remark{\stag{0.24} Remark}  \nl
1) Claim \scite{0.23} enables us to translate results of the form:
two cardinal with omitting types theorem in $\lambda_2$ implies the existence 
of one in $\lambda_1$, provided that types are local in the sense that 
$p \in {\Cal S}(M)$
is determined by $\langle p \restriction N:N \le_{\frak K} M,\|N\| \le
\lambda \rangle$. \newline
2) This enables us to prove implications between cases of
$\lambda$-categoricity, \underbar{if} we have a nice enough theory of types as
in \cite[VIII,\S4]{Sh:c}; if we have in $\lambda_2$ a saturated model,
categoricity in $\lambda_1$ implies categoricity in $\lambda_2$.  Also (if
we know a little more) categoricity in $\lambda_2$ is equivalent to the
non-existence of a non-saturated model in $\lambda_2$. 
\endremark
\bigskip

\proclaim{\stag{0.28} Claim}  1) Assume $M_n \le_{\frak K} M_{n+1},M_n \in
K_\lambda,{\frak K}$ has amalgamation in $\lambda$.  If $p_n \in {\Cal S}
({\Cal M}_n),p_n \le p_{n+1}$ (i.e. $p_n = p_{n+1} \restriction M_n$, see
\scite{0.16}(6)), \underbar{then} there is \nl
$p \in {\Cal S}(\dsize \bigcup_{n < \omega} M_n)$ such that 
$n < \omega \Rightarrow p_n \le p$.
\newline
2) If $\langle M_i:i \le \delta \rangle$ is $\le_{\frak K}$-increasing
continuous, $p_i \in {\Cal S}(M_i),(j < i \Rightarrow p_j \le p_i)$, \nl
$p_i = \text{ tp}(a_i,M_i,N_i)$ and $h_{i,j}$ is a $\le_{\frak K}$-embedding 
of $N_j$ into $N_i$ (for $j < i < \delta$) such that $h_{i,j} \restriction M_j
 = \text{ id}_{M_j},h_{i,j}(a_j) = a_i$, \underbar{then} there is $p_\delta \in
{\Cal S}(M_\delta),i \le \delta \Rightarrow p_i \le p_\delta$. \newline
3) If ${\frak K}$ has amalgamation in $\lambda$ and is stable in $\lambda$ 
(i.e.
$M \in K_\lambda \Rightarrow |{\Cal S}(M)| \le \lambda$), \underbar{then}
\medskip
\roster
\item "{$(a)$}"  every $M \in K_\lambda$ has a universal extension;
\sn
\item "{$(b)$}"  for every $M \in K_\lambda$ and regular $\theta \le \lambda$
there is $N \in K_\lambda$ which is $(\lambda,\theta)$-saturated over $M$.
\endroster
\endproclaim
\newpage

\head {\S1 Weak Diamond} \endhead  \resetall
\bigskip

\definition{\stag{1.1} Definition}  Fix $\lambda$ regular and uncountable.

$$
\align
\text{WDmTId}(\lambda,S,\bar \chi) = \biggl\{ A:&A \subseteq 
\dsize \prod_{\alpha \in S} \chi_\alpha,\text{ and for some function } F 
\text{ with domain} \tag"{$1)$}" \\
  &\dsize \bigcup_{\alpha < \lambda} {}^\alpha(2^{< \lambda}) \text{ mapping }
{}^\alpha(2^{< \lambda}) \text{ into } \chi_\alpha, \\
  &\text{ for every } \eta \in A, \text{ for some } f \in 
{}^\lambda(2^{< \lambda}) \text{ the set} \\
  &\{ \delta \in S:\eta (\delta) = F(f \restriction \delta)\} 
\text{ is not stationary} \biggr\}.
\endalign
$$
\medskip

\noindent
(Note: WDmTId stand for weak diamond target 
ideal) \footnote{in \cite[AP,\S1]{Sh:b}, \cite[AP,\S1]{Sh:f} we express 
cov$_{\text{Wdmt}}(\lambda,S) > \mu^*$ by allowing $f(0) \in \mu^* < \mu$}.
\nl
Here we can replace $2^{< \lambda}$ by any set of this cardinality, and so
we can replace $f \in {}^\lambda(2^{< \lambda})$ by $f_1,\dotsc,f_n \in
{}^\lambda(2^{< \lambda})$ with $F$ being $n$-place.

$$
\text{cov}_{\text{wdmt}}(\lambda,S,\bar \chi) = 
\text{ Min}\biggl\{ |{\Cal P}|:{\Cal P} \subseteq \text{ WDmTId}(\lambda,S) 
\text{ and } {}^S 2 \subseteq \dsize \bigcup_{A \in {\Cal P}} A \biggr\} 
\tag"{$2)$}"
$$

$$
\align
\text{WDmTId}_{< \mu}(\lambda,S,\bar \chi) = \biggl\{ A:&\text{ for some } 
i^* < \mu \text{ and } A_i \in \text{ WDmTId}(\lambda,S) \text{ for} 
\tag"{$3)$}" \\
  &\,i < i^* \text{ we have } A \subseteq \dsize \bigcup_{i < i^*} A_i
\biggr\}
\endalign
$$

$$
\text{WDmId}_{< \mu}(\lambda,\bar \chi) = \biggl\{ S \subseteq \lambda:
\text{cov}_{\text{wdmt}}(\lambda,S) < \mu \biggr\}. \tag "{$4)$}"
$$
\medskip

\noindent
5) Instead of ``$< \mu^+$" we write $\mu$, if we omit $\mu$ we mean
$(2^{< \lambda})$.  If $\bar \chi$ is constantly $2$ we may omit it, if
$\chi_\alpha = 2^{|\alpha|}$ we may write pow instead of $\bar \chi$.
\medskip

$$
\text{Let } \mu_{\text{wd}}(\lambda,\bar \chi) = 
\text{ cov}_{\text{wdmt}}(\lambda,\lambda,\bar \chi). \tag "{$6)$}"
$$
\mn
7)  We say that the weak diamond holds on $\lambda$ if $\lambda \notin
\text{ WDmId}(\lambda)$.
\enddefinition
\bigskip

\noindent
By \cite{DvSh:65}, \cite[XIV,1.18(2),1.8]{Sh:b} 
(presented better in \cite[AP,\S1]{Sh:f}, note: \scite{1.2}(4) below rely on 
\cite{Sh:460}) we have:

\proclaim{\stag{1.2} Theorem} \newline
1)  If $\lambda = \aleph_1,2^{\aleph_0} < 2^{\aleph_1},\mu \le (2^{\aleph_0})$
or even $2^\theta = 2^{< \lambda} < 2^\lambda,\mu = (2^\theta)^+$, or just:
for some $\theta,2^\theta = 2^{< \lambda} < 2^\lambda,\mu \le
2^\lambda$, and $\chi^{< \lambda} < \mu$ for $\chi < \mu$, \ub{then} 
$\lambda \notin \text{ WDmId}_{< \mu}(\lambda)$.  If in addition
$(*)_{< \mu,\lambda}$ below holds, \underbar{then}
$\lambda \notin \text{WDmId}_{< \mu}(\lambda,\text{pow})$, where:
\medskip

$(*)_{\mu,\lambda}$ there are no $A_i \in [\mu]^{\lambda^+}$ for $i < 
2^\lambda$ such that $i \ne j \Rightarrow |A_i \cap A_j| < \aleph_0$
\mn
and $(*)_{< \mu,\lambda}$ means $(*)_{\chi,\lambda}$ holds for $\chi < \mu$.
\smallskip

\noindent
2)  If $\mu \le \lambda^+$ or cf$([\mu_1]^{\le \lambda},\subseteq) < \mu$
for $\mu_1 < \mu$ or $\mu = \aleph_0$ \underbar{then}
WDmId$_{< \mu}(\lambda,\bar \chi)$ is a normal ideal on $\lambda$.  
If this ideal is not trivial, \ub{then} $\lambda = \text{ cf}(\lambda) > 
\aleph_0,2^{< \lambda} < 2^\lambda$. \newline
\smallskip

\noindent
3)  A sufficient condition for $(*)_{< \mu,\lambda}$ is:
\medskip
\roster
\item "{$(a)$}"  $\mu \le 2^\lambda \and (\forall \alpha < \mu)
(|\alpha|^{\aleph_0} < 2^\lambda)$.
\endroster
\medskip

\noindent
4) Another sufficient condition for $(*)_{< \mu,\lambda}$ is:
\medskip
\roster 
\item "{$(b)$}"  $\mu \le 2^\lambda \and \lambda \ge \beth_\omega$.
\endroster
\endproclaim
\bigskip

\remark{\stag{1.2A} Remark}  1) So if cf$(2^\lambda) < \mu$ (which holds if
$2^\lambda$ is singular and 
$\mu = 2^\lambda$) then $(*)_{< \mu,\lambda}$ implies that there is 
$A \subseteq {}^\lambda 2,|A| < 2^\lambda,A \notin
\text{ WDmTId}(\lambda)$. \newline
2) Some related definitions appear in \scite{1.7}; we use them below (mainly
DfWD$_{< \mu}(\lambda)$), but as in a first reading it is recommended to 
ignore them, the definition is given later. \newline
3) We did not look again at the case $(\forall \sigma < \lambda)(2^\sigma <
2^{< \lambda} < 2^\lambda)$.
\endremark
\bigskip

\noindent
As in \cite[3.5]{Sh:88}, (\cite[2.7]{Sh:87a}, \cite[6.3]{Sh:87b}):
\proclaim{\stag{1.3} Claim}  Assume $\lambda \notin 
\text{ WDmId}_{< \mu}(\lambda)$ or at least DfWD$_{< \mu}(\lambda)$ \nl
(where $\lambda = \text{ cf}(\lambda) > \aleph_0$) and ${\frak K}$ is an 
abstract elementary class. \newline
1)  Assume ${\frak K}$ is categorical in $\chi,\lambda = \chi^+$, and 
${\frak K}$ has a model in $\lambda$ (if LS$({\frak K}) \le \chi$ 
this is equivalent to: the model $M \in K_\chi$ is not 
$\le_{\frak K}$-maximal).  Assume further ${\frak K}$ does not have the 
$\lambda$-amalgamation property in $\chi$. \underbar{Then} for any $M_i \in
{\frak K}_\lambda$ for $i < i^* < \mu$, there is $N \in {\frak K}_\lambda$ not
$\le_{\frak K}$-embeddable into any $M_i$ (and the assumptions of part 
(2) below holds). \newline
2)  Assume $M_\eta \in K_{< \lambda}$ for $\eta \in {}^{\lambda >}2,
M_\eta = \dsize \bigcup_{\alpha < \ell g(\eta)} M_{\eta \restriction
(\alpha + 1)},\nu \triangleleft \eta \Rightarrow M_\nu \le_{\frak K} 
M_\eta$, and $M_{\eta \char 94 \langle 0 \rangle},
M_{\eta \char 94 \langle 1 \rangle}$ cannot be amalgamated over $M_\eta$ 
(hence $M_\eta \ne M_{\eta \char 94 \langle \ell \rangle})$.  Set 
$M_\eta =: \dsize \bigcup_{\alpha < \lambda} M_{\eta \restriction \alpha}$ 
for $\eta \in {}^\lambda 2$.  Clearly $M_\eta$ belongs to $K_\lambda$.
For the DfWD$_{< \mu}(\lambda)$ version assume also
\medskip
\roster
\item "{$(*)$}"  $\langle M_\eta:\eta \in {}^{\lambda >}2 \rangle$ is
definable (even just by ${\Cal L}_{\lambda,\lambda}$) in \newline
${\frak B} = ({\Cal H}(\chi),\in,<^*_\chi,{\frak K}_{< \chi},\lambda,\mu)$.
\endroster
\medskip

\noindent
\underbar{Then} for any $N_i \in {\frak K}_\lambda$ for $i < i^* < \mu$, 
there is $\eta \in {}^\lambda 2$ such that: $M_\eta$ is not \newline
$\le_{\frak K}$-embeddable into any $N_i$.
\newline
3)  In part (2), if $LS({\frak K}) \le \lambda$ we can allow 
$N_i \in K_{\kappa_i}$ if $\dsize \sum_{i < i^*} \text{ cf}
([\kappa_i]^\lambda,\subseteq) < \mu$. \newline
4) If we use DfWD$_{< \mu}(\lambda)$: assume only 
$\lambda \notin \text{ WDmId}_{< \mu}(\lambda,\bar \chi)$.
Part (2) holds if $M_\eta$ is defined for $\eta \in \dbcu_{\alpha < \lambda}
\dsize \prod_{i < \alpha} \chi_i$, and $\varepsilon < \zeta < \chi_i,
\eta \in \dsize \prod_{j <i} \chi_j \Rightarrow M_{\eta \char 94 \langle
\varepsilon \rangle},M_{\eta \char 94 \langle \zeta \rangle}$ cannot be
amalgamated over $M_\eta$.  The assumption of Part (4) holds (hence the
conclusion of Part (2)) holds if we assume that for $M \in
K_\chi,i < \lambda$ there are $\chi_i \le_{\frak K}$-extensions of $M$ in
$K_\lambda$, which pairwise cannot be amalgamated over $M$.
\endproclaim
\bigskip

\demo{Proof}  1) It is straightforward to choose 
$M_\eta \in K_\chi$ for $\eta \in
{}^\alpha 2$ by induction on $\alpha$, as required in part (2).  Then use
part (2) to get the desired conclusion. \newline
2) Without loss of generality the universe of $N_i$ is $\lambda$ and the
universe of $M_\eta$ is an ordinal $\gamma_\eta$ such that
$\eta \in {}^{\lambda >}2 \Rightarrow \gamma_\eta < \lambda$ and 
$\eta \in {}^\lambda 2 \Rightarrow \gamma_\eta = \lambda$.
The reader can ignore the ``DfWD$_{< \mu}(\lambda)$" version (ignoring the
$h_\eta$'s,$g$) if he likes.  For $\alpha < \lambda$ and $\eta \in
{}^\alpha 2$ let the function $h_\eta$ be $h_\eta(i) = M_{\eta \restriction
(i+1)}$ for $i < \ell g(\eta)$.  Let $\langle M_\eta:\eta \in {}^{\lambda >}2
\rangle$ be the $<^*_\chi$-first such object.  For each 
$i < i^*$ we define $A_i \subseteq {}^\lambda 2$ by $A_i = \{\eta
\in {}^\lambda 2:M_\eta \text{ can be } \le_{\frak K}$-embedded into
$N_i\}$.  

For $\eta \in A_i$ choose $f_\eta:M_\eta \rightarrow N_i$, a $\le_{\frak K}$-
embedding, hence $f_\eta \in {}^\lambda \lambda$.  We also define a function
$F_i$ from $\dsize \bigcup_{\alpha < \lambda}({}^\alpha 2 \times
{}^\alpha \lambda$) to $\{0,1\}$ by:

$$
\align
F_i(\eta,f) \text{ is}: 0 \quad \text{ \underbar{if} } &f \text{ is a }
\le_{\frak K} \text{-embedding of } M_\eta \text{ into } N_i \\
  &\text{with range } \subseteq \ell g(\eta) \text{ which can be extended to
a } \\
  &\le_{\frak K} \text{-embedding of } M_{\eta \char 94 \langle 0 \rangle}
\text{ into } N_i
\endalign
$$
\bigskip

$\qquad F_i(\eta,f) \text{ is } 1 \quad \text{ otherwise}$.
\medskip

\noindent
Now for any $\eta \in A_i$, the set

$$
E = \{\delta < \lambda:\gamma_{\eta \restriction \delta} = \delta
\text{ and } f_\eta \restriction \delta \text{ is a function from } \delta
\text{ to } \delta\}
$$
\mn
is a club of $\lambda$.  For every $\delta \in E$ 
clearly $F(\eta \restriction \delta,f_\eta \restriction \delta) = 
\eta(\delta)$
(as $M_{\eta \char 94 \langle 0 \rangle},M_{\eta \char 94 \langle 1 \rangle}$
cannot be amalgamated over $M_\eta$).

Hence (for the ``Def" version see Definition \scite{1.7}(2) using 
\scite{1.8}(1), \scite{1.8}(3)) we have 
$A_i \in \text{ WDmTId}^{\text{Def}}(\lambda)$.  
As $i^* < \mu$ clearly $\dsize \bigcup_{i < i^*} A_i \in 
\text{ WDmId}^{\text{Def}}_{< \mu}(\lambda)$ and hence by assumption of 
the claim ${}^\lambda 2 \ne \dbcu_i A_i$.  Take 
$\eta \in {}^\lambda 2 \backslash \dsize \bigcup_{i<i^*} A_i$.  Then
$M_\eta$ is as required. \newline
3) Without loss of generality the universe of $N_i$ is $\kappa_i$.  Let 
${\Cal P}_i \subseteq [\kappa_i]^\lambda$ be a set of minimal
cardinality such that $(\forall B)[B \subseteq \kappa_i \and |B| \le
\lambda \rightarrow (\exists B' \in {\Cal P}_i)(B \subseteq B')]$.  As
$LS({\frak K}) \le \lambda$ we can find for each $A \in {\Cal P}_i$, a model
$N^i_A \le_{\frak K} N_i$ of cardinality $\lambda^+$ whose universe includes
$A$.  Now apply part (2) to $\{N^i_A:i < i^*$ and $A \in {\Cal P}_i\}$. \nl
4) Same proof.  $\hfill\square_{\scite{1.3}}$
\enddemo
\bn
We give three variants of the preceding:
\proclaim{\stag{1.3A} Claim}  1) Assume
\medskip
\roster
\item "{$(*)_1$}"  $\lambda = \text{ cf}(\lambda) > \aleph_0$
{\roster
\itemitem{ (a) }  $M_\eta$ is a $\tau$-model of cardinality $< \lambda$ for
$\eta \in {}^{\lambda >} 2$ and
\sn
\itemitem{ (b) }  for each $\eta \in {}^\lambda 2,\langle M_{\eta \restriction
\alpha}:\alpha < \lambda \rangle$ is $\subseteq$-increasing continuous with
union, called $M_\eta \in K$, of cardinality $\lambda$
\sn
\itemitem { (c) }  if $\eta \in {}^{\lambda >}2,\eta \triangleleft \rho_\ell
\in {}^\lambda 2$ for $\ell = 1,2$ then $M_{\rho_1},M_{\rho_2}$ are not
isomorphic over $M_\eta$
\endroster}
\item "{$(*)_2$}"  $\lambda \notin \text{ WDmId}_{< \mu}(\lambda)$
\endroster
\medskip

\noindent
\underbar{Then}  $I(\lambda,K) \ge \mu$, and in fact we can find $X \subseteq
{}^\lambda 2$ of cardinality $\ge \mu$ such that \newline
$\{M_\rho:\rho \in X\}$ are pairwise non-isomorphic. \newline
2) Assume
\medskip
\roster
\item "{$(*)^d_1$}"  like $(*)_1$ in part (1) but in addition
{\roster
\itemitem{ (d) }  $\langle M_\eta:\eta \in {}^{\lambda >}2 \rangle$ is
definable in ${\frak B} = {\frak B}_\chi$
\endroster}
\item "{$(*)^d_2$}"  $\lambda \notin$ DfWD$_{< \mu}(\lambda)$.
\endroster
\medskip

\noindent
\underbar{Then} the conclusion of part (1) holds. \nl
3) The parallel of \scite{1.3}(4) holds.
\endproclaim
\bigskip

\demo{Proof}  1) Let $\{ N_i:i < i^*\}$ be a maximal subset 
of $\{ M_\rho:\rho \in {}^\lambda 2\}$ consisting of pairwise
non-isomorphic models, and use the proof of \scite{1.3}(2) with
$f_\eta:N_i \simeq M_\eta$. \newline
2), 3)  Left to the reader.  \hfill$\square_{\scite{1.3A}}$
\enddemo
\bigskip

\proclaim{\stag{1.4} Claim}  1)  Assume
\medskip
\roster
\item "{$(*)_1$}"   $M_\eta \in K_{< \lambda}$ for $\eta \in {}^{\lambda >}2,
\langle M_{\eta \restriction \alpha}:\alpha \le \ell g(\eta) \rangle$ is
$\le_{\frak K}$-increasing continuous, and 
$M_{\eta \char 94 \langle 1 \rangle}$ cannot be $\le_{\frak K}$-embedded 
into $M_\nu$ over $M_\eta$ when \nl
$\eta \char 94 \langle 0 \rangle \trianglelefteq \nu \in {}^{\lambda >}2$ 
\sn
\item "{$(*)_2$}"   $\lambda = \text{ cf}(\lambda) > \aleph_0$, 
$\lambda \notin \text{ WDmId}_{< \mu}(\lambda)$, and $\lambda$ is a successor 
cardinal, or at least there is no $\lambda$-saturated normal ideal 
on $\lambda$, or at least WDmId$(\lambda)$ is not $\lambda$-saturated 
(which holds if for some $\theta < \lambda,\{ \delta < \lambda:
\text{cf}(\delta) = \theta\} \notin \text{WDmId}(\lambda)$).
\endroster
\medskip

\noindent
\underbar{Then} there is $A \subseteq {}^\lambda 2,|A| = 2^\lambda$ 
such that: if $\eta_1 \ne \eta_2$ are in $A$ then (taking
$M_\eta = \dbcu_{\alpha < \lambda} M_{\eta \restriction \alpha}$)
\medskip
\roster
\item "{$(\alpha)$}"  $M_{\eta_1} \ncong M_{\eta_2}$ for $\eta_1 \ne
\eta_2 \in A$ and
\sn
\item "{$(\beta)$}"  if $(2^\chi)^+ < 2^\lambda$ for $\chi < \lambda$
\underbar{then} we can also achieve:
$M_{\eta_1}$ cannot be $\le_{\frak K}$-embedded into $M_{\eta_2}$.
\endroster
\medskip
\noindent 
2)  Under the assumptions of \scite{1.3}(1) we can find 
$\langle M_\eta:\eta \in {}^{\lambda >}2 \rangle$ as in the 
assumption of \scite{1.3}(2). \newline
3) Under the assumption of \scite{1.3}(2) the assumption of \scite{1.4}(1) 
holds. \newline
4) Under the assumption of \scite{1.4}(1) we have $I(\lambda,{\frak K}) =
2^\lambda$ and if $(2^\chi)^+ < 2^\lambda$ then $IE(\lambda,{\frak K})
= 2^\lambda$. \newline
5) The parallel of \scite{1.3}(4) holds.
\endproclaim
\bigskip

\demo{Proof of \scite{1.4}}  1) The proof of \cite[3.5]{Sh:88} works (see 
the implications preceding it).  More elaborately, we divide the proof into
cases according to the answer to the following:
\mn
\underbar{Question}:  Is there $\eta^* \in {}^{\lambda >}2$ such that for
every $\nu$ satisfying $\eta^* \trianglelefteq \nu \in {}^{\lambda >}2$ 
there are
$\rho_0,\rho_1 \in {}^{\lambda >}2$ such that: $\nu \triangleleft \rho_0,
\nu \trianglelefteq \rho_1$, and $M_{\rho_0},M_{\rho_1}$ cannot be
amalgamated over $M_{\eta^*}$?
\medskip

\noindent
We can find a function $h:{}^{\lambda >} 2 \rightarrow {}^{\lambda >}2$,
such that:
\medskip
\roster
\widestnumber\item{$(b)_{\text{yes}}$}
\item "{$(a)$}"  the function $h$ is one-to-one, 
preserving $\triangleleft$ and $(h(\nu)) \char 94
\langle \ell \rangle \trianglelefteq h(\nu \char 94 \langle \ell \rangle)$
\sn
\item "{$(b)_{\text{yes}}$}"  \underbar{when} the answer to the question is 
yes, it is exemplified by $\eta^* = h(\langle \rangle)$ and
$M_{h(\nu \char 94 \langle 0 \rangle)},M_{h(\nu \char 94 \langle 1 \rangle)}$
cannot be amalgamated over $M_{\eta^*}$ (for every $\nu \in {}^{\lambda >}2)$
\sn
\item "{$(b)_{\text{no}}$}"  when the answer to the question above is no,
$h(\langle \rangle) = \langle \rangle$ and if $\nu \char 94 \langle 0 \rangle 
\trianglelefteq \rho_0,\nu \char 94 \langle 1 \rangle \triangleleft \rho_1$ 
then $M_{h(\rho_0)},M_{h(\rho_1)}$ can
be amalgamated over $M_{h(\nu)}$.
\endroster
\medskip

\noindent
Without loss of generality $h$ is the identity, by renaming (and we can 
preserve $(*)^d_1$ of \scite{1.3A}(2) in the relevant case).  Also clearly
$M_{\eta \char 94 \langle \ell \rangle} \ne M_\eta$ (by the non-amalgamation
assumption). 
\enddemo
\bigskip

\noindent
\underbar{Case 1}  The answer is yes.   We do not use the non 
$\lambda$-saturation of WDmId$(\lambda)$ in this case. \newline
For any $\eta \in {}^\lambda 2$ and $\le_{\frak K}$-embedding $g$ of 
$M_{\langle \rangle}$ into $M_\eta =: 
\dsize \bigcup_{\alpha < \lambda} M_{\eta \restriction \alpha}$, let

$$
A_{\eta,g} =: \{ \nu \in {}^\lambda 2:\text{there is a } \le_{\frak K}
\text{-embedding of } M_\nu \text{ into } M_\eta \text{ extending } g\}
$$

$$
A_\eta =: \{ \nu \in {}^\lambda 2:\text{there is a } \le_{\frak K}
\text{-embedding of } M_\nu \text{ into } M_\eta \}.
$$
\medskip

\noindent
So: $|A_{\eta,g}| \le 1$ and $\eta \in A_\eta$ 
(as if $\nu_1,\nu_2 \in A_{\eta,g}$
are distinct then for some ordinal $\alpha < \lambda$ and $\nu \in {}^\alpha
2$ we have $\nu =: \nu_0 \restriction \alpha = \nu_1 \restriction
\alpha,\nu_0(\alpha) \ne \nu_1(\alpha)$ and use the choice of
$h(\nu \char 94 \langle \ell \rangle)$).
\sn
Since $A_\eta = \cup \{ A_{\nu,g}:g \text{ is a } \le_{\frak K}
\text{-embedding of } M_{\langle \rangle} \text{ into } M_\eta \}$, we have
$|A_\eta| \le 2^{< \lambda}$. \nl
Hence we can choose by induction on $\zeta < 2^\lambda,\eta_\zeta \in
{}^\lambda 2 \backslash \dsize \bigcup_{\xi < \zeta} A_{\eta_\xi}$ (existing
by cardinality considerations as $2^{< \lambda} < \lambda$).  Then 
$\xi < \zeta \Rightarrow M_{\eta_\xi} \ncong M_{\eta_\zeta}$ so we have
proved clause $(\alpha)$. \newline
If $(2^{< \lambda})^+ < 2^\lambda$, then use the Hajnal free subset theorem
(\cite{Ha61}) to choose distinct
$\{ \eta_\zeta:\zeta < 2^\lambda\} \subseteq {}^\lambda 2$ such that 
$\zeta \ne \xi \Rightarrow \eta_\zeta \notin A_{\eta_\xi}$, so we have 
proved clause $(\beta)$, too.
\bigskip

\noindent
\underbar{Case 2}:  The answer is no. \newline
Again, without loss of generality $M_\eta$ has as universe the ordinal
$\gamma_\eta$. \newline
Let $\langle S_i:i < \lambda \rangle$ be a partition of $\lambda$ to sets,
none of which is in WDmId$(\lambda)$.  For each $i$ we define a function 
$F_i$ as follows:
\medskip
\roster
\item "{{}}"  if $\delta \in S_i,\eta,\nu \in {}^\delta 2,\gamma_\eta = 
\gamma_\nu = \delta$, and $f:\delta \rightarrow \delta$ then
$$
\align
F_i(\eta,\nu,f) = 0 &\text{ if } f \text{ can be extended to an embedding
of} \\
  &M_{\nu \char 94 \langle 0 \rangle} \text{ into some } M_\rho \text{ with }
\eta \char 94 \langle 0 \rangle \trianglelefteq \rho.
\endalign
$$
\endroster
\medskip

$\qquad \qquad \quad F_i(\eta,\nu,f) = 1$ otherwise.
\bigskip

\noindent
So as $S_i \notin \text{ WDmId}(\lambda)$, for some 
$\eta^*_i \in {}^\lambda 2$ we have:
\medskip
\roster
\item "{$(*)$}"  for every $\eta \in {}^\lambda 2,\nu \in {}^\lambda 2$ and
$f \in {}^\lambda \lambda$ the following set of ordinals $i < \lambda$ is 
stationary:
$$
\{ \delta \in S_i:F_i(\eta \restriction \delta,\nu \restriction \delta,
f \restriction \delta) = \eta^*_i(\delta)\}.
$$
\endroster
\medskip

\noindent
Now for any $X \subseteq \lambda$ let $\eta_X,\rho_X \in {}^\lambda 2$ be 
defined by:
\medskip
\roster
\item "{{}}"  if $\alpha \in S_i$ then $i \in X \Rightarrow \eta_X(\alpha) 
= 1 - \eta^*_i (\alpha),i \notin X \Rightarrow \eta_X(\alpha) = 0 
\text{ and}$ \newline 
$\rho_X = \eta_{\{2i:i \in X\} \cup \{2i+1:i \notin X\}}$.
\endroster
\medskip

\noindent
Now we show
\medskip
\roster
\item "{$(*)$}"  if $X,Y \subseteq \lambda$, and $X \ne Y$ then $M_{\rho_X}$
cannot be $\le_{\frak K}$-embedded into $M_{\rho_Y}$.
\endroster
\medskip

\noindent
Clearly $(*)$ will suffice for finishing the proof.

Assume toward a contradiction that $f$ is a 
$\le_{\frak K}$-embedding of $M_{\rho_X}$ into
$M_{\rho_Y}$; as $X \ne Y$ there is $i$ such that $i \in X \Leftrightarrow
i \notin Y$ so there is $j \in \{2i,2i+1\}$ such that $\rho_X \restriction
S_j = \langle 1-\eta^*_j(\alpha):\alpha \in S_j \rangle$ and 
$\rho_Y \restriction S_j$ is identically zero.  Clearly $E = \{ \delta:f
\text{ maps } \delta \text{ into } \delta\}$ is a club of $\lambda$ and
hence $S_j \cap E \ne \emptyset$.

So if $\delta \in S_j \cap E$ then $f \restriction M_{(\rho_Y \restriction 
(\delta + 1))} = f \restriction M_{(\rho_Y \restriction \delta) \char 94 
\langle 0 \rangle}$ extend $f \restriction M_{\rho_Y \restriction \delta}$ 
and is a $\le_{\frak K}$-embedding of it into some 
$M_{(\rho_X \restriction \alpha)},\alpha < \lambda$ large enough.

Now by the choice of $h$ we get

$$
\delta \in S_j \cap E \Rightarrow F(\rho_X \restriction \delta,\rho_Y 
\restriction \delta,f \restriction \delta) = \rho_X(\delta) = 
1 - \eta^*_j(\delta).
$$

But this contradicts the choice of $\eta^*_j$. \newline
2) This corresponds to Case I, where we do not need the weak diamond, but just
$2^{< \lambda} < 2^\lambda$.  \nl
3),4),5).  Check, similarly.  \hfill$\square_{\scite{1.4}}$
\bigskip

\demo{\stag{1.6A} Conclusion}  1) Assume
\medskip
\roster
\item "{$(*)_1$}"  for $\eta \in {}^{\lambda >}2,M_\eta \in K_{< \lambda}$
and $\langle M_{\eta \restriction \alpha}:\alpha \le \ell g(\eta) \rangle$ is
$\le_{\frak K}$-increasing continuous and $M_{\eta \char 94 \langle 1
\rangle}$ cannot be $\le_{\frak K}$-embedded into $M_\nu$ over $M_\eta$ when
$\eta \char 94 \langle 0 \rangle \trianglelefteq \nu \in {}^{\lambda >}2$ and
$\langle M_\eta:\eta \in {}^{\lambda >}2 \rangle$ is definable in ${\frak B}$
\sn
\item "{$(*)^d_2$}"  WDmId$^{\text{Def}}(\lambda)$ or DfWD$^+(\lambda)$ is
not $\lambda$-saturated (which holds if there is no normal $\lambda$-saturated
ideals on $\lambda$ (which holds for non Mahlo $\lambda$) and holds if for
some $\theta,\{ \delta < \lambda:\text{cf}(\delta) = \theta\}$ is not in the
ideal).
\endroster
\mn
\ub{Then} the conclusion of \scite{1.4} holds. 
\enddemo
\bn

\relax From the Definition below, we use here mainly ``superlimit"
\definition{\stag{1.5} Definition}  1)  $M \in {\frak K}_\lambda$ is a
{\it superlimit} if
\medskip
\roster
\item "{$(a)$}"  for every $N \in {\frak K}_\lambda$ such that $M 
\le_{\frak K} N$ there is $M' \in K_\lambda$ such that $N \le_{\frak K} M',
N \ne M'$, and $M \cong M'$;
\sn
\item "{$(b)$}"  if $\delta < \lambda^+$ is limit, $\langle M_i:i < \delta
\rangle$ is $\le_{\frak K}$-increasing and $M_i \cong M$ (for $i < \delta$) 
then $\dsize \bigcup_{i < \delta} M_i \cong M$.
\endroster
\medskip

\noindent
2) For $\Theta \subseteq \{ \mu:\aleph_0 \le \mu \le \lambda,\mu$ regular$\}$ 
we say $M \in {\frak K}_\lambda$ is a $(\lambda,\Theta)$-superlimit if:
\medskip
\roster
\item "{$(a)$}"  clause $(\alpha)$ from part (1) holds and
\sn
\item "{$(b)$}"  if $M_i \cong M$ is ($\le_{\frak K}$)-increasing for 
$i < \mu \in \Theta$ then $\dsize \bigcup_{i < \mu} M_i \cong M$.
\endroster
\medskip

\noindent
3) For $S \subset \lambda^+$ we say $M \in {\frak K}_\lambda$ is a
$(\lambda,S)$-strong limit if:
\medskip
\roster
\item "{$(a)$}"  clause $(\alpha)$ from part (1) holds
\sn
\item "{$(b)$}"  there is a function $F$ from 
$\dsize \bigcup_{\alpha < \kappa}
{}^\alpha(K_\lambda)$ to $K_\lambda$ such that:
{\roster
\itemitem{ $(\alpha)$ }  for any sequence $\langle M_i:
i < \alpha \rangle$ if $\alpha < \kappa,M_0 = M,M_i$ is 
$\le_{\frak K}$-increasing, and $M_i \in {\frak K}_\lambda$, \underbar{then}
$j < \alpha \Rightarrow M_j \le_{\frak K} F(\langle M_i:i < \alpha 
\rangle)$
\sn
\itemitem{ $(\beta)$ }  if $\langle M_i:i < \lambda^+ \rangle$ is 
$\le_{\frak K}$-increasing, $M_0 = M,M_i \in {\frak K}_\lambda$, and for 
\newline
$i < \kappa,M_{i+1} \le_{\frak K} F(\langle M_j:j \le i + 1 \rangle) 
\le_{\frak K} M_{i+2}$ \underbar{then} \nl
$\{\delta \in S|\dsize \bigcup_{i < \delta} M_i \ncong M\}$ is not stationary.
\endroster}
\endroster
\medskip

\noindent
4) $M$ is a $(\lambda,\kappa)$-limit if there is a function $F$ as in
3b$(\alpha)$ such that:
\medskip
\roster
\item "{$(a)$}"  if $\langle M_i:i < \kappa \rangle$ is a 
$<_{\frak K}$-increasing continuous sequence in $K_\lambda$, \newline
$F(\bar M \restriction (i+1)) \le_{\frak K}
M_{i+1}$ then $\dsize \bigcup_{i < \kappa} M_i \cong M$
\item "{$(b)$}"  there is at least one such sequence.
\endroster
\medskip

\noindent
5) $M$ is a $(\lambda,\kappa)$-superlimit is defined similarly, but with $F$ 
omitted and $M_{i+1} \cong M$.
\enddefinition
\bigskip

\proclaim{\stag{1.5A} Claim}
1)  In \scite{1.3}(1) we can replace the categoricity of 
${\frak K}$ in $\chi$ by ``${\frak K}$ has a super limit model in $\chi$" 
which is not an amalgamation base (see Definition \scite{1.5}).  In
this case the assumption of \scite{1.3}(2), and of \scite{1.4}(1) holds.
\newline
2) We can weaken (in \scite{1.4}(5) the existence of superlimit to)
``for some $\kappa = \text{ cf}(\kappa) \le \chi$ there 
is a $(\chi,\{ \kappa\})$-super limit model which is not an amalgamation
base"; provided that we add $\{ \delta < \lambda:\text{cf}(\delta) = 
\kappa\} \notin \text{ WDmId}_{< \mu}(\lambda)$ (but for \scite{2.6}(1) we
need now ``WDmId$_{< \mu}(\lambda) + \{\delta < \lambda:\text{ cf}(\delta) =
\kappa\}$ is not $\lambda$-saturated".  If there is $S \subseteq \{ \delta < 
\lambda:\text{cf}(\delta) = \kappa\}$, which belongs to $I[\lambda]$ but not
to WDmId$_{< \mu}(\lambda)$ we can weaken the model theoretic requirement 
to: there is a $(\chi,\{\kappa\})$-medium limit (see 
\cite[Definition \S3]{Sh:88}) but not used here. \nl
\endproclaim
\bigskip

\proclaim{\stag{1.6} Claim}  Assume $2^\lambda < 2^{\lambda^+}$. \newline
0) If ${\frak K}$ (an abstract elementary class) is categorical in
$\lambda$, LS$({\frak K}) \le \lambda,I(\lambda^+,{\frak K}) <
2^{\lambda^+}$, \ub{then} ${\frak K}_\lambda$ has amalgamation. \nl
1) If ${\frak K}$ (an abstract elementary class) is categorical in $\lambda$ 
and $1 \le IE(\lambda^+,K) < 2^{\lambda^+}$ but $K_{\lambda^{++}} = 
\emptyset$ \underbar{then} ${\frak K}$ has amalgamation in $\lambda$. \newline
2) Assume ${\frak K}$ has amalgamation in $\lambda,LS({\frak K}) \le \lambda,
K_{\lambda^+} \ne \emptyset$ and $K_{\lambda^{++}} = \emptyset$. 
\underbar{Then} there is $M \in K_{\lambda^+}$ saturated above $\lambda$.
\newline
3) If $M$ is $\mu$-saturated above $\lambda,LS({\frak K}) \le \lambda_0 <
\lambda$ and ${\frak K}$ has amalgamation in every $\lambda'_0 \in [\lambda_0,
\lambda)$ \underbar{then} $M$ is $\mu$-saturated above $\lambda_0$. 
\endproclaim
\bigskip

\remark{\stag{1.6B} Remark}  If $I(\lambda^+,{\frak K}) < 
2^{\lambda^+}$, then the assumption ${\frak K}_{\lambda^{++}} = \emptyset$ 
is not used in part (1) of \scite{1.6}; this is \scite{1.4}(1) + (2).  Also
if $(2^\lambda)^+ < 2^{\lambda^+}$ then the assumption $K_{\lambda^{++}} =
\emptyset$ is not needed in part (1) of \scite{1.6}; by \scite{1.4}(1) + (2)
(note (b) of \scite{1.4}(1)).
\endremark
\bigskip

\demo{Proof}  0) By \scite{1.4}(1) applied to $\lambda^+$. \nl
1) If not, we can choose for $\eta \in {}^{\lambda^+ >}2$ a model $M_\eta \in 
{\frak K}_\lambda$ such that \newline 
$[\nu \triangleleft \eta \Rightarrow M_\nu \le_{\frak K} M_\eta]$,
and $M_{\eta \char 94 \langle 0 \rangle},M_{\eta \char 94 \langle 1 \rangle}$
cannot be amalgamated over $M_\eta$.  If \newline
$(2^\lambda)^+ < 2^{\lambda^+}$
we are done by \scite{1.4}(1), so assume $(2^\lambda)^+ = 2^{\lambda^+}$.  
For each \newline
$\eta \in {}^{\lambda^+}2$ let $M_\eta =: \dsize \bigcup_{\alpha < \lambda^+}
M_{\eta \restriction \alpha}$, and let $N_\eta \in {\frak K}_{\lambda^+}$ be 
such that \newline
$M_\eta \le_{\frak K} N_\eta,N_\eta$ is $\le_{\frak K}$-maximal (exists as 
${\frak K}_{\lambda^{++}} = \emptyset$).  Now we choose by 
induction $\zeta < (2^\lambda)^+,\eta_\zeta
\in {}^{\lambda^+}2$ such that $M_{\eta_\zeta}$ is not 
$\le_{\frak K}$-embeddable
into $N_{\eta_\xi}$ for $\xi < \zeta$ (exists by \scite{1.3}(2)).  So 
necessarily for $\xi < \zeta,N_{\eta_\zeta}$ is not 
$\le_{\frak K}$-embeddable into $N_{\eta_\xi}$ (as $M_{\eta_\zeta} \le 
N_{\eta_\zeta}$).  Also, for 
$\xi < \zeta,N_{\eta_\xi}$ is not $\le_{\frak K}$-embeddable into
$N_{\eta_\zeta}$ as otherwise by the maximality of $N_\xi$ this implies 
$N_\xi \cong N_\zeta$.  So 
$\{ N_\zeta:\zeta < 2^{\lambda^+}\}$ exemplifies $IE(\lambda^+,{\frak K}) 
= 2^{\lambda^+}$, contradicting an assumption. \newline
2)  A maximal model in $K_{\lambda^+}$ will do. \nl
3)  Easy.  \hfill$\square_{\scite{1.6}}$
\enddemo
\bigskip

\noindent
\centerline {$* \qquad * \qquad *$}
\bigskip

\noindent
\underbar{\stag{1.6C} Discussion} 
Instead of Weak Diamond we now discuss Definable Weak Diamond, which 
is weaker but suffices. \newline
Compare \cite{MkSh:313}, where many
Cohen subsets are added to $\lambda$ and a combinatorial principle about
amalgamation of configurations $\langle M_s:s \subseteq n,s \ne n \rangle$ is
gotten.

We are interested here in the case $n=1$ (ordinary amalgamation); in \S3,
also $n=2$.  Even more definability can be required.

This is particularly interesting when we look at results under some other
set theory, when combining $2^\lambda = 2^{\lambda^+}$ with definable weak
diamond on $\lambda^+$ is helpful.  This played a major role in the 
preliminary form of this work. 
\bigskip

\definition{\stag{1.7} Definition}  1) In Definition \scite{1.1} we 
add the superscript ${\Cal F}$ if we restrict ourselves to functions
$F \in {\Cal F}$. \newline
2) Fix a model ${\frak B}$ whose universe includes $\lambda$ and has a 
definable pairing function on $\lambda$, and a logic ${\Cal L}$ 
closed under first order operations and substitution; also allow 
``$M \in {\frak K}$" and ``$M \le_{\frak K} N$" in the formulas, if it is
not said otherwise.  Let

$$
\align
{\Cal F}^{\text{Def}}_{{\frak B},{\Cal L}} = \biggl\{ F:&\text{ for some }
g \in {}^\lambda \lambda \text{ and } \bar h = \langle h_\eta:\eta \in
{}^{\lambda >}2 \rangle \text{ where} \\
  &\,h_\eta:\ell g(\eta) \rightarrow \lambda \text{ and }
h_{\eta \restriction \beta} \subseteq h_\eta \text{ for } \beta <
\ell g(\eta) \\
  &\text{ and for some sequence }
\bar \psi = \langle \psi_\alpha:\alpha < \lambda \rangle, \text{ with} \\
  &\,\psi_\alpha \in {\Cal L} \text{ for } \alpha < \lambda \text{ the
following holds for every } \alpha < \lambda \text{ and} \\
  &\,f \in {}^\alpha(2^{< \lambda}):F(f) = 1 \text{ iff } 
({\frak B},\alpha,g,h_f) \models \psi_\alpha \text{ iff } F(f) \ne 0 \biggr\}.
\endalign
$$
\medskip

\noindent
3) The version of weak diamond from \scite{1.1}, restricted to the class 
${\Cal F}$ of \scite{1.7}(2) is called the 
$({\frak B},{\Cal L})$-{\it definitional} version.  If ${\Cal L}$
is ${\Cal L}_{\lambda,\lambda}$ we may omit it.  If ${\frak B}$ has
the form $({\Cal H}(\chi),\in,<^*_\chi,\lambda)$ we write
${\Cal F}^{\text{Def}[(\chi)]}_{\Cal L}$ or ${\Cal F}^{\text{Def}(\chi)}$
instead of ${\Cal F}^{\text{Def}}_{{\frak B},{\Cal L}}$ or
${\Cal F}^{\text{Def}}_{\frak B}$ respectively.  If we omit $\chi$, we mean
$\chi = (2^\lambda)^+$ and we may put Def$(\chi)$ or Def instead of
${\Cal F}^{\text{Def}(\chi)}$ or ${\Cal F}^{\text{Def}}$ in the superscript.
Having the definitional version or the definable weak diamond for $\lambda$
means $\lambda \notin \text{ WDmId}^{\text{Def}}(\lambda)$. \newline
4) Let DfWD$_{< \mu}(\lambda)$ mean that with ${\frak B} = ({\Cal H}(\chi),
\in,<^*_\chi,\lambda,\mu)$ we have $\lambda \notin \text{ WDmId}^{\text{Def}}
_{< \lambda}(\lambda)$.  Instead of ``$< \mu^+$" we write ``$\mu$" and
instead of ``$< 2^{{< \lambda}^+}$" we may write nothing. \newline
5) Let DfWD$^+_{< \mu}(\lambda)$ mean DfWD$_{< \mu}(\lambda)$ together with
the principle $\otimes_\lambda$ below; we adopt the same conventions as
in (4) concerning $\mu$:
\medskip
\roster
\item "{$\bigotimes_\lambda$}"  if for $\eta \in {}^{\lambda >}2,M_\eta$ is
a $\tau_{\frak K}$-model of cardinality $< \lambda,\langle 
M_{\eta \restriction \alpha}:\alpha \le \ell g(\eta) \rangle$ is \nl
$\le$-increasing continuous, for $\eta \in {}^\lambda 2$ we let $M_\eta =
\dsize \bigcup_{\alpha < \lambda} M_{\eta \restriction \alpha}$ and for \nl
$\eta \ne \nu \in {}^\lambda 2,M_\eta$ and $M_\nu$ are not isomorphic over
$M_{\langle \rangle}$, \ub{then} $\{M_\eta / \cong:\eta \in {}^\lambda 2\}$ 
has cardinality $2^\lambda$ \newline
(note that $2^{< \lambda} < 2^\lambda$ implies that).
\endroster
\enddefinition
\bigskip

\proclaim{\stag{1.8} Claim} \nl
1) Assume ${\Cal L}$ first order or at least definable enriching first order.

In the definition of ${\Cal F}^{\text{Def}}_{{\frak B},{\Cal L}}$, we can
replace ``for every $\alpha < \lambda$" by ``for a club of $\alpha <
\lambda$".  In the definition of
${\Cal F}^{\text{Def}}$ we can let $g \in {}^\lambda({}^{\lambda >} 2)$
and $h_\eta:\ell g(\eta) \rightarrow {}^{\lambda >}2$.  In any case 
WDmId$^{\Cal F}_{< \mu}(\lambda)$ increases with ${\Cal F}$ and is 
$\subseteq$ WDmId$_{< \mu}(\lambda)$,
similarly for WDmTId$^{\Cal F}_{< \mu}(\lambda)$. \newline
2)  If ${\Cal F} = {\Cal F}^{\text{Def}}_{\frak B}$ and cf$(\mu)
> \lambda$ \underbar{then} WDmId$_{< \mu}(\lambda)$ is a normal ideal
(but possibly is equal to ${\Cal P}(\lambda)$). \newline
3) Assume $V \models ``\lambda = \chi^+,\chi^{< \chi} =
\chi,\mu > \lambda"$ and $P$ is the forcing notion of adding $\mu$ Cohen
subsets to $\chi$ (i.e. \newline
$\{g:g \text{ a partial function from } \mu \text{ to } \{0\} \text{ with
domain of cardinality } <\chi\}$). \newline
\underbar{Then} in $V^P$ we have WDmId$^{\text{Def}}_{< \mu}(\lambda)$ is 
the ideal of non-stationary subsets of $\lambda$; i.e. with \newline
${\frak B} = ({\Cal H}(\chi),\in,<^*_\mu)^{V^P}$ for \underbar{any} $\chi$.
Also $\otimes_\lambda$ of Definition \scite{1.7}(5) holds.
\endproclaim
\bigskip

\remark{Remark}  In \scite{1.8}(1) we use the assumption on ${\Cal L}$; anyhow
not serious: reread the definition \scite{1.1}(1).

\endremark

\demo{Proof}  1); 2)  By manipulating the $h$'s (using the pairing function
on $\lambda$). \newline
3) See \cite{MkSh:313} or think (the point being that we can break the
forcing, first adding $\bar \psi$ and $g$ (or the $< \mu$ ones) and then
(read \scite{1.1}(1)) choose $\eta \in {}^\lambda 2$ as \newline
$\underset\sim {}\to g \restriction [\gamma,\gamma + \lambda)$ not ``used
before".  Now for any candidate $f \in {}^\lambda({}^{\lambda >} 2)$ for a
club of $\delta < \lambda,\eta(\delta) = g(\gamma + \delta)$ is not used
in the definition of $f \restriction \delta,h_f \restriction \delta$ so
stationarily often $\eta(\delta)$ ``guesses" rightly.) 
\hfill$\square_{\scite{1.8}}$
\enddemo
\bigskip

\proclaim{\stag{1.8A} Claim}  1) If $2^\theta = 2^{< \lambda} < 2^\lambda$
then DfWD$^+(\lambda)$. \newline
2) DfWD$_{< \mu}(\lambda)$ holds when $\lambda \notin \text{ WDmId}_{< \mu}
(\lambda)$ (see \scite{1.2} for sufficient conditions).
\endproclaim
\bigskip

\noindent
\underbar{\stag{1.8B} Discussion}:  
We hope to get successful ``guessing" not just on a stationary set, but on a
positive set for the same ideal for which have guessed; i.e. 
there is $I$ a normal ideal on $\lambda$ such that for $A \in I^+$ 
there is $\eta \in {}^\zeta 2$ guessing $I$-positively; this is 
connected to questions on $\lambda^+$-saturation.  For more see \cite{Sh:638}.

We phrased the following notion originally in the hope of later eliminating 
$\mu_{\text{wd}}(\lambda)$ (i.e. using $2^\lambda$ instead of
$\mu_{\text{wd}}(\lambda)$).
\bigskip

\definition{\stag{1.9} Definition}  1)

$$
\align
\text{UDmId}^{\Cal F}_{< \mu}(\lambda) = \biggl\{ S \subseteq \lambda
:&\text{ for some } i^* < \mu \text{ and } F_i \in {\Cal F} \text{ (for }
i < i^*) \\
  &\text{ for every } \eta \in {}^S 2 \text{ there is } f \in 
{}^\lambda(2^{< \lambda}) \text{ and } i < i^* \\
  &\text{ and club } E \text{ of } \lambda \text{ such that:} \\
  &\text{ for every } \delta \in E \text{ we have:} \\
  &\,\delta \in S \Rightarrow \eta(\delta) = F_i(f \restriction \delta) \\
  &\,\delta \in \lambda \backslash S \Rightarrow 0 = F_i(f \restriction 
\delta) \biggr\}.
\endalign
$$
\medskip

\noindent
2) We omit $\mu$ if $\mu = 1$. \newline
3) $BA^{\Cal F}{(\lambda)}$ is defined as the family of $S \subseteq \lambda$
such that for some $F \in {\Cal F}$ and $\eta = O_S$ the condition above
holds.
\enddefinition
\bigskip

\proclaim{\stag{1.10} Claim}  Assume ${\Cal F} = 
{\Cal F}^{\text{Def}}_{\frak B}$. \newline
1) In the definition \scite{1.9}(1) we can replace 
$f \in {}^\lambda(2^{< \lambda})$ by
$f \in {}^\lambda 2$ or $f \in {}^\lambda({}^{\lambda >} 2)$. \newline
2) UDmId$^{\Cal F}(\lambda)$ is a normal ideal on $\lambda$ (but possibly
is ${\Cal P}(\lambda)$). \newline
3) $BA^{\Cal F}(\lambda)$ is a Boolean algebra of subsets of $\lambda$
including all non-stationary subsets of $\lambda$ and even UDmId$^{\Cal F}
(\lambda)$, and is closed under unions of $< \lambda$ sets and even under
diagonal union. \newline
4) If $S \in BA^{\Cal F}(\lambda)$ and $F \in {\Cal F}$ then for some
$\eta \in {}^S 2$ we have:
\medskip
\roster
\item "{$(*)$}"  for every $f \in {}^\lambda({}^{\lambda >}2)$ we have
$$
\{ \delta \in S:\eta(\delta) = F(f \restriction \delta)\} \ne \emptyset
\text{ mod UDmId}^{\Cal F}(\lambda).
$$
\endroster
\medskip

\noindent
5) UDmId$^{\Cal F}_{< \mu}(\lambda) \subseteq$ WDmId$^{\Cal F}_{< \mu}
(\lambda) \subseteq$ WDmId$_{< \mu}(\lambda)$ and they increase with
${\Cal F}$ and $\lambda \in \text{ UDmId}^{\Cal F}_{< \mu}(\lambda)
\Leftrightarrow \lambda \in \text{ WDmId}^{\Cal F}_{< \mu}(\lambda)$. 
\endproclaim
\bigskip

\demo{Proof} Straightforward.  \hfill$\square_{\scite{1.10}}$
\enddemo
\bigskip

\noindent
\underbar{\stag{1.12} Discussion}  Remember
\medskip
\roster
\item "{$(*)_1$}"  If $V = ``\chi = \chi^{< \chi} \and 2^\chi = \chi^+"$,
$P$ is the forcing notion of adding $\mu > \chi^+$ Cohen subsets to $\chi$ 
\underbar{then} in $V^\chi$, any equivalence relation on ${\Cal P}(\lambda)$ 
definable with parameters $X \subseteq \chi$ and ordinals which has at least 
$\chi^{++}$ equivalence classes has at least $\mu$ equivalence classes
\ermn
and (see \cite[XVI,\S2]{Sh:f})
\mr
\item "{$(*)_2$}"  ZFC is consistent with $CH+$ for some stationary, 
costationary $S \subseteq \omega_1$ we have
\medskip
{\roster
\itemitem{ $(a)$ }  WDmId$(\aleph_1) = \{A \subseteq \omega_1:A \backslash S$ 
is not stationary$\}$
\sn
\itemitem{ $(b)$ }  ${\Cal D}_{\omega_1} + S$ is $\aleph_2$-saturated.
\endroster}
\item "{$(*)_3$}"  $ZFC + GCH$ is consistent with $\{ \delta < \aleph_2:
\text{cf}(\delta) = \aleph_1\} \in \text{ WDmId}(\aleph_2)$
\sn
\item "{$(*)_4$}"  $ZFC + 2^{\aleph_1} < 2^{\aleph_2}$ is consistent with
$\{ \delta < \aleph_2:\text{cf}(\delta) = \aleph_0\} \in \text{ WDmId}
(\aleph_2)$.
\endroster
\newpage

\head {\S2 First attempts} \endhead  \resetall
\bigskip

Given amalgamation in $K_\lambda$ (cf. 2.2) we try to define and analyze 
types $p \in {\Cal S}(M)$ for $M \in K_\lambda$.  But types here 
(as in \cite{Sh:300}) are not sets of formulas.  They may instead be
represented by triples $(M,N,a)$ with $M \le_{\frak K} N$ and 
$a \in N \backslash M$.  We look for nice types (i.e. triples) and try 
to prove mainly the density of the set of minimal types).

To simplify matters we allow uses of stronger assumptions than are
ultimately desired
(e.g. $2^{\lambda^+} > \lambda^{++}$ and/or $K_{\lambda^{+3}} =
\emptyset$).  These will later be eliminated.  However the first extra
assumption is still a ``mild set theoretic assumption", and the second is
harmless if we think only of proving our main theorem \scite{0.B} and not 
on subsequent continuations.

So the aim of this section is to show that we can start to analyze such 
classes.
\bigskip

\demo{\stag{2.1} Hypothesis}  ${\frak K}$ is an abstract elementary class.
\enddemo
\bigskip

\proclaim{\stag{2.2} Claim}  Assume
\medskip
\roster
\item "{$(*)^2_\lambda$}"  $K$ is categorical in $\lambda;1 \le 
I(\lambda^+,K) < 2^{\lambda^+};LS({\frak K}) \le \lambda$ and:
$2^\lambda < 2^{\lambda^+}$, or at least the definable weak diamond holds for
$\lambda^+$ holds.
\endroster
\medskip

\noindent
\underbar{Then} \newline
1)  ${\frak K}_\lambda$ has amalgamation. \newline
2)  If $I(\lambda^{++},K) = 0$ \underbar{then} ${\frak K}$ has a 
model in $\lambda^+$ which is universal homogeneous above $\lambda$, hence
saturated above $\lambda$ (see \scite{0.18}(2)). \newline
3) If $I(\lambda^{++},K) = 0$ then $M \in K_\lambda \Rightarrow 
|{\Cal S}(M)| \le \lambda^+$.
\endproclaim
\bigskip

\demo{Proof}  1)  If amalgamation fails in $K_\lambda$ and $2^\lambda <
2^{\lambda^+}$, \underbar{then} the assumptions of \scite{1.3}(1) hold with
$\lambda^+$ in place of $\lambda$.  Hence by \scite{1.4}(2) the statement
$(*)_1$ of \scite{1.4}(1) (see there) holds and easily also $(*)_2$ of
\scite{1.4}(1), hence by \scite{1.4}(4) we have $I(\lambda^+,K) = 
2^{\lambda^+}$, a contradiction.
If $2^\lambda = 2^{\lambda^+}$, we are using the variants from \scite{1.7}.
\newline
2)  As $I(\lambda^{++},K) = 0 < I(\lambda^+,K)$, there is $M \in
K_{\lambda^+}$ which is maximal.  If $M$ is not universal
homogeneous above $\lambda$ then there are $N_0,N_1 \in K_\lambda$ with
$N_0 \le_{\frak K} M$ and $N_0 \le_{\frak K} N_1$ such that $N_1$ cannot 
be $\le_{\frak K}$-embedded into $M$ over $N_0$.
Use \scite{0.17}(2) to get a contradiction. \newline
3) Follows from (2).  \hfill$\square_{\scite{2.2}}$
\enddemo
\bigskip

\definition{\stag{2.3A} Definition}  \newline
1)(a) $K^3_\lambda = \{(M_0,M_1,a):M_0 \le_{\frak K} M_1 \text{ are both in } 
K_\lambda \text{ and } a \in M_1 \backslash M_0\}$.

$$
(M_0,M_1,a) \le (M'_0,M'_1,a') \text{ \underbar{if} } 
a = a',M_0 \le_{\frak K} M'_0,M_1 \le_{\frak K} M'_1 \tag"{$(b)$}"
$$

$$
\align
(M_0,M_1,a) \le_h &(M'_0,M'_1,a') \text{ \underbar{if} } 
h(a) = a', \text{ and for } \ell=0,1 \text{ we have:} \tag "{$(c)$}" \\
  &h \restriction M_\ell
\text{ is a } \le_{\frak K} \text{-embedding of } 
M_\ell \text{ into } M'_\ell.
\endalign
$$

$$
\align
(M_0,M_1,a) < &(M'_0,M'_1,a') \text{ \underbar{if} } 
(M_0,M_1,a) \tag"{$(d)$}" \\
  &\le (M'_0,M'_1,a) \text{ and } M_0 \ne M'_0.
\endalign
$$

$$
\text{ similarly } <_h \tag"{$(e)$}"
$$
\medskip

\noindent
2)  $(M_0,M_1,a) \in K^3_\lambda$ has the weak extension property \ub{if}
there is $(M'_0,M'_1,a) \in K^3_\lambda$ such that $(M_0,M_1,a) \le 
(M'_0,M'_1,a)$ and $M_0 \ne M_1$. \newline
3) $(M_0,M_1,a) \in K^3_\lambda$ has the extension property \underbar{if}: for
every $N_0 \in {\frak K}_\lambda$ and \newline
$\le_{\frak K}$-embedding $f$ of $M_0$ into $N_0$ there are $N_1,b$ and $g$ 
such that: $(M_0,M_1,a) \le_g (N_0,N_1,b) \in K^3_\lambda$ and $g \supseteq f$
(so $g(a) = b$ and $g$ is a $\le_{\frak K}$-embedding of $M_1$ into $N_1$).
\enddefinition
\bigskip

\proclaim{\stag{2.3} Claim}  Assume
\medskip
\roster
\item "{$(*)^3_\lambda$}"  $LS({\frak K}) \le \lambda,K$ is categorical in 
$\lambda$ and in $\lambda^+$, and $1 \le I(\lambda^{++},K)$.
\endroster
\medskip

\noindent
\underbar{Then} every $(M_0,M_1,a) \in K^3_\lambda$ has the weak extension
property that is: \nl
If $M_0 \le_{\frak K} M_1$ are in $K_\lambda$ and $a \in M_1 \backslash M_0$,
then we can find $M'_1$ in ${\frak K}_\lambda$ such that: $M_0 <_{\frak K}
M_1$ hence $M_0 \ne M'_0$ and $(M_0,M_1,a) \le (M'_0,M'_1,a)$.
\endproclaim
\bigskip

\demo{Proof}   We can choose $\langle N_i,a_i:i < \lambda^+ \rangle$ 
such that:
\medskip
\roster
\item "{$(a)$}"  $N_i \in K_\lambda \text{ is } \le_{\frak K} 
\text{-increasing continuous in }i$;
\sn
\item "{$(b)$}"  $h_i \text{ is an isomorphism from } M_1 \text{ onto } 
N_{i+1} \text{ such that}$ \newline
$h_i(M_0) = N_i,h_i(a) = a_i$.
\endroster
\medskip

\noindent
Now as $a \in M_1 \backslash M_0$ clearly $i < j < \lambda^+ \Rightarrow 
a_i \in N_{i+1} \le_{\frak K} N_j \and a_j \notin N_{i+1}$ hence 
$\dsize \bigcup_{i < \lambda^+} N_i \in K_{\lambda^+}$.

By \scite{0.17}(1) applied to $\lambda^+$ there are $M'_0 \le_{\frak K} M'_1$ 
in $K_{\lambda^+},M'_0 \ne M'_1$, and there is $b \in M'_1 \backslash M'_0$.  
As $K$ is categorical in $\lambda^+$,
without loss of generality $M'_0 = \dsize \bigcup_{i < \lambda^+} N'_i$. \nl
Let $\chi$ be large enough and ${\frak B} \prec ({\Cal H}(\chi) \in,<^*_\chi)$
be such that $\lambda \subseteq {\frak B},\|{\frak B}\| = \lambda$ and 
$\langle N_i,a_i:i < \lambda^+ \rangle,M'_0,M'_1,b$ and 
the definition of ${\frak K}$ belong to ${\frak B}$.
\newline
Let $\delta = {\frak B} \cap \lambda^+$, so $\delta \in (\lambda,\lambda^+)$
is a limit ordinal and

$$
N_\delta \le_{\frak K} N_{\delta + 1} \le_{\frak K} M'_1,
$$

$$
N_\delta \le_{\frak K} (M'_1 \cap {\frak B}) \le_{\frak K} M'_1
$$

$$
{\frak B} \cap M'_0 = N_\delta
$$

$$
N_{\delta + 1} \cap (M'_1 \cap {\frak B}) =  N_\delta,
$$
\medskip

\noindent
so for some $N$ we have:
$$
N \in {\frak K}_\lambda,N \le_{\frak K} M'_1, \text{ and }
(N_{\delta+1} \cup (M'_1 \cap {\frak B})) \subseteq N
$$
\medskip

\noindent
so (see Definition \scite{2.3A}(1) above)

$$
(N_\delta,N_{\delta + 1},a_\delta) \le_{\frak K} 
(M'_1 \cap {\frak B},N,a_\delta),
$$
\medskip

\noindent
and $b$ witnesses $N_\delta \ne M'_1 \cap {\frak B}$. \newline
As $(M_0,M_1,a) \cong (N_\delta,N_{\delta +1},a_\delta)$, the result
follows).  \hfill$\square_{\scite{2.3}}$
\enddemo 
\bigskip

\definition{\stag{2.3B} Definition}  1) $(M_0,M_1,a) \in K^3_\lambda$ is 
\ub{minimal} when:

$$
\text{if }(M_0,M_1,a) \le_{h_\ell} (M'_0,M^\ell_1,a_\ell) \in K^3_\lambda
\text{ for } \ell = 1,2,
$$

$$
\text{\underbar{then} } \text{tp}(a_1,M'_0,M^1_1) = \text{ tp}(a_2,M'_0,
M^2_1).
$$
\mn
2)  $(M_0,M_1,a) \in K^3_\lambda$ is \ub{reduced} when:

$$
\text{if }(M_0,M_1,a) \le (M'_0,M'_1,a) \in {\frak K}^3_\lambda \text{ then }
M'_0 \cap M_1 = M_0.
$$
\medskip

\noindent
3) We say $p \in {\Cal S}(M_0)$ is \ub{minimal}, where 
$M_0 \in K_\lambda$, if for some $a,M_1$ we have: 
$p = \text{ tp}(a,M_0,M_1)$ and $(M_0,M_1,a) \in
K^3_\lambda$ is minimal. \newline
4) We say $p \in {\Cal S}(M_0)$ is \ub{reduced} where $M_0 \in K_\lambda$, if 
for some $a,M_1$ we have \newline
$p = \text{ tp}(a,M_0,M_1)$ and $(M_0,M_1,a) \in K^3_\lambda$ is reduced.
\enddefinition
\bigskip

\demo{\stag{2.4} Fact}  1) For every $(M_0,M_1,a) \in K^3_\lambda$ there is a
reduced $(M'_0,M'_1,a)$ such that: $(M_0,M_1,a) \le (M'_0,M'_1,a) \in 
K^3_\lambda$. \newline
2)  Assume $\langle (M_{0,\alpha},M_{1,\alpha},a):\alpha < \delta \rangle$
is an increasing sequence of members of $K^3_\lambda$
\mr
\item "{$(a)$}"  if $\delta < \lambda^+$ \underbar{then}
$(M_{0,\alpha},M_{1,\alpha},a) \le (\dsize \bigcup_{\beta < \delta} 
M_{0,\beta}, \dsize \bigcup_{\beta < \delta} M_{1,\beta},a) \in K^3_\lambda$.
\sn
\item "{$(b)$}"  If $\delta = \lambda^+$ the
result may be in $K^3_{\lambda^+}$: if $\{ \alpha < \delta:M_{0,\alpha} \ne
M_{0,\alpha + 1}\}$ is cofinal, this holds. 
\sn
\item "{$(c)$}"  If $\delta < \lambda^+$ and each $(M_{0,\alpha},
M_{1,\alpha},a)$ is reduced \underbar{then} so is \newline
$(\dsize \bigcup_{\beta < \delta} M_{0,\beta}, \dsize \bigcup_{\beta < \delta}
M_{1,\beta},a)$.
\ermn
3) If $(M_0,M_1,a) \le (M'_0,M'_1,a)$ are in $K^3_\lambda$ and the 
first triple is minimal \underbar{then} so is the second. \newline
4) If $(M_0,M_1,a) \le (M'_0,M'_1,a)$ are in $K^3_\lambda$ \underbar{then}
$\text{tp}(a,M_0,M_1) \le \text{ tp}(a,M'_0,M'_1)$; 
(Definition \scite{0.14}(1)). \newline
5) If $K_\lambda$ has amalgamation, then: $(M_0,M_1,a) \in K^3_\lambda$ is
minimal if and only if:
\medskip
\roster
\item "{$(*)$}"  If $(M_0,M_1,a) \le_{h_\ell} (M'_0,M'_1,a_\ell) \in
K^3_\lambda$ for $\ell = 1,2$ and $h_1 \restriction M_0 = h_2 \restriction
M_0$ \underbar{then} tp$(a_1,M'_0,M'_1) = \text{ tp}(a_2,M'_0,M'_1)$.
\endroster
\medskip

\noindent
6) If there is no maximal member \footnote{will be applied for $\lambda^+$}
of $K^3_\lambda$ and there are $N_0 <_{\frak K} N_1$ in $K_\lambda$,
\ub{then} there are $N^0 <_{\frak K} N^1$ in ${\frak K}_{\lambda^+}$. \newline
7) If every triple in $K^3_\lambda$ has the weak extension property, and there
are $N_0 <_{\frak K} N_1$ in ${\frak K}_\lambda$ \underbar{then} there are
$N^0 <_{\frak K} N^1$ in ${\frak K}_{\lambda^+}$. \newline
8) If LS$({\frak K}) \le \lambda$ and every triple in $K^3_\lambda$ has the
extension property and $K^3_\lambda \ne \emptyset$ \underbar{then} no 
$M \in K_{\lambda^+}$ is
$<_{\frak K}$-maximal hence $K_{\lambda^{++}} \ne \emptyset$. \newline
9) If $LS({\frak K}) \le \lambda$ and $K_{\lambda^+} \ne \emptyset$,
\underbar{then} $K^3_\lambda \ne \emptyset$.
\enddemo
\bigskip

\demo{Proof}  Easy.  Note that part (7) is \scite{2.3}.
\enddemo
\bigskip

\proclaim{\stag{2.5} Claim}  Assume $(*)^2_{\lambda^+} + (*)^3_\lambda$
(i.e. the hypothesis of \scite{2.2} and \scite{2.3}) and $2^{\lambda^+} > 
\lambda^{++}$, and $K_{\lambda^{+3}} = \emptyset$. \newline

\ub{Then} in ${\frak K}^3_\lambda$ the minimal triples are dense 
(i.e. above every triple in $K^3_\lambda$ there is a minimal one).
\endproclaim
\bigskip

\remark{Remark}  We do not intend to adopt the hypotheses 
``$2^{\lambda^+} > \lambda^{++}",K_{\lambda^{+3}} = \emptyset$ indefinitely.
They will be eliminated in \S3.
\endremark
\bigskip

\demo{Proof}  If not, we can choose by induction on $\alpha < \lambda^+$,
for $\eta \in {}^\alpha 2$ a triple $(M^0_\eta,M^1_\eta,a_\eta)$ and
$h_{\eta,\nu}$ for $\nu \trianglelefteq \eta$ such that:
\medskip
\roster
\widestnumber\item{ (iii) }
\item "{$(i)$}"  $(M^0_\eta,M^1_\eta,a_\eta) \in K^3_\lambda$
\sn
\item "{$(ii)$}"  $\nu \triangleleft \eta \Rightarrow (M^0_\nu,M^1_\nu,
a_\nu) \le_{h_{\eta,\nu}} (M^0_\eta,M^1_\eta,a_\eta)$
\sn
\item "{$(iii)$}"  $\nu_0 \triangleleft \nu_1 \triangleleft \nu_2
\Rightarrow h_{\nu_2,\nu_0} = h_{\nu_2,\nu_1} \circ h_{\nu_1,\nu_0}$
\sn
\item "{$(iv)$}"  $(M^0_{\eta \char 94 \langle \ell \rangle},M^1
_{\eta \char 94 \langle \ell \rangle},h_{\eta \char 94 \langle \ell \rangle,
\eta} \restriction M^0_\eta)$ for $\ell = 0,1$ are equal
\sn
\item "{$(v)$}"  $\text{tp}(a_{\eta \char 94 \langle 0 \rangle},
M^0_{\eta \char 94 \langle 0 \rangle},M^1_{\eta \char 94 \langle 0 
\rangle}) \ne \text{tp}(a_{\eta \char 94 \langle 1 \rangle},
M^0_{\eta \char 94 \langle 1 \rangle},M^1_{\eta \char 94 \langle 1 
\rangle})$; \newline
this makes sense as $M^0_{\eta \char 94 \langle 0 \rangle} = 
M^0_{\eta \char 94 \langle 1 \rangle}$
\sn
\item "{$(vi)$}"  if $\eta \in {}^\delta 2$ and 
$\delta < \lambda^+$ is a limit ordinal, \underbar{then}
$M^\ell_\eta = \dsize \bigcup_{\alpha < \delta} h_{\eta,\eta \restriction
\alpha} (M^\ell_{\eta \restriction \alpha})$ \newline
for $\ell = 0,1$,
\sn
\item "{$(vii)$}"  $(M^0_{<>},M^1_{<>},a_{<>}) \in K^3_\lambda$ is a triple
above which there is no minimal one.
\endroster
\medskip

\noindent
This is straightforward: for $\alpha = 0$ choose a triple in $K^3_\lambda$
above which supposedly there is no minimal triple; in limit $\alpha$ take 
limits of diagrams (chasing the $h$'s); in 
successor $\alpha$, use non-minimality and \scite{2.4}(5). 

Let $M^* \in {\frak K}_{\lambda^{++}}$ be saturated above $\lambda$
(exists by \scite{2.2}(2) so it is necessarily homogeneous universal above
$\lambda^+$, hence above $\lambda$; 
note: $\lambda$ there stands for $\lambda^+$ here).

We choose by induction on $\alpha < \lambda^+$ for $\eta \in {}^\alpha 2$, a
$\le_{\frak K}$-embedding $g_\eta$ of $M^0_\eta$ into $M^*$ such that:

$$
\nu \triangleleft \eta \Rightarrow g_\nu = g_\eta \circ h_{\eta,\nu}.
$$

$$
g_{\eta \char 94 \langle 0 \rangle} = g_{\eta \char 94 \langle 1 \rangle}.
$$
\medskip

\noindent
This is clearly possible. 
Let $N^0_\eta = M^* \restriction \text{ Rang}(g_\eta)$.  
For $\eta \in {}^{\lambda^+}2$ let $N^0_\eta = M^* \restriction
\dbcu_{\alpha < \lambda^+} \text{ Rang}(g_{\eta \restriction \alpha})$ and
let $g_\eta = \dbcu_{\alpha < \lambda^+} g_{\eta \restriction \alpha}$.  
Chasing arrows we can find for $\eta \in {}^{\lambda^+}2$ a limit
to $(\langle M^0_{\eta \restriction \alpha},M^1_{\eta \restriction \alpha},
a_{\eta \restriction \alpha} \rangle,h_{\eta \restriction \beta,
\eta \restriction \alpha}:\alpha < \beta < \lambda^+)$ say
$(M^0_\eta,M^1_\eta,a_\eta) \in K^3_{\lambda^+}$ and $h_{\eta,\nu}$ for
$\nu \triangleleft \eta$ as usual.  Let $f_\eta$ be the function from
$M^0_\eta$ into $M^*$ such that for $\alpha < \lambda^+$ we have
$f_\eta \circ h_{\eta,\eta \restriction \alpha} = g_\eta$.  So $f_\eta$ is
a $\le_{\frak K}$-embedding of $M^0_\eta$ into $M^*$.  
So we can extend $f_\eta$ to $f^+_\eta$,
a $\le_{\frak K}$-embedding of $M^1_\eta$ into $M^*$. \newline

Let $a^*_\eta = f^+_\eta(a_\eta)$ for $\eta \in {}^{\lambda^+}2$. \newline

As $2^{\lambda^+} > \lambda^{++}$ for some $\eta_0 \ne \eta_1$ we have
$a^*_{\eta_0} = a^*_{\eta_1}$.  So for some $\alpha < \lambda^+$, \newline
$\eta_0 
\restriction \alpha = \eta_1 \restriction \alpha$ but $\eta_0(\alpha) \ne
\eta_1(\alpha)$, without loss of generality $\eta_\ell(\alpha) = \ell$ and
by clause $(v)$ above we get a contradiction. \hfill$\square_{\scite{2.5}}$
\enddemo
\bigskip

\proclaim{\stag{2.5A} Claim}  1) Assume $(*)^2_\lambda$ or just
\mr
\item "{$(*)^{2^-}_\lambda$}"  ${\frak K}$ has amalgamation in $\lambda$ and
$LS({\frak K}) \le \lambda$.
\ermn
If $M_0 \le_{\frak K} N_0 \in K_{\lambda^+}$ and $(M_0,M_1,a) \in 
K^3_\lambda$ \underbar{then} there is $N \in K_{\le \lambda^+}$ such that: 
$N_0 \le_{\frak K} N$ and for every
$c \in N$ satisfying tp$(c,M_0,N) = \text{ tp}(a,M_0,M_1)$, there is a
$\le_{\frak K}$-embedding $h$ of $M_1$ into $N$ extending id$_{M_0}$ such that
$h(a) = c$ and $N \notin K_{\lambda^+} \Rightarrow N$ is a 
$<_{\frak K}$-maximal member of $K_\lambda$. \nl
2) Assume $M_0 \le_{\frak K} N_0 \in K_\lambda$ and $(M_0,M_1,a) \in
K^3_\lambda$ has the weak extension property. \ub{Then} there is $N \in
K_{\lambda^+}$ such that: $N_0 \le_{\frak K} N$ and for every $c \in N'$
either for some $N' \in K_\lambda$ we have $N_0 \cup \{c\} \subseteq
N \le_{\frak K} N$ and $c$ does not strongly realize tp$(a,M_0,M_1)$ \ub{or}
there is an $\le_{\frak K}$-embedding $h$ of $M_1$ into $N$ extending
id$_{M_0}$ such that $h(a) = c$.
\endproclaim
\bigskip

\demo{Proof}  1) We choose by induction on $\alpha < \lambda^+$, a model
$N_\alpha \in K_\lambda$ increasing (by $\le_{\frak K}$) continuous such
that: for $\alpha$ even $N_\alpha \ne N_{\alpha +1}$ if $N_\alpha$ is not
$\le_{\frak K}$-maximal, and for $\alpha$ odd let $\beta_\alpha = \text{ Min}
\{\beta:\beta = \alpha +1 \text{ or } \beta \le \alpha$ and there is
$c \in N_\beta$ such that there is no $\le_{\frak K}$-embedding $h$ 
of $M_1$ into
$N_\alpha$ extending id$_{M_0}$ such that $h(a)=c$ but for some $N \in
K_\lambda,N_\alpha \le_{\frak K} N$ and there is a $\le_{\frak K}$-embedding
$h$ of $M_1$ into $N$ extending id$_{M_0}$ such that $h(a)=c\}$, and if
$\beta_\alpha \le \alpha$ then 
choose $N$ exemplifying this and let $N_{\alpha +1} = N$.  
By the definition of type we are done. \nl
2) Same proof, note that the non $\le_{\frak K}$-maximality of $N_\alpha$
(and hence $N$) follows by a weak extension property.
\enddemo
\bigskip

\proclaim{\stag{2.6} Claim}  Assume $(*)^2_\lambda$ or just:
\medskip
\roster
\item "{$(*)^{2^-}_\lambda$}"  ${\frak K}$ has amalgamation for $\lambda$ 
and LS$({\frak K}) \le \lambda$.
\endroster
\medskip

\noindent
1)  Assume that above $(M_0,M_1,a) \in K^3_\lambda$ there is no minimal 
member of $K^3_\lambda$ \underbar{then} $(M_0,M_1,a)$ itself has the 
extension property.
\newline
2) If $(M_0,M_1,a) \in K^3_\lambda,M_0 \le_{\frak K} N \in K$ and the 
number of $c \in N$ such that tp$(c,M_0,N) = \text{ tp}(a,M_0,M_1)$ 
is $> \lambda$ \underbar{then} $(M_0,M_1,a)$ has the extension property. 
\newline
3) Assume above $(M_0,M_1,a) \in K^3_\lambda$ there is no minimal member of
$K^3_\lambda$ \underbar{then}
\medskip
\roster
\item"{$(*)$}"   for some $N$ we have: $M_0 \le_{\frak K} N$ and $N$ 
is as required in part (2).
\endroster
\endproclaim
\bigskip

\remark{\stag{2.6A} Remark}  1) See \scite{3.16A}. \newline
2) Note that $(*)^2_\lambda$ is from \scite{2.2} and $(*)^2_\lambda
\Rightarrow (*)^{2^-}_\lambda$ by \scite{2.2}.
\endremark
\bigskip

\demo{Proof}  1) Follows by part (2) and (3) (part (b)). \newline
2) Without loss of generality $N$ has cardinality $\lambda^+$ and also
is as in \scite{2.5A}. \newline
By \scite{0.17}(2) for any $M'_0$ such that $M_0 \le_{\frak K} M'_0 \in 
K_\lambda$ there is $N_1,N \le_{\frak K} N_1 \in K_{\lambda^+}$ and a 
$\le_{\frak K}$-embedding $h$ of $M'_0$ into $N_1$ extending id$_{M_0}$ with
$h_c(a) = c$.  
Now some $c \in N \backslash h(M'_0)$ realizes tp$(a,M_0,M_1)$ and there is 
an embedding $h_c$ of $M_1$ into $N$ extending id$_{M_0}$ such that $h_c(a) =
c$.  Lastly let $N'_1 \le_{\frak K} N_1$ be of cardinality $\lambda$ and 
include Rang$(h_c) \cup
\text{ Rang}(h)$ (send $a$ to $c$ via $h_c$). \newline
3) We first prove
\mr
\item "{$(*)_0$}"  For some $M^+_0,M_0 \le_{\frak K} M^+_0 \in 
K_\lambda$ and tp$(a,M_0,M_1)$ has $> \lambda$ extensions in 
${\Cal S}(M^+_0)$ (in fact $\ge \text{ min}\{2^\mu:2^\mu > \lambda\}$).
\endroster
\enddemo
\bigskip

\demo{Proof of $(a)$}  Let $M^0_\eta,M^1_\eta,a_\eta,h_{\eta,\nu}\,
(\eta \in {}^{\lambda^+ >}2$ and $\nu \trianglelefteq \eta)$ be as 
in the proof of \scite{2.5} (i.e. satisfy $(i)-(vi)$ there) and $M^0_{<>} =
M^0,M^1_{<>} = M^1$.  
Let $\mu = \text{ Min}\{ \mu:2^\mu > \lambda\}$, so 
${}^{\mu >} 2$ has cardinality $\le \lambda$ and $\mu \le \lambda$.  Let 
${}^{\mu >} 2 =
\{ \eta_\zeta:\zeta < \zeta^*\}$ be such that $\eta_\xi \triangleleft
\eta_\zeta \Rightarrow \xi < \zeta$ and so $\zeta^* < \lambda^+$ and without
loss of generality is a limit ordinal.  Now we can choose by induction on 
$\zeta \le \zeta^*$ a model $M^*_\zeta \in K_\lambda$ and, if $\zeta <
\zeta^*$ also a function $g_{\eta_\zeta}$ such that:
\medskip
\roster
\item "{$(\alpha)$}"  $M^*_\zeta$ is $\le_{\frak K}$-increasing continuous 
in $\zeta$
\sn
\item "{$(\beta)$}"  $M^*_0 = M_{\langle \rangle}$
\sn
\item "{$(\gamma)$}"  $g_{\eta_\zeta}$ is a $\le_{\frak K}$-embedding of 
$M^0_{\eta_\zeta}$ into $M^*_{\zeta +1}$
\sn
\item "{$(\delta)$}"  if $\eta_\xi \triangleleft \eta_\zeta$ then
$g_{\eta_\zeta} \circ h_{\eta_\zeta,\eta_\xi} = g_{\eta_\xi}$
\sn
\item "{$(\varepsilon)$}"  if $\xi < \xi_0,\xi < \xi_1,\eta_{\xi_0} =
\eta_\xi \char 94 \langle 0 \rangle,\eta_{\xi_1} = \eta_\xi \char 94
\langle 1 \rangle$ then \newline
($M^*_{\xi_0} = M^*_{\xi_1}$ and) $g_{\eta_{\xi_0}} = g_{\eta_{\xi_1}}$.
\endroster
\medskip

\noindent
So for $\eta \in {}^\mu 2$ we can find $g_\eta$, 
a $\le_{\frak K}$-embedding of $M^0_\eta$ into $M^*_{\zeta^*}$,
such that \newline
$g_\eta \circ h_{\eta,\eta
\restriction \alpha} = g_{\eta \restriction \alpha}$ for every 
$\alpha < \mu$.  We also can let \newline
$p^0_\eta = g_\eta[\text{tp}(a_\eta,M^0_\eta,
M^1_\eta)] \in {\Cal S}(M^*_{\zeta^*} \restriction \text{Rang}(g_\eta))$, 
and find $p_\eta$ such that $p^0_\eta \le p_\eta \in {\Cal S}(M^*_{\zeta^*})$ 
(possible as ${\frak K}_\lambda$ has amalgamation by \scite{2.2}(1) if
$(*)^2_\lambda$ holds and by $(*)^{2^-}_\lambda$ otherwise).  

For $\eta \in {}^\mu 2$ and $\alpha \le \mu$ let
$N^0_{\eta \restriction \alpha} = M^*_{\zeta^*} \restriction$ Rang 
$(g_{\eta \restriction \alpha})$.  Clearly for \nl
$\eta \in {}^{\mu \ge}2,N^0_\eta$ is well defined; 
$\eta \triangleleft \nu \in {}^{\mu \ge} 2 
\Rightarrow N^0_\eta \le_{\frak K} N^0_\nu$; and $N^0_{\eta \char 94 
\langle 0 \rangle} =
N^0_{\eta \char 94 \langle 1 \rangle}$.  Also letting for $\eta \in {}^\mu 2$,
and $\alpha \le \mu$, the type $p^0_{\eta \restriction \alpha}$ be
$p^0_\eta \restriction N^0_{\eta \restriction \alpha}$ we have: 
$p^0_\eta \in {\Cal S}(N^0_\eta)$ is well defined, $\eta \triangleleft \nu
\in {}^{\mu \ge} 2 \Rightarrow p^0_\eta \le p^0_\nu$ and $p^0_{\eta
\char 94 \langle 0 \rangle} \ne p^0_{\eta \char 94 \langle 1 \rangle}$.
Hence for $\eta_0 \ne \eta_1$ from ${}^\mu 2$ we have 
$p_{\eta_0} \ne p_{\eta_1}$.  So $|{\Cal S}(M^*_{\zeta^*})| \ge 2^\mu >
\lambda$.
\enddemo
\bigskip

\demo{Proof of $(*)$}  We choose by induction on $i < \lambda^+,N_i \in
K_\lambda$ which is $\le_{\frak K}$-increasing continuous, $N_0 =
M^+_0$ ($M^+_0$ is from $(*)_0$ above) and for each $i$ some $c_i \in 
N_{i+1}$ realizes over $N_0$ a
(complete) extension of $p = \text{ tp}(a,M_0,M_1)$ not realized in $N_i$.
There is such type by clause $(*)$ above and there is such an $N_{i+1}$ as
${\frak K}$ has amalgamation in $\lambda$.  Clearly $c_i \notin N_i$ and so
$\dsize \bigcup_{i < \lambda^+} N_i$ is as required. \newline
4) By using \scite{0.17}(2) repeatedly $\lambda^+$ times. 
\hfill$\square_{\scite{2.6}}$
\enddemo
\bigskip

\proclaim{\stag{2.7} Claim}  Assume $(*)^{2^-}_{\lambda^+}$ (from
\scite{2.6}; that is ${\frak K}$ has amalgamation in $\lambda$ and
$LS({\frak K}) \le \lambda$).

If $(M_0,M_1,a) \le (M'_0,M'_1,a)$ are from $K^3_\lambda$, and the second has
the extension property, \underbar{then} so does the first.
\endproclaim
\bigskip

\demo{Proof}  Use amalgamation over $M_0$: if $M_0 \le_{\frak K} N_0 \in 
K_\lambda$ we can find $N'_0$ such that \newline
$M'_0 \le_{\frak K} N'_0 \in K_\lambda$ and there is a 
$\le_{\frak K}$-embedding of $N_0$ into $N'_0$ over $M_0$.  Now use 
``$(M'_0,M'_1,a)$ has the extension property" for $N'_0$. 
\hfill$\square_{\scite{2.7}}$
\enddemo
\bn
Now we introduce
\definition{\stag{2.8} Definition}  For any models
$M,M_0 \in K_\lambda$, any type $p \in {\Cal S}(M_0)$ and \newline
$f_0:M_0 \overset\text{onto}\to{\underset\text{iso}\to
\longrightarrow} M$ we let ${\Cal S}_p(M) = {\Cal S}^p_M = \{ f_0(f(p)):f \in 
\text{ AUT}(M_0)\}$.  Note: ${\Cal S}^p_M$ does not depend on $f_0$.
If $K$ is categorical in $\lambda,{\Cal S}_p(M)$ is well defined for every
$M \in K_\lambda$.  We write also ${\Cal S}_{\text{tp}(a,M_0,M_1)}(M)$ or 
${\Cal S}_{(M_0,M_1,a)}(M)$ when $(M_0,M_1,a) \in K^3_\lambda$.
\enddefinition
\bigskip

\proclaim{\stag{2.9} Claim}  Assume $(*)^{2^-}_\lambda + 
(*)^2_{\lambda^+} + (*)^3_\lambda + 2^{\lambda^+} < 2^{\lambda^{++}} + 
K_{\lambda^{+3}} = \emptyset$. \newline
If $(M_0,M_1,a) \in K^3_\lambda$ is minimal \underbar{then} it has the 
extension property.
\endproclaim
\bigskip

\remark{Remark}  Instead of $2^{\lambda^+} < 2^{\lambda^{++}}$, we can just
demand the definable weak diamond.
\endremark
\bigskip

\demo{Proof}  Assume not.  By the previous two claims (\scite{2.6}(1),
\scite{2.7}) we may assume that $(M_0,M_1,a)$ is minimal.
As ${\frak K}$ has amalgamation in $\lambda^+$ by \scite{2.2}(1), there 
is $M^* \in K_{\lambda^{++}}$ which is saturated 
above $\lambda^+$ (as $K_{\lambda^{+3}} = \emptyset$) hence $M^*$ is
saturated above $\lambda$ (by \scite{1.6}(3)).
By \scite{2.4}(1) + \scite{2.7}, without loss of generality
\medskip
\roster
\item "{$\bigotimes_0$}"  $(M_0,M_1,a)$ is reduced.
\endroster
\medskip

\noindent
Let $h:M_0 \rightarrow M^* \in K_{\lambda^+}$ be a 
$\le_{\frak K}$-embedding and let 
$p = \text{ tp}(a,M_0,M_1)$.  If $h(p)$ is realized in $M^*$ by 
$\ge \lambda^+$ elements we are done by \scite{2.6}(2).  So assume
\medskip
\roster
\item "{$\bigotimes_1$}"  $h(p)$ is realized by $\le \lambda$ members of
$M^*$.
\ermn
Similarly
\mr
\item "{$\bigotimes^+_1$}"   $g$ is realized by $\le \lambda$ members of
$M^*$ for $q = g(\text{tp}(a,M'_0,M'_1))$ if $g$ is a
$\le_{\frak K}$-embedding of $M'_0$ into $M^*$ and
$(M'_0,M'_1,a) \ge (M_0,M_1,a)$.
\ermn
Next we prove
\medskip
\roster
\item "{$\bigotimes_2$}"  for some reduced $(M'_0,M'_1,a) \ge (M_0,M_1,a)$
from $K^3_\lambda$ we have \newline
$|{\Cal S}_{(M'_0,M'_1,a)}(M'_0)| > \lambda^+$.
\endroster
\enddemo
\bigskip

\demo{Proof of $\bigotimes_2$}  If not, we build two non-isomorphic members
of $K_{\lambda^+}$ as follows.
\medskip

\noindent
\underbar{First}:  Choose by induction on $i < \lambda^+,(N_{0,i},N_{1,i},a)
\in K^3_\lambda$ reduced (see \scite{2.4}(1)), increasing continuously (see 
\scite{2.4}(2)), with $N_{0,i} \ne N_{0,i+1},(N_{0,0},N_{1,0},a) = 
(M_0,M_1,a)$; this is possible as
$(N_{0,i},N_{1,i},a) \in K^3_\lambda$ has the weak extension property 
(by \scite{2.3} see \scite{2.3A}(1)).  
Let $N^1 = \dsize \bigcup_{i < \lambda^+} N_{0,i}$.
\medskip

\noindent
\underbar{Second}:  Choose by induction on $i < \lambda^+,N^0_i 
\le_{\frak K} M^*,\| N^0_i\| = \lambda,N^0_i$ strictly increasing continuous 
such that: 
\medskip
\roster
\item "{$(*)$}"  for every $\beta < \lambda^+$, and
$q \in {\Cal S}_{(M'_0,M'_1,a)}(N^0_\beta)$ for some $\gamma \in
(\beta,\lambda^+)$ there are no $N',N^0_\gamma \le_{\frak K} N' \in 
K_\lambda$ and $c \in N' \backslash N^0_\gamma$ such that $c$ realize $q$.
\endroster
\mn
This is straightforward by $\otimes_1 + \otimes_2$ and bookkeeping.
Let $N^0 = \dsize \bigcup_{i < \lambda^+} N^0_i$.
By categoricity of $K_{\lambda^+}$ there is an isomorphism $g$ from $N^1$
onto $N^0$, so \newline
$E = \{ \delta < \lambda^+:g \text{ maps } N_{0,\delta} \text{ onto }
N^0_\delta\}$ is a club of $\lambda^+$.  Now let $\delta^* \in E$, and 
apply $(*)$ for $\beta = \delta^*,q = g(\text{tp}(a,N_{0,\delta^*},
N_{1,\delta^*}))$ to get $\gamma$.  Choose $\delta \in E$ which is $> \gamma$.
Now $N_{1,\delta}$ gives a contradiction.
\hfill$\square_{\otimes_2}$
\medskip

\noindent
Without loss of generality
\medskip
\roster
\item "{$\bigotimes_3$}"  $|{\Cal S}_{(M_0,M_1,a)}(M_0)| > \lambda^+$ and
$p = \text{ tp}(a,M_0,M_1)$.
\endroster
\medskip

\noindent
Next we claim
\roster
\item "{$\bigotimes_4$}"  If $M \in K_\lambda,M \le_{\frak K} M^*,\Gamma
\subseteq \bigcup\{{\Cal S}_p(M'):M' \le_{\frak K} M,\| M'\| = \lambda\}$, 
\newline
$| \Gamma | \le \lambda^+$, then
$$
\align
\Gamma^* =: \{q \in {\Cal S}_p(M):&\text{there is } M',
M \le_{\frak K} M',\| M' \| = \lambda, \\
  &M' \text{ realizes } q \text{ but no } r \in \Gamma\}
\endalign
$$
has cardinality $\lambda^{++}$, in fact $|{\Cal S}_p(M) \backslash 
\Gamma^*| \le \lambda^+$.
\endroster
\enddemo
\bigskip

\demo{Proof of $\otimes_4$}  Without loss of generality $|M^*| = \lambda
^{++}$ (i.e. the universe of $M^*$ is $\lambda^{++}$).  For every 
$q \in {\Cal S}_p(M)$ there is a triple $(M_0,M_{1,q},a_q)$ isomorphic to
$(M_0,M_1,a)$ (hence reduced) such that tp$(a_q,M_0,M_{1,q}) = q$.  As
$M^*$ is saturated above $\lambda$, by \scite{0.19} without loss of generality
$M_{1,q} \le_{\frak K} M^*$.

Without loss of generality $\delta < \lambda^{++} \and (\lambda^+$ divides
$\delta$) $\Rightarrow M_\delta =: M^* \restriction \delta \le_{\frak K} M^*$.
Now
\medskip
\roster
\item "{$(*)_0$}"  $q_1 \ne q_2 \Rightarrow a_{q_1} \ne q_{q_2}$ and
\sn
\item "{$(*)_1$}"  $a_q \notin \delta \and \delta < \lambda^{++} 
\and \lambda^+ \text{ divides } \delta \Rightarrow M_{1,q} \cap M_\delta = M$.
\endroster
\medskip

\noindent
[Why?  As $(M,M_{1,q},a_q)$ is reduced]. \newline
Now if $r \in \Gamma$, say $r \in {\Cal S}_p(M'')$ then by
$\otimes_1$ we know $A_r = \{ c \in M^*:c \text{ realizes } r\}$
has cardinality $\le \lambda$ and hence $A = \bigcup\{ A_r:r \in \Gamma\}$ has
cardinality $\le \lambda^+$, so we can find $\delta < \lambda^{++}$
divisible by $\lambda^+$ such that $A \subseteq \delta$.
So as (by $\otimes_3$) we have $|{\Cal S}_p(M)| > \lambda^+$ hence we can find
$q[\delta] \in {\Cal S}_p(M)$ such that $a_{q(\delta)} \notin \delta$ hence
$(M,M_{1,q[\delta]},a_{q[\delta]})$, exemplifies the conclusion of 
$\otimes_4$. \hfill$\square_{\otimes_4}$
\enddemo
\bigskip

\noindent
\underbar{Final contradiction}:  By $\otimes_4$ we can construct
$2^{\lambda^+}$ non-isomorphic members of $K_{\lambda^+}$ using 
\scite{1.4}(1) as follows.  We choose by induction on $\alpha < \lambda^+$, 
for every
$\eta \in {}^\alpha 2,M_\eta,p^0_\eta,p^1_\eta$ such that:
\medskip
\roster
\item "{$(a)$}"  $M_{<>} = M_0$
\sn
\item "{$(b)$}"  $M_\eta \in K_\lambda$
\sn
\item "{$(c)$}"  $\langle M_{\eta \restriction \beta}:\beta \le \ell g(\eta)
\rangle$ is $<_{\frak K}$-increasing continuous
\sn
\item "{$(d)$}"  $p_\eta \in {\Cal S}_p(M_{\eta \restriction \beta})$
\sn
\item "{$(e)$}"  for $\beta \le \alpha$, we have 
$p^0_\eta,p^1_\eta \in {\Cal S}(M_\eta)$ and: $M_\eta$ realizes 
$p^\ell_{\eta \restriction \beta}$ iff $\beta < \alpha \and \ell = 
\eta(\beta)$.
\endroster
\medskip

\noindent
If $\alpha = 0$ or $\alpha$ is a limit, there is no problem to define 
$M_\eta$ for $\eta \in {}^\alpha 2$.  If $M_\eta$ is defined, we can 
choose by induction on $i < \lambda^{++},(N_{\eta,i},a_{\eta,i})$ such that
$(M_\eta,N_{\eta,i},a_{\eta,i}) \in K^3_\lambda$, tp$(a_{\eta,i},M_\eta,
N_{\eta,i}) \in {\Cal S}_p(M_\eta)$ and $N_{\eta,i}$ omits any \newline
$q \in \{p^\ell_{\eta \restriction \beta}:\beta < \ell g(\eta),\ell \ne
\eta(\beta)\} \cup \{\text{tp}(b,M_\eta,N_{\eta,j}):j < i \text{ and}$
\newline
$b \in N_{\eta,j} \text{ and tp}(b,M_\eta,N_{\eta,j}) \in 
{\Cal S}_p(M_\eta)\}$.
By $\otimes_4$ we can choose $(N_{\eta,j},a_{\eta,j})$.  \newline
Hence $|W_{\eta,i}| \le \lambda$ where
$W_{\eta,i} = \{j < \lambda^{++}:\text{for some } b \in N_{\eta,i}$ 
we have tp$(b,M_\eta,N_{\eta,i}) = \text{tp}(a_{\eta,j},M_\eta,
N_{\eta,j})\}$. \newline
Hence we can find $i < j < \lambda^{++}$ such that $i \notin W_{\eta,j} 
\and j \notin W_{\eta,i}$.  Let $M_{\eta \char 94 \langle 0
\rangle} = N_{\eta,i},M_{\eta \char 94 \langle 1 \rangle} = N_{\eta,j},
p^0_\eta = \text{ tp}(a_{\eta,i},M_\eta,N_{\eta,i}),p^1_\eta = \text{ tp}
(a_{\eta,j},M_\eta,N_{\eta,j})$. \newline
Let for $\eta \in {}^{\lambda^+}2,
M_\eta = \dsize \bigcup_{\alpha < \lambda^+} M_{\eta \restriction \alpha}$,
and apply \scite{1.4}(1). \hfill$\square_{\scite{2.8}}$
\bigskip

\demo{\stag{2.10} Conclusion}  $[(*)^{2^-}_\lambda + (*)^2_{\lambda^+} + 
(*)^3_\lambda + 2^{\lambda^+} < 2^{\lambda^{++}} + K_{\lambda^{+3}} = 
\emptyset]$. \newline

Every $(M_0,M_1,a) \in K^3_\lambda$ has the extension property.
\enddemo
\bigskip

\demo{Proof}  By \scite{2.7} and \scite{2.6} + \scite{2.8}. 
\hfill$\square_{\scite{2.10}}$
\enddemo
\bigskip

\remark{\stag{2.10A} Remark}  Conclusion \scite{2.10} says in other words: \nl
\ub{if}
\mr
\item "{$(a)$}"  $LS({\frak K}) \le \lambda$
\sn
\item "{$(b)$}"  $K$ is categorical in $\lambda$ and in $\lambda^+$
\sn
\item "{$(c)$}"  $1 \le I(\lambda^{++},K) < 2^{\lambda^{++}}$
\sn
\item "{$(d)$}"  $K_{\lambda^{+3}}$ is empty
\sn
\item  "{$(e)$}"  $2^{\lambda^+} < 2^{\lambda^{++}}$ (or just definable weak
diamond)
\ermn
\ub{then} every triple 
$(M_0,M_1,a)$ in $K^3_\lambda$ has the extension property.
\endremark
\bigskip

\proclaim{\stag{2.11} Claim}  [$(*)^3_\lambda$, in other words
$LS({\frak K}) \le \lambda;{\frak K}$ categorical in $\lambda$ and in
$\lambda^+$; and $1 \le I(\lambda^{++},K) < 2^{\lambda^{++}}$]. \newline

If $M_0 \le_{\frak K} M_1$ are in $K_\lambda$ \underbar{then} we can 
find $\alpha < \lambda^+$ and
$\langle N_i:i \le \alpha \rangle$ which is $\le_{\frak K}$-increasing 
continuous, $N_i \in K_\lambda,(N_i,N_{i+1},a_i) \in K^3_\lambda$ is reduced, 
$M_0 = N_0$, and $M_1 \le_{\frak K} N_\alpha$.
\endproclaim
\bigskip

\demo{Proof}  If not, we can contradict categoricity in $K_{\lambda^+}$ 
(similar to the proof of $\otimes_2$ during the proof of \scite{2.8}).

Without loss of generality $M_0 \ne M_1$.  We choose by induction on $i <
\lambda^+$, \newline
$N^0_i \in K_\lambda,\le_{\frak K}$-increasing continuous such 
that $(N^0_i,N^0_{i+1}) \cong (M_0,M_1)$ (possible by \scite{2.4}(9) and
the categoricity of $K$ in $\lambda$).  Let
$N^0 = \dsize \bigcup_{i < \lambda^+} N^0_i$.

We choose by induction on $i < \lambda^+,N^1_i \in K_\lambda,
\le_{\frak K}$-increasing continuous and $a_i$ such that
$(N^1_i,N^1_{i+1},a_i) \in K^3_\lambda$ is reduced and let $N^1 =
\dsize \bigcup_{i < \lambda^+} N^1_i$ (possible by \scite{2.4}(1) and the
categoricity of $K$ in $\lambda$).
So by the categoricity in $\lambda^+$ without loss of generality $N^1 = N^0$,
hence for some $\delta_1 < \delta_2 < \lambda^+$ we have

$$
N^0_{\delta_1} = N^1_{\delta_1},N^0_{\delta_2} = N^1_{\delta_2}.
$$
\medskip

\noindent
By changing names $(N^0_{\delta_1},N^0_{\delta_1+1}) = (M_0,M_1)$ and so
$\langle N_{\delta_2+i}:i \le \delta_2 - \delta_1 \rangle$ is as required.
\hfill$\square_{\scite{2.11}}$
\enddemo
\bigskip

\demo{\stag{2.12} Conclusion}  $[(*)^{2^-}_\lambda + 
(*)^2_{\lambda^+} + (*)^3_\lambda + 2^{\lambda^+} < 2^{\lambda^{++}} +
K_{\lambda^{+3}} = \emptyset$, i.e. the \nl
assumption of \scite{2.10}]. \newline
$K_\lambda$ has disjoint amalgamation ($M_2,M_1$ are
in disjoint amalgamation over $M_0$ in $M_3$ if $M_0 \le_{\frak K} M_\ell 
\le_{\frak K} M_3,M_1 \cap M_2 = M_0$).
\enddemo
\bigskip

\demo{Proof}  By \scite{2.11} and iterated applications of \scite{2.10}.
\hfill$\square_{\scite{2.12}}$
\enddemo
\newpage

\head {\S3 Non-structure} \endhead  \resetall
\bigskip

The first major aim of this section is to prove the density of minimal 
types using as set theoretic assumptions only $2^\lambda <
2^{\lambda^+} < 2^{\lambda^{++}}$ from cardinal arithmetic.  The second aim
is to prepare for a proof of a weak form of uniqueness of amalgamation in
${\frak K}_\lambda$.  Our aim is also to explain various methods. 
The proofs are similar to the ones in \cite[\S6]{Sh:87b}. 
\sn
The immediate role of this section is to get many models in $\lambda^{++}$
from the assumption ``the minimal triples in $K^3_\lambda$ are not dense":
in \scite{3.16A} we get this under some additional assumptions, and in
\scite{3.27} we get it using only the additional assumption $I(\lambda,
K^{+3})=0$, which suffices for our main theorem (this does not suffice
for the theorem of \cite{Sh:600}, hence is eliminated there).

But the section is done in a more general fashion, so let us first explain
two general results concerning the construction of many models based on
repeated ``failures of amalgamation" or ``nonminimality of types".

In \scite{3.12}, we give a construction assuming the ideal of small subsets 
of $\lambda^+$ (that is 
WDmId$(\lambda^+)$) is not $\lambda^{++}$-saturated, as exemplified by
$\langle S_\alpha:\alpha < \lambda^{++} \rangle$.  We build for 
$\eta \in {}^{(\lambda^{+2})>}2$ models $M_\eta \in K_{\lambda^+}$ such
that $M_\eta = \dsize \bigcup_{\alpha < \lambda^+} M_{\eta,\alpha},
|M_\eta| = \lambda \times (1 + \ell g \eta)$ and
$\nu \triangleleft \eta \Rightarrow M_\nu \le_{\frak K} M_\eta$.  
Building $M_{\eta \char 94 \langle \ell \rangle}$ manufacture
$M_{\eta \char 94 \langle \ell \rangle,\alpha + 1}$ as a limit of models
$\langle M_{\eta \char 94 \langle \ell \rangle,\alpha}:\alpha < \lambda^+ 
\rangle$, a representation of $M_{\eta \char 94 \langle \ell \rangle}$,
usually in a way predetermined simply, except when 
$\alpha \in S_{\ell g(\eta)}$ and $\ell = 1$, and then we consult a 
weak diamond sequence.  This is like \scite{1.4}(1), but there we use 
our understanding of models in $K_\lambda$ to build many models in 
$K_{\lambda^+}$ while here we build models in $K_{\lambda^{++}}$, thus 
getting $2^{\lambda^{+2}}$ models in $\lambda^{+2}$.  We even get 
$2^{(\lambda^{+2})}$ models in $K_{\lambda^{+2}}$ with none 
$\le_{\frak K}$-embeddable into any other.
\medskip

A second proof \scite{3.17B} is like \scite{1.3}(1) in the sense that 
we get only close to $2^{\lambda^{++}}$ models.  It is similar to 
\cite[6.4]{Sh:87b}, and the parallel
to \cite[6.3]{Sh:87b} holds here.  So we have to find an analog of 
\cite[definition 6.5,6.7]{Sh:87b}.  But there we use fullness on the side 
(meaning: $M \in K_\lambda$ is full over $N \in K_\lambda$ if
$N \le_{\frak K} M$, and $(M,c)_{c \in N}$ is saturated), but we do not have
this yet.

We still have not explained the framework of this section.  In
\scite{3.0} - \scite{3.2} we present construction frameworks $\bold C$, 
which involve sequences of models of length $\le \lambda$ each of cardinality 
$< \lambda$ and in particular, define local and nice $\bold C$.
In our applications here $\lambda^+$ plays the role of $\lambda$ (and
$< \lambda^+$ is specialized to $\lambda$).  \nl
Then in \scite{3.3} - \scite{3.4} we present examples of such frameworks.  Our
intention is to use the limit of a sequence 
$\langle M_\alpha:\alpha < \lambda \rangle$ as an approximation to a model of 
cardinality $\lambda^+$.  For this we define in \scite{3.5} - \scite{3.6} 
a successor relation (next approximation), modulo a
``$< \lambda$-amalgamation choice function"; this is denoted 
$\bar M^1 \le^{\text{at}}_{F_1} \bar M^2$.  Iterating it we get the quasi 
order $\le_F$ (see \scite{3.7}).  In \scite{3.8} we define the key coding 
properties, (of an amalgamation choice function $F$ for the framework
$\bold C$).  The intention is that these coding properties suffice to build 
many non-isomorphic models in $\lambda^+$.  In \scite{3.10} we give the 
``atomic step" for this construction.

In \scite{3.11} we prove the existence of $2^{\lambda^+}$ non-isomorphic
models, using the $\lambda$-coding property.  As we do not have this
in some applications we have in mind, we next turn to the weak 
$\lambda$-coding property in \scite{3.12} as well as the weak (local) 
$\lambda$-coding property and corresponding properties of $F$ (all
in Definition \scite{3.13}, \scite{3.17A}), connect them \scite{3.14}, and
prove there are many models in \scite{3.17B}. \nl
Lastly, \scite{3.16A}, \scite{3.27} deal with our concrete case: if the minimal
triples in ${\frak K}^3_\lambda$ are not dense then in most cases failures
of amalgamation lead to the $\lambda^+$-coding property and hence to many
models in cardinality $\lambda^{++}$.

Note generally that we mainly axiomatize the construction of models in 
$\lambda^+$, not how we get $\bar M',\bar M \le^{\text{at}}_{F,a} \bar M' \in
\bold S eq_\lambda$ that is coding properties; for the last point, see the
examples just cited.

Later, in \scite{6.5B}, we shall need again to use the machinery from this
section, in trying to prove that there are enough cases of disjoint
amalgamation in ${\frak K}_\lambda$. \nl
We may want to turn the framework presented here into a more general
one; see more in \cite{Sh:600}.
\bigskip

\demo{\stag{3.0} Context} \nl
1) ${\frak K}$ is an abstract elementary class. \newline
2) But $=_M$ or $=_{\frak K}$ is just an equivalence relation, i.e. for 
$M \in K,=^M$ is
an equivalence relation on $|M|$, moreover a congruence relation relative
to all relations and functions in $\tau(M)$ that is for $R \in \tau(M)$ an 
$n$-ary relation, we have

$$
\dsize \bigwedge_{\ell < n} a_i =^M b_i \Rightarrow \langle a_0,\dotsc,
a_{n-1} \rangle \in R^M \equiv \langle b_0,\dotsc,b_{n-1} \rangle \in R^M).
$$
\medskip

\noindent
We let $\|M\| = |(M)/ =^M|$ and

$$
K_{\lambda,\mu} = \{M:|M| \text{ has } \mu \text{ elements and }
|M|/ =^M \text{ has } \lambda \text{ elements}\}.
$$
\medskip

\noindent
3) Now the meaning of $\le_{\frak K}$ should be clear but 
$M <_{\frak K} N$ means
$(M \in K,N \in K$ and) $M \subseteq N$ and $M/=^M \le_{\frak K} N/=^N$ and
$\exists a \in N[a/=^N \notin (M/=^N)]$ i.e. $(\exists a \in N)(\forall b \in
M)(\neg a =^N b)$. \newline
4) $K^3_\lambda = \{(M,N,a):M \le_{\frak K} N$ are from 
$K_{\lambda,\lambda}$ and $a \in N,(a/=^N) \notin (M/=^N)\}$. \nl
5) In this content ``$R$ is an isomorphism relation from $M_1$ onto $M_2$
if
\mr
\item "{$(a)$}"  $R \subseteq M_1 \times M_2$
\sn
\item "{$(b)$}"  $a_1 =^{M_1} b_1 \and a_2 =^{M_2} b_2 \Rightarrow a_1 R a_2
\leftrightarrow b_1 R b_2$
\sn
\item "{$(c)$}"  $(\forall x \in M_1)(\exists y \in M_2) x Ry$
\sn
\item "{$(d)$}"  $(\forall y \in M_2)(\exists x \in M_1) x Ry$
\sn
\item "{$(e)$}"  if $Q \in \tau(M_1) = \tau(M_2)$ is an $n$-place relation
and $a_\ell Q b_\ell$ for $\ell = 0,\dotsc,n-1$ then $(a_0,\dotsc,a_{n-1})
\in Q^{M_1} \rightarrow (b_0,\dotsc,b_{n-1}) \in Q^{M_2}$.
\endroster
\enddemo
\bigskip

\demo{\stag{3.0A} Explanation}  The need of \scite{3.0}(2) is 
just to deal with amalgamations which are not necessarily disjoint.  
If we use disjoint amalgamation, we can
omit \scite{3.0}(2) below on Definition \scite{3.5}, $a$ disappears so $F$ 
is four place and use $K_\lambda$ instead of $K_{\lambda,\lambda}$.  This is
continued in \cite[2.17]{Sh:600}.  May be better understood after reading
\scite{3.5}, after clause (c).
\enddemo
\bigskip

\definition{\stag{3.OB} Definition}  Let $\lambda$ be regular uncountable
and ${\frak K}$ an abstract elementary class.
\newline
1) A $\lambda$-construction framework $\bold C = \bold C({\frak K}^+,
{\bold Seq},\le^*)$ means (we shall use it below with $\lambda^+$ playing the
role of $\lambda$):
\medskip
\roster
\item "{$(a)$}"  $\tau^+ = \tau^+({\frak K}^+)$ is a vocabulary extending
$\tau$.  ${\frak K}^+$ is an abstract elementary class satisfying 
axioms I,II,III from \scite{0.2} and $M \le_{{\frak K}^+} N \Rightarrow M 
\restriction \tau \le_{\frak K} N \restriction \tau$.  Furthermore
${\frak K}^+ = {\frak K}^+_{< \lambda}$.  As above, equality (in $\tau$) is
a congruence relation.
\sn
\item "{$(b)$}"  ${\bold Seq} = \dsize \bigcup_{\alpha \le \lambda}
{\bold Seq}_\alpha$ where for $\alpha \le \lambda,{\bold Seq}_\alpha$ is a
subset of \newline
$\{\bar M:\bar M = \langle M_i:i < \alpha \rangle,M_i \in {\frak K}^+$ is 
$\le_{{\frak K}^+}$-increasing continuous$\}$. \newline
For $\alpha = \lambda$ we require further that $M/=^M$ has cardinality 
$\lambda$, where \nl
$M = \dsize \bigcup_{i < \lambda} M_i$.
\sn
\item "{$(c)$}"  $\le^*$ is a relation on triples $x,y,t$ written 
$x \le^*_t y$ for $x,y \in {\bold Seq}$ and $t$ a set of pairwise 
disjoint closed intervals of $\ell g(x)$.
\ermn
We require:
\medskip
\roster
\item "{$(d)$}"  ${\bold Seq}$ is closed under isomorphism and initial
segments.
\sn
\item "{$(e)$}"  If $\bar M^1 \le^*_t \bar M^2$ and $\gamma \in \cup t$
then $M^1_\gamma \le_{{\frak K}^+} M^2_\gamma$ and hence 
$M^1_\gamma \restriction \tau \le_{\frak K} M^2_\gamma \restriction \tau$.
\sn
\item "{$(f)$}"  If $\bar M^1 \le^*_t \bar M^2,s \subseteq t$, and
$\bar M^2 \trianglelefteq \bar M^3 \in {\bold Seq}$ then $\bar M^1 \le^*_s
\bar M^3$.
\sn
\item "{$(g)$}"  If $t$ is a set of closed pairwise disjoint intervals
of $\ell g(\bar M)$ and $\bar M \in {\bold Seq}$ then $\bar M \le^*_t \bar M$.
\endroster
\enddefinition
\bigskip

\noindent
\underbar{\stag{3.1} Convention/Definition}: \nl
1) From now on $\bold C$ will be a $\lambda$-construction framework. \newline
2) If $\bar M \in {\bold Seq}_\lambda$ then we let $\bar M = \langle M_i:
i < \lambda \rangle$ and $M =: \dsize \bigcup_{i < \lambda} M_i$; similarly
with $\bar M^x = \langle M^x_i:i < \lambda \rangle$. \newline
3) $K^{\text{qr}}_\lambda = \{(\bar M,\bold f):\bar M \in 
{\bold Seq}_\lambda \text{ and }
\bold f:\lambda \rightarrow \lambda\}$. \newline
4) If $(\bar M^\ell,\bold f^\ell) \in K^{\text{qr}}_\lambda$ for
$\ell = 1,2$ then $(\bar M^1,\bold f^1) \le (\bar M^2,\bold f^2)$ means
that: for some club $E$ of $\lambda$, we have 
\mr
\item "{$(a)$}"  $\delta \in E \Rightarrow \bold f^1(\delta) \le 
\bold f^2(\delta)$ and 
\sn
\item "{$(b)$}"  $\bar M^1 \le^*_t \bar M^2$ where 
$t = t_{E,\bold f^1} = \{[\delta,\delta + \bold f^1(\delta)]:\delta \in E\}$.
\endroster
\medskip

\noindent
5) ${\bold Seq}^s = \{ \bar M \in {\bold Seq}:
\dsize \bigcup_i |M_i| \text{ is a set of ordinals } < \lambda^+\}$; 
similarly for ${\bold Seq}_\alpha$. \newline
6) $K^{\text{qs}}_\lambda = \{(\bar M,\bold f) \in 
K^{\text{qr}}_\lambda:\bar M \in {\bold Seq}^s_\lambda\}$. 
\medskip

\noindent
7) $\bold C$ is \ub{local} 
(respectively, revised local) if (a-c) following hold:
\medskip
\roster
\item "{$(a)$}"  $\bar M = \langle M_i:i < \alpha \rangle \in {\bold Seq}
_\alpha$ iff:
{\roster
\itemitem{ $(\alpha)$ }  $\bar M$ is $\le_{{\frak K}^+}$-increasing continuous
in ${\frak K}^+_{< \lambda}$
\sn
\itemitem{ $(\beta)$ }  $i+1 < \alpha \Rightarrow \langle M_i,M_{i+1} \rangle
\in {\bold Seq}_2$
\sn
\itemitem{ $(\gamma)$ }  if $\alpha = \lambda$ then $|M/=^M| = \lambda$,
(recall $M = \dsize \bigcup_{i < \lambda} M_i$)
\endroster}
\item "{$(b)$}"  for $\bar M^1,\bar M^2 \in {\bold Seq}$ and $t$ a set of
pairwise disjoint closed intervals contained in $\ell g(\bar M^1)$ we have:
$$
\bar M^1 \le^*_t \bar M^2 \text{ \underbar{iff} } [\gamma_1,\gamma_2] \in t
\text{ implies}
$$
{\roster
\itemitem{ $(\alpha)$ } $\gamma \in [\gamma_1,\gamma_2] \Rightarrow 
M^1_\gamma \le_{{\frak K}^+} \bar M^2_\gamma$
\sn
\itemitem{ $(\beta)$ }  in the local case: $\gamma \in [\gamma_1,\gamma_2)
\Rightarrow \langle M^1_\gamma,M^1_{\gamma + 1} \rangle \le^*_{\{[0,1]\}} 
\langle M^2_\gamma,M^2_{\gamma + 1} \rangle$; \nl
in the revised local case:  if $\ell g(\bar M^1),\ell g(\bar M^2) < \lambda$
then \nl
$\gamma \in [\gamma_1,\gamma_2]
\Rightarrow \langle M^1_\gamma,M^1 \rangle \le^*_{\{[0,1]\}} \langle
M^2_\gamma,M^2 \rangle$, and generally for some club $E$ of $\lambda,
\gamma \in [\gamma_1,\gamma_2] \and \gamma < \delta \in E \Rightarrow \langle
M^1_\gamma,\dbcu_{\beta < \delta} M^1_\beta \rangle \le^*_{\{[0,1]\}}
\langle M^2_\gamma,\dbcu_{\beta < \delta} M^2_\beta \rangle$ (and if
$\ell g(\bar M^\ell) = \alpha_\ell < \delta,
\dbcu_{\beta < \delta} M^\ell_\beta$ means 
$\dbcu_{\beta < \alpha^\ell} M^\ell_\beta)$.
\endroster}
\item "{$(c)$}"  if $\langle M^\zeta_0,M^\zeta_1 \rangle \in
{\bold Seq}_2$ for $\zeta < \zeta^* < \lambda,\langle M^\zeta_0:\zeta \le 
\zeta^* \rangle$ and $\langle M^\zeta_1:\zeta \le \zeta^* \rangle$ are 
$\le_{{\frak K}^+}$-increasing continuous, and $\zeta <
\zeta^* \Rightarrow \langle M^\zeta_0,M^\zeta_1 \rangle \le^*_{\{[0,1]\}}
\langle M^{\zeta +1}_0,M^{\zeta +1}_1 \rangle$ \underbar{then} 
$\langle M^0_0,M^0_1 \rangle \le^*_{\{[0,1]\}} \langle M^{\zeta^*}_0,
M^{\zeta^*}_1 \rangle \in \bold S eq$.
\ermn
So intervals $[\alpha,\alpha] \in t$ are irrelevant for the local version. 
In the revised local version it is natural to add monotonicity for
$\le_{\{[0,1]\}}$.
\medskip

\noindent
8) For $\alpha \le \lambda$ we say $\bold C$ is \ub{closed} for $\alpha$ if:
\medskip
\roster
\item "{$(\alpha)$}"  $\bar M = \langle M_i:i < \alpha \rangle \in
{\bold Seq}$ \underbar{iff} $\beta < \alpha \Rightarrow \bar M 
\restriction (\beta +1) \in {\bold Seq}$
\sn
\item "{$(\beta)$}"  if $\bar M^\ell = 
\langle M^\ell_i:i < \alpha_\ell \rangle \in {\bold Seq}$ for 
$\ell =1,2$ and $\alpha = \alpha_1 < \alpha_2$, \ub{then}
$$
\bar M^1 \le^*_t \bar M^2 \Leftrightarrow \bar M^1 \le^*_t \bar M^2 
\restriction \alpha.
$$
\endroster
\medskip

\noindent
9) $\bold C$ is \ub{disjoint} if: $\bar M^1 \le^*_t \bar M^2,
[\gamma_1,\gamma_2] \in t,\gamma \in [\gamma_1,\gamma_2)$ implies 
$M^1_\gamma = M^1_{\gamma +1} \cap M^2_\gamma$. \nl
$\bold C$ is \ub{truly disjoint} if: 
$\bar M^1 \le^*_t \bar M^2,[\gamma_1,\gamma_2] \in t,
\gamma \in [\gamma_1,\gamma_2]$ implies $M^1_\gamma = M^1 \cap 
M^2_\gamma$.
\medskip

\noindent
10) In $K^{\text{qr}}_\lambda$, we say $(\bar M,\bold f)$ is a \ub{m.u.b.} 
(minimal upper bound) of $\langle (\bar M^\xi,\bold f^\xi):\xi < \delta 
\rangle$ if
\mr
\item "{$(a)$}"  $\xi < \delta \Rightarrow (\bar M^\xi,\bold f^\xi) \le 
(\bar M,\bold f)$ and 
\sn
\item "{$(b)$}"  for any $(\bar M,\bold f')$ satisfying (a), for some club 
$E$ of $\lambda$ we have: if $\alpha \in E$ and $j \le \bold f(\alpha)$ then
$\bold f(\alpha) \le \bold f'(\alpha)$ and $M_{\alpha + j} \le_{{\frak K}^+} 
M'_{\alpha + j}$.
\ermn
When we require an increasing sequence in $K^{\text{qr}}_\lambda$ to be
continuous we mean that a m.u.b. is used at limits.
\medskip

\noindent
11) We say $\bold C$ is \ub{explicitly local} if it is local and
\medskip
\roster
\item "{$(d)$}"  if $\zeta^* < \lambda$ is a limit ordinal,
$\langle M^\zeta_0,M^\zeta_1 \rangle \in \bold S eq_2$
for $\zeta \le \zeta^*$ and for $\ell = 0,1$ the sequence
$\langle M^\zeta_\ell:\zeta < \zeta^* \rangle$
is $\le_{\frak K}$-increasing continuous, $M^\zeta_\ell \le_{{\frak K}^+}
M^{\zeta^*}_\ell$, and $\zeta < \xi \le \zeta^* \Rightarrow \langle
M^\zeta_0,M^\zeta_1 \rangle \le^*_{\{[0,1]\}} \langle M^\xi_0,M^\xi_1 \rangle$
\underbar{then}
$\langle \dsize \bigcup_{\zeta < \zeta^*} M^\zeta_0,\dsize \bigcup_{\zeta <
\zeta^*} M^\zeta_1 \rangle$ is $\le_{\{[0,1]\}} (M^{\zeta^*}_0,M^{\zeta^*}_1)$.\endroster
\medskip

\noindent
12) $\bold C$ is \ub{closed} if it is closed for every ordinal $\le \lambda$.
\medskip
\noindent
13) $\bold C$ is semi (respectively almost) \ub{closed}, as witnessed 
by $G$, if:
\medskip
\roster
\item "{$(a)$}"  $\bold C$ is closed for every limit ordinal $\delta <
\lambda$;
\sn
\item "{$(b)$}"  $G$ is a function from $\bold S eq_{< \lambda}$ to
$\bold S eq_{< \lambda}$ such that \newline
$\bar M \triangleleft G(\bar M)$;
\sn
\item "{$(c)$}"  $\bar M = \langle M_\alpha:\alpha < \lambda \rangle$ belongs
to $\bold S eq_\lambda$ if $\bar M$ obeys $G$, which means: $\beta < \lambda
\Rightarrow \bar M \restriction \beta \in \bold S eq_\beta$ and $\{ \alpha <
\lambda:G(\bar M \restriction \alpha) \triangleleft \bar M\}$ 
is unbounded in $\lambda$
\sn
\item "{$(d)$}"  in the almost closed version, we add: $G(\bar M)$ 
depends on $\bigcup \bar M$ only.
\endroster
\medskip

\noindent
14) $\bold C$ is \ub{$\lambda$-nice} if
\medskip
\roster
\item "{$(a)$}"  $\le$ is a transitive on $K^{qr}_\lambda$;
\sn
\item "{$(b)$}"  any increasing continuous sequence in 
$K^{qr}_\lambda$ of length $< \lambda^+$ has a m.u.b.; (see part (10))
(not necessarily unique)
\sn
\item "{$(c)$}"  $\bold C$ is closed (see part (12)).
\endroster
\medskip
 
\noindent
15) $\bold C$ is almost $\lambda$-nice (as witnessed by $G$) is 
defined similarly, replacing ``closed" by ``almost closed" (witnessed by $G$).
\bigskip

\proclaim{\stag{3.2} Claim}  Let $\bold C$ be a local (or revised local)
$\lambda$-construction framework. \newline
1) If $(\bar M^\ell,\bold f^\ell) \in K^{\text{qr}}_\lambda$ for 
$\ell = 1,2$ and $E$ is a club of $\lambda^+$ and $(*)$ below \ub{then} \nl
$(\bar M^1,\bold f^1) \le (\bar M^2,\bold f^2)$ when
\medskip
\roster
\item "{$(*)$}"  if $\delta \in E$ then $M^1_{\delta + i} = M^2_{\delta+i}$
for $i \le \bold f^1(\delta)$ and $\bold f^1(\delta) \le \bold f^2(\delta)$;
in the ``revised local" version assume in addition that $M^1 = M^2$.
\endroster
\medskip

\noindent
2)  $\le$ is transitive and a reflexive relation on $K^{qr}_\lambda$. \newline
3)  Any increasing continuous sequence of pairs from $K^{qr}_\lambda$
of length $< \lambda^+$ has a minimal upper bound. \newline
4) If in addition $\bold C$ is explicitly local (see Definition 
\scite{3.1}(11))
\underbar{then} any increasing sequence in $K^{qr}_\lambda$ of length
$< \lambda^+$ has a lub. \newline
5) $\bold C$ is $\lambda$-nice (hence in particular, $\lambda$-closed).
\endproclaim
\bigskip

\demo{Proof}  1)  Check clause (b) of Definition \scite{3.1}(7) and
Definition \scite{3.1}(a). \newline
2)  Use clauses (b,c) of Definition \scite{3.1}(7)).  In (c) take
$\zeta^* = 2$ \newline
3)  Without loss of generality the elements of the sequence are 
$(\bar M^\xi,\bold f^\xi) \in K^{\text{qr}}_\lambda$ for \newline
$\xi < \mu$, where $\mu$ is a regular cardinal $\le \lambda$.  
For $\xi < \zeta < \mu$, 
let $E_{\xi,\zeta}$ be a closed unbounded subset of $\lambda$ exemplifying 
Definition \scite{3.1}(4) for 
$(\bar M^\xi,\bold f^\xi) \leqq (\bar M^\zeta,\bold f^\zeta)$. 
First assume $\mu < \lambda$.  Let $E \subseteq \dsize \bigcap_{\xi <
\zeta < \mu} E_{\xi,\zeta} \subset \lambda$ 
be a closed unbounded subset of $\lambda$,
such that: $\alpha \in E \Rightarrow \alpha + (\underset {\xi < \mu} {}\to 
{\text{sup} } \bold f^\xi(\alpha)) + 1 < \text{ Min}(E \backslash 
(\alpha + 1))$.  
Let $E = \{ \alpha_i:i < \lambda\}$ with $\alpha_i$ increasing continuously
with $i$.  Notice that for every $i,\xi < \zeta
< \mu \Rightarrow \bold f^\xi(\alpha_i) \leqq \bold f^\zeta(\alpha_i)$.
Let $E^* = \{i:\alpha_i = i\}$.
\sn
We now define $\bar M = \langle M_j:j < \lambda \rangle$ by defining $M_j$ 
by induction on $j$.  If 
$j = \alpha_j \in E^*$ let $M_j = \dbcu_{\xi < \mu} M^\xi_j$. If
$\alpha < j \le \alpha + \bold f^\xi(\alpha)$ for some 
$\alpha \in E^*$ and some
$\xi < \mu$, let $\xi' = \sup(\xi:\alpha + \bold f^\xi(\alpha) \le j)$ and set
$M_j = \dbcu_{\xi > \xi'} M^\xi_j$.  If $j = \underset\xi {}\to \sup(\alpha
+ \bold f^\xi(\alpha))$ for some $\alpha \in E$ and $j > \alpha 
+ \bold f^\xi(\alpha)$ for each $\xi < \mu$.  Let $M_j = \dbcu_{\beta < j} 
M_\beta$.  Finally, if $j < \lambda$ does not fall under any of the previous 
cases, let $M_j =  \dbcu_{\xi < \mu} M^\xi_{\alpha_j}$.

We claim that $\bar M \in {\bold S}eq_\lambda$.  
One checks that $\bar M$ is
continuous and increasing, the main point being that if $\alpha \in E^*$ and
$\alpha < j_1 < \alpha + f^{\xi_1}(\alpha) \le j_2 < \alpha + f^{\xi_2}
(\alpha)$ with $\xi_1 < \xi_2 < \mu$, then $M^{\xi_1}_{j_1} \le
M^{\xi_2}_{j_1} \le M^{\xi_2}_{j_2}$.  One must also check that $\langle
M_j,M_{j+1} \rangle \in {\bold S}eq_2$ for all $j$.  This 
follows from clause (c) of Definition \scite{3.1}(7).
\sn
Let $\bold f$ be defined by 
$\bold f(\alpha_i) = \sup\{\bold f^\xi(\alpha_i):\xi < \mu\}$ if 
$i \in E^*$ and $\bold f(\alpha_i) = 0$ otherwise.  
Clearly $(\bar M^\xi,\bold f^\xi) \le (\bar M,\bold f)$ for $\xi < \mu$.

What about being a $\le$-m.u.b.?  Assume that 
$(\bar M',\bold f') \in K^{qr}_\lambda$ and $\xi < \mu
\Rightarrow (\bar M^\xi,\bold f^\xi) \le (\bar M',\bold f')$.  So for each
$\xi < \mu$ some club $E'_\xi$ of $\lambda$ exemplifies Definition
\scite{3.1}(4), and let $E' =: \dsize \bigcap_{\xi < \mu} E'_\xi \cap E^*$,
a club of $\lambda$.

Now for $\delta \in E'$ we have $(\forall \xi < \mu)(\bold f^\xi(\delta) \le
\bold f'(\delta))$ hence $\bold f(\delta) = 
\underset {\xi < \mu} {}\to \sup \bold f^\xi(\delta) \le \bold f'(\delta)$, 
so $\delta \in E' \Rightarrow \bold f(\delta) \le \bold f'(\delta)$.
Similarly $\delta \in E' \and j \le \bold f(i) \Rightarrow M_{\delta +j}
\le_{\frak K} M'_{\delta +j}$.  So clearly $(\bar M,\bold f) \le
(\bar M',\bold f')$ and $(\bar M,\bold f)$ is 
a minimal u.b. (see Definition \scite{3.1}(10)). \newline
4)  As in the proof of part (3), let $\langle(\bar M^\xi,\bold f^\xi):\xi <
\mu \rangle$ be as therein and let $(\bar M,\bold f)$ be constructed as above.
For proving it is a lub, let $(M^\xi,\bold f^\xi) \le (\bar M',\bold f')$ 
for $\xi < \mu$, and define $E'$ as there.  For $\delta \in E',j < \bold f
(\delta)$ we have $\langle (M^\xi_{\delta +j},M^\xi_{\delta + j+1}):\xi \in 
(\xi_{\delta,j+1},\mu) \rangle$ is $\le^*_{\{[0,1]\}}$-increasing continuous 
and $\langle M^\xi_{\delta +j},M^\xi_{\delta +j+1} \rangle \le^*_{\{[0,1]\}}
(M'_{\delta +j},M'_{\delta +j+1})$ 
for $\xi \in (\xi_{\delta,j+1},\mu)$, so
as $\bold C$ is explicitly local by clause (d) in Definition \scite{3.1}(11) 
we have \newline
$\langle M_{\delta +j},M_{\delta +j+1} \rangle = 
\left < \dsize \bigcup_{\xi \in (\xi_{\delta,j+1},\mu)} M^\xi_{\delta +j},
\dsize \bigcup_{\xi \in (\xi_{\delta,j+1},\mu)} M^\xi_{\delta + j+1}
\right> \le_{\{[0,1]\}} (M'_{\delta+j},M'_{\delta+j+1})$ is as 
required. \newline
The proof for the case $\mu = \lambda$ is similar, using diagonal
intersection. \nl
5)  Left to the reader. \hfill$\square_{\scite{3.2}}$
\enddemo
\bigskip

\centerline {$* \qquad * \qquad *$}
\bigskip

It may clarify matters 
if we introduce some natural cases of $\bold C$.  We shall use
the forthcoming $\bold C^0$ in our construction of many models in
${\frak K}_{\lambda^{+2}}$.
\definition{\stag{3.3} Definition}  For $\ell \in \{0,1,2\}$ and
$\lambda = \text{ cf}(\lambda) > LS({\frak K})$, let
$\bold C = \bold C^\ell_{{\frak K},\lambda}$ consist of
\medskip
\roster
\item "{$(a)$}"  $\tau^+ = \tau$, \newline
${\frak K}^+ = \{M \in {\frak K}_{< \lambda}:\text{ if } \lambda 
\text{ is a successor cardinal then } \|M\|^+ = \lambda\}$ \newline
(with $=^M$ being equality)
\sn
\item "{$(b)$}"  ${\bold Seq}_\alpha = \{ \bar M:\bar M =
\langle M_i:i < \alpha \rangle$ is a $\le_{\frak K}$-increasing continuous \nl
sequence of members of ${\frak K}_{< \lambda} \text{ and if } \alpha =
\lambda \text{ then } \dbcu_{i < \alpha} M_i \text{ has cardinality }
\lambda\}$
\sn
\item "{$(c)$}"  $\bar M <^*_t \bar N$ when:
{\roster
\itemitem{ $(\alpha)$ }  $\bar M = \langle M_i:i < \alpha^* \rangle,\bar N =
\langle N_i:i < \beta^* \rangle$ are from ${\bold Seq}$
\sn
\itemitem{ $(\beta)$ }  if $[\gamma_1,\gamma_2] \in t$ then: \newline

$\qquad (i) \quad \,\, \gamma \in [\gamma_1,\gamma_2] \Rightarrow 
M_\gamma \le_{\frak K} N_\gamma$ \newline

$\qquad (ii) \quad \,\,$ if $\ell =1$, then in addition 
$\gamma \in [\gamma_1,\gamma_2) \Rightarrow M_\gamma = M_{\gamma +1} \cap 
N_\gamma$ \nl

$\qquad (iii) \quad$ if $\ell = 2$, then in addition \nl

$\qquad \qquad \qquad \gamma \in [\gamma_1,\gamma_2) \Rightarrow M_\gamma =
M \cap N_\gamma$ where $M = \dbcu_{i < \alpha^*} M_i$.
\endroster}
\endroster
\enddefinition
\bn
Note that $\bold C^1_{{\frak K},\lambda}$ is interesting when we have disjoint
amalgamation in the appropriate cases.
\demo{\stag{3.3A} Fact} \nl
1)  If $\ell = 0$ or $1$, then 
$\bold C^\ell_{{\frak K},\lambda}$ is an explicitly local 
$\lambda$-construction framework, (hence $\lambda$-nice by 
\scite{3.2}(5)) and ${\frak K}^+$ satisfies axioms I-VI. \nl
2)  If $\ell =0$ or $2$, then $\bold C^\ell_{{\frak K},\lambda}$ 
is an explicitly revised local $\lambda$-construction framework 
(hence $\lambda$-nice by \scite{3.2}(5)) and ${\frak K}^+$ satisfies 
axioms I-VI.
\enddemo
\bigskip

\definition{\stag{3.4} Definition}  If 
$\lambda = \text{ cf}(\lambda) > LS({\frak K})$ then
$\bold C^3_{{\frak K},\lambda}$ consists of
\medskip
\roster
\item "{$(a)$}"  $\tau^+ = \tau \cup \{P,<\},{\frak K}^+$ is the set of
$(M,P^M,<^M)$ where $M \in {\frak K}_{< \lambda},P^M \subseteq M,<^M$ 
a linear ordering of $P^M$ (but $=^M$ may be as in \scite{3.0}(2)) and
$M_1 \le_{{\frak K}^+} M_2$ iff $(M_1 \restriction \tau) \le_{\frak K} (M_2
\restriction \tau)$ and $M_1 \subseteq M_2$
\sn
\item "{$(b)$}"  ${\bold Seq}_\alpha = \{\bar M:\bar M = \langle M_i:i \le
\alpha \rangle$ is an increasing continuous sequence of members of
${\frak K}^+$ and $\langle M_i \restriction \tau:i \le \alpha \rangle$ is
$\le_{\frak K}$-increasing, and for \newline
$i < j < \alpha:P^{M_i}$ is a proper initial segment of 
$(P^{M_j},<^{M_j})$ and there is a first element in the difference$\}$ \nl
we denote the $<^{M_{i+1}}$-first element of $P^{M_{i+1}} \backslash P^{M_i}$,
by $a_i[\bar M]$
\sn
\item "{$(d)$}"  $\bar M <^*_t \bar N$ \underbar{iff} \newline
$\bar M = \langle M_i:i < \alpha^* \rangle,\bar N = \langle N_i:
i < \alpha^{**} \rangle$ are from ${\bold Seq},t$ is a set of pairwise 
disjoint closed intervals of $\alpha^*$ and for any $[\alpha,\beta] \in t$ 
we have $(\beta < \alpha^*$ and): 
\sn
$\gamma \in [\alpha,\beta] \Rightarrow M_\gamma \le_{\frak K} N_\gamma \and
a_\gamma[\bar M] \notin N_\gamma$, moreover \newline
$a_\gamma[\bar M] = a_\gamma[\bar N]$.
\endroster
\enddefinition
\bigskip

\demo{\stag{3.4A} Fact}  $\bold C^3_{{\frak K},\lambda}$ is an explicitly
local and revised local $\lambda$-construction framework (hence 
$\lambda$-nice by \scite{3.2}(5)) and ${\frak K}^+$ satisfies axioms I-VI.
\enddemo 
\bn
We now introduce amalgamation choice functions.
The use of ``$F$ a $\lambda$-amalgamation choice function" is to help use the
weak diamond, by taking out most of the freedom in choosing amalgams.
This gives possibilities for coding (\scite{3.8}, \scite{3.10}).
\definition{\stag{3.5} Definition}  1) We say that $F$ is a 
$\lambda$-amalgamation choice function for the construction framework
$\bold C$ \underbar{if} $F$ is a five place function satisfying:
\medskip
\roster
\item "{$(a)$}"  if $M_\ell \in K^+_{< \lambda}$ for $\ell < 3,M_0
\le_{{\frak K}^+} M_1,M_0 \le_{{\frak K}^+} M_2,M_1 \cap M_2 = M_0$ (before
dividing by $=^{M_\ell});a \in M_2$ and $(\forall b \in M_0)[\neg a 
=^{M_1} b]$; and $A$ is a set such that $A \cup |M_1| \cup |M_2|$ is a 
set of ordinals \underbar{then} \newline
$F(M_0,M_1,M_2,A,a)$, if defined, is a member $N$ of ${\frak K}^+$ with 
universe \newline
$A \cup |M_1| \cup |M_2|$, which $\le_{{\frak K}^+}$-extends $M_1$ and 
$M_2$ and $(a/=^N) \notin (M_1/=^N)$;
\sn
\item "{$(b)$}"  [uniqueness] \newline
if $(M_0,M_1,M_2,A,a)$ and $(M'_0,M'_1,M'_2,A',a')$ are as above and in the
domain of $F,f$ is an order preserving mapping from $A \cup |M_1| \cup |M_2|$
onto $A' \cup |M'_1| \cup |M'_2|$ such that $f \restriction M_\ell$ is an
isomorphism from $M_\ell$ onto $M'_\ell$ for $\ell = 0,1,2$ (so preserving
$=^{M_\ell}$, and its negation) and
$f(a) = a'$ \underbar{then} $f$ is an isomorphism from $F(M_0,M_1,M_2,A,a)$
onto $F(M'_0,M'_1,M'_2,A',a')$;
\sn
\item "{$(c)$}"  if $F(x_0,x_1,x_2,x_3,x_4)$ is well defined then
$x_0,x_1,x_2,x_3,x_4$ are as in part (a).
\endroster
\medskip
\noindent
Observe that as $M_1 \cap M_2 = M_0$ in (a), if we do not have disjoint
amalgamation then we are {\it forced} to allow $=_N$ to be a nontrivial
congruence. \nl
2) If $F$ is defined whenever the conditions in part (a) hold
and $A \backslash M_1 \backslash M_2$ has large enough cardinality 
\underbar{then} we say $F$ is full (if $\lambda = \mu^+$, it suffices to
demand $A \backslash M_1 \backslash M_2$ has cardinality $\mu$). \nl
3) We say $F$ has strong uniqueness \ub{if}
\mr
\item "{$(d)$}"  if $(M_0,M_1,M_2,A,a)$ and $(M'_0,M'_1,M'_2,A',a')$ are as
above and in the domain of $F$ and for $\ell = 0,1,2$ we have $R_\ell$ is
a isomorphism relation $R_\ell$ from $M_\ell$ onto $M'_\ell$ such that
$R_0 = R_1 \cap (M_0 \times M'_0) = R_2 \cap (M_0 \times M'_0)$ and $|A| = 
|A'|$, \ub{then} there is an isomorphism relation $R$ from $M = F(M_0,M_1,
M_2,A,a)$ onto $M' = F(M'_0,M'_1,M'_2,A',a)$ such that $R_\ell = R \cap
(M_\ell \times M'_\ell)$.
\endroster
\enddefinition
\bigskip

\definition{\stag{3.6} Definition}  Assume $\bold C$ is a 
$\lambda$-construction framework and $F$ is a $\lambda$-amalgamation choice
function for $\bold C$.  Let $(\bar M^\ell,\bold f^\ell) \in 
K^{\text{qr}}_\lambda$ for $\ell = 1,2$. \newline
1) $(\bar M^1,\bold f^1) <^{\text{at}}_{F,a}(\bar M^2,\bold f^2)$ (if we
omit $a$, this means for some $a$; ``at" stands for {\it atomic} extension;
we may write $\le^{\text{at}}_{F,a}$ instead of $<^{\text{at}}_{F,a}$) 
means that:
\medskip
\roster
\item "{$(a)$}"  $(\bar M^1,\bold f^1) \le (\bar M^2,\bold f^2)$
\sn
\item "{$(b)$}"  for some club $E$ of $\lambda$, for every $\delta \in E$
taking $e_\delta =: [\delta,\delta + \bold f^1(\delta)]$ we have
\smallskip
\noindent
{\roster
\itemitem{ $(*)$ }  if $\beta < \gamma$ are successive 
members of $e_\delta$ then: \newline
$M^2_\gamma = F(M^1_\beta,M^1_\gamma,M^2_\beta,|M^2_\gamma|,a)$
\sn
\itemitem{ $(**)$ }   $|M^2_\gamma| = |M^1_\gamma| \cup |M^2_\beta| \cup
\{i:i$ an ordinal not in $|M^1_\gamma| \cup |M^2_\beta|$ and otp$(|M^2_\gamma|
\cap i \backslash |M^1_\gamma| \cup |M^2_\beta|) < \|M^2_\gamma\|$
\sn
\itemitem{ $(***)$ }  $\|M^2_\gamma\| = \text{ Min}\{\|N\|:N = F(M^1_\beta,
M^1_\delta,M^2_\beta,|N|,a)\}$.
\endroster}
\ermn
A suitable club $E$ may be called a {\it witness} for the relation. 
Implicit in clause (b) in $a \in M^2$ and $\neg(\exists b)(b \in M^1 \and
a =^{M^2} b)$. \nl
2) $(\bar M',\bold f') \le_F (\bar M'',\bold f'')$ means that: there is
a sequence $\langle (\bar M^\zeta,\bold f^\zeta):\zeta \le \xi \rangle$ such
that:
\medskip
\roster
\item "{$(a)$}"  $\xi < \lambda$
\sn
\item "{$(b)$}"  $(\bar M^\zeta,\bold f^\zeta) \in K^{\text{qr}}_\lambda$ is
$\le$-increasing continuous in $\zeta$ \newline
(remember Definition \scite{3.1}(10))
\sn
\item "{$(c)$}"  for each $\zeta < \xi$ we have $(\bar M^\zeta,\bold f^\zeta)
\le^{\text{at}}_F (\bar M^{\zeta + 1},\bold f^{\zeta +1})$
\sn
\item "{$(d)$}"  $(\bar M',\bold f') = (\bar M,\bold f^0)$ and
$(\bar M'',\bold f'') = (\bar M^\xi,\bold f^\xi)$.
\ermn
A club $E$ which witnesses all the relations in (c) is called a witness for
the relation $\le_F$. \nl
3) $<^{\text{at},*}_{F,a}$ is defined similarly to part (1) but we demand
in clause (b) only that $e_\delta \subseteq [\delta,\delta + \bold f^2
(\delta)]$ is closed and
$\{\delta,\delta + \bold f^2(\delta)\} \subseteq e_\delta$; the 
requirement from clause (a) is unchanged and we require also:

$$
\text{ if } \beta \in [\delta,\delta + \bold f^2(\delta)] \text{ then }
M^2_\beta = M^2_{\text{max}(e_\delta \cap (\beta + 1))}.
$$
\medskip

\noindent
Then define $\le^*_F$ by iterating $\le^{at,*}_F$. \newline
4) We may replace $F$ by $\bold F$, a family of such functions.  Then in 
each case in \scite{3.6}(2)(c) we use one such that $F$'s.
$\bold F^*$ is the family of all such $F$'s. \newline
5) $(\bar M^1,\bold f^1) <_{F,a} (\bar M^2,\bold f^2)$ means that for some
$(\bar M,\bold f)$ we have $(\bar M^1,\bold f^1) \le^{at}_{F,a} (\bar M,
\bar f) \le_F (\bar M^2,\bold f^2)$. \newline
6) We define mub as in \scite{3.1}(10).
\enddefinition
\bigskip

\remark{\stag{3.6A} Remark}  1) What we prove below on 
$<^{\text{at}}_{F,a},\le_F$ also holds for $<^{\text{at},*}_{F,a},\le^*_F$.
\nl
2) Note: using $\bold F$ rather than $F$ may help in proving cases of 
Definition \scite{3.13}, but we can use one $F$ which codes all members of
$\bold F$  by asking on $A \backslash M_1 \backslash M_2$, though 
artificially. \nl
3) We can replace $F$ by $\langle F_\eta:\eta \in \lambda$ 
a sequence of ordinals of length $< \lambda,\eta(1+i) < 2,\eta(0) <
2^{< \lambda} \rangle$, each $F_\eta$ with uniqueness \scite{3.5}(3) and
$(**)$ of \scite{3.6}(1)(b) is replaced by $M^2_\gamma = 
f_{\eta \restriction \delta}(M^1_\beta,M^1_\gamma,M^2_\beta,|M^2_\gamma|,a)$,
and omit $(***)$ there. 
\endremark
\bigskip

\proclaim{\stag{3.7} Claim}  If $\bold C$ is nice, \ub{then} on 
$K^{\text{qr}}_\lambda,\le_F$ is a quasi-order,
and every increasing continuous sequence of length less than $\lambda^+$ has
a mub.
\endproclaim   
\bigskip

\demo{Proof}  Check.  \hfill$\square_{\scite{3.7}}$
\enddemo
\bigskip

\definition{\stag{3.8} Definition}  1) We say a $\lambda$-amalgamation choice
function $F$ for ${\frak K}^+$ has the $\lambda$-{\it coding property} 
for $\bold C$ if: ${\bold Seq}_\lambda \ne \emptyset$, and for every 
$\bar M^1 \in {\bold Seq}_\lambda$, function
$\bold f^1:\lambda \rightarrow \lambda$, and $S \subseteq \lambda$
we can find $\bar M^{2,\eta} \in {\bold Seq}_\lambda$ for $\eta \in
{}^\lambda 2$ with $\eta$ extending $0_{\lambda \backslash S}$ that is
$\eta \restriction (\lambda \backslash S)$ being constantly zero, a function
$\bold f^2:\lambda \rightarrow \lambda$ such that $\bold f^1 
\le_{{\Cal D}_\lambda} \bold f^2,\bold f^2 \restriction (\lambda \backslash
S) = \bold f^1 \restriction (\lambda \backslash S)$ and an element $a_\eta$
of $M^{2,\eta}$ (usually $a_\eta = a$) such that:
\medskip
\roster
\item "{$(*)_1$}"  $(\bar M^1,\bold f^1) \le^{\text{at}}_{F,a}
(\bar M^{2,\eta},\bold f^1)$ for all $\eta$ extending $0_{\lambda
\backslash S}$, and \nl
$\eta \restriction \alpha = \nu \restriction \alpha \Rightarrow M^{2,\eta}
\restriction \alpha = M^{2,\nu} \restriction \alpha$;
\sn
\item "{$(*)_2$}"  for some club $E$ of $\lambda$ the following is
impossible: for some $\eta_3,\eta_4 \in {}^\lambda 2$ extending
$0_{\lambda \backslash S}$, for $\ell =3,4$ we have 
$(\bar M^1,\bold f^2) \le_{F,a^\ell}(\bar M^\ell,\bold f^\ell)$ witnessed 
by a club $E_\ell$ we have (abusing our notation we are dividing by the
equality congrugence) $f_\ell$ a $\le_{\frak K}$-embedding of 
$M^{2,\eta_\ell} \restriction \tau$ into $M^\ell \restriction \tau$ 
over $M^1_\lambda$, and for some $\delta \in E \cap E_3 \cap E_4 \cap S$ 
we have $\bar M^3 \restriction \delta = \bar M^4 \restriction \delta,
a^3 \in M^{2,\eta_3}_\delta,a^4 \in M^{2,\eta_4}_\delta,f_3(a^3) =
f_4(a^4),\bold f^3 \restriction \delta = 
\bold f^4 \restriction \delta,f_3 \restriction M^{2,\eta_3}_\delta = f_4
\restriction M^{2,\eta_4}_\delta$, and $\ell \in \{3,4\} \Rightarrow 
\text{ Rang}(f_\ell \restriction M^{2,\eta_\ell}) \subseteq M^\ell_\delta,
\eta_3 \restriction \delta = \eta_4 \restriction \delta,\eta_3(\delta) \ne
\eta_4(\delta)$.
\endroster
\medskip

\noindent
2) We say that $F$ has the weak $\lambda$-coding property if above we
restrict ourselves to the cases $\bold f^1 \restriction S = 0_S$.  We can 
even restrict ourselves to
the cases $\bold f^1 \in {\Cal F} \subseteq {}^\lambda \lambda$
provided that $0_\lambda \in {\Cal F}$ and we demand $\bold f^2 \in
{\Cal F}$. \newline
3) If we replace $F$ by $\bold F$, a family of such functions it means we use
Definition \scite{3.6}(4). \nl
4) We say ${\frak K}$ has a coding property if some $\lambda$-amalgamation
choice function $F$ has this property.  Typically the actual choice of $F$ is
irrelevant as long as its domain is sufficiently rich.
\enddefinition
\bigskip

\demo{\stag{3.8A} Observation}  The restriction above to $\eta$ such that
$\eta$ extends $0_{\lambda \backslash S},S$ is natural but inessential, 
as we can extend the definition
of $M^{2,\eta}$ to all $\eta$ in ${}^\lambda 2$ by defining $M^{2,\eta} =
M^{2,\eta'}$ where $\eta' \restriction S = \eta \restriction S$ and $\eta'
\restriction (\lambda \backslash S)$ is constantly zero.  Then the
same properties will hold.
\enddemo
\bigskip

\remark{\stag{3.9} Remark} 

1) For a local construction framework $\bold C$ in \scite{3.8}(1) the 
conditions $(*)_1$ and $(*)_2$ can be replaced by local requirements.  
For example, in condition $(*)_1$, we may take \scite{3.1}(7b) into account.
\nl
2) In \scite{3.8}$(*)_2$ a sufficient condition for the impossibility of the
stated conditions on $E,\eta^3,\eta^4$ and $\delta$, where $\eta^3 
\restriction \delta = \eta^4 \restriction \delta = \eta$, say, is that there
is $\bar a \in M^{2,\eta}_\delta$ so that:
\mr
\item "{$(*)_3$}"  tp$(\bar a,M^1_{\delta + \bold f^2(\delta)},
M^{2,\eta^3}) \ne \text{ tp}(\bar a,M^1_{\delta + \bold f^2(\delta)},
M^{2,\eta^4})$.
\ermn
We can even allow $\bar a$ to be infinite here, say a full listing of
$M^{2,\eta}_\delta$.

To see that this suffices, suppose that we also have the conditions of
\scite{3.8}$(*)_2$.  Then for $\ell =3,4$ as $M^{2,\eta_\ell}_{\delta +
j+1}$ is given by

$$
F(M^1_{\delta +j},M^1_{\delta +j+1},M^{2,\eta_\ell}_{\delta +j},
|M^{2,\eta_\ell}_{\delta +j+1}|,a^\ell)
$$
\mn
where $j < \bold f^2(\delta)$, we find $M^{2,\eta_3}_{\delta + \bold f^2
(\delta)} = M^{2,\eta_4}_{\delta + \bold f^2(\delta)} = M^*$, say. \nl
Since $M^{2,\eta_3}$ and $M^{2,\eta_4}$ can be amalgamated over $M^*$, we
have

$$
\text{tp}(f_3(\bar a),M^1_{\delta + \bold f^2(\delta)},M^{3,\eta_3}) =
\text{tp}(f_4(\bar a),M^1_{\delta + \bold f^2(\delta)},M^{4,\eta_4}).
$$
\mn
On the other hand by $(*)_2$ we have

$$
\text{tp}(\bar a,M^1_{\delta + \bold f^2(\delta)},M^{2,\eta_\ell}) =
\text{tp}(f_\ell(\bar a),M^1_{\delta + \bold f^2(\delta)},M^{\ell,\eta_\ell})
$$
\mn
and this gives a contradiction. \nl
3) To understand Definitions \scite{3.8}(1,2) you may look at the places
they are verified, such as \scite{3.16A} and \scite{3.16C}.  Also see 
\scite{3.14}(3,4). \nl
4) The next lemma deduces from the criterion in \scite{3.8}(1) another one 
which is natural for use in a non-structure theorem. \newline
5) Note that \scite{3.8}(1) implies: for every $(\bar M^1,\bold f^1)$ there
is $(\bar M^2,\bold f^2)$ such that \newline
$(\bar M^1,\bold f^1) \le^{at}_{F,a} (\bar M^2,\bold f^2)$ and 
$M^1_\lambda \ne M^2_\lambda$. \newline
6) In Definition \scite{3.8}(1) we can replace $<^{at}_{F,a}$ by $<_F$ or
$<_{F,a}$ with no harm as $<_{F,a}$ satisfies the requirement on
$<^{at}_{F,a}$ and starting from it we again get $<_F$. \newline
7) In $(*)_2$ of Definition \scite{3.8}(1) for some function $H$ depending 
on $(\bar M^1,f^1)$ we may add the further restriction
$\ell \in \{3,4\}$ and $\alpha < \delta$ implies $\delta \in H(\bar M^\ell
\restriction \alpha,f_\ell \restriction M^\ell_\alpha)$.  I.e. this weakening
of the demand does not change the desired conclusions. \nl
8) We can weaken the demand in $(*)_2$ of \scite{3.8}(1) to 
extensions which actually
arise but this seems more cumbersome. While the adaptation is straightforward,
we have no application in mind. \nl
9) In \scite{3.8}(1), $(*)_2$ we may strengthen the requirement by excluding
the case where the club $E$ is allowed to depend on $\eta$.  That is, we
consider quadruples $(E^\eta,\bar M^\eta,\bold f^\eta,f^\eta)$ for
$\eta \in {}^\lambda 2$ such that $(\bar M^1,\bold f^2) \le_F (\bar M^\eta,
\bold f^\eta)$ is witnessed by a club $E^\eta$ in $\lambda$ and $f^\eta$ is
a $\le_{\frak K}$-embedding of $M^{2,\eta} \restriction \tau$ into
$M^\eta \restriction \tau$ over $M^1$.  We require:
\mr
\item "{$(*)'_2$}"  for no $\eta^3,\eta^4 \in {}^\lambda 2$ and \nl
$\delta \in E^{\eta^3} \cap E^{\eta^4} \cap E \cap S$ do we have:

$$
\eta_3 \restriction \delta = \eta_4 \restriction \delta,\eta_3(\delta) \ne
\eta_4(\delta),\bold f^3 \restriction \delta = \bold f^4 \restriction
\delta;
$$

$$
\bar M^{2,\eta_3} \restriction [\delta,\delta + \bold f^2(\delta)] =
\bar M^{2,\eta_4} \restriction [\delta,\delta + \bold f^3(\delta)];
$$

$$
f^{\eta_3} \restriction M^{2,\eta_3}_\delta = f^{\eta_4} \restriction
M^{2,\eta_4}_\delta;
$$

$$
f^{\eta_\ell}[M^{2,\eta_\ell}_\delta] \subseteq M^{\eta_\ell}_\delta
\text{ for } \ell = 3,4.
$$
\ermn
10) In \scite{3.8} we can also require the models $M^1,M^{2,\eta}$ to have
universes $\lambda(1 + \alpha)$ and $\lambda(1 + \alpha +1)$ respectively
for some $\alpha$, with $\lambda(1 + \alpha) \in M^{2,\eta}_0$.  This will 
not change much.
\endremark
\bigskip

\proclaim{\stag{3.10} Lemma}  Assume $(\exists \mu)(\lambda = \mu^+ 
\and 2^\mu < 2^{\mu^+})$ or at least the definitional weak diamond 
on $\lambda$ holds.  Assume $\bold C$ is $\lambda$-nice, $J =
\text{WDmId}^{\text{Def}}(\lambda)$. \newline
If the $\lambda$-amalgamation choice function $F$ has the $\lambda$-coding 
property, \ub{then} it has the explicit $(\lambda,J)$-coding property, which 
means: if $(\bar M^1,\bold {\bar f}^1) \in K^{\text{qs}}_\lambda$ and 
$S \subseteq \lambda$ satisfies $S \notin J$ \underbar{then} we can find
$(\bar M^2,\bold f^2) \in K^{\text{qs}}_\lambda$ such that:
\medskip
\roster
\item "{$(a)$}"  $(\bar M^1,\bold f^1) \le^{\text{at}}_F 
(\bar M^2,\bold f^1)$ and
$\bold f^1 \restriction (\lambda \backslash S) = \bold f^2 \restriction
(\lambda \backslash S)$
\sn
\item "{$(b)$}"  if $(\bar M^1,\bold f^2) \le_F (\bar M^3,\bold f^3) \in 
K^{\text{qs}}_\lambda$ \underbar{then} $M^2 \restriction \tau$ cannot be 
$\le_{\frak K}$-embedded into $M^3 \restriction \tau$ over $M^1$.
\endroster
\endproclaim
\bigskip

\demo{Proof}  The proof is straight forward once 
you digest the meaning of weak diamond.

Let $S = \lambda$.  
Suppose $\langle \bar M^{2,\eta}:\eta \in {}^\lambda 2 \rangle$ and
$\bold f^2$ are as in \scite{3.8}(1), taking note of \scite{3.8A}.  Then we
claim there is $\nu \in {}^\lambda 2$ for which \scite{3.10} holds on taking
$\bar M^2$ to be $\bar M^{2,\nu}$.  Assume toward a contradiction that this
fails for each $\nu$.  Then clause (b) fails, and for each $\nu \in
{}^\lambda 2$ we have some $\bar M^{3,\nu}$ and $\bold f^{3,\nu}$ with:

$$
(\bar M^1,\bold f^2) \le_F (\bar M^{3,\nu},\bold f^{3,\nu})
\text{ witnessed by a club } E^\nu;
$$

$$
f_\nu:M^{2,\nu} \restriction \tau \rightarrow M^{3,\nu} \restriction \tau
\text{ over } M^1 \restriction \tau \text{ a } \le_{\frak K}
\text{-embedding}.
$$
\mn
Now by the definition of WDmId$^{\text{Def}}(\lambda)$ we can find
$\eta_3 \ne \eta_4$ in ${}^\lambda 2$ and $\delta \in E^{\eta_3} \cap
E^{\eta_4} \cap E^* \cap S$ with $E^* = \{\alpha < \lambda$ limit: $\beta <
\alpha$ implies $\beta + \bold f^2(\delta),\beta + \bold f^1(\delta) <
\alpha\}$, as forbidden in $(*)_2$ of \scite{3.8}(1). 
\hfill$\square_{\scite{3.10}}$
\enddemo
\bigskip

\noindent
Now we can give a reasonable non-structure theorem.
\proclaim{\stag{3.11} Theorem}  Assume $\bold C$ is 
$\lambda$-nice, $(\exists \mu)(\lambda = \mu^+ \and 2^\mu < 2^{\mu^+}) \and
2^\lambda < 2^{\lambda^+}$, or at least DWD$(\lambda)$, and
DWD$^+(\lambda^+)$.  Let $J = \text{WDmId}^{\text{Def}}(\lambda)$. 

If $F$ has the $(\lambda,J)$-coding property, then 
$I(\lambda^+,{\frak K}) \ge 2^{\lambda^+}$; also if 
$2^{\lambda^+} > (2^\lambda)^+$
then IE$(\lambda^+,{\frak K}) \ge 2^{\lambda^+}$.
\endproclaim
\bigskip

\demo{Proof}  We choose by induction on $\alpha < \lambda^+$, for every
$\eta \in {}^\alpha 2$, a pair $(\bar M^\eta,\bold f^\eta)$ such that:
\medskip
\roster
\item "{$(a)$}"  $(\bar M^\eta,\bold f^\eta) \in K^{\text{qs}}_\lambda$
\sn
\item "{$(b)$}"  if $\nu \triangleleft \eta$ then $(\bar M^\nu,\bold f^\nu)
\le_F (\bar M^\eta,\bold {\bar f}^\eta)$
\sn
\item "{$(c)$}"  $(\bar M^\eta,\bold f^\eta) 
\le^{\text{at}}_{F,a_{\eta \char 94
\langle 0 \rangle}}(\bar M^{\eta \char 94 \langle 0 \rangle},
\bold f^{\eta \char 94 \langle 0 \rangle})$
\sn
\item "{$(d)$}"  $(\bar M^\eta,\bold{\bar f}^\eta) 
\le^{\text{at}}_{F,a_{\eta \char 94 \langle 1 \rangle}}(\bar M^{\eta 
\char 94 \langle 1 \rangle},\bold f^{\eta \char 94 \langle 1 \rangle})$
\sn
\item "{$(e)$}"  if $\ell g(\eta)$ is a limit ordinal and
$(\bar M^\eta,\bar f^{\eta \char 94 \langle 0 \rangle})
\le_F (\bar M',\bold f')$ \underbar{then} \newline
$M^{\eta \char 94 \langle 1
\rangle} \restriction \tau$ cannot be $\le_{\frak K}$-embedded into 
$M' \restriction \tau$ over $M^\eta$
\sn
\item "{$(f)$}"  if $\alpha$ is limit ordinal, then 
$(\bar M^\eta,\bold f^\eta)$ is a $\le_F$-mub of $\langle 
(\bar M^{\eta \restriction i},\bold f^{\eta \restriction i}):i < \alpha 
\rangle$.
\endroster
\medskip

\noindent
For $\alpha = 0$ note that as ${\bold Seq}_\lambda \ne \emptyset$ also
${\bold Seq}_\lambda^s \ne \emptyset$ hence $K^{\text{qs}}_\lambda 
\ne \emptyset$.
For $\alpha$ limit use \scite{3.3}, \scite{3.5}.  
For $\alpha = \beta + 1,\beta$ a limit ordinal and $\nu \in
{}^\beta 2$, define $(\bar M^{\nu \char 94 \langle \ell \rangle},
\bold f^{\nu \char 94 \langle \ell \rangle})$ for $\ell = 0,1$ by 
\scite{3.10}.  If $\alpha = \beta + 1,\beta$ non-limit use \scite{3.9}(3).
Let $M^\eta = \dsize \bigcup_{\alpha < \lambda^+} 
M^{\eta \restriction \alpha} \restriction \tau$ for $\eta \in 
{}^{(\lambda^+)} 2$.   Now note $\{a_{\eta \restriction
(i+1)}/=^{M^\eta}:i < \lambda^+\} \subseteq M^\eta/=^{M^\eta}$ are
pairwise distinct so $M^\eta \in K_{\lambda^+}$.  Now we
can apply \scite{1.4} (with $\lambda^+$ here standing for $\lambda$
there). \hfill$\square_{\scite{3.11}}$ 
\enddemo
\bigskip

\centerline {$* \qquad * \qquad *$}
\bigskip

\noindent
Unfortunately, in some interesting cases we get only weak coding.
\proclaim{\stag{3.12} Theorem}  Assume $\bold C$ is $\lambda$-nice,
$(\exists \mu)(\lambda = \mu^+$ and $2^\mu < 2^\lambda < 2^{\lambda^+}$ and 
WDmId$(\lambda^+)$ is not $\lambda^+$-saturated (or at 
least DfWD$(\lambda)$, DfWD$^+(\lambda^+)$ and 
WDmId$^{\text{Def}}(\lambda)$ is not $\lambda^+$-saturated)).
\newline
If $F$ has the weak $\lambda$-coding property (see Definition \scite{3.8}(2)),
or at least the parallel of the conclusion of \scite{3.10}, 
\underbar{then} $I(\lambda^+,K) \ge 2^{\lambda^+}$.
\endproclaim
\bigskip

\demo{Proof}  We can find $\langle S^*_\alpha:\alpha < \lambda^+ \rangle$
such that:

$$
S^*_\alpha \subseteq \lambda
$$

$$
[\alpha < \beta \Rightarrow |S^*_\alpha
\backslash S^*_\beta| < \lambda]
$$
\medskip

\noindent
and
$$
S'_\alpha =: S^*_{\alpha+1} \backslash S^*_\alpha \notin \text{ WDmId}
^{\text{def}}(\lambda).
$$
\enddemo
\bigskip

\noindent
We again choose by induction on $\alpha < \lambda^+$ for every $\eta \in
{}^\alpha 2$ a pair $(\bar M^\eta,\bold f^\eta)$ such that:
\medskip
\roster
\item "{$(a)$}"  $(\bar M^\eta,\bold f^\eta) \in K^{\text{qs}}_\lambda$
\sn
\item "{$(b)$}"  if $\nu \triangleleft \eta$, \ub{then} 
$(\bar M^\nu,\bold f^\nu) \le_F (\bar M^\eta,\bold f^\eta)$
\sn
\item "{$(c)$}"  $f^\eta \restriction (\lambda \backslash S_{\ell g(\eta)})
= 0_{\lambda \backslash S_{\ell g(\eta)}}$
\sn
\item "{$(d)$}"  $(\bar M^\eta,\bold f^\eta) <^{\text{at}}_{F,a_{\eta
\char 94 \langle 0 \rangle}}(\bar M^{\eta \char 94 \langle 0 \rangle},
\bold f^{\eta \char 94 \langle 0 \rangle})$
\sn
\item "{$(e)$}"  $(\bar M^\eta,\bold f^\eta) <^{\text{at}}
_{F,a_{\eta \char 94 \langle 1 \rangle}}(\bar M^{\eta \char 94 \langle 1
\rangle},\bold f^\eta)$
\sn
\item "{$(f)$}"  if $(\bar M^\eta,\bold f^{\eta \char 94 \langle 0 \rangle})
\le_F (\bar M',\bold f')$ \underbar{then} $M^{\eta \char 94 \langle 1 \rangle}
\restriction \tau$ cannot be $\le_{\frak K}$-embedded \newline
into $M' \restriction \tau$ over $M^\eta$
\sn
\item "{$(g)$}"  if $\alpha$ is a limit ordinal, \ub{then} 
$(\bar M^\eta,\bold f^\eta)$ is a $\le_F$-mub of $\langle (\bar M^{\eta
\restriction i},\bold f^{\eta \restriction i}):i < \alpha \rangle$.
\endroster
\medskip

\noindent
Again there are no problems (the difference is in clause (c)). \newline
Then we apply \scite{1.4}(1) (or \scite{1.6A}(1)). 
\hfill$\square_{\scite{3.12}}$
\bigskip

\definition{\stag{3.13} Definition}  
1) We say a $\lambda$-amalgamation choice function $F$ for $\bold C$
has the weak \footnote{there is no clear relation between ``weak local"
and ``local" $\lambda$-coding inspite of the name} (or weakly) local
$\lambda$-coding property for $\bold C$ if: 
\mr
\item "{$(*)_1$}"  Assume $\langle M_0,M_1 \rangle \in {\bold Seq}_2,M_0
\le_{{\frak K}^+} N_0 \in {\frak K}^+_{< \lambda},M_1 \cup N_0 \subseteq
\lambda^+$, and $|M_1| \cap |N_0| = |M_0|$ and 
$a \in N_0,a/=^{N_0} \notin (M_0/ 
=^{N_0})$  (i.e. $(\forall b \in M_0)(\neg a =^{N_0} b)$). \newline
\underbar{Then} we can find $N^1,N^2 \in {\frak K}^+_{< \lambda}$ such that:
{\roster
\itemitem{ (a) }  $N^1 = F(M_0,M_1,N_0,|N^1|,a)$
\sn
\itemitem{ (b) }  $\langle N_0,N^1 \rangle \in {\bold Seq}_2$ and
$\langle N_0,N^2 \rangle \in {\bold Seq}_2$
\sn
\itemitem{ (c) }  $M_1 \le_{{\frak K}^+} N^1$ and $M_1 \le_{{\frak K}^+} 
N^2$
\sn
\itemitem{ (d) }  $N^1 \restriction \tau,N^2 \restriction \tau$ are 
contradicting \footnote{Of course, we may consider only ones in ``legal"
extensions.  We can also note that for the intended use, the disjointness
is automatic (so \scite{3.0}(2) not needed)}
amalgamations of $M_1 \restriction \tau,N_0 \restriction \tau$ over $M_0
\restriction \tau$; i.e. for no $N',h$ do we have: $(N^1 \restriction \tau)
\le_{\frak K} N' \in {\frak K}_{< \lambda}$ and $h$ is a 
$\le_{\frak K}$-embedding of $N^2$ into $N'$ over $M_1 \cup N_0$; or at least
\sn
\itemitem { (d)$^-$ }  $(N^1,N^2)$ is a $\tau$-contradicting pair of
amalgamations of $M_1,N_0$ over $M_0$ which just says: if $N^1
\le_{{\frak K}^+} N \in {\frak K}^1_{< \lambda}$ then there is no 
$\le_{\frak K}$-embedding $h$ of $N^2 \restriction \tau$ into $N 
\restriction \tau$ over $M_1 \cup N_0$ (i.e. is the identity on $M_1$ and 
on $N_0$)
\endroster}
\endroster
\medskip

\noindent
(note: this is not necessarily symmetric; and we use just the $\tau$-reducts
of $N^2,M_0$, \newline
$M_1,N_0$ so we can replace them by \newline
$N^2 \restriction \tau,M_0 \restriction \tau,M_1 \restriction \tau,N_0
\restriction \tau$ respectively). \newline
2) We say that a $\lambda$-amalgamation choice function $F$ for $\bold C$
has the {\it local} $\lambda$-coding property if: 
\medskip
\roster
\item "{$(*)_2$}"  if $\bar M = \langle M_j:j < \lambda \rangle \in 
{\bold Seq}_\lambda,\bar N = \langle N_j:j \le \delta + i \rangle \in
{\bold Seq}_{\delta +i},a \in N_0$, \nl
$(a/=^{N_{\delta +i}}) \notin M_{\delta +i}/=^{M_{\delta +i}}$ and \nl
$\bar M \restriction (\delta +i+1)\le^*_{\{[\delta,\delta +i]\}} \bar N$, and
\nl
$N_{\delta +j+1} = F(M_{\delta +j},M_{\delta +j+1},N_{\delta +j},
|N_{\delta +j+1}|,a)$ for $j < i$ \underbar{then} for some $i_1,i_2 
\in (i,\lambda)$ and
$\bar N^\ell = \langle N^\ell_\alpha:\alpha < \delta + i_\ell \rangle
\in {\bold Seq}_{\delta +i_\ell}$ for $\ell =1,2$ we have:
{\roster
\itemitem{ (a) } $\bar N^\ell \restriction (\delta+i+1) = \bar N$ for
$\ell = 1,2$
\sn
\itemitem{ (b) } for $j \in [\delta +i,\delta +i_1)$ we have \newline
$N^1_{j+1} = F(M_j,M_{j+1},N^1_j,|N^1_{j+1}|,a)$
\sn
\itemitem{ (c) }  $M_{\delta +i_2} \le_{{\frak K}^+} N^2_{\delta +i_2}$
\sn
\itemitem{ (d) }  $N^1_{\delta +i_1} \restriction \tau,N^2_{\delta +i_2}
\restriction \tau$ are contradictory amalgamations of \newline
$M_{\delta +i_1} \restriction \tau$ and $N_\delta \restriction 
\tau$ over $M_\delta \restriction \tau$ or at least
\sn
\itemitem{ (d)$^-$ }  $(N^1_{\delta + i_1},N^2_{\delta + i_2})$ are
$\tau$-contradictory amalgamations of $M_{\delta + i(*)}$ and $N_\delta$ over
$M_\delta$ where $i(*) = \text{ Min}\{i_1,i_2\}$.
\endroster}
\ermn
(So if $i=0$, this gives us a possibility to amalgamate, helpful for \nl
$i \in \lambda \backslash \dbcu_{\delta \in E} [\delta,\delta + \bold f
(\delta)))$.
\medskip

\noindent
3) We say that a $\lambda$-amalgamation choice function $F$ for $\bold C$
has the {\it weaker} local $\lambda$-coding property for
$\bold C$ if:
\medskip
\roster
\item "{$(*)_3$}"  as in part (2) but $i=0$.
\endroster
\medskip

\noindent
4) In 1), 2), 3) above we say $\bold C$ has weak local or the local or 
the weaker local coding
property respectively, if we omit the mention of $F$ meaning for some $F$
(clause (a) in $(*)_1$, clause (b) in $(*)_2$.)
\enddefinition
\bigskip

\proclaim{\stag{3.14} Claim}  1) If $F$ has the weak local $\lambda$-coding
property for $\bold C$ or $F$ has the local $\lambda$-coding property for
$\bold C$ \underbar{then} $F$ has the weaker local $\lambda$-coding property
for $\bold C$. \newline
2) Assume
\medskip
\roster
\item "{$(a)$}"  $\bold C$ has a local $\lambda$-construction framework
\sn
\item "{$(b)$}"  $F$ has the $\lambda$-coding (or weaker $\lambda$-coding)
(or the weak $\lambda$-coding) property for \nl
$\bold C$.
\endroster
\medskip

\noindent
\ub{Then} for some $F'$ we have:
\roster
\item "{$(\alpha)$}"  $F'$, too, is a $\lambda$-amalgamation choice function
for $\bold C$
\sn
\item "{$(\beta)$}"  if $F(N_0,N_1,N_2,A,a)$ is well defined and its
$\tau$-reduct is $<_{\frak K} M \in K_{< \lambda}$ and $A \subseteq A'
\subseteq \lambda^+,|A'| < \lambda$, \ub{then} for some $A'',A' \subseteq A''
\subseteq \lambda^+,|A''| < \lambda$, and $F'(N_0,N_1,N_2,A',a)$
is well defined and $F(N_0,N_1,N_2,A,a) \le_{{\frak K}^+} 
F'(N_0,N_1,N_2,A'',a)$ and $M \le_{\frak K} F'(N_0,N_1,N_2,A'',a) \restriction
\tau$
\sn
\item "{$(\gamma)$}"  $F'$ has the local or weak local or weaker local
\newline
(respectively as in (b)) $\lambda$-coding property for $\bold C$. 
\endroster
\medskip

\noindent
3) If $\bold C$ is local ($\lambda$-construction framework), $F$ a 
$\lambda$-amalgamation choice function, with the weak (or just weaker) local 
$\lambda$-coding property, \underbar{then} $F$ has the weak $\lambda$-coding 
property (hence under the set theoretic assumptions of \scite{3.12}, 
$I(\lambda^+,{\frak K}) \ge 2^{\lambda^+}$). \newline
4)  If $\bold C$ is local ($\lambda$-construction framework), $F$ a 
$\lambda$-amalgamation choice function, with the local 
$\lambda$-coding property, \underbar{then} $F$ has the $\lambda$-coding 
property (hence under the set theoretic assumptions of \scite{3.11}, we have
$I(\lambda^+,{\frak K}) \ge 2^{\lambda^+}$ and $(2^\lambda)^+ <
2^{\lambda^+} \Rightarrow IE(\lambda^+,{\frak K}) \ge 2^{\lambda^+}$).
\endproclaim
\bigskip

\remark{Remark}  The parallel of part (2) holds for local and weaker local
property \underbar{if} $F$ acts on sequences.  See \scite{3.17A} below.
\endremark
\bigskip

\demo{Proof}  1) Check. \newline
2) Here we use clause $(d)^-$ rather than $(d)$. \newline
3) Similar (or read the proof of \scite{3.16A}). \newline
4) Check. \hfill$\square_{\scite{3.14}}$
\enddemo
\bn
In this context we may consider
\definition{\stag{3.17A} Definition}  We say that $F$ is a 
$\lambda$-amalgamation choice function for sequences, for $\bold C$ if:
\medskip
\roster
\item "{$(a)$}"  if $x = F(x_1,x_2,x_3,x_4,x_5,x_6)$ is defined then for
some $\alpha_1,\alpha_2,\alpha_3,\alpha < \lambda$ we have
$x_\ell = \bar M^\ell \in \bold S eq^s_{\alpha_\ell}$ for $\ell < 3,
\bar M^1 \triangleleft \bar M^2,t = x_4$ is a set of pairwise disjoint
intervals $\subseteq \alpha_i,\bar M^1 <^*_t \bar M^3,A = x_5$ a set of
$< \lambda$ ordinals $< \lambda^+$, \newline
$x = \bar M \in \bold S eq^s_\alpha,
M_\alpha$ has universe $A,\bar M^1 \le^*_{t \cup \{[\alpha_1,\alpha_2]\}}
\bar M$
\sn
\item "{$(b)$}"  [uniqueness] \newline
as in Definition \scite{3.5}.
\endroster
\enddefinition
\bigskip

\proclaim{\stag{3.17B} Claim}  Assume:
\medskip
\roster
\item "{$(a)$}"  $(\exists \mu < \lambda)(2^\mu = 2^{< \mu} < 2^\lambda),
2^\lambda < 2^{\lambda^+}$ or at least we have the definable weak diamond 
for $\lambda$ and $\lambda^+$; and $2^{< \lambda} \ge \lambda^+$.
\sn
\item "{$(b)$}"  $\bold C$ is a nice ($\lambda$-construction framework).
\sn
\item "{$(c)$}"  $F$ is a $\lambda$-amalgamation choice function for
$\bold C$.
\sn
\item "{$(d)$}"  $F$ has the weaker $\lambda$-coding property for $\lambda$.
\endroster
\medskip

\noindent
\underbar{Then} $I(\lambda^+,{\frak K}) \ge \mu_{\text{wd}}(\lambda^+)$.
\endproclaim
\bigskip

\demo{Proof}  Straight forward using \scite{1.3}. \newline

Let $\langle M_i:i < i^* \rangle$ list the models in ${\frak K}_{\lambda^+}$
up to isomorphisms and assume toward contradiction that
$i^* < \mu_{\text{wd}}(\lambda^+)$.  It is enough
to choose by induction on $\alpha < \lambda^+$ a sequence
$\bar M^\eta_0$ for $\eta \in {}^\alpha 2$ such that:
\medskip
\roster
\item "{$(a)$}"  $\bar M^\eta \in \bold S eq_\lambda$
\sn
\item "{$(b)$}"  $M^\eta_\lambda$ has universe $\gamma_\eta < \lambda^+$
\sn
\item "{$(c)$}"  $\nu \triangleleft \eta \Rightarrow \bar M^\nu \le_F
\bar M^\eta$
\sn
\item "{$(d)$}"  if $\delta = \ell g(\eta)$ is a limit ordinal then 
$\bar M^\eta$ is a mub of $\langle \bar M^{\eta \restriction \alpha}:\alpha
< \delta \rangle$
\sn
\item "{$(e)$}"  if $\eta \char 94 \langle 0 \rangle \trianglelefteq \nu
\in {}^{\lambda^+} 2$ then $M^{\eta \char 94 \langle 1 \rangle}_\lambda 
\restriction \tau$ cannot be $\le_{\frak K}$-embedded into \newline
$M^\nu_\lambda \restriction \tau$ over $M^\eta_\lambda \restriction \tau$.
\endroster
\medskip

\noindent
This is possible by \scite{3.14}(3). \newline
Having the $M_\eta,\eta \in {}^{\lambda >} 2$ we get the conclusion by
\scite{1.3A}.  \hfill$\square_{\scite{3.17B}}$
\enddemo
\bigskip

\remark{\stag{3.13A} Remark}  1) If we are just interested in $I(\lambda^+,
K)$ rather than also in $IE(\lambda^+,K)$, then we can change the definition
of $\tau$-contradictory to:
\medskip
\roster
\item "{${}$}"  $N^1,N^2$ are $\tau$-contradictory amalgamations of
$M_1,N$ over $M_0$ \ub{if} 
$M_0 \le_{{\frak K}^+} M_1 \le_{{\frak K}^+} N^\ell,
M_0 \le_{{\frak K}^+} N_0 \le_{{\frak K}^+} N^\ell,M_0 = M_1 \cap N_0$ and
there are no $N^1_*,N^2_* \in {\frak K}^+_{< \lambda}$ such that: $N^\ell
\le_{{\frak K}^+} N^\ell_*$ and $N^1_* \restriction \tau,N^2_* \restriction
\tau$ are isomorphic over $M_1 \cup N_0$.
\endroster
\medskip

\noindent
2) Note that it is unreasonable to assume that we will always use the local
versions: e.g. if we have a superlimit model in ${\frak K}_{\lambda^+}$
and we want to have $\bar M \in \bold S eq_\lambda \Rightarrow M_\lambda$
superlimit, we have to add some global condition 
(see \cite{Sh:600}).  Also possibly we
will have in the construction (i.e. in \scite{3.11} or \scite{3.12}) that
$\bar M_\eta$ has two ``contradictory" extensions 
$\bar M_{\eta \char 94 \langle 0 \rangle},\bar M_{\eta \char 94 \langle 1 
\rangle}$, (see clauses (e) and (f) in their proof) only when 
cf$(\ell g(\eta)) = \theta$, where $S = \{ \delta:\delta < \lambda^+$ and 
cf$(\delta) = \theta\} \notin \text{ WDmId}(\lambda)$; or even 
$\ell g(\eta) \in S$, for a given $S \in \text{ WDmId}(\lambda)^+$.  
We shall deal with such cases when needed.
\endremark
\bigskip

\remark{Remark}  We intend to continue this elsewhere.
\endremark
\bigskip

\proclaim{\stag{3.16A} Lemma}  Let ${\frak K}$ be an abstract elementary class
with LS$({\frak K}) \le \lambda$ which is categorical in $\lambda$ and in
$\lambda^+$, with $1 \le I(\lambda^{++},K) < 2^{\lambda^{++}}$. Assume that
$2^\lambda < 2^{\lambda^+} < 2^{\lambda^{++}}$, or at least that the 
definitional weak diamond holds for both $\lambda^+$ and $\lambda^{++}$.

If there is a model in $K_{\lambda^+}$ which is saturated over $\lambda$, 
\ub{then} the minimal triples are dense in $K^3_\lambda$.
\endproclaim
\bigskip

\demo{Proof}  Let $\bold C$ be $\bold C^0_{{\frak K},\lambda^+}$ (see
Definition \scite{3.3}) hence $\bold C$ is explicitly local 
$\lambda$-construction framework (by \scite{3.3A}(1)).  Suppose toward 
contradiction that above $(M^*,N^*,a^*) \in 
K^3_\lambda$, there is no minimal triple.  We claim in this case that there
is a $\lambda^+$-amalgamation choice function $F$ for $\bold C$ with the
$\lambda^+$-coding property, with domain the quadruples $(M_0,M_1,M_2,A,b)$
such that: $(M^*,N^*,a)$ embeds in $(M_0,M_2,b)$; $A$ and the universes of
$M_1,M_2$ are contained in $\lambda^{++}$; $|A \backslash (|M_1| \cup |M_2|)|
= \lambda$; and $M_0 \le_{\frak K} M_1,M_2$.  Then applying \scite{3.10} and
\scite{3.11} we get a contradiction.

We first make two observations concerning triples $(M,N,b)$ lying above
$(M^*,N^*,b)$.  Any such triple has the extension property by \scite{2.7} 
(or just \scite{2.6}(1)) and
hence we can manufacture a $\lambda^+$-amalgamation choice function with the
specified domain.  Furthermore, there is $M' \in K_{\lambda^+}$ with
$M \le_{\frak K} M'$ such that tp$(a,M,N)$ has more than one extension to
$M'$, by the failure of minimality.

Let us show that any $\lambda^+$-amalgamation choice function $F$ with the
specified domain has the $\lambda^+$-coding property on $\bold C =
\bold C^0_{{\frak K},\lambda^+}$.

Let $\bar M^1 \in \bold S$eq$_{\lambda^+}$ and let $\bold f^1:\lambda
\rightarrow \lambda$.  For any set $S$ we must find sequences 
$\bar M^{2,\eta}$ (depending on $S$) as in \scite{3.14}(1).  The approach will
be to first build suitable $M^{2,\eta}$ for all $\eta \in {}^\lambda 2$,
independent of $S$, then restrict appropriately given $S$.

Let $M^1 = \dbcu_i M^1_i$.  Then $M^1 \in {\frak K}_{\lambda^+}$ and by our
assumptions $M^1$ is therefore saturated over $\lambda$.  Hence we may
suppose $M^* \le_{\frak K} M^1$.  We may also suppose that the universe of
$M^1$ is an ordinal, and we may choose a subset $A^2$ of $\lambda^{++}$ which
is the union of an increasing continuous sequence $A^2_\alpha$ (for $\alpha <
\lambda^+$) so that: $A^2_\alpha \cap |M^1| = |M^1_\alpha|$ and $A^2_0
\backslash |M^1_0|$ and $A^2_{\alpha +1} \backslash (A^2_\alpha \cup
|M^1_{\alpha +1}))$ have cardinality $\lambda$.  Let $E$ be the club:

$$
\{\delta < \lambda:\text{ for } \alpha < \delta, \text{ we have } \alpha +
\bold f^1(\alpha) + 1 < \delta\}.
$$
\mn
We now define triples $(M^*_\eta,N^*_\eta,a^*)$ for all $\eta \in {}^i 2$,
by induction on $i$, together with ordinals $\alpha_\eta$ satisfying the
following conditions:
\mr
\item "{$(a)$}"  for any $\eta \in {}^\lambda 2$, the sequence
$(M^*_{\eta \restriction j},N^*_{\eta \restriction j},a^*)(j \le i)$ is
increasing and continuous; and similarly the $\alpha_\eta$ are increasing and
continuous.
\sn
\item "{$(b)$}"  $(M^*_{<>},N^*_{<>}) = (M^*,N^*)$.
\sn
\item "{$(c)$}"  $M^*_\eta = M^1_{\alpha_\eta}$ and the universe of
$N^*_\eta$ is $A^2_{\alpha_\eta}$.
\sn
\item "{$(d)$}"  If $\delta \in E,\eta \in {}^\delta 2$ and $\alpha_\eta =
\delta$, then:
{\roster
\itemitem{ $(d1)$ }  for $i < \bold f^1(\delta)$ and $\eta \triangleleft
\nu \in {}^{\delta +i+1}2$, the model $N^*_\nu$ is given by $F$ applied to
amalgamate $M^1_{\delta+i+1}$ and $N^*_{\nu \restriction (\delta+i)}$ over
$M^1_{\delta +i}$, using $A^2_{\delta +i+1}$ and keeping $a^*$ out of
$M^*_{\delta + i+1}$
\sn
\itemitem{ $(d2)$ }  for all $\nu,\nu' \in {}^{\delta + \bold f^1(\delta)+1}2$
extending $\eta$ if $\nu' \ne \nu$, then for some $\beta$
tp$(a^*,M^*_\beta,N^*_\nu) \ne \text{ tp}(a^*,M^*_\beta,N^*_\nu)$
\sn
\itemitem{ $(d3)$ }  for non-zero $i \le \bold f^1(\delta)$ and 
$\nu$ such that $\eta \triangleleft \nu \in {}^{\delta + i}2$ we have 
$M^*_\nu = M^1_{\alpha_\nu}$ and $\alpha_\nu = \alpha_{\ell g(\nu)} = 
\alpha_\eta + i = \alpha_\delta +i$ (so $M^*_\nu = M^1_{\delta +i}$)
\sn
\itemitem{ $(d4)$ }  for non-zero $i \le \bold f^1(\delta)$ and $\nu,\rho \in
{}^{\delta +i}2$ such that $\eta \triangleleft \nu \and \eta \triangleleft
\rho$ we have $N^*_\nu = N^*_\rho$ call it $N^1_\eta$.
\endroster}
\ermn
During carrying the definition the main point is guaranteeing clause $(d)$.
So let $\eta \in {}^\delta 2$ and assume that $\alpha_\eta = \delta \in E^*$.
First we define by induction on $i \le \bold f^1(\delta)$, a model
$N^1_{\eta,i}$ such that $N^1_i$ has universe $A^2_{\delta +i},N^1_{\eta,i}$
is $\le_{\frak K}$-increasing, contradiction, $M^1_{\delta +i} \le_{\frak K}
N^1_{\eta,i}$ and $N^1_{\eta,0} = N^*_\eta$ and $f,i < \bold f^1(\delta)$
then 

$$
N^1_{\eta,i+1} = F(M^1_{\delta +i},M^1_{\delta +i+1},N^1_{\eta,i},
A^2_{\delta +i+1},a^*).
$$
\mn
Next we choose by induction on $i \le \bold f^1(\delta)$ for each $\rho \in
{}^i 2$, an ordinal $\beta_{\eta,\ell} \in [\delta + \bold f^1(\delta),
\lambda^+)$ and model $N^1_{\eta,\rho}$ such that:
\mr
\widestnumber\item{${iii}$}
\item "{$(i)$}" $N^1_{\eta,\rho} \in K_\lambda$ has universe $A^2_{\beta_{
\eta,\rho}}$
\sn
\item "{$(ii)$}"  $M^1_{\beta_{\eta,\rho}} \le_{\frak K} N^1_{\eta,\rho}$
and $(a^* / =^{N^1_{\eta,\rho}}) \notin M^1_{\beta_{\eta,\rho}}/=^{N^1
_{\eta,\rho}}$
\sn
\item "{$(iii)$}"  $\beta_{\eta,<>} = \delta + \bold f^1(\delta)$
\sn
\item "{$(iv)$}"  $\rho_1 \triangleleft \rho_2 \Rightarrow \beta_{\eta,\rho_1}
< \beta_{\eta,\rho_2} \and N^1_{\eta,\rho_1} \le_{\frak K} N^1_{\eta,\rho_2}$
\sn
\item "{$(v)$}"  $i$ limit $\Rightarrow N^1_{\eta,\rho} = \dbcu_{\zeta < i}
N^1_{\eta,\rho \restriction \zeta}$
\sn
\item "{$(vi)$}"  $\beta_{\eta,\rho \char 94 <0>} =
\beta_{\eta \char 94 \rho \char 94 <1>}$ and \nl
tp$(a^*,M^1_{\beta_{\eta,\rho \char 94 <0>}},
N^1_{\eta,\rho \char 94 <0>}) \ne \text{ tp}(a^*,
M^1_{\beta_{\eta,\rho \char 94 <1>}},N^1_{\eta,\rho \char 94 <1>})$.
\ermn
There is no problem to carry the definition.  Now let

$$
\alpha_{\eta \char 94 \rho} = \delta + i = \ell g(\eta \char 94 \rho)
\text{ if } \rho \in {}^i 2 \text{ and } i \le \bold f^1(\delta)
$$

$$
\alpha_{\eta \char 94 \rho} = \beta_{\eta,\rho} \text{ if }
\rho \in {}^{\bold f^1(\delta) +1} 2
$$

$$
M^*_{\eta \char 94 \rho} = M^1_{\alpha_{\eta \char 94 \rho}} \text{ if }
\rho \in {}^i 2,i \le \bold f^1(\delta) +1
$$

$$
N^*_{\eta \char 94 \rho} = N^1_{\eta,i} \text{ if }
\rho \in {}^i 2 \text{ and } i \le \bold f^1(\delta)
$$

$$
N^*_{\eta \char 94 \rho} = N^1_{\eta,\rho} \text{ if }
\rho \in {}^{\bold f^1(\delta) +1} 2.
$$
\mn
Now check.
\mn
Having carried the induction for $\eta \in {}^{\lambda^+}2$ we let
$\bar M^{2,\eta} = \langle N^*_{\eta \restriction \alpha}:\alpha <
\lambda^+ \rangle$ and $\bold f^2 = \bold f^1 + 1$.  We have to check that
the demand in Definition \scite{3.8} holds.
\mn
Note that this is essentially the proof mentioned in \scite{3.9}(9).

By our initial remarks there is little difficulty in carrying out this
induction. \nl
We then set $\bar M^{2,\eta} = \langle N^*_{\eta \restriction i}:i <
\lambda^+ \rangle$ for $\eta \in {}^{\lambda^+} 2$.  Given a set $S \subseteq
\lambda^+$ we consider $\bar M^{2,\eta}$ for $\eta \in {}^{\lambda^+} 2$ 
extending $0_{\lambda \backslash S}$, together with the function 
$\bold f^2$ equal to $\bold f^1 + 1$ on $S$ and to $\bold f^1$ on $S$.  
We claim that the two conditions of Definition \scite{3.14}(1) are met.

The first of these is a condition on the type of construction allowed, and,
of course, it has been obeyed, notably in $(d1)$ above:
\mr
\item "{$(*)_1$}"  $(\bar M^1,\bold f^1) \le^{\text{at}}_{F,a}
(\bar M^{2,\eta},\bold f^1)$; and $M^{2,\eta} \restriction \alpha$ is
determined by $\eta \restriction \alpha$.
\ermn
The second condition referred to a club $E$, which can be the intersection
of the club we have defined above with $\{\delta:\alpha_\delta = \delta\}$.
This condition goes as follows:
\mr
\item "{$(*)_2$}"  it is impossible to find sequences $\eta^3,\eta^4$
(extending $0_{\lambda \backslash S}$), extensions $(\bar M^1,\bold f^2) 
\le_F (\bar M^3,\bold f^3),
(\bar M^4,\bold f^4)$ witnessed by clubs $E^3,E^4$ (i.e. $E^\ell$ is the
intersection of the clubs which witness the atomic relations 
$\le^{\text{at}}_F$ implicit in $\le_F$), and embeddings $f_\ell:
M^{2,\eta^\ell} \rightarrow M^\ell$ ($\ell = 3$ or 4) over $M^1$ such that
for some $\delta \in E \cap E^3 \cap E^4 \cap S$ we have:
{\roster
\itemitem{ $(i)$ }  $\bar M^3 \restriction (\delta + \bold f^2(\delta) +1) =
\bar M^4 \restriction (\delta + \bold f^2(\delta) + 1)$;
\sn
\itemitem{ $(ii)$ }  $\bold f^3 \restriction \delta = \bold f^4 \restriction
\delta$;
\sn
\itemitem{ $(iii)$ }  $\eta^3 
\restriction \delta = \eta^4 \restriction \delta$
(call the restriction $\eta$) and $\eta^3(\delta) \ne \eta^4(\delta)$;
\sn
\itemitem{ $(iv)$ }  $f_3,f_4$ are equal on $M^{2,\eta}_\delta$; and
\sn
\itemitem{ $(v)$ }  for $\ell = 3,4,f_\ell$ maps $M^{2,\eta}_\delta$ into
$M^\ell_\delta$.
\endroster}
\ermn
Suppose on the contrary we have $\eta^3,\eta^4$ (extending
$0_{\lambda \backslash S}$),$(\bar M^1,\bold f^2) \le_F (\bar M^3,\bold f^3)$,
\nl
$(\bar M^4,\bold f^4),E^3,E^4,f_3,f_4$ and $\delta$ as above.  

Let $\hat M = M^3_{\delta + \bold f^2(\delta)}$.  It follows from condition
(i) and the fact that $\delta$ belongs to the witnessing clubs $E^3,E^4$
that $\hat M = M^4_{\delta + \bold f^2(\delta)}$.  Then $f_3,f_4$ provide
embeddings of $N^*_{\eta^3 \restriction (\delta + \bold f^1(\delta) +1)}$
and $N^*_{\eta^3 \restriction (\delta + \bold f^1(\delta)+1)}$ into $\hat M$
which agrees on $N^*_\eta$ (hence on $a^*$) and on $M^1_\eta$.  By
\scite{3.9}(2) we are done.  \hfill$\square_{\scite{3.16A}}$
\enddemo
\bigskip

\proclaim{\stag{3.16B} Lemma}  Let ${\frak K}$ be an abstract elementary
class with LS$({\frak K}) \le \lambda$ which is categorical in $\lambda$ and
in $\lambda^+$, with $1 \le I(\lambda^{++},K) < 2^{\lambda^{++}}$.  Assume
that $2^\lambda < 2^{\lambda^+} < 2^{\lambda^{++}}$, or at least that the
definitional weak diamond holds for both $\lambda^+$ and $\lambda^{++}$.
\nl
Then:
\mr
\item "{$(*)$}"  for any $M \in K_{\lambda^+}$ and any triple $(M^0,N^0,
a^0)$ in $K^3_\lambda$ with $M^0 \le_{\frak K} M$, we can find sequences
$\bar M = \langle M_i:i < \lambda^+ \rangle,\bar N = \langle N_i:i <
\lambda^+ \rangle$ such that
\sn
{\roster
\itemitem{ (a) }   $(M^0,N^0,a^0) = (M_0,N_0,a^0)$;
\sn
\itemitem{ (b) }  $(M_i,N_i,a)$ is increasing and continuous in $K^3_\lambda$;
\sn
\itemitem{ (c) }  the union of the $M_i$ is $M$;
\sn
\itemitem{ (d) }  the set $S(\bar M,\bar N,a)$ of $\delta < \lambda^+$ such
that for some $j > \delta$ for all $i \ge j$
\block {if we have $(M_j,N_j,a) \le_\ell (M_i,N^\ell,a)$ for
$\ell =1,2$ then we can amalgamate $(M_i,N^1,a)$ and $(M_i,N^2,a)$ over
$(M_i,N_\delta,a)$} \endblock
is stationary in $\lambda^+$.
\endroster}
\endroster
\endproclaim
\bigskip

\demo{Proof}  Otherwise, we claim that any full $\lambda^+$-amalgamation
choice function will have the $\lambda^+$-coding property.

Let $\bar M \in \bold S \text{eq}_{\lambda^+},\bold f:\lambda^+ \rightarrow
\lambda^+$ and $S \subseteq \lambda^+$ be given.  Then as $K$ is categorical 
in $\lambda^+$, we may suppose that $M$ is the union of the $M_i$.  If (c)
fails, then (d) fails on a club of $\delta$, providing enough failures of
amalgamation to carry the proof as in \scite{3.6A}.
\hfill$\square_{\scite{3.6A}}$
\enddemo
\bigskip

\proclaim{\stag{3.16C} Lemma}  Let ${\frak K}$ be an abstract elementary class
with LS$({\frak K}) \le \lambda$ which is categorical in $\lambda,\lambda^+$
and $\lambda^{++}$ and with no model in cardinality $\lambda^{+3}$.

Suppose that there is no model in $K_{\lambda^+}$ saturated above $\lambda$
and that:
\mr
\item "{$(*')$}"  for any $M \in K_{\lambda^+}$ and any triple $(M^0,N^0,
a^0)$ in $K^3_\lambda$ with $M^0 \le_{\frak K} M$, we can find sequences
$\bar M = \langle M_i:i < \lambda^+ \rangle,\bar N = \langle N_i:i <
\lambda^+ \rangle$ such that
\sn
{\roster
\itemitem{ (a) }   $(M^0,N^0,a^0) = (M_0,N_0,a^0)$;
\sn
\itemitem{ (b) }  $(M_i,N_i,a)$ is increasing and continuous in $K^3_\lambda$;
\sn
\itemitem{ (c) }  the union of the $M_i$ is $M$;
\sn
\itemitem{ (d) }  the set $S(\bar M,\bar N,a)$ of $\delta < \lambda^+$ such
that for some $j > \delta$ for all $i \ge j$
\block {if we have $(M_j,N_j,a) \le_{h_\ell} (M_i,N^\ell,a)$ for
$\ell =1,2$ then we can amalgamate $(M_i,N^1,a)$ and $(M_i,N^2,a)$ over
$(M_i,N_\delta,a)$} \endblock
is stationary in $\lambda^+$.
\endroster}
\ermn
\ub{Then} the minimal triples are dense in $K^3_\lambda$.
\endproclaim
\bigskip

\demo{Proof}  Suppose that there is no minimal triple above $(M^*,N^*,a^*)$.
It suffices to show that there is no maximal model in $K_{\lambda^{++}}$ and
as $K_{\lambda^{++}}$ is categorical, this will follow from the existence of
a single pair of models $(M',N')$ in $K_{\lambda^{++}}$ with $M' <_{\frak K}
N'$.  So it suffices to show:
\mn
\centerline {every triple $(M,N,a)$ in $K^3_{\lambda^+}$ has a proper
extension in $K^3_{\lambda^+}$} 
\mn
as the desired pair $(M',N')$ can then be built as the limit of an
increasing continuous chain.

Fix $(M,N,a)$ in $K^3_{\lambda^+}$.  As there is no model saturated over
$\lambda$ in $K_{\lambda^+}$, there is some $M_0$ in $K_\lambda$ over which
there are more than $\lambda^+$ types.  By $\lambda$-categoricity we may
suppose $M_0 \le_{\frak K} M$.  Fix a triple $(M_0,N_0,b)$ in $K^3_\lambda$
for which tp$(b,M_0,N_0)$ is not realized in $M$.

Apply $(*')$ to $M$ and $(M_0,N_0,b)$ to get sequences $\bar M,\bar N$ of
length $\lambda^+$ as in $(*')$.  Let $S = S(\bar M,\bar N,b)$.  In
particular $M = \bigcup M_i$.  Let $M^1 = \bigcup N_i$.  As ${\frak K}$ has
amalgamation in $\lambda^+$ we may suppose $M^1,N \le_{\frak K} N^1$ with
$N^1 \in {\frak K}_{\lambda^+}$.

We can also choose an increasing continuous sequence $(M'_i,N'_i,a^*)$
beginning with $(M^*,N^*,a^*)$ such that each $(M'_i,N'_i,a^*)$ is reduced
and tp$(a^*,M'_i,N'_i)$ has more than one extension in ${\Cal S}(M'_{i+1})$, 
using the failure of minimality and \scite{2.5}(1).  By categoricity we may 
suppose $M = \bigcup M'_i$.  Set $M^2 = \bigcup N'_i$.  By amalgamation we 
may suppose $M^2,N^1 \le_{\frak K} N^2 \in K_{\lambda^+}$.

We claim that one of the triples $(M^1,N^1,a)$ or $(M^2,N^2,a)$ is a proper
extension of $(M,N,a)$.  Suppose on the contrary that $a$ belongs to both
$M^1$ and $M^2$.

Represent $N^2$ as the union of a continuous $\le_{\frak K}$-increasing 
chain $\langle N^*_i:i < \lambda^+ \rangle$ of models in $K_\lambda$. \nl
Let $E$ be

$$
\{i < \lambda^+:M_i = M'_i;N^*_i \cap M = M_i;N^*_i \cap M^1 = N_i;N^*_i \cap
M^2 = N'_i\},
$$
\mn
a club in $\lambda^+$.

Fix $\delta \in E \cap S$ such that $a$ is in $N_\delta$ and $N'_\delta$.
We show now that $a^* \in N_\delta$.  Now $(M'_\delta,N'_\delta,a^*)$ is
reduced.  If $a^* \notin N_\delta$ then $(N_\delta,N^2,a^*)$ lies over
$(M'_\delta,N'_\delta,a^*)$ and hence $N_\delta \cap N'_\delta \subseteq
M'_\delta$; but the element $a$ witnesses the failure of this condition.
So $a^* \in N_\delta$.

Let $j > \delta$ be chosen in accordance with the definition of $S(\bar M,
\bar N,b)$ and let $i > j' > j_a,k,j' \in E$.  As tp$(a^*,M'_{j'},N'_{j'})$
has more than one extension to $M'_{j' +1}$, the same applies to $M'_i$.
However, $M'_i = M_i$ and $M'_{j'} = M_{j'}$ and thus tp$(a^*,M_{j'},
N^2_{j'}) = \text{ tp}(a^*,M_{j'},N^*_{j'}) = \text{ tp}(a^*,M_{j'},N_{j'})$
has more than one extension over $M_i$.  Thus, $M_i$ and $N_{j'}$ may be
amalgamated in two incompatible ways over $M_{j'}$, getting $N^+$ and $N^-$,
say (i.e. the $N^+$ and $N^-$ cannot be amalgamated over $M_i$ preserving
the images of $a^*$).  Furthermore $(M_i,N^\pm,b)$ lies above $(M_j,N_j,b)$
in $K^3_{\lambda^+}$, that is, $b$ is not mapped into $M_i$, because $M$
does not realize tp$(b,M_0,N_0)$.  However, this contradicts the definition
of $S$, as the triples $(M_i,N^+,b)$ and $(M_i,N^-,b)$; cannot be amalgamated
over $(M_i,N_\delta)$ since $a^*$ belongs to $N_\delta$.
\hfill$\square_{\scite{3.16C}}$
\enddemo
\bigskip

\proclaim{\stag{3.27} Theorem}  Let ${\frak K}$ be an abstract elementary
class with LS$({\frak K}) \le \lambda$ which is categorical in $\lambda$
and in $\lambda^+$ with $1 \le I(\lambda^{++},K) < 2^{\lambda^{++}}$. Assume
that $2^\lambda < 2^{\lambda^+} < 2^{\lambda^{++}}$, or at least that the
definitional weak diamond holds for both $\lambda^+$ and $\lambda^{++}$.

Then under either of the following assumptions, the minimal triples are
dense in $K^3_\lambda$:
\mr
\item "{$(A)$}"  $K$ is categorical in $\lambda^{++}$ and has no model in
cardinality $\lambda^{+3}$;
\sn
\item "{$(B)$}"  there is a model saturated above $\lambda$ in cardinality
$\lambda^+$.
\endroster
\endproclaim
\bigskip

\demo{Proof}  By the previous lemmas \scite{3.16A}, \scite{3.16B},
\scite{3.16C}.  Note that $(*')$ is exactly the negation of $(*)$.
\hfill$\square_{\scite{3.27}}$
\enddemo
\bigskip

\remark{\stag{3.28} Remark}  This will be proved without the additional
assumptions $(A,B)$ in \cite{Sh:600}.  In any case this does not affect the
proof of Theorem \scite{0.1}
\endremark
\bigskip

\proclaim{\stag{3.29} Claim}  Let ${\frak K}$ be an abstract elementary class
with LS$({\frak K}) \le \lambda$ which is categorical in $\lambda$ and in
$\lambda^+$, with $1 \le I(\lambda^{++},K) < 2^{\lambda^{++}}$, and with no
model in cardinality $\lambda^{+3}$.  Assume that $2^\lambda < 2^{\lambda^+}
< 2^{\lambda^{++}}$.

\ub{Then} the minimal triples are dense in $K^3_\lambda$.
\endproclaim
\bigskip

\demo{Proof}  If $2^{\lambda^+} > \lambda^{++}$ we get the conclusion 
by \scite{2.5}. If
$2^{\lambda^+} = \lambda^{++}$ then as $2^\lambda < 2^{\lambda^+}$ we have
$2^\lambda = \lambda^+$.  Thus, there is a model in $K_{\lambda^+}$ which is
saturated above $\lambda$, and Lemma \scite{3.27} applies.
\hfill$\square_{\scite{3.29}}$
\enddemo
\newpage

\head {\S4 Minimal types} \endhead  \resetall
\bn
We return to the analysis of minimal types initiated in \scite{2.8}. \newline
We use from \S2 only \scite{2.1}, \scite{2.4}, \scite{2.6} 
so there are repetitions.
\bigskip

\demo{\stag{4.0} Hypothesis}
\medskip
\roster
\item "{$(a)$}"   ${\frak K}$ is an abstract elementary class with
$LS({\frak K}) \le \lambda$ (for simplicity $K_{< \lambda} = \emptyset$).
\sn
\item "{$(b)$}"  ${\frak K}$ is categorical in $\lambda,\lambda^+$
with $K_{\lambda^{+2}} \ne \emptyset$ (note: $(a) + (b) = (*)^3_\lambda$ of
\scite{2.3}).  
\sn
\item "{$(c)$}"  ${\frak K}$ has amalgamation in $\lambda$ (\scite{2.2}(1)),
\nl
so by $(a) + (c)$, we have $(*)^2_\lambda$ from \scite{2.6} hence
${\frak K}$ satisfies the model theoretic properties 
which were deduced in \scite{2.3}-\scite{2.4},\scite{2.6} in particular:
{\roster
\itemitem{ (i) }   every $(M,N,a) \in K^3_\lambda$ has the weak extension
property \scite{2.3};
\sn
\itemitem{ (ii) }  criteria for the extension property \scite{2.6};
\sn
\itemitem{ (iii) }  basic Definitions and properties 
\scite{2.3A}, \scite{2.4}.
\endroster}  
\endroster
\enddemo
\bigskip

\definition{\stag{4.1A} Definition}  1) If 
$p \in {\Cal S}(N),N \in K_\lambda$ and $N' \in K_\lambda$ remember
(from Definition \scite{2.9})

$$
{\Cal S}_p(N') = \{f(p):f \text{ is an isomorphism from } N \text{ onto }
N'\}
$$
\medskip

\noindent
and let

$$
\align
{\Cal S}_{\ge p}(N') = \biggl\{ q \in {\Cal S}(N'):&\,q \text{ not algebraic 
(i.e. not realized by any} \\
  &\,c \in N') \text{ and, for some } N'' \in K_\lambda, \\
  &\,N'' \le_{\frak K} N', \text{ we have }
q \restriction N'' \in {\Cal S}_p(N'') \biggr\}.
\endalign
$$
\medskip

\noindent
2) We say the type $p \in {\Cal S}(N)$ is $\lambda$-algebraic if 
$\|N\| \le \lambda$ and for every
$M$ such that $N \le_{\frak K} M$ we have: 
$\lambda \ge |\{c \in M:\text{tp}(c,N,M) = p\}|$. 
\enddefinition
\bigskip

\proclaim{\stag{4.1} Claim}  If $(M_0,M_1,a) \in K^3_\lambda$ is minimal, 
\underbar{then} it has the extension property.
\endproclaim
\bigskip

\demo{Proof}  Let $p^* = \text{ tp}(a,M_0,M_1)$, and assume it is a
counter-example.   We note:
\medskip
\roster
\item "{$\bigotimes_1$}"  for some $M^*$ we have $M_0 \le_{\frak K} 
M^* \in K_\lambda$ but for no $M^+$ and $b$ do we have $M^* \le_{\frak K} 
M^+ \in K_\lambda,b \in M^+ \backslash
M^*$ and $b$ realizes $p^*$. \newline
[Why?  If not for every $N,M_0 \le_{\frak K} N \in K_\lambda$, we can
find $N_1$, \newline
$N \le_{\frak K} N_1 \in K_\lambda$ and $b \in N_1 \backslash N$ which
realizes $p^*$.  Hence (as ${\frak K}_\lambda$ has amalgamation in $\lambda$)
we can find $N_2$ such that $N_1 \le_{\frak K} N_2 \in 
K_\lambda$, and $g$ a $\le_{\frak K}$-embedding of $M_1$ into $N_2$ 
extending id$_{M_0}$ such that $g(a) = b$.
\newline
This proves the extension property].
\medskip
\item "{$\bigotimes_2$}"  if $p \in {\Cal S}_{\ge p^*}(N)$ and 
$N \in K_\lambda$ and $N \le_{\frak K} N^* \in K$, \underbar{then} the 
set of elements of $b \in N^*$ 
realizing $p$ has cardinality $\le \lambda$. \newline
[Why?  by \scite{4.0}(c)(ii); so indirectly \scite{2.6}(2)].
\medskip
\item "{$\bigotimes_3$}"  if $N \in K_\lambda$, then $|{\Cal S}_{\ge p^*}(N)| 
> \lambda^+$.
\endroster
\enddemo
\bigskip

\demo{Proof of $\bigotimes_3$}  If $N$ forms a counterexample, as $K$ is
categorical in $\lambda$ and using $\otimes_2$ we can find 
$\langle N_i:i < \lambda^+ \rangle$, $\le_{\frak K}$-increasing continuous 
sequence of members of $K_\lambda$ such that:
\medskip
\roster
\item "{$(*)$}"  for every $\alpha < \lambda^+$ and $q \in {\Cal S}_{\ge p^*}
(N_\alpha)$, for some $\beta = \beta_q < \lambda^+$ we have: for no
$N',b$ do we have $N_\beta \le_{\frak K} N' \in K_\lambda,b \in N' \backslash
N_\beta$ and $b$ realizes $q$.
\endroster
\medskip

\noindent
So $N_{\lambda^+} = \dsize \bigcup_{i < \lambda^+} N_i$ has the property
\medskip
\roster
\item "{$(**)$}"  if $\bar N' = \langle N'_\alpha:\alpha < \lambda^+
\rangle$ is a representation of $N_{\lambda^+}$ then for a club of $\delta <
\lambda^+$ for every $q \in {\Cal S}_{\ge p^*}(N'_\delta)$ for a club of
$\beta \in (\delta,\lambda^+)$, for no $N',b$ do we have: $N_\beta 
\le_{\frak K} N' \in K_\lambda,b \in N' \backslash N_\beta$ and $b$ 
realizes $q$.
\endroster
\enddemo
\bigskip

\noindent
On the other hand, we can choose by induction on $\alpha < \lambda^+$ a
triple \newline
$(N_{0,\alpha},N_{1,\alpha},a) \in K^3_\lambda$ increasing continuous
in $\alpha$ such that \newline
$(N_{00},N_{10},\alpha) = (M_0,M_1,a)$ and $N_{0,\alpha} \ne
N_{0,\alpha + 1}$ (existence by the weak extension property; i.e.
\scite{2.3} = \scite{4.0}(c)(i)).  \newline
Now $N_{0,\lambda^+} = \dsize
\bigcup_{\alpha < \lambda^+} N_{0,\alpha} \in K_{\lambda^+}$ does not
satisfy the statement $(**)$: \newline
$\langle N_{0,\alpha}:\alpha < \lambda^+
\rangle$ is a representation of $N_{0,\lambda^+}$, and for every \newline
$\alpha,\text{tp}(a,N_{0,\alpha},N_{1,\alpha})$ extend
tp$(a,N_{0,0},N_{1,0}) = \text{ tp}(a,M_0,M_1) = p^*$ hence \newline
tp$(a,N_{0,\alpha},N_{1,\alpha}) \in {\Cal S}_{\ge p^*} (N_{0,\alpha})$ 
satisfies: for every $\beta \in (\alpha,\lambda^+)$, there is $N',N_{0,\beta} 
\le_{\frak K} N'$ and some $b \in N' \backslash N_{0,\beta}$ realizes 
$\text{tp}(a,N_{0,\alpha},N_{1,\alpha})$; simply
choose $(N',b) = (N_{1,\beta},a)$.  So $N_{0,\lambda^+},N_{\lambda^+}$
cannot be isomorphic (as one satisfies $(**)$ the other not).  But both are
in $K_{\lambda^+}$, contradicting the categoricity of ${\frak K}$ in
$\lambda^+$. \hfill$\square_{\otimes_3}$
\bigskip

\noindent
To finish the proof of \scite{4.1} it is enough to prove 
\proclaim{\stag{4.2} Claim}  If $p^* \in {\Cal S}(M_0)$ is minimal, 
$M_0 \in K_\lambda$,
\underbar{then} $N \in K_\lambda \Rightarrow |{\Cal S}_{\ge p^*}(N)| 
\le \lambda^+$.
\endproclaim
\bigskip

\demo{Proof}  By \scite{4.2B}, \scite{4.2C} below.  Note that
${\Cal S}_{\ge p^*}(N)$ has the same cardinality for every $N \in K_\lambda$.
\enddemo
\bigskip

\proclaim{\stag{4.2A} Claim}  1) If $N_1 \le_{\frak K} N_2$ are in 
$K_\lambda$ and $p_1 \in {\Cal S}(N_1)$ is minimal and is omitted by $N_2$ 
\underbar{then} $p_1$ has a unique extension in ${\Cal S}(N_2)$, call it 
$p_2$, and \newline
$p_1 \in {\Cal S}_{\ge p^*}(N_1) \Rightarrow [p_2 \in {\Cal S}
_{\ge p^*}(N_2)$ and $p_2$ is minimal]. \newline
2)  If $N_1 \le_{\frak K} N_2$ are in $K_\lambda,p_1 \in {\Cal S}(N_1)$ 
minimal, \underbar{then} $p_1$ has a unique non-algebraic extension in 
${\Cal S}(N_2)$ called $p_2$, it is minimal and $p_1 \in {\Cal S}_{\ge p^*}
(N_1) \Rightarrow p_2 \in {\Cal S}_{\ge p^*}(N_2)$. \newline
3)  (Continuity) if $\langle N_i:i \le \alpha \rangle$ is a
$\le_{\frak K}$-increasing continuous sequence of members of 
$K_\lambda,p_i \in {\Cal S}(N_i),p_0$ minimal, $p_i \in {\Cal S}(N_i)$ 
extends $p_0$ and is non-algebraic then $\langle p_i:i \le \alpha \rangle$ 
is increasing continuously.
\endproclaim
\bigskip

\demo{Proof of \scite{4.2A}}  Easy.  E.g., \newline
3) If $i < j \le \alpha$ then $p_j \restriction N_i$ is well defined, it
belongs to ${\Cal S}(N_i)$, also it is non-algebraic and extending $p_0$ hence
by the uniqueness (=\scite{4.2A}(1)) we have $p_i = p_j \restriction N_i$.
If $\delta \le \alpha,p_\delta \in {\Cal S}(N_\delta)$ extends $p_i$ for
$i < \delta$; if $p'_\delta \in {\Cal S}(N_\delta)$ extends each $p_i(i <
\delta)$ then it extends $p_0$ and is non-algebraic hence by uniqueness
$p'_\delta = p_\delta$. 
\enddemo
\bigskip

\proclaim{\stag{4.2B} Claim}  If $N \in K_\lambda,{\Cal S} \subseteq 
{\Cal S}(N)$ and $|{\Cal S}| > \lambda^+$, \underbar{then} we can find $N^*,
N_i$ in $K_\lambda$, (for $i < \lambda^{++}$) such that:
\medskip
\roster
\item "{$(\alpha)$}"   $N \le_{\frak K} N^* <_{\frak K} N_i$
\sn
\item "{$(\beta)$}"   for no $i_0 < i_1 < \lambda^{++}$ and 
$c_\ell \in N_{i_\ell} \backslash N^*$ (for $\ell = 0,1$) do we have \newline
$\text{tp}(c_0,N^*,N_{i_0}) = \text{ tp}(c_1,N^*,N_{i_1})$
\sn
\item "{$(\gamma)$}"  there are $a_i \in N_i$ (for $i < \lambda^{++}$) 
such that tp$(a_i,N,N_i) \in {\Cal S}$ (and they are pairwise distinct).
\endroster
\endproclaim
\bigskip

\remark{Remark}  We use here less than Hypothesis \scite{4.0}:
\medskip
\roster
\item "{$(*)$}"  ${\frak K}$ is abstract elementary class with amalgamation
in $\lambda$, categorical in $\lambda$, \newline
$K_{\lambda^+} \ne \emptyset$.  
\endroster
\medskip

\noindent
The same applies to \scite{4.2C}.
\endremark
\bigskip

\demo{Proof}  Without loss of generality $|N| = \lambda$; 
now choose by induction on $\alpha < \lambda^{++},\bar N^\alpha,N_\alpha,
a_\alpha$ such that:
\medskip
\roster
\item "{$(A)$}"  $N_\alpha \in K_{\lambda^+}$ has a set of elements  
$\lambda \times (1 + \alpha)$ and $N_\alpha$ is $\le_{\frak K}$-increasing 
continuous in $\alpha$
\smallskip
\noindent
\item "{$(B)$}"  $\bar N^\alpha = \langle N^\alpha_i:i < \lambda^+
\rangle$ is a representation of $N_\alpha$ (i.e. is 
$\le_{\frak K}$-increasing continuous, $\|N^\alpha_i\| \le \lambda$ and
$N_\alpha = \dsize \bigcup_{i < \lambda^+} N^\alpha_i$)
\smallskip
\noindent
\item "{$(C)$}"  for $\alpha < \lambda^{++}$ successor, if $i < j <
\lambda^+,p \in {\Cal S}(N^\alpha_i)$ is realized in $N^\alpha_j$ 
and is $\lambda$-algebraic (see Definition \scite{4.1A}(2)) \underbar{then} 
for no $N',b$ do we have $N^\alpha_j \le_{\frak K} N' \in K_\lambda$ 
and $b \in N' \backslash N^\alpha_j$ realizes $p$ (actually not needed)
\smallskip
\noindent
\item "{$(D)$}"  $N = N^0_0 \le_{\frak K} N_0$ and $a_\alpha \in 
N_{\alpha + 1} \backslash N_\alpha$ realizes some $p_\alpha \in {\Cal S}$ 
not realized in $N_\alpha$.
\smallskip
\noindent
\item "{$(E)$}"  If $\aleph_0 < \text{ cf} (\alpha) \le \lambda$ 
\underbar{then} let for $j < \lambda$: \newline
$M^\alpha_j = \cap\{ \dsize \bigcup_{\beta \in C}
N^\beta_j:C \text{ a club of } \alpha\}$ whenever possible, i.e. the result
is in $K_\lambda$ and $\le_{\frak K} N^\alpha$
\smallskip
\noindent
\item "{$(F)$}"  for each $\alpha < \lambda^{++}$, for a club $E^0_\alpha$
of ordinals $i < \lambda^+$ we have $(N^\alpha_i,N^{\alpha + 1}_i,a_\alpha)$ 
is reduced; \newline
hence (as tp$(a_\alpha,N^\alpha_i,N^{\alpha +1}_i)$ extends
tp$(a_\alpha,N,N^{\alpha +1}_i)$ which is not realized in $N_\alpha$):
\smallskip
\noindent
\item "{$(G)$}"  for every ($i \in E^0_\alpha$ and) 
$b \in N^{\alpha +1}_i \backslash 
N^\alpha_i$ the type $\text{tp}(b,N^\alpha_i,N^{\alpha + 1}_i)$ is 
not realized in $N_\alpha$ \newline
(a key point).
\endroster
\medskip

\noindent
There is no problem to carry out 
the construction (${\frak K}$ has amalgamation in $\lambda$). \newline
Let $w^\alpha_i =: \{ \beta:N^{\alpha + 1}_i \cap N_{\beta + 1}
\nsubseteq N_\beta\}$, so necessarily $|w^\alpha_i| \le \lambda,w^\alpha_i$
is increasing continuous in $i < \lambda^+$ and $\alpha = \dsize
\bigcup_{i < \lambda^+} w^\alpha_i$ and for $\beta < \alpha$ let \newline
$\bold i(\beta,\alpha) = \text{ Min}\{i:\beta \in w^\alpha_i\}$. \newline
Now for every $\alpha \in S^* =: \{ \delta < \lambda^{++}:\text{cf}(\delta)
= \lambda^+\}$, the set 

$$
\align
E_\alpha =: \biggl\{ i < \lambda^+:&\,i \text{ limit }, N \le_{\frak K} 
N^{\alpha + 1}_i,a_\alpha \in N^{\alpha + 1}_i, \\
  &\text{ and for every } \beta < \alpha \text{ if } \\
  &\,\beta \in w^\alpha_i \text{ then} \\
  &\,N^\beta_i = N^\alpha_i \cap N_\beta \text{ and for } j < i \\
  &\text{ the closure of } w^\alpha_j \text{ (in } \alpha) 
\text{ is included in } w^\alpha_i \\
  &\text{ and } \beta_1 < \beta_2 \and \beta_1 \in w^\alpha_j \and \beta_2
\in w^\alpha_j \Rightarrow \bold i(\beta_1,\beta_2) < i \biggl\}
\endalign
$$
is a club of $\lambda^+$.
\medskip

\noindent
As we can assume $\lambda > \aleph_0$ (ignoring $\lambda = \aleph_0$ as 
was treated earlier in \cite{Sh:88} though for a PC$_{\aleph_0}$ class, or
see \cite{Sh:600},\S2), we can choose $j_\alpha \in E_\alpha$ such that
cf$(j_\alpha) = \aleph_1$ and let $\delta_\alpha = \sup(w^\alpha_{j_\alpha})$,
now $w^\alpha_{j_\alpha}$ is closed under $\omega$-limits (as
$\langle w^\alpha_j:j \le \alpha \rangle$ is increasing continuous, $j <
\alpha \Rightarrow \text{ closure}(w^\alpha_j) \subseteq w^\alpha_{j+1})$ 
and $\aleph_1 = \text{ cf(otp } w^\alpha_{j_\alpha})$
so there is $\langle \beta_\varepsilon:\varepsilon < \omega_1 \rangle$ 
increasing continuous with limit $\delta_\alpha,\beta_\varepsilon \in 
w^\alpha_{j_\alpha}$ so $\varepsilon < \zeta < \omega_1 \Rightarrow 
N^{\beta_\varepsilon}_j = N^{\beta_\zeta}_j \cap N_{\beta_\varepsilon}$ 
and easily $N^\alpha_j \cap N_{\beta_\varepsilon} = N^{\beta_\varepsilon}_j$ 
hence
\medskip
\roster
\item "{$\bigoplus$}"  $M^{\delta_\alpha}_{j_\alpha} = \cap
\biggl\{ \dsize \bigcup_{j \in C} N^\beta_j:C \text{ a club of }
\delta_\alpha \biggr\}$ \newline
\smallskip
\noindent
so $N^\alpha_{j_\alpha} = M^{\delta_\alpha}_{j_\alpha}$ (see
\cite[\S4]{Sh:351}).
\endroster
\medskip

\noindent
By Fodor lemma for some $j^*,\alpha^*$ and stationary 
$S \subseteq S^*,\alpha \in S^* \Rightarrow j_\alpha = j^* \and 
\delta_\alpha = \delta^*$.  So for all
$\alpha \in S,N^\alpha_{j_\alpha}$ are the same say $N^*$.  So $N^* \in
K_\lambda$, for $\alpha \in S$, \nl
$q_\alpha = \text{ tp}(a_\alpha,N^*,N^{\alpha+1}_{j_\alpha})$ 
extend $p_\alpha(\in {\Cal S})$.  Also if $r \in {\Cal S}
(N^*)$ is realized in $N^{\alpha +1}_{j_\alpha}$ say by $b$ (for some
$\alpha \in S$) then no member of $\bigcup \{ N^{\beta +1}_{j^*} \backslash
N^\beta_{j^*}:\beta \in S \cap \alpha\}$ realizes it (holds by clause 
$(G)$, see clause $(F)$).  So the sets 
$\Gamma_\alpha = \{\text{tp}(b,N^{\delta^*}_{j^*},N^{\alpha +1}_{j^*}):
b \in N^{\alpha +1}_{j^*} \backslash N^*\}$ are
pairwise disjoint and each has a member extending $p_\alpha \in {\Cal S}$ (as 
exemplified by $a_\alpha$ and $p_\alpha$ is not extended by any
$p \in \dsize \bigcup_{\beta < \alpha} \Gamma_\beta$ (as $p_\alpha$ is not
realized in $N_\alpha$)).  \hfill$\square_{\scite{4.2B}}$
\enddemo
\bigskip

\proclaim{\stag{4.2C} Claim}  Assume $p^*$ is a counterexample to \scite{4.2}.
\newline
1) If $N \in K_\lambda,\Gamma \subseteq {\Cal S}_{\ge p^*}(N),|\Gamma| \le 
\lambda^+$ \underbar{then}

$$
\align
\{ p \in {\Cal S}_{\ge p^*}(N):&\text{ for some } N',N \le_{\frak K} 
N' \in K_\lambda \text{ and some } b \in N' \\
  &\text{ realizes } p \text{ but no } b \in N' \text{ realizes any }
q \in \Gamma \}
\endalign
$$
\medskip

\noindent
has cardinality $\ge \lambda^{++}$. \newline
2) We can find $N \in K_{\lambda^+}$, and $N_i,
N <_{\frak K} N_i \in K_{\lambda^+}$ for $i < 2^{\lambda^+}$ such that the
set $\Gamma_i = \{\text{tp}(a,N,N_i):a \in N_i\}$ are pairwise distinct, in
fact, no one embeddable into another (so we get $I(\lambda^+,{\frak K}) =
2^{\lambda^+}$ and if $(2^\lambda)^+ < 2^{\lambda^+}$ then $IE(\lambda^+,
{\frak K}) = 2^{\lambda^+}$ thus contradicting the categoricity in
$\lambda^+$ from the assumptions).
\endproclaim
\bigskip

\demo{Proof}  1) Apply \scite{4.2B} with $N$ and ${\Cal S}_{\ge p^*}(N)$ 
here standing for $N$ and ${\Cal S}$ there, so we get $N^*,
N_i(i < \lambda^{++})$ and $M^*$ such that $N^* \le_{\frak K} 
N_i \in {\frak K}_\lambda$, and $\Gamma_i = \{\text{tp}(a,N^*,N_i):a \in 
N_i \backslash N^*\}$ are pairwise disjoint and there are $p_i \in 
\Gamma_i,p_i \restriction N \in {\Cal S}_{\ge p^*}(N)$ pairwise 
distinct; now $p_i$ is not algebraic hence $p_i \in {\Cal S}_{\ge p^*}(N^*)$.
As $K$ is categorical in $\lambda$,
without loss of generality $N^* = N$, so all but $\le \lambda^+$ of the models
$N_i$ can serve as the required $N'$. \newline
2) Now by part (1) of \scite{4.2C} we can choose by induction on 
$i < \lambda^+$, \newline
$\langle (N_\eta,\Gamma_\eta):\eta \in {}^i 2 \rangle$ 
such that:
\medskip
\roster
\item "{$(a)$}"  $N_\eta \in K_\lambda$ and $\Gamma_\eta \subseteq
\dsize \bigcup_{j < i} {\Cal S}_{\ge p^*}(N_{\eta \restriction j})$ and
$|\Gamma_\eta| \le \lambda$
\sn
\item "{$(b)$}"  if $\nu \triangleleft \eta$ then $N_\nu \le_{\frak K} 
N_\eta$ and $\Gamma_\eta \subseteq \Gamma_\nu$
\sn
\item "{$(c)$}"  some $p \in \Gamma_{\eta \char 94 \langle 0 \rangle}$ is
from ${\Cal S}(N_\eta)$ and is realized in $N_{\eta \char 94 \langle 1 
\rangle}$ \newline
(similarly for $\Gamma_{\eta \char 94 \langle 1 \rangle},N_{\eta \char 94 
\langle 0 \rangle}$)
\sn
\item "{$(d)$}"  if $i$ is a limit ordinal, then 
$N_\eta = \dsize \bigcup_{j < i} N_{\eta \restriction j}$ and
$\Gamma_\eta = \dsize \bigcup_{j < i} \Gamma
_{\eta \restriction j}$.
\endroster
\medskip

\noindent
The successor case is done by \scite{4.2C}(1) (you may object that the type 
in $\Gamma_\eta$
are not from ${\Cal S}_{\ge p^*}(N_\eta)$ but from $\dsize \bigcup_{j < i}
{\Cal S}_{\ge p^*}(N_{\eta \restriction j})$?  However they are minimal - see
\scite{4.2A}(1)). \newline
For $\eta \in {}^{\lambda^+}2$ let $N_\eta = \dsize \bigcup_{i < \lambda^+}
N_{\eta \restriction i}$.  Now by 1.4(1), $\{ N_\eta / \cong:\eta \in
{}^{\lambda^+}2\}$ has cardinality $2^{\lambda^+}$, contradiction.
\hfill$\square_{\scite{4.2C}},\square_{\scite{4.2}},\square_{\scite{4.1}}$
\enddemo
\bigskip

\proclaim{\stag{4.3} Claim}  1) If $M_0 \in K_\lambda,M_1 \in 
K_{\lambda^+}$, and $M_0 \le_{\frak K} M_1$ \underbar{then} every 
minimal $p \in {\Cal S}(M_0)$ is realized in $M$. \newline
2)  Every $M_1 \in K_{\lambda^+}$ is saturated at least for minimal types
(i.e. if $M_0 \le_{\frak K} M_1$, \newline
$M_0 \in K_\lambda$ and $M_1 \in
K_{\lambda^+}$ \underbar{then} every minimal $p \in {\Cal S}(M_0)$ is 
realized in $M_1$). \newline
3)  If $M \in K_\lambda$ then $\{p \in {\Cal S}(M):p \text{ minimal}\}$ has
cardinality $\le \lambda^+$.
\endproclaim
\bigskip

\demo{Proof}  1) Let $\bar N = \langle N_\alpha:\alpha < \lambda^+ \rangle$
be a representation of $N_{\lambda^+} \in K_{\lambda^+}$.  Let \newline
$N \in K_\lambda,p \in {\Cal S}(N)$ be minimal we ask
\medskip
\roster
\item "{$(*)_p$}"  is there a club of $\alpha < \lambda^+$ such that every
$q \in {\Cal S}_{\ge p}(N_\alpha)$ is realized in $\lambda^+$?
\endroster
\medskip

\noindent
By \scite{4.2} there is $N'_{\lambda^+} \in K_{\lambda^+}$ for 
which the answer is yes, hence,
as ${\frak K}$ is categorical in $\lambda^+$, this holds for $N_{\lambda^+}$.
So this holds for every minimal $p$.  Now if $N' \le_{\frak K} 
N_{\lambda^+},N' \in K_\lambda$ and $p \in {\Cal S}(N')$ is minimal then 
for some $\alpha,N' \le_{\frak K} N_\alpha$ and for every $\beta \in 
[\alpha,\lambda^+)$, 
$p$ has a unique non-algebraic extension $p_\beta \in {\Cal S}(N_\beta)$ 
(which necessarily is minimal).
Now $p_\beta \in {\Cal S}_{\ge p}(N_\beta)$ 
hence for a club of $\beta < \lambda^+,
p_\beta$ is realized in $N_{\lambda^+}$, so we have finished the proof of
part (1). \newline
2)  By part (1).  \newline
3)  Follows by part (1).  \hfill$\square_{\scite{4.3}}$
\enddemo
\bigskip

\noindent
\relax From \scite{4.1}, \scite{2.6}(1) we can conclude
\demo{\stag{4.4} Conclusion}  Every $(M_0,M_1,a) \in K^3_\lambda$ has 
the extension property.
\enddemo
\newpage

\head {\S5 Inevitable types and stability in $\lambda$} \endhead  \resetall
\bigskip

\demo{\stag{5.0} Hypothesis}  Assume the model theoretic assumptions 
from \scite{4.0} and
\medskip
\roster
\item "{$(d)$}"  there is a minimal member of $K^3_\lambda$ (follows from
the conclusion of \scite{3.27}).
\endroster
\enddemo
\bigskip

\definition{\stag{5.1} Definition}  We call $p \in {\Cal S}(N)$ inevitable 
if: \newline
$N \le_{\frak K} M \and N \ne M \Rightarrow \text{ some }
c \in M \text{ realizes } p$. \nl
We call $(M,N,a) \in K^3_\lambda$ inevitable if tp$(a,M,N)$ is inevitable.
\enddefinition 
\bigskip

\noindent
So by \scite{4.1} - \scite{4.3} we shall deduce
\proclaim{\stag{5.2} Claim}  1) If there is a minimal triple in $K^3_\lambda$,
\underbar{then} there is an inevitable \newline
$p = \text{ tp}(a,N,N_1)$ with $(N,N_1,a) \in K^3_\lambda$ minimal. \newline
2) Moreover, if $p_0 \in {\Cal S}(N_0)$ is minimal, $N_0 \in K_\lambda$
\underbar{then} we can find $N_1, N_0 \le_{\frak K} N_1 \in K_\lambda$ 
such that the unique non-algebraic extension $p_1$ of $p_0$ in 
${\Cal S}(N_1)$ is inevitable.
\endproclaim
\bigskip

\demo{Proof of \scite{5.2}}  1) Follows by part (2). \newline
2)  Let $(M_0,M_1,a) \in K^3_\lambda$ be minimal and 
$p_0 = \text{ tp}(a,M_0,M_1)$.  We try to choose by induction
on $i$ a model $N_i$ such that: $N_0 = M_0,N_i \in K_\lambda$ is 
$\le_{\frak K}$-increasing continuously and $N_i$ omits $p_0,N_i \ne N_{i+1}$.
If we succeed, $\dsize \bigcup_{i < \lambda^+} N_i$ is a member of
$K_{\lambda^+}$ which is non-saturated for minimal types, 
contradicting \scite{4.3}(2).
As for $i=0$, $i$ limit we can define, necessarily for some $i$ we have 
$N_i$ but not $N_{i+1}$.  Now $p_0$ has a unique extension 
in ${\Cal S}(N_i)$ which we call $p_i$ and $p_0$ has no 
algebraic extension in ${\Cal S}(N_i)$. \newline
[why?  as $N_i$ omits $p_0$].  So $p_i$ is the unique extension of $p_0$
in ${\Cal S}(N_i)$ [by \scite{4.2A}(1)], and so
\medskip
\roster
\item "{$(*)$}"  if $N_i \le_{\frak K} N' \in K_\lambda$ and $N' \ne N$,
\underbar{then} $p_i$ is realized in $N'$.
\endroster
\medskip

\noindent
By L.S. we can omit ``$N' \in K_\lambda$", so $(N_i,p_i)$ are as required.
\hfill$\square_{\scite{5.2}}$
\enddemo
\bigskip

\demo{\stag{5.2A} Fact}  Inevitable types have few $(\le \lambda)$
conjugates (i.e. for $p \in {\Cal S}(M_0)$ inevitable $M_0 \in K_\lambda,
M_1 \in K_\lambda$ we have $|{\Cal S}_p(M_1)| \le \lambda$), moreover
$|\{p \in {\Cal S}(N):p \text{ inevitable} \}| \le \lambda$
for $N \in K_\lambda$.  
\enddemo
\bigskip

\demo{Proof}  Easy. 
\enddemo
\bn
The following construction plays a central role in what remains of this paper.
\proclaim{\stag{5.3} Claim}  For any limit $\alpha < \lambda^+$, we can find
$\langle N_i:i \le \alpha \rangle$ and $\langle p_i:i \le \alpha \rangle$
such that:
\medskip
\roster
\widestnumber\item{(viii)}
\item "{$(i)$}"  $N_i \in K_\lambda$,
\sn
\item "{$(ii)$}"  $N_i$ is $\le_{\frak K}$-increasing continuous
\sn
\item "{$(iii)$}"  $p_i \in {\Cal S}(N_i)$ is minimal, 
\sn
\item "{$(iv)$}"  $p_i$ increases continuously (see \scite{4.2A}(3))
\sn
\item "{$(v)$}"  $p_0$ is inevitable
\sn
\item "{$(vi)$}"  $p_\alpha$ is inevitable
\sn
\item "{$(vii)$}"   $N_i \ne N_{i+1}$ moreover some 
$c \in N_{i+1} \backslash N_i$ realizes $p_0$ (hence $p_i$).
\endroster
\endproclaim
\bigskip

\remark{Remark}  Why not just build a non-saturated model in order to prove
\scite{5.3}?  Works, too.
\endremark
\bigskip

\demo{Proof}  Choose $N^0 <_{\frak K} N^1$ in $K_{\lambda^+}$ (so
$N^0 \ne N^1$), such a pair exists as 
$K_{\lambda^{+2}} \ne \emptyset$.  Let 
$N^\ell = \dsize \bigcup_{i < \lambda^+} N^\ell_i$ with
$N^\ell_i \in K_\lambda$ being $\le_{\frak K}$-increasing 
continuously in $i$.   Now $E_0 = \{ \delta < \lambda^+:
N^1_\delta \ne N^0_\delta \text{ and } N^1_\delta \cap N^0 = N^0_\delta\}$
is a club of $\lambda^+$. 
Without loss of generality $E_0 = \lambda^+$.
\medskip

For each $c \in N^1 \backslash N^0$, the set

$$
X_c =: \{ i < \lambda^+:c \in N^1_i \text{ and } (N^0_i,N^1_i,c) 
\text{ is minimal}\}
$$
\medskip

\noindent
is empty or an end segment of $\lambda^+$ hence

$$
\align
E_1 = \biggl\{ \delta < \lambda^+:&\text{(i)   } \delta \text{ limit}\\
  &\text{(ii)  } \text{if } i < \delta \text{ and } p \in {\Cal S}(N^0_i) 
\text{ is minimal inevitable} \\
  &\qquad \text{and realized in } N^0 \backslash N^0_\delta 
\text{ then it is}\\
  &\qquad \text{realized in } N^0_\delta \backslash N^0_i \text{ (actually
automatic)} \\
  &\text{(iii)  } \text{if } c \in N^1_\delta \backslash N^0 \text{ (hence }
\exists i < \delta,c \in N^1_i) \text{ and } X_c \\
  &\qquad \text{is non-empty \underbar{then} } \delta \in X_c \text{ and}
\text{ min}(X_c) < \delta \biggr\}
\endalign
$$
\medskip

\noindent
is a club of $\lambda^+$ (see \scite{5.2A}).
\medskip

Now for $\delta \in E_1$, we have $N^0_\delta <_{\frak K} N^1_\delta$, 
so by \scite{5.2}(1) there is $c_\delta \in N^1_\delta \backslash 
N^0_\delta$ such that:

$$
\gather
(N^0_\delta,N^1_\delta,c_\delta) \text{ is minimal}\\
\text{tp}(c_\delta,N^0_\delta,N^1_\delta) \text{ is inevitable}
\endgather
$$
\bigskip

As $\delta$ is limit, for some $i < \delta,c \in N^1_i$, also $\delta \in X_c$
hence there is $j$ such that: $i < j < \delta \and j \in X_c$ hence 
$(N^0_j,N^1_j,c)$ is 
minimal; choose such $j_\delta,c_\delta$.  Let $\kappa = \text{ cf}(\kappa)
= \text{ cf}(\alpha) \le \lambda$, so for some $j^*,c^*$ we have

$$
S = \{ \delta \in E_1:\text{cf}(\delta) = \kappa,j_\delta = j^*,
c_\delta = c^*\}
$$
\medskip

\noindent
is stationary in $\lambda^+$.
\medskip

Choose $e$ closed $\subseteq E_1$ of order type $\alpha + 1$
with first element and last element in $S$; for $\zeta \in [j^*,\lambda^+)$ 
let $p_\zeta = \text{ tp}(c^*,N^0_\zeta,N^1_\zeta)$.  (In fact, we could
have: all non-accumulation member of $e$ are in $S$; no real help.)
\medskip

Now $\langle N^0_\zeta,p_\zeta:\zeta \in e \rangle$ is as required (up to
re-indexing)(clause (viii) holds by clause (ii) in the definition of $E_1$). 
\hfill$\square_{\scite{5.3}}$
\enddemo
\bigskip

\proclaim{\stag{5.4} Claim}  Assume $\langle N_i,p_i:i \le \alpha \rangle$ 
is as in \scite{5.3}, $\alpha < \lambda$ divisible by $\lambda$.  \ub{Then} 
any $p \in {\Cal S}(N_0)$ is realized in $N_\alpha$, moreover $N_\alpha$ is 
universal in $K_\lambda$ over $N_0$.
\endproclaim
\bigskip

\demo{Proof}  (Similar to the proof of \scite{0.15}; which is 
\cite[II,\S3]{Sh:300}).

Let $N_0 \le_{\frak K} M_0 \in K_\lambda,a \in M_0 \backslash N_0$ 
we shall show that $\text{tp}(a,N_0,M_0)$ is realized in $N_\alpha$.

Let $\alpha = \dsize \bigcup_{i < \lambda} S_i, \langle S_i:i < \lambda
\rangle$ pairwise disjoint, each $S_i$ unbounded in $\alpha$,
$\lambda$ divides $\text{otp}(S_i)$ and $\text{Min}(S_i) \ge i$.  
We choose by induction on $i \le \alpha$ the following:

$$
N^1_i,M^1_i,h_i,\langle a_\zeta:\zeta \in S_i \rangle \quad
(\text{the last one only if } i < \alpha)
$$
\medskip

\noindent
such that:
\medskip
\roster
\item "{$(a)$}"  $N^1_i \le_{\frak K} M^1_i$ are in $K_\lambda$
\sn
\item "{$(b)$}"  $N^1_i$ is $\le_{\frak K}$-increasing continuous in $i$
\sn
\item "{$(c)$}"  $M^1_i$ is $\le_{\frak K}$-increasing continuous in $i$
\sn
\item "{$(d)$}"  $(N^1_0,M^1_0) = (N_0,M_0)$
\sn
\item "{$(e)$}"  $\langle a_\zeta:\zeta \in S_i \rangle$ is a list of
$\{ c \in M^1_i:c \text{ realizes } p_0\}$
\sn
\item "{$(f)$}"  $h_i$ is an isomorphism from $N_i$ onto $N^1_i$
\sn
\item "{$(g)$}"  $j < i \Rightarrow h_j \subseteq h_i$ and
$h_0 = \text{ id}_{N_0}$
\sn
\item "{$(h)$}"  $a_i \in N^1_{i+1}$ (note: $M^1_i \cap N^1_{i+1} \ne
N^1_i$ in general).
\endroster
\enddemo
\bigskip

\noindent
\underbar{For $i=0$}:  See clauses $(d),(g)$

$$
N^1_0 = N_0, \, M^1_0 = M_0, \, h_0 = \text{ id}_{N_0}.
$$
\bigskip

\noindent
\underbar{For $i=$ limit}:  Let $N^1_i = \dsize \bigcup_{j < i} N^1_j$ and
$M^1_i = \dsize \bigcup_{j < i} M^1_j$ and $h_i = \dsize 
\bigcup_{j < i} h_j$ and lastly
choose $\langle a_\zeta:\zeta \in S_i \rangle$ by clause $(e)$.
\bigskip

\noindent
\underbar{For $i=j+1$}:  Note $a_j$ is already defined, it belongs to $M^1_j$
and it realizes $p_0$.
\bigskip

\noindent
\underbar{Case 1}:  $a_j \in N^1_j$ (so clause $(h)$ is no problem).
\newline
Use amalgamation on $N_j,N_i,M^1_j$ and the mapping $\text{id}_{N_j},h_i$,i.e.

$$
\CD
N_i             @>>>         M^1_i \\
@Aid_{N_j}AA               @AAid_{N^1_0}A \\
N_j             @>h_i>>        N^1_j
\endCD
$$
\bigskip

\noindent
\underbar{Case 2}:  $a_j \notin N^1_j$.

Then $\text{tp}(a_j,N^1_j,M^1_j)$ is not algebraic, extending the minimal
type $p_0 \in {\Cal S}(N_0)$. \newline
Also by clause $(viii)$ of \scite{5.3} there is 
$c \in N_i \backslash N_j$ which realizes $p_0$.  As 
$p_0 \in {\Cal S}(N)$ is minimal

$$
h_j \left( \text{tp}(c,N_j,N_i) \right) = \text{ tp}(a_j,N^1_j,M^1_j)
$$
\medskip

\noindent
so acting as in Case 1 we can also guarantee $h_i(c) = a_j$, so $a_j \in
\text{ Rang}(h_i) = N^1_i$ as required.
\bigskip

\centerline {$* \qquad * \qquad *$}
\bigskip

In the end we have $N^1_\alpha \le_{\frak K} M^1_\alpha$.  If $N^1_\alpha = 
M^1_\alpha$, then $h^{-1}_\alpha \restriction M_0 = 
h^{-1}_\alpha \restriction M^1_0 = h^{-1}_\alpha \restriction N^1_0$ 
show that $M_0$ can be embedded into $N_\alpha$ over $N_0$ as required.  So
assume $N^1_\alpha <_{\frak K} M^1_\alpha$.  Now $p_\alpha \in {\Cal S}
(N_\alpha)$ is inevitable hence $h_\alpha(p_\alpha) \in {\Cal S}
(N^1_\alpha)$ is inevitable. \newline
Hence some $d \in M^1_\alpha \backslash N^1_\alpha$ realizes $h_\alpha
(p_\alpha)$ hence $d$ realizes $h_\alpha(p_\alpha) \restriction N^1_0 =
p_0$; also $\alpha$ is a limit ordinal so for some $i < \alpha,d \in M^1_i$
hence for some $\zeta \in S_i$ we have $a_\zeta = d$, hence

$$
d = a_\zeta \in N^1_{\zeta + 1} \subseteq N^1_\alpha,
$$
\medskip

\noindent
contradicting the choice of $d$. \newline
So we are done. \hfill$\square_{\scite{5.4}}$
\bigskip

\demo{\stag{5.5} Conclusion}  If $N \in K_\lambda$ then:
\medskip
\roster
\item "{$(a)$}"  $|{\Cal S}(N)| = \lambda$
\sn
\item "{$(b)$}"  there is $N_1, N <_{\frak K} N_1 \in K_\lambda$ such that
$N_1$ is universal over $N$ in $K_\lambda$
\sn
\item "{$(c)$}"  for any regular $\kappa \le \lambda$ we can demand that
$(N_1,c)_{c \in N}$ is $(\lambda,\kappa)$-saturated (see \scite{0.2}(1))
\endroster
\enddemo
\bigskip

\remark{\stag{5.5R} Remark}  In fact amalgamation in $\lambda$ and 
stability in $\lambda$ (i.e. $(a)$ of \scite{5.5}) implies $(b)$ and $(c)$ 
of \scite{5.5}.
\endremark
\bigskip

\demo{\stag{5.5A} Conclusion}  The $N \in K_{\lambda^+}$ is saturated above
$\lambda$ (i.e. over models in $K_\lambda$!).
\enddemo 
\bigskip

\proclaim{\stag{5.6} Claim}  Assume $\kappa = \text{ cf}(\kappa) \le \lambda$.
\newline
There are $N_0,N_1,a,N^+_0,N^+_1$ such that
\medskip
\roster
\widestnumber\item{(iii)}
\item "{$(i)$}"  $(N_0,N_1,a) \in K^3_\lambda$ and
\sn
\item "{$(ii)$}"  $(N_0,N_1,a) \le (N^+_0,N^+_1,a) \in K^3_\lambda$ and 
\sn
\item "{$(iii)$}"  $(N^+_0,c)_{c \in N_0}$ is $(\lambda,\kappa)$-saturated,
\sn
\item "{$(iv)$}"  $\text{tp}(a,N_0,N_1) \text{ is minimal inevitable and}$
\sn
\item "{$(v)$}"  $\text{tp}(a,N^+_0,N^+_1) \text{ is minimal inevitable}$.
\endroster
\endproclaim
\bigskip

\demo{Proof}  As in the proof of \scite{5.3} as

$$
\align
E_2 = \biggl\{ \delta:&\text{ for every } i < \delta,(N^0_\delta,c)_{c \in
N_i} \text{ is saturated}\\
  &\text{ of cofinality } \text{cf}(\delta) \biggr\}
\endalign
$$
\medskip

\noindent
is a club of $\lambda^+$.  \hfill$\square_{\scite{5.6}}$
\enddemo
\bigskip

\proclaim{\stag{5.7} Claim}  1) In $K_\lambda$ we have disjoint amalgamations.
\newline
2) If $M \le_{\frak K} N$ are in $K_\lambda$ and $p \in {\Cal S}(M)$
non-algebraic then for some $N',c$ we have: $N \le_{\frak K} N' \in K_\lambda$
and $c \in N' \backslash N$ realizes $p$.
\endproclaim
\bigskip

\demo{Proof}  1) First note
\medskip
\roster
\item "{$\bigotimes$}"  if $M \le_{\frak K} N$ in $K_\lambda$ we can find 
$\alpha < \lambda^+$,and sequence $\langle M_i:i \le \alpha \rangle$ which is 
$\le_{\frak K}$-increasing continuous, and
$\langle a_i:i < \alpha \rangle$ such that $(M_i,M_{i+1},a_i)$ is minimal
and reduced and $N \le_{\frak K} M_\alpha,M = M_0$.
\newline
[Why?  There is a minimal reduced pair, hence we can find $\langle M_i:i <
\lambda^+ \rangle \le_{\frak K}$-increasing continuous, 
$(M_i,M_{i+1},a_i)$ minimal reduced and $M = M_0$.  So by \scite{5.5A} we know
$\dsize \bigcup_{i < \lambda^+} M_i \in
K_{\lambda^+}$ is saturated, hence we can embed $N$ into $\dsize \bigcup
_{i < \lambda^+} M_i$ over $N$ so this embedding is into some $M_\alpha,
\alpha < \lambda^+$.]
\endroster
\medskip

So given $M \le_{\frak K} M^1,M^2$, without 
loss of generality $M^\ell = M^\ell_{\alpha_\ell},\langle (M^\ell_i,a^\ell_i):
i \le \alpha_\ell \rangle$ as above,
and start to amalgamate using the extension property and ``reduced". \newline
${}$ \hfill$\square_{\scite{5.7}}$
\enddemo
\bigskip

\remark{\stag{5.8} Remark}  We could prove \scite{5.7} earlier 
using ``reduced triples".  I.e. note that
for some $\langle M^1_i:i < \lambda^+ \rangle \in {\bold Seq}_{\lambda^+}
[\bold C^1_{{\frak K},\lambda^+}]$, for each $i$ for some $a$ the triple
$(M^1_i,M^1_{i+1},a) \in K^3_\lambda$ is reduced.  Hence if $M
\le_{\frak K} N$ from $K_\lambda$, for some $\bar M = \langle M_i:i \le
\alpha \rangle,\le_{\frak K}$-increasing continuous, $\langle M_i,M_{i+1},
b_i) \in K^3_\lambda$ is reduced, $M_0 = M,N \le M_\alpha \in K_\lambda$
(otherwise find $\langle M^2_i:i < \lambda^+ \rangle \in {\bold Seq}
_{\lambda^+}
[\bold C^1_{{\frak K},\lambda^+}]$ with $(M_i,M_{i+1}) \cong (M,N)$, hence 
$M^1 = \dsize \bigcup_{i < \lambda^+} M^1_i,M^2 = \dsize
\bigcup_{i < \lambda^+} M^2_i$ are non-isomorphic members of $K_{\lambda^+}$,
contradiction).  Now prove by induction on $\beta \le \alpha$ that if $M 
\le_{\frak K} N_0 \in K_\lambda$ then $N_0,M_\beta$ has disjoint amalgamation
over $M_0 = M$ (i.e. we need to decompose only one side).
\endremark
\bn
\ub{\stag{5.9} Question}:  If $M \in K_\lambda,p \in {\Cal S}(M)$ is minimal,
is it reduced?  Or at least, if $M_0 \le_{\frak K} M_1$ are in $K_\lambda,
p_1 \in {\Cal S}(M_\ell)$, no algebraic $p_0 = p_1 \restriction M_0,p_0$ is
minimal and reduced is also $p_1$ reduced? \nl
Probably true, but not needed for our aims (here and in \cite{Sh:600},
\scite{7.4B}(2)) is an approximation.  Can prove it if $\lambda <
\lambda^{\aleph_0}$ or there are E.M. models. 
\newpage

\head {\S6 Proving for ${\frak K}$ categorical in $\lambda^{+2}$} \endhead  \resetall
\bigskip

\demo{\stag{6.0} Hypothesis}  Assume the model theoretic assumptions from 
\scite{4.0} + \scite{5.0} and so the further model theoretic properties 
deduced in \S4 + \S5. \newline
We use heavily \scite{4.3}.
\enddemo
\bigskip

\definition{\stag{6.1} Definition}  1) We say $(M_0,M_1,M_2) \in 
K^{\text{uniq}}_\lambda$ (has unique amalgamation in $K_\lambda$) when
\medskip
\roster
\item "{$(a)$}"  $M_0,M_1,M_2 \in K_\lambda$
\sn
\item "{$(b)$}"  $M_0 \le_{\frak K} M_1$ and $M_0 \le_{\frak K} M_2$
\sn
\item "{$(c)$}"  \underbar{if} for $i = 1,2$ we have $g^i_\ell:M_\ell
\rightarrow N_i \in K_\lambda$ such that:
{\roster
\itemitem{ (i) }  $g^i_\ell$ a $\le_{\frak K}$-embedding,
\sn
\itemitem{ (ii) }  $g^i_0 \subseteq g^i_1$ and $g^i_0 \subseteq g^i_2$ for
$i =1,2$
\sn
\itemitem{ (iii) }  $\text{Rang}(g^i_1) \cap \text{ Rang}(g^i_2) = 
\text{ Rang}(g^i_0)$ \newline
(disjoint amalgamation) for $i=1,2$
\endroster}
\underbar{then} we can find $N \in K_\lambda$ and $\le_{\frak K}$-embeddings
$$
f^i:N_i \rightarrow N \text{ for } i =1,2
$$
\underbar{such that}
$$
\dsize \bigwedge_{\ell < 3} f^1 \circ g^1_\ell = f^2 \circ g^2_\ell.
$$
\endroster
\medskip

\noindent
2) Let $K^{2,\text{uq}}_\lambda$ is the class of pairs $(M_0,M_2)$ such that
$M_0 \le_{\frak K} M_2$ are both in $K_\lambda$ and 
$[M_0 \le_{\frak K} M_1 \in K_\lambda \Rightarrow (M_0,M_1,M_2) \in 
K^{\text{uniq}}_\lambda]$ and
let $K^{3,\text{uq}}_\lambda$ be the class of pairs \newline
$(M_0,M_2) \in K^{2,\text{uq}}_\lambda$ satisfying $M_0 \ne M_2$.
\enddefinition
\bigskip

\proclaim{\stag{6.2} Claim}  1) If $(M_0,M_1,M_2) \in 
K^{\text{uniq}}_\lambda$ \underbar{then}
\medskip
\roster
\item "{$(a)$}"   $(M_0,M_2,M_1) \in K^{\text{uniq}}_\lambda$
\sn
\item "{$(b)$}"  if $M_0 \le_{\frak K} M'_2 \le_{\frak K} M_2$ \underbar{then}
$(M_0,M_1,M'_2) \in K^{\text{uniq}}_\lambda$.
\endroster
\medskip

\noindent
2) Assume $M_0 \le_{\frak K} M_2$ are from $K_\lambda$ and 
$M_1 \in K_\lambda$ is universal over $M_0$.  Then \newline
$(M_0,M_2) \in K^{2,\text{uq}}_\lambda \Leftrightarrow (M_0,M_1,M_2) \in 
K^{\text{uniq}}_\lambda$.
\endproclaim
\bigskip

\demo{Proof} 1)a)  Trivial. \newline
b)  Chase arrows (using disjoint amalgamation; i.e. \scite{5.7}). \newline
2) Follows by \scite{6.2}(1)(a)+(b) and the definition. 
\hfill$\square_{\scite{6.2}}$
\enddemo
\bigskip

\proclaim{\stag{6.3} Lemma}  Suppose
\medskip
\roster
\item "{$\bigotimes$}"  there is $(M_0,M_1,M_2) \in K^{\text{uniq}}_\lambda$
such that $M_0 \ne M_2$ and $M_1$ is universal over $M_0$. 
\endroster
\medskip

\noindent
\underbar{Then}:  there are $N^0 <_{\frak K} N^1$ in 
$K_{\lambda^+}$ such that:
\medskip
\roster
\item "{$(a)$}"  $N^0 \ne N^1$
\sn
\item "{$(b)$}"  for every $c \in N^1 \backslash N^0$ there is $M = M_c$
satisfying $N^0 \le_{\frak K} M \le_{\frak K} N^1$ and $N^0 \ne M$ 
and $c \in N^1 \backslash M$.
\endroster
\endproclaim
\bigskip

\demo{Proof}  Choose $\langle N^0_i:i < \lambda^+ \rangle$, a sequence of
members of $K_\lambda$ which is $\le_{\frak K}$-increasing continuous, 
such that:

$$
(N^0_i,N^0_{i+1}) \cong (M_0,M_2).
$$
\medskip

So $N^0_i \ne N^0_{i+1}$ hence $N_0 = \dsize \bigcup_{i < \lambda^+} N^0_i
\in K_{\lambda^+}$ and without loss of generality $|N_0| = \lambda^+$. \nl
We now choose by induction $i < \lambda^+,N^1_i$ and $M_{i,c}$ for
$c \in N^1_i \backslash N^0_i$ such that:
\medskip
\roster
\item "{$(a)$}"  $N^0_i \le_{\frak K} N^1_i \in K_\lambda$ 
and $N^0_i \ne N^1_i$
\sn
\item "{$(b)$}"  $N^1_i$ is $\le_{\frak K}$-increasing continuous in $i$
\sn
\item "{$(c)$}"  $j < i \Rightarrow N^1_j \cap N^0_i = N^0_j$; moreover
$N^1_i \cap |N_0| = N^0_i$
\sn
\item "{$(d)$}"  $N^0_i \le_{\frak K} M_{i,c} \le_{\frak K} N^1_i$
\sn
\item "{$(e)$}"  $c \notin M_{i,c}$
\sn
\item "{$(f)$}"  $N^0_i \ne M_{i,c}$
\sn
\item "{$(g)$}"  if $j < i$ and $c \in N^1_j \backslash N^0_j$ \underbar{then}
$M_{i,c} \cap N^1_j = M_{j,c}$.
\endroster
\enddemo
\bigskip

\noindent
\underbar{For $i = 0$}:  Choose $N^1_i$ such that \newline
$N^0_i \le_{\frak K} N^1_i,(N^1_i,c)_{c \in N^0_i}$ saturated 
(any cofinality will do)
then by disjoint amalgamation easy to define the $M_{0,c}$ (remembering
clause $(c)$).
\bigskip

\noindent
\underbar{For $i$ limit}:  Straightforward.
\bigskip

\noindent
\underbar{For $i = j+1$}:  First we disjointly amalgamate getting 
$N'_i \in K_\lambda$ such that \newline
$N^0_i \le N'_i,N^1_j \le_{\frak K} N'_i$ and $|N'_i| \cap |N_0| =
|N^0_i|$ (as set of elements). \newline
Let $N^1_i$ be such that:

$$
\gather
N'_i \le_{\frak K} N^1_i \in K_\lambda \\
(N^1_i,c)_{c \in N'_i} \text{ is saturated (any cofinality will do)}\\
|N^1_i| \cap |N_0| = |N^0_i|.
\endgather
$$
\medskip

Lastly we shall find the $M_{i,c}$'s, 
the point is that $(N^0_j,N^0_i,N^1_j) \in K^{\text{uniq}}_\lambda$ 
(by \scite{6.2}(2)). \newline
By $\bigotimes$ and Claim \scite{6.2} we could have done the 
amalgamation in two steps and use uniqueness.  Then by 
uniqueness of saturated extensions embed the result inside $N^1_i$ and 
similarly deal with new $c$'s.
\medskip

Now let $N_1 =: \dsize \bigcup_{1 < \lambda^+} N^1_i$ and for
$c \in N_1 \backslash N_0$ let $M_c = \bigcup \{ M_{i,c}:c \in N'_i\}$ they
are as required. \hfill$\square_{\scite{6.3}}$
\bigskip

\remark{Remark}  The proof of \scite{6.4} below is like 
\cite[proof of 3.8 stage(c)]{Sh:88}.  The aim is to contradict that under
$I(\lambda^+,{\frak K}) = 0$ there \ub{are} maximal triples.
\endremark
\bigskip

\demo{\stag{6.4} Conclusion}  Assume ${\frak K}$ has amalgamation in
$\lambda^+$.  If $\bigotimes$ of \scite{6.3}, 
\underbar{then} there is no maximal triple $(M,N,a)$ in $K^3_{\lambda^+}$.
\enddemo
\bigskip

\demo{Proof}  We can get by \scite{6.3} a contradiction. \newline
[Why?  Assume $(N_0,N_2,a) \in K^3_\lambda$ maximal, $(N^0,N^1)$ as in the
conclusion (i.e. (a) + (b)) of \scite{6.3}; by categoricity in $\lambda^+$
without loss of generality $N_0 = N^0$ and let $N_1 = N^1$.  Now ${\frak K}$ 
has amalgamation for $\lambda^+$ so there are $N \in K_{\lambda^+}$ and
$f$ such that $f:N_2 \rightarrow N$ is a $\le_{\frak K}$-embedding of $N_2$
into $N$ over $N_0$ and $N_1 \le_{\frak K} N$.  If $f(a) \notin N_1$, 
\underbar{then} $(N_0,N_2,a) <_f (N_1,N,f(a))$ contradict
maximality.  If $f(a) \in N_1$, then $M_{f(a)}$ is well defined (see 
\scite{6.3}) and $(N_0,N_2,a) <_f (M_{f(a)},N,f(a)$) contradicts maximality.]
\hfill$\square_{\scite{6.4}}$
\enddemo
\bigskip

\remark{\stag{6.4A} Remark}  1) Another proof is to replace the assumption
``${\frak K}$ has amalgamation in $\lambda^+$" by $I(\lambda^+,K) <
2^{\lambda^{+2}}$.  We start with $N_0,N_1,N_2,a$ as
above and building, for every $S \subseteq \lambda^{+2}$, a sequence
$\langle M^S_\alpha:\alpha < \lambda^{+2} \rangle$ of members of
$K_{\lambda^+}$, which is $\le_{\frak K}$-increasing continuous, and 
$\alpha \in S \Rightarrow (M^S_\alpha,M^S_{\alpha + 1},a^S_\alpha) \cong 
(N_0,N_2,a)$, and $\alpha \in \lambda^{+2} \backslash S \Rightarrow 
(M^S_\alpha,M^S_{\alpha+1}) \cong (N_0,N_1)$ which are as in (a) + (b) of 
\scite{6.3}.  Let 
$M^S = \dsize \bigcup_{\alpha < \lambda^{+2}} M^S_\alpha \in K_{\lambda^{+2}}$
and from $M^S/\cong$ we can reconstruct $S/{\Cal D}_{\lambda^{+2}}$.  So here
we use $I(\lambda^{+2},K) < 2^{\lambda^{+2}}$ but not the definitional weak
diamond for $\lambda^{++}$. \newline
2) Note that if $2^{\lambda^+} < 2^{\lambda^+}$ then the assumption of
\scite{6.4A}(1) implies the assumption of \scite{6.4}.
\endremark
\bigskip

\proclaim{\stag{6.5} Claim}  Assume
\medskip
\roster
\item "{$(*)$}"  $2^\lambda < 2^{\lambda^+} < 2^{\lambda^{++}}$ \newline
(or at least the definitional weak diamond for $\lambda^+,\lambda^{++}$).
\endroster
\medskip

\noindent
If $\bigotimes$ 
of \scite{6.3} fails, we get $I(\lambda^{+2},K) \ge \mu_{\text{wd}}(\lambda)$.
\endproclaim
\bigskip

\demo{Proof}  By \scite{3.17B} and \scite{6.5A} below.
\enddemo
\bigskip

\noindent
The following serves to prove \scite{6.5}
\proclaim{\stag{6.5A} Claim}  Assume $M \in K_\lambda \Rightarrow |{\Cal S}
(M)| \le \lambda$. \newline

If $K^{3,\text{uq}}_\lambda = \emptyset$ (see Definition \scite{6.1}(2)),
\underbar{then} there is $F$, an amalgamation choice function $F$ for
$\bold C = \bold C^0_{{\frak K},\lambda^+}$ with the weak $\lambda^+$-coding
property.
\endproclaim
\bigskip

\demo{Proof}  The point is that if $\bar M = \langle M_\alpha:\alpha <
\lambda^+ \rangle \in \bold S eq_{\lambda^+}[\bold C]$ and 
$\alpha < \lambda^+,M_\alpha <_{\frak K} N \in K_\lambda$, \ub{then} 
for some $\beta \in (\alpha,\lambda^+)$ we have:
\medskip
\roster
\item "{${}$}"  $M_\beta$ is universal over $M_\alpha$, so as
$K^{3,\text{uq}}_\lambda = \emptyset$ necessarily \newline
$(M_\alpha,M_\beta,N) \notin K^{\text{uniq}}_\lambda$ and rest should be
clear.
\ermn
Of course, we use the extension property. \hfill$\square_{\scite{6.5A}}$
\enddemo
\bigskip

\remark{\stag{6.5Aa} Remark}  We can work in the context of \S3, we need the
existence of a saturated (equivalent by super limit) $M \in K_{\lambda^+}$.
We now say how to replace $\mu_{\text{wd}}(\lambda^{+2})$ by 
$2^{\lambda^{+2}}$.
\endremark
\bigskip

\proclaim{\stag{6.5B} Claim}  1) Assume each $M \in K_{\lambda^+}$ is
saturated above $\lambda$.

If $(M,N,a) \in K^3_\lambda$ it and every $(M',N',a) \in K^3_\lambda$ above it
has the extension property but for every $(M'',N'',a) \ge (N,N,a)$ (all in
$K^3_{\lambda^+}$ for some $M^* \ge_{\frak K} M''$ from $K_{\lambda^+}$, in
amalgamation $(M^*,N^*,a) \ge (M'',N'',a)$ the type of $M^* \cup N$ inside
$N^*$ is not determined \underbar{then} some $F$ (actually $\bold F^*$) has
the $\lambda^+$-coding. \newline
2) If above we just require that the type of $M^* \cup N''$ inside $N^*$ is
not determined, \underbar{then} some $F$ (actually $\bold F^*$) has weak
$\lambda^+$-coding. \newline
3) We can restrict ourselves to disjoint embedding; i.e. use 
$(K^3_{\lambda^+},\le_{\text{dj}})$.
\endproclaim
\bigskip

\demo{\stag{6.5C} Discussion}  1) We get $IE(\lambda^{+2},{\frak K}) = 
2^{\lambda^{+2}}$ when $(2^{\lambda^+})^+ < 2^{\lambda^{+2}}$.  See more in
\cite{Sh:600}.
\enddemo
\bigskip

\proclaim{\stag{6.6} Theorem}  Assume $I(\lambda^{+2},K) < \mu_{\text{wd}}
(\lambda^{+2})$ and $(*)$ of \scite{6.5} (or at least the conclusion of
\scite{6.5}).  Then  
$I(\lambda^{+2},K) = 1 \Rightarrow I(\lambda^{+3},K) > 0$.
\endproclaim
\bigskip

\remark{Remark}  As in \cite[\S3]{Sh:88}.
\endremark
\bigskip

\demo{Proof}  By \scite{0.17}(1) it is enough to show that for some $M \in 
K_{\lambda^{++}}$ there is $M'$, \newline
$M \le_{\frak K} M' \in K_{\lambda^{++}},M \ne M'$. \newline
[Why?  As then we can choose by induction on $i < \lambda^{+3}$ models
$M_i \in K_{\lambda^{+2}}$, \newline
$\le_{\frak K}$-increasing continuous, $M_i \ne M_{i+1}$, for
$i = 0$ use $K_{\lambda^{+2}} \ne \emptyset$, for $i$ limit take union, for
$i = j+1$ use the previous sentence; so $M_{\lambda^{+3}} = \cup \{ M_i:i <
\lambda^{+3}\} \in K_{\lambda^{+3}}$ as required.] \newline
By \scite{6.5}, the statement $\bigotimes$ of \scite{6.3} holds so 
we can find $(N^0,N^1)$ as there so by \scite{6.4} there is in 
$K^3_{\lambda^+}$ no 
maximal member.  This implies (easy, see \scite{2.4}(6)) that there are 
$M^* \le_{\frak K} N^*$ from $K_{\lambda^{+2}}$ such that $M^* \ne N^*$ 
which as said above (by categoricity in $\lambda^+$), suffices.  
\hfill$\square_{\scite{6.6}}$
\enddemo 
\newpage

\head {\S7 Extensions and conjugacy} \endhead  \resetall
\bigskip

\demo{\stag{7.0} Hypothesis}  Assume the model theoretic 
assumptions from \scite{4.0} + \scite{5.1} and the further model theoretic 
properties deduced since then (but not in \scite{6.5},\scite{6.6}) or just
\medskip
\roster
\item "{$(a)$}"  ${\frak K}$ abstract elementary class
\sn
\item "{$(b)$}"  ${\frak K}$ has amalgamation in $\lambda$
\sn
\item "{$(c)$}"  ${\frak K}$ is categorical in $\lambda$ (can be weakened)
\sn
\item "{$(d)$}"  ${\frak K}$ is stable in $\lambda$ (see \scite{5.5}, clause
$(a)$)
\sn
\item "{$(e)$}"  there is an inevitable $p \in {\Cal S}(N)$ for 
$N \in K_\lambda$ (holds by \scite{5.2})
\sn
\item "{$(f)$}"  the basic properties in type theory.
\endroster
\enddemo
\bigskip

We now continue toward eliminating the use of $I(\lambda^{++},K) = 1$ 
(in \scite{6.6}), and give more information.
We first deal with the nice types in ${\Cal S}(N),N \in K_\lambda$ 
in particular the realize/materialize problem which is here: if $N_1
\le_{\frak K} N_2$ are in $K_\lambda,p_\ell \in {\Cal S}(N_\ell)$ is minimal,
$p_1 \le p_2$ are they conjugate? (i.e. $p_2 \in {\Cal S}_{p_1}(N_2)$).
\bigskip

\proclaim{\stag{7.1} Claim}  If $N \in K_\lambda$ and $p \in {\Cal S}(N)$ 
is minimal and reduced or just $p$ is reduced (see Definition 
\scite{2.3A}(7)), \underbar{then} $p$ is inevitable.
\endproclaim
\bigskip

\demo{Proof}  Suppose $N,p$ form a counterexample.  We can then find $N_1$ and
$a$ such that $N \le_{\frak K} N_1 \in K_\lambda,a \in N_1 \backslash N$ 
and $p = \text{ tp}(a,N,N_1)$ and $(N,N_1,a)$ is reduced.  As $p$ is not 
inevitable, there is
$N_2$ such that: $N \le_{\frak K} N_2 \in K_\lambda,N \ne N_2$ but no element of
$N_2$ realizes $p$.  By amalgamation in $K_\lambda$, without loss of
generality there is $N_3 \in K_\lambda$ such that $\ell \in \{1,2\}
\Rightarrow N_\ell \le_{\frak K} N_3$.  By \scite{5.2} (i.e.\scite{7.0}(e)) 
there is
$q \in {\Cal S}(N)$ which is inevitable so there are $c_\ell \in N_\ell$ with
$q = \text{ tp}(c_\ell,N,N_\ell)$ for $\ell \in \{1,2\}$.
By the equality of types (and amalgamation in $K_\lambda$) there is
$N^+ \in K$, a $\le_{\frak K}$-extension of $N_1$ and a 
$\le_{\frak K}$-embedding of $N_2$ into $N^+$ over $N$ such that
$f(c_2) = c_1$; so without loss of generality $N^+ = N_3$ and $f$ is the
identity, hence $c_1 = c_2$.  Now $a \notin N_2$ as $p = \text{ tp}(a,N,
N_1)$ is not realized in $N_2$.  So $(N,N_1,a) \le (N_2,N_3,a)$ and
$N_2 \cap N_1 \backslash N \ne \emptyset$ contradicting ``$(N,N_1,a)$ 
is reduced".  \hfill$\square_{\scite{7.1}}$
\enddemo
\bigskip

\proclaim{\stag{7.2} Claim}  1) If $\kappa = \text{ cf}(\kappa) \le \lambda$ 
and $\bar N = \langle N_i:i \le \omega \kappa \rangle$ is an
$\le_{\frak K}$-increasingly continuous sequence,
$N_i \in K_\lambda,N_{i+1}$ universal over $N_i$, and $p \in {\Cal S}
(N_{\omega \kappa})$ is minimal reduced (or minimal inevitable) 
\underbar{then} for some $i < \omega \kappa$ we have
$p \restriction N_i \in {\Cal S}(N_i)$ is minimal (so $p$ is the unique, 
non-algebraic extension of $p \restriction N_i$ in ${\Cal S}
(N_{\omega \kappa})$ (and of course, there is one)). \newline
2) If $\lambda \ge \kappa = \text{ cf}(\kappa),\bar N = \langle N_i:i \le
\kappa \rangle$ is $\le_{\frak K}$-increasing continuous in $K_\lambda$ 
and \newline
$p \in {\Cal S}(N_\kappa)$ is minimal and reduced and 
the set $Y =: \{i < \kappa:N_{i+1} \text{ is }
(\lambda,\kappa) \text{-saturated over } N_i\}$ is unbounded in $\kappa$
\underbar{then} for every large enough $i \in Y$ there is an isomorphism $f$
from $N_{i+1}$ onto $N_\kappa$ which is the identity on $N_i$ and
\footnote{in fact, this implies $(*)$}
\medskip
\roster
\item "{$(*)$}"  $f$ maps $p \restriction N_{i+1} \in {\Cal S}(N_{i+1})$ to
$p \in {\Cal S}(N_\kappa)$.
\ermn
Hence as $p$ is minimal reduced, so is $p \restriction N_{i+1}$.
\endproclaim
\bigskip

\demo{Proof}  1) We can choose $(N^0_i,N^1_i,a) \in K^3_\lambda$ for
$i < \lambda^+$ reduced, $\le$-increasing continuous such 
that $N^0_i \ne N^0_{i+1}$.  
Let $N_\ell = \dsize
\bigcup_{i < \lambda^+} N^\ell_i$.  As in the proof of \scite{5.3} for 
$c \in N_1 \backslash N_0$

$$
I^*_c = \{ j < \lambda^+:c \in N^1_j \text{ and } \text{tp}(c,N^0_j,N^1_j)
\text{ is minimal}\}
$$
\medskip

\noindent
is empty or is an end segment of $\lambda^+$ and

$$
\align
E = \biggl\{ \delta < \lambda^+:&\text{ \underbar{if} } c \in N^1_\delta
\text{ and } I^*_c \ne \emptyset \text{ \underbar{then} } I^*_c \cap \delta\\
  &\text{ is an unbounded subset of } \delta;\text{and if } \alpha < \delta \\
  &\text{ then for some } \beta \in (\alpha,\delta),N^0_\beta \text{ is
universal over } N^0_\alpha \\
  &\text{ and if } Pr \text{ is one of the properties reduced and/or} \\
  &\text{ inevitable and/or minimal and there is } i \ge \delta 
\text{ such that}\\
  &\,(N^0_i,N^1_i,c) \text{ has } Pr, \text{ then there are arbitrarily} \\
  &\text{ large such } i < \delta \biggr\}
\endalign
$$
\medskip

\noindent
is a club of $\lambda^+$; for the universality demand in the definition 
of $E$ use
categoricity in $\lambda^+$.  Let $\delta \in \text{ acc}\text{(acc}(E)),
\text{cf}(\delta) = \kappa$,
let $\langle \alpha_\zeta:\zeta < \omega \kappa \rangle$ be an 
increasing continuous sequence of ordinals from $E$ with limit $\delta$, 
now set $\alpha_{\omega \kappa} = \delta$ and 
$N'_\zeta =: N^0_{\alpha_\zeta}$.
\medskip

So there is an isomorphism $f$ from $N_{\omega \kappa}$ onto 
$N^0_{\alpha_{\omega \kappa}}$ such that for every $\zeta < \omega \kappa$ 
we have $N^0_{\alpha_{2 \zeta}} \le_{\frak K} f(N_{\alpha_{2 \zeta}}) 
\le_{\frak K} N^0_{\alpha_{2 \zeta +1}}$ (so if
$\zeta$ is a limit ordinal, then $N^0_{\alpha_\zeta} = N^0_{\alpha_{2 \zeta}} 
= f(N_\zeta)$), so without loss of generality $f$ is the identity.
As $p \in {\Cal S}(N_{\omega \kappa})$ is inevitable (by assumption or by
\scite{7.1}) and $N_{\omega \kappa} = N^0_{\alpha_{\omega \kappa}} 
<_{\frak K} N^1_{\alpha_{\omega \kappa}}$, for some $c \in 
N^1_{\alpha_{\omega \kappa}} \backslash N^0_{\alpha_{\omega \kappa}}$ we have 
$p = \text{ tp}(c,N^0_{\alpha_{\omega \kappa}},
N^1_{\alpha_{\omega \kappa}})$, so for some $\beta < \alpha_{\omega \kappa}$ 
we have $c \in N^1_\beta$.  As $p$ is minimal (by assumption) 
clearly $\delta \in I_c$, but
$\delta \in E$ so $\text{Min}(I_c) < \delta$, but $I_c$ is an end segment of
$\lambda^+$ hence without loss of generality for some $\zeta < \omega \kappa$ 
we have $\beta = \alpha_\zeta \in I_c$.  So for $\xi \in 
(\zeta,\omega \kappa)$, both
$p \in {\Cal S}(N_{\omega \kappa})$ and $p \restriction N_\xi \in {\Cal S}
(N_\xi)$ are non-algebraic extensions of the minimal $p \restriction 
N^0_{\alpha_\zeta} \in {\Cal S}(N^0_{\alpha_\zeta})$ and 
$N^0_{\alpha_\zeta} \le_{\frak K} N_\xi \le_{\frak K} N_{\omega \kappa}$, 
all in $K_\lambda$, so we have proved part (1).
\medskip

\noindent
2) Without loss of generality every $\zeta < \kappa$ is in $Y$.  We can find
$\langle N'_\zeta:\zeta \le \omega \kappa \rangle$ as in part (1), moreover
satisfying ``$N'_{\zeta +1}$ is $(\lambda,\kappa)$-saturated over
$N_\zeta$" and such that: for every $\zeta \le \kappa$ we have $N_\zeta =
N'_{\omega \zeta}$.  So again choose $\zeta < \kappa$ as there, we 
set $\beta = \alpha_{\omega \zeta} \in I^*_c$.
If $\xi \in Y \and \omega \xi > \zeta$ clearly by the uniqueness of
$(\lambda,\kappa)$-saturated models there is an isomorphic $f$ from 
$N_{\xi +1} = N'_{\omega(\xi + 1)}$ onto $N_\kappa = N'_{\omega \kappa}$ 
over $N_\xi = N'_{\omega \kappa}$, and $f(p \restriction N_{\xi + 1}) = p$ 
is proved as above. \hfill$\square_{\scite{7.2}}$
\enddemo
\bigskip

\proclaim{\stag{7.3} Claim}  1) If $M_0 \le_{\frak K} M_1$ are in 
$K_\lambda$ and the types $p_\ell \in {\Cal S}(M_\ell)$ are minimal 
reduced, for $\ell = 0,1$ 
and $p_0 = p_1 \restriction M_0$ \underbar{then}
$p_0,p_1$ are conjugate; (i.e. there is an isomorphic $f$ from $M_0$ 
onto $M_1$ such that $f(p_0) = p_1$). \newline
2) If in addition $M \le_{\frak K} M_0$ and $M_0,M_1$ are $(\lambda,
\kappa)$-saturated over $M$, \ub{then} $p_0,p_1$ are conjugate over $N$.
\endproclaim
\bigskip

\remark{Remark}  Note that $p$ minimal (or reduced) implies that $p$ is not
algebraic.
\endremark
\bigskip

\demo{Proof}  1) Let $\langle (N^0_i,N^1_i,a):i < \lambda^+ \rangle$ and $E$
be as in the proof of \scite{7.2} and \newline
$\kappa = \text{ cf}(\kappa) \le \lambda$.  
For each $\delta \in S_\kappa =: \{ \alpha < \lambda^+:\alpha \in E 
\text{ and } \text{cf}(\alpha) = \kappa\}$, and minimal reduced $p 
\in {\Cal S}(N^0_\delta)$, we know that for some $i_p < \delta,p 
\restriction N^0_{i_p}$ is minimal reduced [why?  by \scite{7.2}(1),(2)] 
and some $q_p \in {\Cal S}(N^0_{i_p})$ is conjugate to $p$.  For \newline
$\kappa = \text{ cf}(\kappa) \le \lambda,q \in {\Cal S}(N^0_i),i < \lambda^+,
r \in {\Cal S}(N^0_i)$ minimal let
\medskip

\noindent
$A^{\kappa,i}_{q,r} = \bigl\{ \delta < \lambda^+:\text{there is a type } 
p \text{ such that }
r \subseteq p \in {\Cal S}(N^0_\delta),p$ non-algebraic \newline

$\qquad \qquad \text{(this determines } p$), $p$ minimal reduced,
$i_p = i,q_p = q$ \newline

$\qquad \qquad \text{(and clearly } p \restriction N^0_i = r) \}$.  \newline
\smallskip

\noindent
Next let

$$
\align
E_1 = \biggl\{ \delta < \lambda^+:&\text{ for every } \kappa 
= \text{ cf}(\kappa) \le \lambda, \\
  &\,\,r,q \in {\Cal S}(N^0_i) \text{ and } i < \delta, 
\text{ \underbar{if} } 
A^{\kappa,i}_{q,r} \text{ is well defined and} \\
  &\text{ unbounded in } \lambda^+ \text{ \underbar{then} it is unbounded in }
 \delta \biggr\}.
\endalign
$$
\medskip

\noindent
So if $\delta_1 \in E_1,\kappa = \text{ cf}(\delta_1),p_1 \in {\Cal S}
(N^0_\delta)$ is minimal reduced, \underbar{then} we can find \newline
$\delta_0 < \delta_1,\text{cf}(\delta_0) = \kappa$, and 
$p_0 \in {\Cal S}(N^0_{\delta_0})$ minimal reduced with
\newline
$q_{p_1} = q_{p_0},i_{p_1} = i_{p_0},p_0 \restriction N^0_{i_{p_0}} = p_1
\restriction N^0_{i_{p_1}}$ call it $r$, it is necessarily minimal.
\medskip

As $p_1,p_0$ extend $r,N^0_{i_{p_0}} = N^0_{i_{p_1}} \le_{\frak K} 
N^0_{\delta_0} \le_{\frak K} N^0_{\delta_1}$, necessarily 
$p_1 = p_0 \restriction N^0_{\delta_0}$,
and also they are both conjugate to $q_{p_0} = q_{p_1}$ hence they are
conjugate.
\newline
Next we prove 
\medskip
\roster
\item "{$(*)$}"  if $M_0 <_{\frak K} M_1$ are in $K_\lambda,M_1$ is 
$(\lambda,\kappa)$-limit over $M_0$, $p'_0 \in {\Cal S}(M_0)$ is minimal 
reduced and $p'_0 \le p'_1 \in {\Cal S}(M_1),p'_1$ non-algebraic, 
\underbar{then} $p'_0,p'_1$ are conjugate.
\endroster
\medskip

\noindent
Above we have a good amount of free choice in choosing $p_1 \in {\Cal S}
(N^0_{\delta_1})$ (it should be minimal and reduced) so we could have chosen
$p_1$ to be conjugate to $p'_0$ (i.e. is in 
${\Cal S}_{p'_0}(N^0_{\delta_1})$; now also the corresponding 
$p_0$ is conjugate to $p_1$ hence
$p_0$ is conjugate to $p'_0$, hence we can find an isomorphism $f_0$ from 
$M_0$ onto $N^0_{\delta_0},
f_0(p'_0) = p_0$, and extend it to an isomorphism $f_1$ from $M_1$ onto
$N^0_{\delta_1}$, so necessarily $f(p'_1) = p_1$ (as $p_1$ is the unique 
non-algebraic extension).  As $p_0,p_1$ are conjugate through $f_1$ also 
$p'_0,p'_1$ are conjugate.  So $(*)$ holds.
\medskip

\noindent
Now assume just $M_0 \le_{\frak K} M_1$ are in $K_\lambda,
p_0 \in {\Cal S}(M_0)$ minimal
reduced, $p_1 \in {\Cal S}(M_1)$ the unique non-algebraic extension of $p_0$
and it is reduced (and necessarily minimal).
There is $M_2,M_1 \le_{\frak K} M_2 \in K_\lambda,M_2$ is 
$(\lambda,\kappa)$-limit over $M_1$ hence also over $M_0$ and let $p_2$ be
the unique non-algebraic extension of $p_1$ in ${\Cal S}(M_2)$ hence
$p_2$ is also the unique non-algebraic extension of $p_0$ in ${\Cal S}(M_2)$.
\newline
Using $(*)$ on $(M_0,M_2,p_0,p_2)$ and on $(M_1,M_2,p_1,p_2)$ and get that
$p_0,p_2$ are conjugate and that $p_1,p_2$ are conjugate hence $p_1,p_2$ are 
conjugate, the required result. \newline
2) Similar proof.  \hfill$\square_{\scite{7.3}}$
\enddemo
\bigskip

\proclaim{\stag{7.4A} Claim}  1) Assume $M_1 \le_{\frak K} M_2$ are in 
$K_\lambda$ and $M_2$ is $(\lambda,\kappa)$-saturated over $M_1$.   If
$p_1 \in {\Cal S}(M_1)$ is minimal and reduced, \ub{then} $p_2$, the unique
non-algebraic extension of $p_1$ in ${\Cal S}(M_2)$, is reduced (and, of
course, minimal). \nl
2) There is no need to assume ``$p_1$ reduced.
\endproclaim
\bigskip

\demo{Proof}  1) We can find $\langle N_i:i \le \kappa \rangle$, an
$\le_{\frak K}$-increasingly continuous sequence in $K_\lambda$ such that
$N_{i+1}$ is $(\lambda,\kappa)$-saturated over $N_i$ and $N_\kappa = M_1$.
So by \scite{7.2}(1),(2), for some 
$\zeta < \kappa$ we have: $p_2 \restriction 
N_\zeta$ is minimal and for some isomorphism $f$ from $N_{\zeta +1}$ onto
$N_\kappa$ we have $f(p_1 \restriction N_\zeta) = p_1$ and $f \restriction
N_\zeta = \text{ id}_{N_\zeta}$.  Also $M_1,N{\zeta+1}$ are isomorphic over
$N_\zeta$ (as both are $(\lambda,\kappa)$-saturated over it) hence there is
an isomorphism $g$ from $N_{\zeta+1}$ onto $M_2$ over $N_\zeta$.  Now
$p_1 = f(p_1 \restriction N_{\zeta +1})$ and $f_2 =: g(p_1 \restriction
N_{\zeta +1})$ are non-algebraic extensions of $p_1 \restriction N_\zeta$ 
which is minimal, hence $p_1 = p_2 \restriction M_1$ and $p_2$ is as mentioned
in \scite{7.4A}.  Now $g \circ f^{-1}$ show that $p_1,p_2$ are conjugate so
as $p_1$ is reduced also $p_2$ is reduced. \nl
2) Easy as we can find $N,M_1 \le_K N,q \in {\Cal S}(N)$ minimal reduced;
without loss of generality 
$N \le_{\frak K} M_2$ and $M_2$ is $(\lambda,\kappa)$-saturated
over $N$, and apply part (1).  \hfill$\square_{\scite{7.4A}}$
\enddemo
\bigskip

\proclaim{\stag{7.4B} Claim}  Assume
\mr
\item "{$(a)$}"  $N_{i,j} \in K_\lambda$ for $i \le \delta_1,j \le \delta_2$
\sn
\item "{$(b)$}"  $\langle N_{i,j}:j \le \delta_2 \rangle$ is
$\le_{\frak K}$-increasingly continuous for each $i \le \delta_1$
\sn
\item "{$(c)$}"  $\langle N_{i,j}:i < \delta_1 \rangle$ is
$\le_{\frak K}$-increasingly continuous for each $j \le \delta_2$
\sn
\item "{$(d)$}"  $\langle N_{i,j}:i \le \delta_1,j \in \delta_2 \rangle$ is
smooth, i.e. $N_{i_1,j_1} \cap N_{i_2,j_2} = N_{\text{min}\{i_1,i_2\}} \cap
N_{\text{min}\{j_1,j_2\}}$
\sn
\item "{$(e)$}"  $N_{i+1,j+1}$ is universal over $N_{i,j+1} \cup N_{i+1,j}$
(i.e. $N_{i+1,j+1}$ is universal over some $N'_{i+1,j+1}$ where
$N_{i,j+1} \cup N_{i+1,j} \subseteq N'_{i+1,j+1}$,
\sn
\item "{$(f)$}"  $\delta_1$ is divisible by 
cf$(\delta_2) \times \lambda \times \omega$ \nl
(and even easier if $\delta_1 = 1$!).
\ermn
\ub{Then} $N_{\delta_1,\delta_2}$ is $(\lambda,\text{cf}(\delta_1))$-saturated
over $N_{i,\delta_1}$ for $i < \delta_1$.
\endproclaim
\bigskip

\demo{Proof}  Without loss of generality $\delta_2 = \text{ cf}(\delta_2)$.
(Why?  let $\langle \alpha_\varepsilon:\varepsilon \le \text{ cf}(\delta_1)
\rangle$ be increasingly continuous with limit $\delta_1$ such that
[$\varepsilon$ limit $\leftrightarrow \alpha_\varepsilon$ limit], and use
$N'_{i,\varepsilon} = N_{i,\alpha_\varepsilon}$). \nl
For $i < \delta_1,j < \delta_2$ let $M_{i,j},M'_{i,j}$ be such that
$M_{i+1,j} \cup M_{i,j+1} \subseteq M_{i,j} \le_{\frak K} M'_{i,j} \le
M_{i+1,j+1}$ and $M'_{i,j}$ is $(\lambda,\text{cf}(\lambda))$-saturated over
$M_{i,j}$.
\sn
Now let $p \in {\Cal S}(N)$ be minimal and reduced; and for $i \le \delta_1,
j \le \delta_2$ let $p_{i,j} \in {\Cal S}(N_{i,j})$ be the unique
non-algebraic extension of $p$ in ${\Cal S}(N_{i,j})$ so it is minimal.  Now
for $i < \delta_1$, note that $\langle M'_{i + \varepsilon,\varepsilon}:
\varepsilon < \delta_2 \rangle$ is $\le_{\frak K}$-increasing 
(not continuous!) and $M'_{i + \varepsilon +1,\varepsilon +1}$ is $(\lambda,
\text{cf}(\lambda))$-saturated over $M'_{i + \varepsilon,\varepsilon}$ and
$\dbcu_{\varepsilon < \delta} M'_{i + \varepsilon,\varepsilon} = 
\dbcu_{\varepsilon < \delta_1} N_{i+1,\varepsilon +1} = N_{i + \delta_2,
\delta_2}$, and $N_{0,0} \le_{\frak k} M'_{1,0}$, hence by \scite{7.4A} we
know that $p_{i + \delta_2,\delta_2}$ is reduced (and minimal).  In fact,
similarly $\alpha < \delta_1 \and \text{ cf}(\alpha) = \text{ cf}(\delta_2)
\Rightarrow p_{\alpha,\delta_2}$ is reduced.  As
$N_{i+1,j+1} \ne N_{i+1,j} \cup N_{i,j+1}$ and clause (d) (smoothness)
necessarily $N_{i + \delta_2,\delta_2} < N_{i + \delta_2 +1,\delta_2}$, hence
some $c \in N_{i + \delta_2 + i,\delta_2} \backslash 
N_{i + \delta_2,\delta_2}$ realizes $p_{i + \delta_2,\delta_2}$.  So if
$\alpha \le \delta_1$ is divisible by $\delta_2 \times \lambda$ and has
cofinality cf$(\delta_2)$ and $\beta < \alpha$, then by \scite{5.4}
$N_{\alpha,\delta_2}$ is universal over $N_{\beta,\delta_2}$.  As
$\delta_1$ is divisible by cf$(\delta_2) \times \lambda \times \omega$ we
are done.  \hfill$\square_{\scite{7.4B}}$
\enddemo
\bigskip

\proclaim{\stag{7.5} Lemma}  1) For every $N \in K_{\lambda^+}$ we can find a
representation \newline
$\bar N = \langle N_i:i < \lambda^+ \rangle$, with
$N_{i+1}$ being $(\lambda,\text{cf}(\lambda))$-saturated over $N_i$.
\newline
2)  If for $\ell = 1,2$ we have $N^\ell = \langle N^\ell_i:i < \lambda^+ 
\rangle$ as in part (1) \underbar{then} there is an isomorphism $f$ from 
$N^1$ onto
$N^2$ mapping $N^1_i$ onto $N^2_i$ for each $i < \lambda^+$.  Moreover
for any $i < \lambda^+$ and isomorphism $g$ from $N^1_i$ onto $N^2_i$ we
can find an isomorphism $f$ from $N^1$ onto $N^2$ extending $g$ and mapping
$N^1_j$ onto $N^2_j$ for each $j \in [i,\lambda^+)$. \newline
3)  If $N^0 \le_{\frak K} N^1$ are in $K_{\lambda^+}$ then we can 
find representation
$\bar N^\ell$ of $N^\ell$ as in (1) with $N^0_i = N^0 \cap N^1_i$, (so
$N^0_i \le_{\frak K} N^1_i)$. \newline
4)  For any strictly increasing function $\bold f:\lambda^+ \rightarrow
\lambda^+$, we can find $N_{i,\varepsilon}$ for \newline
$i < \lambda^+,\varepsilon \le \lambda \times (1 + \bold f(i))$ such that:
\medskip
\roster
\item "{$(a)$}"  $N_{i,\varepsilon} \in K_\lambda$
\sn
\item "{$(b)$}"  $\langle N_{i,\varepsilon}:\varepsilon \le \lambda \times
(1 + \bold f(i)) \rangle$ is strictly $\le_{\frak K}$-increasing continuous
\sn
\item "{$(c)$}"  for each $\varepsilon,\langle N_{i,\varepsilon}:i \in
[i_\varepsilon,\lambda^+) \rangle$ is a representation 
as in (1) where \newline
$i_\varepsilon = \text{ Min}\{i:\varepsilon \le \lambda \times (1 + \bold f
(i))\}$
\sn
\item "{$(d)$}"  if $\varepsilon < \lambda \times (1 + f(i))$ and 
$i < j < \lambda^+$ then $N_{j,\varepsilon} \cap N_{i,\lambda \times (1 +
\bold f(i))} = N_{i,\varepsilon}$
\sn
\item "{$(e)$}"  $N_{i+1,\varepsilon +1}$ is $(\lambda,\aleph_0)$-saturated
over $N_{i+1,\varepsilon} \cup N_{i,\varepsilon +1}$.
\endroster
\endproclaim
\bigskip

\demo{Proof}  Straight. \nl
4) First use $\bold f':\lambda^+ \rightarrow \lambda^+$ which is
$\bold f'(i) = \lambda^\omega \times \bold f(i)$.  Then define the
$N_{i,\varepsilon} \, (\varepsilon < \lambda \times (1 + \bold f'(i);i <
\lambda^+)$. ``Forget" about ``$N_{i+1,\varepsilon}$ is 
$(\lambda,\text{cf}(\lambda))$-saturated over $N_{i,\varepsilon}$",  
remember we have disjoint amalgamation by \scite{5.7}.  Now by \scite{7.4B}, 
even for $\varepsilon$ limit divisible by $\lambda^3$ we get 
$N_{i + \lambda,\varepsilon}$ is $(\lambda,\text{cf}(\lambda))$-saturated 
over $N_{i,\varepsilon}$, so renaming all is O.K.). 
\hfill$\square_{\scite{7.5}}$
\enddemo
\bigskip

\noindent
We can deduce from \scite{7.4A}, but to keep the door open to other 
uses we shall not use
\proclaim{\stag{7.6} Claim}  If $\kappa_\ell = \text{ cf}(\kappa_\ell) \le
\lambda$, and $N_\ell$ is $(\lambda,\kappa_\ell)$-saturated over $N$ for
$\ell =1,2$ \underbar{then} $N_1,N_2$ are isomorphic over $N$.
\endproclaim
\bigskip

\demo{Proof}  We can define by induction on $i \le \lambda \times \kappa_1$,
and then by induction on \newline
$j \le \lambda \times \kappa_2,M_{i,j}$ such that:
\medskip
\roster
\item "{$(a)$}"  $M_{i,j} \in K_\lambda$
\sn
\item "{$(b)$}"  $M_{0,0} = N$
\sn
\item "{$(c)$}"  $i_1 \le i \and j_1 \le j \Rightarrow M_{i_1,j_1}
\le_{\frak K} M_{i,j}$
\sn
\item "{$(d)$}"  $M_{i_1,j_1} \cap M_{i_2,j_2} = M_{\text{min}\{i_1,i_2\},
\text{min}\{j_1,j_2\}}$
\sn
\item "{$(e)$}"  $M_{i,j}$ is $\le_{\frak K}$-increasing continuous in $i$
\sn
\item "{$(f)$}"  $M_{i,j}$ is $<_{\frak K}$-increasing continuous in $j$
\sn
\item "{$(g)$}"  $M_{0,j} \ne M_{0,j+1}$
\sn
\item "{$(h)$}"   $M_{i+1,0} \ne M_{i,0}$
\sn
\item "{$(i)$}"  $M_{i+1,j+1}$ is universal over $M_{i+1,j} \cup M_{i,j+1}$.
\endroster
\medskip

\noindent
There is no problem by \scite{5.7}(1) (using the existence of disjoint 
amalgamation). \newline
Now $M_{\lambda \times \kappa_1,\lambda \times \kappa_2}$ is the union of the
strictly $<_{\frak K}$-increasing sequence \newline
$\langle M_{0,0} \rangle \char 94
\langle M_{\lambda \times i,\lambda \times \kappa_2}:i < \kappa_1 \rangle$
hence by \scite{7.4B} is $(\lambda,\kappa_1)$-saturated over $M_{0,0} = N$ 
hence $M_{\lambda \times \kappa_1,\lambda \times \kappa_2} \cong_N N_1$.
Similarly $M_{\lambda \times \kappa_1,\lambda \times \kappa_2}$ is the union
of the strictly $\le_{\frak K}$-increasing sequence $\langle M_{0,0} \rangle
\char 94 \langle M_{\lambda \times \kappa_1,\lambda \times j}:j < \kappa_2
\rangle$ hence is $(\lambda,\kappa_2)$-saturated over $M_{0,0} = N$, 
hence
$M_{\lambda \times \kappa_1,\lambda \times \kappa_2} \cong_N N_2$.  Together
$N_1,N_2$ are isomorphic over $N$. \nl
${{}}$  \hfill$\square_{\scite{7.6}}$
\enddemo
\newpage

\head {\S8 Uniqueness of amalgamation in ${\frak K}_\lambda$} \endhead  \resetall
\bigskip

We deal in this section only with $K_\lambda$. \newline
We want to, at least, approximate unique amalgamation using as starting
point $\bigotimes$ of \scite{6.3} (see also \scite{6.5}), i.e. 
$K^{3,\text{uq}}_\lambda \ne \emptyset$.
\bigskip

\demo{\stag{8.0} Hypothesis}  1) Assume hypothesis \scite{7.0}, so
\medskip
\roster
\item "{$(a)$}"  ${\frak K}$ abstract elementary class
\sn
\item "{$(b)$}"  ${\frak K}$ has amalgamation in $\lambda$
\sn
\item "{$(c)$}"  ${\frak K}$ is categorical in $\lambda$ (can be weakened)
\sn
\item "{$(d)$}"  ${\frak K}$ is stable in $\lambda$ (see \scite{5.5},
clause $(a)$)
\sn
\item "{$(e)$}"  there is an inevitable $p \in {\Cal S}(N)$ for 
$N \in K_\lambda$ (holds by \scite{5.2})
\sn
\item "{$(f)$}"  the basic properties in type theory.
\endroster
\medskip

\noindent
2) $(M^*,N^*)$ is some pair in \newline
$K^{3,\text{uq}}_\lambda = \{(M_0,M_2):
M_0 \le_{\frak K} M_2$ are in $K_\lambda$ and for every $M_1,M_0
\le_{\frak K} M_1 \in K_\lambda \Rightarrow (M_0,M_1,M_2) \in K^{\text{uniq}}
_\lambda$; equivalently $\text{for some } M_1,(M_0,M_1,M_2)$ are as in 
$\bigotimes$ of \scite{6.3}$\}$ (eventually the choice does not matter; 
if each time instead of $\cong(M^*,N^*)$ we write 
$\in K^{2,\text{uq}}_\lambda$, see \scite{8.9}; but if we start with this
definition then the uniqueness theorems will be more cumbersome).
\enddemo
\bigskip

\definition{\stag{8.1} Definition}  Assume $\bar \delta = \langle \delta_1,
\delta_2,\delta_3 \rangle,\delta_2$ a limit ordinal $< \lambda^+$ but
$\delta_1,\delta_3$ are $< \lambda^+$ and may be 1.  We say that 
$NF_{\lambda,\bar \delta}(N_0,N_1,N_2,N_3)$ (we say
$N_1,N_2$ are saturated and smoothly amalgamated in $N_3$ over $N_0$ for
$\bar \delta$) when:
\medskip
\roster
\item "{$(a)$}"  $N_\ell \in K_\lambda$ for $\ell \in \{ 0,1,2,3\}$
\sn
\item "{$(b)$}"  $N_0 \le_{\frak K} N_\ell \le_{\frak K} N_3$ for $\ell = 1,2$
\sn
\item "{$(c)$}"  $N_1 \cap N_2 = N_0$ (i.e. in disjoint amalgamation)
\sn
\item "{$(d)$}"  $N_1$ is ($\lambda$,cf$(\delta_1)$)-saturated over $N_0$
\sn
\item "{$(e)$}"  $N_2$ is ($\lambda$,cf$(\delta_2)$)-saturated over $N_0$, if
$\delta_1 =1$ this means just $N_0 \le_{\frak K} N_2$
\sn
\item "{$(f)$}"  there are $N_{1,i},N_{2,i}$ for $i \le \lambda \times
\delta_1$ (called the witness) such that:
{\roster
\itemitem{ $(\alpha)$ }  $N_{1,0} = N_0,N_{1,\lambda \times \delta_1} =
N_1$ 
\sn
\itemitem{ $(\beta)$ }  $N_{2,0} = N_2$
\sn
\itemitem{ $(\gamma)$ }  $\langle N_{\ell,i}:i \le \lambda \times \delta_1 
\rangle$ is $<_{\frak K}$-increasing continuous for $\ell = 1,2$
\sn
\itemitem{ $(\delta)$ }  $(N_{1,i},N_{1,i+1}) \cong (M^*,N^*)$
\sn
\itemitem{ $(\varepsilon)$ }  $N_{2,i} \cap N_1 = N_{1,i}$
\sn
\itemitem{ $(\zeta)$ }  $N_3$ is $(\lambda$,cf$(\delta_3)$)-saturated
over $N_{2,\lambda \times \delta_1}$; if $\delta_3 = 1$ this means just
\newline
$N_{2,\lambda \times \delta_1} \le_{\frak K} N_3$
\endroster}
\endroster
\enddefinition
\bn
\ub{Discussion}:  Why this definition of NF?  We need a nonforking notion with
the usual properties.  We first describe a version depending on $\langle
\delta_0,\delta_1,\delta_2 \rangle$ and get $NF = NF_{\lambda,\bar \delta},
\bar \delta$ works like a scaffold -- eventually $\bar \delta$ disappears.
\bigskip

\definition{\stag{8.2} Definition}  1) We say
$\nonforkin{N_1}{N_2}_{N_0}^{N_3}$ (or
$N_1,N_2$ are smoothly amalgamated over $N_0$ inside $N_3$ or
$NF_\lambda(N_0,N_1,N_2,N_3)$) \underbar{if} we can find $M_\ell \in
K_\lambda$ (for $\ell < 4$) such that:
\medskip
\roster
\item "{$(a)$}"  $NF_{\lambda,\langle \lambda,\lambda,\lambda \rangle}
(M_0,M_1,M_2,M_3)$
\sn
\item "{$(b)$}"  $N_\ell \le_{\frak K} M_\ell$ for $\ell < 4$
\sn
\item "{$(c)$}"  $N_0 = M_0$
\sn
\item "{$(d)$}"  $M_1,M_2$ are $(\lambda,\text{cf}(\lambda))$-saturated 
over $N_0$ (follows by (a) see clauses (d), (e) of \scite{8.1}).
\endroster
\enddefinition
\bigskip

\proclaim{\stag{8.3} Claim}  1) If $\bar \delta = \langle \delta_1,\delta_2,
\delta_3 \rangle,\delta_\ell$ a limit ordinal $< \lambda^+$ and
$N_\ell \in K_\lambda$ for $\ell < 3$, and $N_1$ is $(\lambda$,
cf$(\delta_1)$)-
saturated over $N_0$ and $N_2$ is $(\lambda$,cf$(\delta_2)$)-saturated
over $N_0$ and $N_0 \le_{\frak K} N_1,N_0 \le_{\frak K} N_2$ 
and for simplicity
$N_1 \cap N_2 = N_0$.  \underbar{Then} we can find $N_3$ such that
$NF_{\lambda,\bar \delta}(N_0,N_1,N_2,N_3)$. \newline
2) Moreover, we can choose any $\langle N_{1,i}:i \le \lambda \times
\delta_1 \rangle$ as in \scite{8.1}(f)$(\alpha),(\gamma),(\delta)$ 
as part of the witness.
\endproclaim
\bigskip

\demo{Proof}  Straight (remembering \scite{7.4A} (and uniqueness of the
$(\lambda,\text{cf}(\delta_1))$-saturated model over $N_0$)).
\hfill$\square_{\scite{8.3}}$
\enddemo
\bigskip

\proclaim{\stag{8.4} Claim}  If Definition \scite{8.1}, if $\delta_3$ is a
limit ordinal, \ub{then} without loss of generality
(even without changing $\langle N_{1,i}:i \le \lambda \times \delta_1
\rangle$)
\medskip
\roster
\item "{$(g)$}"  $N_{2,i+1}$ is $(\lambda$,cf$(\delta_2)$)-saturated over
$N^1_{i+1} \cup N^2_i$ (which means it is \newline
$(\lambda$,cf$(\delta_2)$)-saturated over some $N$, 
where $N^1_{i+1} \cup N^2_i \subseteq N \le_{\frak K} N_{2,i+1}$).
\endroster
\endproclaim
\bigskip

\demo{Proof}  So assume $NF_{\lambda,\bar \delta}(N_0,N_1,N_2,N_3)$ holds
as witnessed by $\langle N_{\ell,i}:i \le \lambda \times \delta_\ell \rangle$
for $\ell = 1,2$.  Now we choose by induction on $i \le \lambda \times
\delta_1$ a model $M_{2,i} \in K_\lambda$ such that:
\medskip
\roster
\widestnumber\item{ (iii) }
\item "{$(i)$}"  $N_{2,i} \le M_{2,i}$
\sn
\item "{$(ii)$}"  $M_{2,0} = N_2$
\sn
\item "{$(iii)$}"  $M_{2,i}$ is $\le_{\frak K}$-increasing continuous
\sn
\item "{$(iv)$}"  $M_{2,i} \cap N_{2,\lambda \times \delta_1} = N_{2,i}$
moreover $M_{2,i} \cap N_3 = N_{2,i}$
\sn
\item "{$(v)$}"  $M_{2,i+1}$ is $(\lambda$,cf$(\delta_2)$)-saturated over
$M_{2,i} \cup N_{2,i+1}$.
\endroster
\medskip

There is no problem to carry the definition.  Let $M_3$ be such that
$M_{2,\lambda \times \delta_1} \le_{\frak K} M_3 \in K_\lambda$ and $M_3$ is
$(\lambda$,cf$(\delta_3)$)-saturated over $M_{2,\lambda \times \delta_1}$.
So both $M_3$ and $N_3$ are $(\lambda$,cf$(\delta_3)$)-saturated over
$N_{2,\lambda \times \delta_1}$, hence they are isomorphic over
$N_{2,\lambda \times \delta_1}$ so let $f$ be an isomorphism from $M_3$
onto $N_3$ which is the identity over $N_{2,\lambda \times \delta_1}$.
\newline
Clearly $\langle N_{1,i}:i \le \lambda \times \delta_1 \rangle,
\langle f(M_{2,i}):i \le \lambda \times \delta_1 \rangle$ are also 
witnesses for \newline
$NF_{\lambda,\bar \delta}(N_0,N_1,N_2,N_3)$ satisfying the extra demand
$(g)$.  \hfill$\square_{\scite{8.4}}$
\enddemo
\bigskip

\proclaim{\stag{8.5} Claim}  (Weak Uniqueness).

Assume that for $x \in \{ a,b\}$, we have $NF_{\lambda,\bar \delta^x}
(N^x_0,N^x_1,N^x_2,N^x_3)$ holds as witnessed by
$\langle N^x_{1,i}:i \le \lambda \times \delta^x_1 \rangle,\langle N^x_{2,i}:
i \le \lambda \times \delta^x_1 \rangle$ and $\delta_1 =: \delta^a_1 = 
\delta^b_1$,
cf$(\delta^a_2) = \text{ cf}(\delta^b_2)$ and \newline
cf$(\delta^a_3) = \text{ cf}(\delta^b_3) \ge \aleph_0$.
\medskip

Suppose further that $f_\ell$ is an isomorphism from $N^a_\ell$ onto
$N^b_\ell$ for $\ell = 0,1,2$, moreover: $f_0 \subseteq f_1,f_0 \subseteq
f_2$ and $f_1(N^a_{1,i}) = N^b_{1,i}$.

\underbar{Then} we can find an isomorphism $f$ from $N^a_3$ onto $N^b_3$
extending $f_1 \cup f_2$.
\endproclaim
\bigskip

\demo{Proof}  Without loss of generality $N^x_{2,i+1}$ is $(\lambda,\text{cf}
(\delta_2))$-saturated over $N^x_{1,i+1} \cup N^x_{2,i}$ (by \scite{8.4}, note
``without changing the $N_{1,i}$'s).  Now we choose by induction on \newline
$i \le \lambda \times \delta_1$ an isomorphism $g_i$ from $N^a_{2,i}$ onto 
$N^b_{2,i}$ such that: $g_i$ is increasing in $i$ and $g_i$ extends 
$(f_1 \restriction N^a_{1,i}) \cup f_2$. \newline
For $i = 0$ choose $g_0 = f_2$ and for $i$ limit let $g_i$ be
$\dsize \bigcup_{j < i} g_j$ and for $i = j+1$ use $(N_{1,i},N_{1,i+1})
\cong (M^*,N^*)$ (see \scite{8.1}) and the extra saturation clause (g).  
Now we can extend $g_{\lambda \times \delta_1}$ to an isomorphism from 
$N^a_3$ onto $N^b_3$ as $N^x_3$ is $(\lambda,\text{cf}(\delta_3))$-saturated 
from $N^x_{2,\lambda \times \delta_1}$ (for $x \in \{ a,b\}$); note that
knowing \scite{8.5} possibly the choice of \newline
$\langle N_{1,i}:i \le \lambda \times \delta_1 \rangle$ matters.  
\hfill$\square_{\scite{8.5}}$
\enddemo
\bigskip

\noindent
Now we prove an ``inverted" uniqueness
\proclaim{\stag{8.6} Claim}  Suppose that 
\medskip
\roster
\item "{$(a)$}"  for $x \in \{ a,b\}$ we have 
$NF_{\lambda,\bar \delta^x}(N^x_0,N^x_1,N^x_2,N^x_3)$
\sn
\item "{$(b)$}"  $\bar \delta^x = \langle \delta^x_1,\delta^x_2,\delta^x_3
\rangle,\delta^a_1 = \delta^b_2$, 
$\delta^a_2 = \delta^b_1$,cf$(\delta^a_3) = \text{ cf}(\delta^b_3)$ all
limit ordinals
\sn
\item "{$(c)$}"  $f_0$ is an isomorphism from $N^a_0$ onto $N^b_0$
\sn
\item "{$(d)$}"  $f_1$ is an isomorphism from $N^a_1$ onto $N^b_2$
\sn
\item "{$(e)$}"  $f_2$ is an isomorphism from $N^a_2$ onto $N^b_1$
\sn
\item "{$(f)$}"  $f_0 \subseteq f_1$ and $f_0 \subseteq f_2$.
\endroster
\medskip

\noindent
\underbar{Then} there is an isomorphism from $N^a_3$ onto $N^b_3$ extending
$f_1 \cup f_2$.
\endproclaim
\bigskip

\noindent
Before proving
\proclaim{\stag{8.6A} Subclaim}  1) For any limit ordinals 
$\delta^a_1,\delta^a_2,\delta^a_3 \le \lambda$ we can find $M_{i,j}$ 
(for \newline
$i \le \lambda \times \delta^a_1 \text{ and } j \le \lambda \times 
\delta^a_2)$ and $M_3$ such that:
\medskip
\roster
\item "{$(A)$}"  $M_{i,j} \in K_\lambda$
\sn
\item "{$(B)$}"  $i_1 \le i_2 \and j_1 \le j_2 \Rightarrow M_{i_1,j_1} 
\le_{\frak K} M_{i_2,j_2}$
\sn
\item "{$(C)$}"  if $i \le \lambda \times \delta_1$ is a limit ordinal and
$j \le \lambda \times \delta_2$ then $M_{i,j} = \dsize \bigcup_{\zeta < i}
M_{\zeta,j}$
\sn
\item "{$(D)$}"  if $i \le \lambda \times \delta_1$ and $j \le 
\lambda \times \delta_2$ is a limit ordinal then $M_{i,j} = \dsize \bigcup
_{\xi < j} M_{i,\xi}$
\sn
\item "{$(E)$}"  $(M_{0,j},M_{0,j+1}) \cong (M^*,N^*)$
\sn
\item "{$(F)$}"  $(M_{i,0},M_{i+1,0}) \cong (M^*,N^*)$
\sn
\item "{$(G)$}"  for $i_1,i_2 \le \lambda \times \delta_1$ and $j_1,j_2 \le
\lambda \times \delta_2$ we have \newline
$M_{i_1,j_1} \cap M_{i_2,j_2} =
M_{\text{min}\{i_1,i_2\},\text{min}\{j_1,j_2\}}$,
\sn
\item "{$(H)$}"  $M_{\lambda \times \delta_1,\lambda \times \delta_2}
\le_{\frak K} M_3 \in K_\lambda$ moreover \newline
$M_3$ is $(\lambda,\text{cf}(\delta_3))$-saturated over $M_{\lambda \times
\delta_1,\lambda \times \delta_2}$.
\endroster
\medskip

\noindent
2) Moreover, it is O.K. if $\langle M_{0,j}:j \le \lambda \times
\delta^a_2 \rangle,\langle M_{i,0}:i \le \lambda \times \delta^a_1 \rangle$
are pregiven as long as both $\le_{\frak K}$-increasing continuous in 
${\frak K}_\lambda$ satisfying (E) + (F) and $M_{0,\lambda \times \delta^a_2} 
\cap M_{\lambda \times \delta^a_2,0} = M_{0,0}$.
\endproclaim
\bigskip

\demo{Proof}  1) This is done by induction on $i$ and for fixed $i$ by 
induction on $j$.  For $i=0, j=0$ let $M_{0,0} \in K_\lambda$ be arbitrary,
for $i=0,j$ limit use clause $(D)$, for $i=0,j = \xi + 1$ use clause $(E)$.
For $i$ limit use clause $(C)$ (and check).  For $i = \zeta + 1$, if
$j=0$ use clause $(F)$, if $j$ limit use clause $(D)$ and if $j = \xi + 1$
use the existence of disjoint amalgamation (i.e. \scite{5.7}).

Lastly, choose $M_3 \in K_\lambda$ which is $(\lambda,\text{cf}
(\delta^a_3))$-saturated over
$M_{\lambda \times \delta^a_1,\lambda \times \delta^a_2}$. \newline
2) Similarly.  \hfill$\square_{\scite{8.6A}}$
\enddemo
\bigskip

\demo{Proof of \scite{8.6}}  Let $M_{i,j},M_3$ be as in \scite{8.6A}.
For $x \in \{ a,b\}$ as 
$NF_{\lambda,\bar \delta^x}(N^x_0,N^x_1,N^x_2,N^x_3)$, we know that 
there are witnesses
$\langle N^x_{1,i}:i \le \lambda \times \delta^x_1 \rangle,\langle N^x_{2,i}:
i \le \lambda \times \delta^x_1 \rangle$ for this, so \newline
$\langle N^x_{1,i}:i
\le \lambda \times \delta^x_1 \rangle$ is 
$\le_{\frak K}$-increasing continuous
and $(N^x_{1,i},N^x_{1,i+1}) \cong (M^*,N^*)$.  \newline
So $\langle N^a_{1,i}:i \le \lambda \times \delta^a_1 \rangle$ is
$\le_{\frak K}$-increasing continuous sequences with each
successive pair isomorphic to $(M^*,N^*)$ hence by \scite{8.6A}(2) without
loss of generality there is an isomorphism $g_1$ from 
$N^a_{1,\lambda \times \delta^a_1}$ onto $M_{\lambda \times 
\delta^a_1}$, mapping $N^a_{1,i}$ onto $M_{i,0}$; remember
$N^a_{1,\lambda \times \delta^a_1} = N^a_1$.  Let $g_0 = g_1 \restriction
N^a_0 = g_1 \restriction N^a_{1,0}$ so $g_0 \circ f^{-1}_0$ is an isomorphism
from $N^b_0$ onto $M_{0,0}$.
\medskip

As $\delta^b_1 = \delta^a_2$, using \scite{8.6A}(2) fully without loss of
generality there is an isomorphism $g_2$ from
$N^b_{1,\lambda \times \delta^a_2}$ onto $M_{0,\lambda \times \delta^a_2}$
mapping $N^b_{1,j}$ onto $M_{0,j}$ (for $j \le \lambda \times \delta^a_2)$
and $g_2$ extends $g_0 \circ f^{-1}_0$. \newline
Now we want to use the weak uniqueness \scite{8.5} and for this note:
\medskip
\roster
\item "{$(\alpha)$}"  $NF_{\lambda,\bar \delta^a}(N^a_0,N^a_1,N^a_2,N^a_3)$
as witnessed by $\langle N^a_{1,i}:i \le \lambda \times \delta^a_1 \rangle$,
\newline
$\langle N^a_{2,i}:i \le \lambda \times \delta_1 \rangle$ \newline
[why?  an assumption]
\sn
\item "{$(\beta)$}"  $NF_{\lambda,\bar \delta^a}(M_{0,0},M_{\lambda \times 
\delta^a_1,0},M_{0,\lambda \times \delta^a_2},M_3)$ as witnessed by the
sequences \newline
$\langle M_{i,0}:i \le \lambda \times \delta_1 \rangle,\langle
M_{i,\lambda \times \delta^a_2}:i \le \lambda \times \delta^a_2 \rangle$
\newline
[why? check]
\sn
\item "{$(\gamma)$}"  $g_0$ is an isomorphism from $N^a_0$ onto $M_{0,0}$
\newline
[why?  see its choice]
\sn
\item "{$(\delta)$}"  $g_1$ is an isomorphism from $N^a_1$ onto
$M_{\lambda \times \delta^a_1,0}$ mapping $N^a_{1,i}$ onto $M_{i,0}$ for
$i \le \lambda \times \delta^a_1$ and extending $g_0$ \newline
[why?  see the choice of $g_1$ and of $g_0$]
\sn
\item "{$(\varepsilon)$}"  $g_2 \circ f_2$ is an isomorphism from
$N^a_2$ onto $M_{0,\lambda \times \delta^a_2}$ extending $g_0$ \newline
[why?  $f_2$ is an isomorphism from $N^a_2$ onto $N^b_1$ and $g_2$ is an
isomorphism from $N^b_1$ onto $M_{0,\lambda \times \delta^a_1}$ extending
$g_0 \circ f^{-1}_0$ and $f_0 \subseteq f_2$].
\endroster
\medskip

So there is by \scite{8.5} an isomorphism $g^a_3$ from $N^a_3$ onto 
$M_3$ extending $g_1$ and $g_2 \circ f_2$.
\medskip

We next want to apply \scite{8.5} to the $N^b_1$'s; so note:
\medskip
\roster
\item "{$(\alpha)'$}"  $NF_{\lambda,\bar \delta^b}(N^b_0,N^b_1,N^b_2,N^b_3)$
as witnessed by $\langle N^b_{1,i}:i \le \lambda \times \delta^a_2 \rangle$,
\newline
$\langle N^b_{2,i}:i \le \lambda \times \delta^a_2 \rangle$
\sn
\item "{$(\beta)'$}"  $NF_{\lambda,\bar \delta^b}(M_{0,0},M_{0,\lambda \times 
\delta^a_2},M_{\lambda \times \delta^a_1,0},M_3)$ as witnessed by the
sequences \newline
$\langle M_{0,j}:j \le \lambda \times \delta^a_2 \rangle,
\langle M_{\lambda \times \delta^a_1,j}:j \le \lambda \times 
\delta^a_1 \rangle$
\sn
\item "{$(\gamma)'$}"  $g_0 \circ (f_0)^{-1}$ is an 
isomorphism from $N^b_0$ onto $M_{0,0}$ \newline
[why?  Check.]
\sn
\item "{$(\delta)'$}"  $g_2$ is an isomorphism from $N^b_1$ onto
$M_{0,\lambda \times \delta^a_2}$ mapping $N^b_{1,j}$ onto $M_{0,j}$ for
$j \le \lambda \times \delta^a_2$ and extending $g_0 \circ (f_1)^{-1}$
\newline
[why?  see the choice of $g_2$: it maps $N^b_{1,j}$ onto $M_{0,j}$]
\sn
\item "{$(\varepsilon)'$}"  $g_1 \circ (f_1)^{-1}$ is an isomorphism from
$N^b_2$ onto $M_{\lambda \times \delta^a_0}$ extending $g_0$ \newline
[why?  remember $f_1$ is an isomorphism from $N^a_1$ onto $N^b_2$ extending
$f_0$ and the choice of $g_1$: it maps $N^a_1$ onto $M_{\lambda \times 
\delta^a_{1,0}}$].
\endroster
\medskip

\noindent
So there is an isomorphism $g^b_3$ form $N^b_3$ onto $M_3$ extending
$g_2,f_1 \circ (f_1)^{-1}$. \newline
Lastly $(g^b_3)^{-1} \circ g^a_3$ is an isomorphism from $N^a_3$ onto
$N^b_3$ (chase arrows). Also

$$
\align
((g^b_b)^{-1} \circ g^a_3) \restriction N^a_1 &= (g^b_3)^{-1}(g^a_3 
\restriction N^a_1) \\
  &= (g^b_3)^{-1} g_1 = ((g^b_3)^{-1} \restriction M_{\lambda \times
\delta^a_{1,0}}) \circ g_1 \\
  &= (g^b_3 \restriction N^b_2)^{-1} \circ g_1 = ((g_1 \circ (f_1)^{-1})^{-1})
\circ g_1 \\
  &= (f_1 \circ (g_1)^{-1}) \circ g_1 = f_1.
\endalign
$$
\medskip

\noindent
Similarly  $((g^b_3)^{-1} \circ g^a_3) \restriction N^a_2 = f_2$.
\newline
So we have finished. \hfill$\square_{\scite{8.6}}$
\enddemo
\bigskip

\proclaim{\stag{8.7} Claim}  [Uniqueness].  Assume for 
$x \in \{ a,b\}$ we have \newline
$NF_{\lambda,\bar \delta^x}(N^x_0,N^x_1,N^x_2,N^x_3)$ and
cf$(\delta^a_1) = \text{ cf}(\delta^b_1),\text{cf}(\delta^a_2) 
= \text{ cf}(\delta^b_2),\text{cf}(\delta^a_3) = \text{ cf}(\delta^b_3)$,
all $\delta^x_\ell$ limit ordinals.

If $f_\ell$ is an isomorphism from $N^a_\ell$ onto $N^b_\ell$ for $\ell < 3$
and $f_0 \subseteq f_1,f_0 \subseteq f_2$ \underbar{then} there is an
isomorphism $f$ from $N^a_3$ onto $N^b_3$ extending $f_1,f_2$.
\endproclaim
\bigskip

\demo{Proof}  Let $\bar \delta^c = \langle \delta^c_1,\delta^c_2,
\delta^c_3 \rangle = \langle \delta^a_2,\delta^a_1,\delta^a_3 \rangle$; 
by \scite{8.3} there are
$N^c_\ell$ (for $\ell \le 3$) such that $NF_{\lambda,\bar \delta^c}(N^c_0,
N^c_1,N^c_2,N^c_3)$.  There is for $x \in \{ a,b\}$ an isomorphism $g^x_0$
from $N^a_0$ onto $N^c_0$ (as $K_\lambda$ is categorical in $\lambda$) and
without loss of generality $g^b_0 = g^a_0 \circ f_0$.  Similarly
for $x \in \{a,b\}$
there is an isomorphism $g^x_1$ from $N^x_1$ onto $N^c_2$ extending $g^x_0$
(as $N^x_1$ is
$(\lambda,\text{cf}(\delta^x_1))$-saturated over $N^x_0$ and also
$N^c_2$ is $(\lambda,\text{cf}(\delta^c_2))$-saturated over $N^c_0$ and
$\text{cf}(\delta^c_2) = \text{cf}(\delta^a_1) = \text{cf}(\delta^x_1))$ 
and without loss of genrality $g^b_1 = g^a_1 \circ f_1$).  Similarly for
$x \in \{a,b\}$ 
there is an isomorphism $g^x_2$ from $N^x_2$ onto $N^c_1$ extending $g^x_0$
(as $N^x_2$ is
$(\lambda,\text{cf}(\delta^x_2))$-saturated over $N^x_0$ and also $N^c_1$
is $(\lambda,\text{cf}(\delta^c_1))$-saturated over $N^c_0$ and
$\text{cf}(\delta^c_1) = \text{cf}(\delta^a_2) = \text{cf}(\delta^x_2))$
and without loss of generality $g^b_2 = g^b_2 \circ f_2$.
\newline
So by \scite{8.6} for $x \in \{a,b\}$ there is an isomorphism 
$g^x_3$ from $N^x_3$ onto $N^c_3$
extending $g^x_1$ and $g^x_2$.  
Now $(g^b_3)^{-1} \circ g^a_3$ is an isomorphism
from $N^a_3$ onto $N^b_3$ extending $f_1,f_2$ as required.
\hfill$\square_{\scite{8.7}}$
\enddemo
\bigskip

\demo{\stag{8.8} Conclusion}  [Symmetry]. \newline

If $NF_{\lambda,\langle \delta_1,\delta_2,\delta_3 \rangle}(N_0,N_1,N_2,N_3)$
\underbar{then} 
$NF_{\lambda,\langle \delta_2,\delta_1,\delta_3 \rangle}(N_0,N_2,N_1,N_3)$.
\enddemo
\bigskip

\demo{Proof}  By \scite{8.6} (and \scite{8.7}).
\enddemo
\bigskip

\proclaim{\stag{8.9} Claim}  1) In Definition \scite{8.1} we can replace
$(N_{1,i},N_{1,i+1}) \cong (M^*,N^*)$ by \newline
$(N_{1,i},N_{1,i+1}) \in K^{3,\text{uq}}$.
\endproclaim
\bigskip

\demo{Proof}  Like the proof of \scite{8.6} (get $(M_{i,0},M_{i+1,0}) \cong
(M^*,N^*)$, \newline
$(M_{0,j},M_{0,j+1}) \in K^{3,\text{uq}}_\lambda)$, but as we shall
not use it, we do not elaborate. \hfill$\square_{\scite{8.9}}$
\medskip

Now we turn to smooth amalgamation (not necessarily saturated, see
Definition \scite{8.2}).
\enddemo
\bigskip

\proclaim{\stag{8.10} Claim}  1) If $NF_{\lambda,\bar \delta}(N_0,N_1,N_2,
N_3)$ each $\delta_\ell$ limit \underbar{then} 
$NF_\lambda(N_0,N_1,N_2,N_3)$ (see Definition \scite{8.2}). \newline
2) In Definition \scite{8.2} we can add:
\medskip
\roster
\item "{$(d)^+$}"  $M_\ell$ is $(\lambda,\text{cf}(\lambda))$-saturated 
over $N_0$ and moreover over $N_\ell$,
\sn
\item "{$(e)$}"    $M_3$ is $(\lambda,\text{cf}(\lambda))$-saturated over 
$M_1 \cup M_2$ (actually this is given by $(f)(\zeta)$ of Definition 
\scite{8.1}). 
\endroster
\endproclaim
\bigskip

\demo{Proof}  1) By \scite{8.6A} we can find $M_{i,j}$ for 
$i \le \lambda \times (\delta_1 + \lambda),j \le \lambda \times 
(\delta_2 + \lambda)$ for $\bar \delta' =: \langle \delta_1 + \lambda,
\delta_2 + \lambda,\delta_3 \rangle$ and choose $M'_3 \in K_\lambda$ which 
is $(\lambda,\text{cf}(\delta_3))$-saturated over $M_{\lambda \times 
\delta_1,\lambda \times \delta_2}$.  So 
$NF_{\lambda,\bar \delta}(M_{0,0},M_{\lambda \times
\delta_1,0},M_{0,\lambda \times \delta_2},M_3)$, hence by \scite{8.7} 
without loss of generality 
$M_{0,0} = N_0,M_{\lambda \times \delta_1,0} = N_1,
M_{0,\lambda \times \delta_2} = N_2$, and $N_3 = M'_3$.  Lastly, let $M_3$
be $(\lambda,\text{cf}(\lambda))$-saturated over $M'_3$.
Now clearly also \newline
$NF_{\lambda,\langle \delta_1 + \lambda,\delta_2 + \lambda,\delta_3 +
\lambda \rangle}(M_{0,0},M_{\lambda \times (\delta_1 + \lambda),0},
M_{0,\lambda \times (\delta_2 + \lambda)},M_3)$ and \newline
$N_0 = M_{0,0},
N_1 = M_{\lambda \times \delta_2,0} \le_{\frak K} 
M_{\lambda \times (\delta_2 + 
\lambda),0},N_2 = M_{0,\lambda \times \delta_2} \le_{\frak K} 
M_{0,\lambda \times
(\delta_2 + \lambda)}$ \newline
and $M_{\lambda \times (\delta_1 + \lambda),0}$ is
$(\lambda,\text{cf}(\lambda))$-saturated over 
$M_{\lambda \times \delta_1,0}$ and
$M_{0,\lambda \times (\delta_2 + \lambda)}$ is \newline
$(\lambda,\text{cf}(\lambda))$-saturated over $M_{0,\lambda \times \delta_2}$ 
and $N_3 = M'_3 \le_{\frak K} M_3$.  So
we get all the requirements for $NF_\lambda(N_0,N_1,N_2,N_3)$ (as witnessed
by $\langle M_{0,0},M_{\lambda \times (\delta_1 + \lambda),0},
M_{0,\lambda \times (\delta_2 + \lambda)},M_3 \rangle$). \newline
2) Similar proof.  \hfill$\square_{\scite{8.10}}$
\enddemo
\bigskip

\proclaim{\stag{8.11} Claim}  [Uniqueness of smooth amalgamation]:  If
$NF_\lambda(N^x_0,N^x_1,N^x_2,N^x_3)$ for $x \in \{ a,b\},f_\ell$ an
isomorphism from $N^a_\ell$ onto $N^b_\ell$ for $\ell < 3$ and
$f_0 \subseteq f_1,f_0 \subseteq f_2$ \underbar{then} $f_1 \cup f_2$ can be
extended to a $\le_{\frak K}$-embedding of $N^a_3$ into some 
$\le_{\frak K}$-extension of $N^b_3$
(so if $N^x_3$ is $(\lambda,\kappa)$-saturated 
over $N^x_1 \cup N^x_2$ for $x =
a,b$, we can extend $f_1 \cup f_2$ to an isomorphism from $N^a_3$ onto
$N^b_3$.
\endproclaim
\bigskip

\demo{Proof}  For $x \in \{ a,b\}$ let the sequence $\langle M^x_\ell:
\ell < 4 \rangle$ be a witness to \newline
$NF_\lambda(N^x_0,N^x_1,N^x_2,N^x_3)$ as in \scite{8.2}, \scite{8.10}(2) so
in particular \newline
$NF_{\lambda,\langle \lambda,\lambda,\lambda \rangle}
(M^x_0,M^x_1,M^x_2,M^x_3)$.  By chasing arrows and uniqueness, i.e. 
\scite{8.6} without loss of generality $M^a_\ell = M^b_\ell$ for $\ell < 4$ 
and $f_0 = \text{ id}_{N^a_0}$.  As $M^a_1$ is 
$(\lambda,\text{cf}(\lambda))$-saturated over $N^a_1$ and also
over $N^b_1$ and $f_1$ is an isomorphism from $N^a_1$ onto $N^b_1$, clearly
there is an automorphism $g_1$ of $M^a_1$ such that $f_1 \subseteq g_1$,
hence also $\text{id}_{N^a_0} = f_0 \subseteq f_1 \subseteq g_1$.  
Similarly there is an
automorphism $g_2$ of $M^a_1$ extending $f_1$ hence $f_0$.  So $g_\ell \in
\text{ AUT}(M^a_\ell)$ for $\ell = 1,2$ and
$g_1 \restriction M^a_0 = f_0 =
g_2 \restriction M^a_0$.  By the uniqueness of $NF_{\lambda,\langle
\lambda,\lambda,\lambda \rangle}$ (i.e. Claim \scite{8.7}) there is an 
automorphism $g_4$ of
$M^a_4$ extending $g_1 \cup g_2$.  This proves the desired conclusion.
\hfill$\square_{\scite{8.11}}$
\enddemo
\bigskip

\proclaim{\stag{8.12} Claim}  Assume
\medskip
\roster
\item "{$(a)$}"  $\bar \delta = \langle \delta_1,\delta_2,\delta_3 \rangle,
\delta_\ell < \lambda^+$ is a limit ordinal for $\ell = 1,2,3$; \newline
$N_0 \le_{\frak K} N_\ell \le_{\frak K} N_3$ for $\ell = 1,2$ and
\sn
\item "{$(b)$}"  $N_\ell$ is $(\lambda,\text{cf}(\delta_\ell))$-saturated
over $N_0$
\sn
\item "{$(c)$}"  $N_3$ is cf$(\delta_3)$-saturated over $N_1 \cup N_2$.
\endroster
\medskip

\noindent
\underbar{Then}  $NF_\lambda(N_0,N_1,N_2,N_3)$ iff $NF_{\lambda,\bar \delta}
(N_0,N_1,N_2,N_3)$.
\endproclaim
\bigskip

\demo{Proof}  The ``if" direction holds by \scite{8.11}.  For the ``only if"
direction, by the proof of \scite{8.10}(1) (and Definition \scite{8.1}, 
\scite{8.2}) we can find $M_\ell(\ell \le 3)$ such that 
$NF_{\lambda,\bar \delta}(M_0,M_1,M_2,M_3)$ and clauses (b), (c), (d) of
Definition \scite{8.2} holds so by \scite{8.10} also
$NF_\lambda(M_0,M_1,M_2,M_3)$.  Easily there are for $\ell < 3$, an
isomorphism $f_\ell$ from $M_\ell$ onto $N_\ell$ such that 
$f_0 = f_\ell
\restriction M_\ell$.  By the uniqueness for smooth amalgamation 
(i.e. \scite{8.11}) we can find an isomorphism $f_3$ from $M_3$ 
onto $N_3$ extending $f_1 \cup f_2$.  So as
$NF_{\lambda,\bar \delta}(M_0,M_1,M_2,M_3)$ holds also $NF_{\lambda,
\bar \delta},(f_0(M_0),f_3(M_1),f_3(M_2),f_3(M_3))$; \newline
i.e. $NF_{\lambda,\bar \delta}(N_0,N_1,N_2,N_3)$ as required. 
\hfill$\square_{\scite{8.12}}$
\enddemo
\bigskip

\proclaim{\stag{8.13} Claim} [Monotonicity]:  If 
$NF_\lambda(N_0,N_1,N_2,N_3)$ and
$N_0 \le_{\frak K} N'_1 \le_{\frak K} N_1$ and $N_0 \le_{\frak K} 
N'_2 \le_{\frak K} N_2$ and $N'_1 \cup N'_2
\subseteq N'_3 \le_{\frak K} N_3$ 
\underbar{then} $NF_\lambda(N_0,N'_1,N'_2,N'_3)$.
\endproclaim
\bigskip

\demo{Proof}  Read Definition \scite{8.2}.
\enddemo
\bigskip

\proclaim{\stag{8.14} Claim}  [Symmetry]:  $NF_\lambda(N_0,N_1,N_2,N_3)$ holds
\underbar{if and only if} \newline
$NF_\lambda(N_0,N_2,N_1,N_3)$ holds.
\endproclaim
\bigskip

\demo{Proof}  By Claim \scite{8.8} (and Definition \scite{8.2}).
\enddemo
\bigskip

\proclaim{\stag{8.15} Claim}  Assume $\alpha < \lambda^+$ is an ordinal and
for $x \in \{a,b,c\}$ the sequence $\langle N^x_i:i \le \alpha \rangle$ is a
$\le_{\frak K}$-increasing sequence of members of $K_\lambda$,
for $x = a,b$ the sequence is $\le_{\frak K}$-increasing continuous for
$i \le \alpha,N^b_i \cap N^a_\alpha = N^a_i,N^c_i \cap N^a_\alpha = N^a_i,
N^a_i \le_{\frak K} N^b_i \le_{\frak K} N^c_i$ and
$N^c_i$ is $(\lambda,\kappa_i)$-saturated over $N^b_i$ and \newline 
$\text{NF}_{\lambda,\bar \delta^i}(N^a_i,N^a_{i+1},N^c_i,N^b_{i+1})$
(so $i < \alpha \Rightarrow N^c_i \le_{\frak K} N^b_{i+1}$) where \newline
$\bar \delta^i = \langle
\delta^i_1,\delta^i_2,\delta^i_3 \rangle$ sequence of limit ordinals,
$i < \alpha \Rightarrow \delta^{i+1}_2 = \delta^3_i$, and for $i < 0$ limit,
cf$(\delta^3_i) = \dsize \sum_{j < i} \delta^3_j,\delta_1 =
\dsize \sum_{\beta < \alpha} \delta^1_\beta$ and $\delta_3 = \kappa_\alpha,
\bar \delta = \langle \delta_1,\delta^0_2,\delta_3 \rangle$.  
\underbar{Then} $\text{NF}_{\lambda,\bar \delta}(N^a_0,N^a_\alpha,N^b_0,
N^c_\alpha)$.
\endproclaim
\bigskip

\demo{Proof}  Use uniqueness of \scite{8.7}, lastly use \scite{8.7} to 
show $N^b_\alpha$ is $(\lambda,\text{cf}(\alpha))$-saturated over 
$N^a_\alpha \cup N^b_0$.  \hfill$\square_{\scite{8.15}}$
\enddemo
\bigskip

\proclaim{\stag{8.16} Claim}  Assume that $\alpha < \lambda^+$ and for $x \in
\{ a,b\}$ we have $\langle N^x_i:i \le \alpha \rangle$ is \newline
$\le_{\frak K}$-increasing continuously sequence of members of 
$K_\lambda$. \newline 
1)  If $NF_\lambda(N^a_i,N^a_{i+1},N^b_i,N^b_{i+1})$ for each $i < \alpha$
\underbar{then} $NF_\lambda(N^a_0,N^a_\alpha,N^b_0,N^b_\alpha)$. \newline
2)  If  $\alpha_1 < \lambda^+,\alpha_2 < \lambda^+$ and 
$M_{i,j} \, (i \le \alpha_1,j \le \alpha_2)$ are as in \scite{8.6A}, 
and for each $i < \alpha_1,j < \alpha_2$ we have:

$$
\nonforkin{M_{i,j+1}}{M_{i+1,j}}_{M_{i,j}}^{M_{i+1,j+1}}
$$
\medskip

$$
\text{\underbar{then}  }
\nonforkin{M_{i,0}}{M_{0,j}}_{M_{0,0}}^{M_{\alpha_1,\alpha_2}} \text{ for }
i \le \alpha_1, j \le \alpha_2.
$$
\endproclaim
\bigskip

\demo{Proof}  1) We first prove special cases and use them to prove more
general cases. \newline
\smallskip

\noindent
\underbar{Case A}:  $N^a_{i+1}$ is $(\lambda,\delta^1_i)$-saturated over
$N^a_i$ and $N^b_{i+1}$ is $(\lambda,\delta^2_i)$-saturated over \newline
$N^a_{i+1} \cup N^b_i$ for $i < \alpha$. \newline
\medskip

We can choose by induction on $i \le \alpha,N^c_i \in K_\lambda$ such that
\medskip
\roster
\item "{$(a)$}"  $N^b_0 \le_{\frak K} N^c_0 \le_{\frak K} N^b_{i+1},
N^c_0$ is $(\lambda,\delta^2_0)$-saturated over $N^b_0$, and \newline
$NF_{\lambda,\langle \delta^0_1,\delta^0_2,
\delta^0_2 \rangle}(N^a_0,N^a_1,N^c_0,N^b_1)$
\sn
\item "{$(b)$}"  $N^c_\alpha \in K_\lambda$ is $(\lambda,
\delta^\alpha_3)$-saturated over $N^b_\alpha$.
\endroster
\medskip

\noindent
(Possible by uniqueness; i.e. \scite{8.11} and monotonicity, i.e.
\scite{8.13}).  Now we can use \scite{8.15}.
\bigskip

\noindent
\underbar{Case B}:  
For each $i < \alpha$ we have: $N^a_{i+1}$ is $(\lambda,\kappa_i)$-saturated
over $N^a_i$.  \newline
Let $\bar \delta^i = (\kappa_i,\lambda,\lambda)$. \newline
We can find a $\le_{\frak K}$-increasing sequence
$\langle M^x_i:i \le \alpha \rangle$ for $x \in \{ a,b,c\}$, continuous for
$x=a,b$ such that $i < \alpha \Rightarrow M^b_i \le_{\frak K} 
M^c_i \le_{\frak K} M^b_{i+1}$ and $M^b_\alpha \le_{\frak K} M^c_\alpha$ and
$NF_{\lambda,\bar \delta^i}(M^a_i,M^a_{i+1},M^c_i,M^b_{i+1})$ by
choosing $M^a_i,M^b_i,M^c_i$ by induction on $i$.  By Case A we know that
$NF_\lambda(M^a_0,M^a_\alpha,M^b_0,M^c_\alpha)$ holds.
\medskip

We can now choose an isomorphism $f^a_0$ from $N^a_0$ onto $M^a_0$ (exists
as $K$ is categorical in $\lambda$) and then a $\le_{\frak K}$-embedding 
of $N^b_0$ into $M^b_0$ extending $f^a_0$.  Next we choose by induction on
$i \le \alpha,f^a_i$ an isomorphism from $N^a_i$ onto $M^a_i$ such that:
$j < i \Rightarrow f^a_j \subseteq f^a_i$, possible by ``uniqueness of the
$(\lambda,\kappa_i)$-saturated model over $M^a_i$". (\scite{0.22}) \medskip

Now we choose by induction on $i \le \alpha$, a $\le_{\frak K}$-embedding 
$f^b_i$ of
$N^b_i$ into $M^b_i$ extending $f^a_i$ and $f^b_j$ for $j < i$.  For
$i = 0$ we have done it, for $i$ limit use $\dsize \bigcup_{j < i} f^b_j$, 
lastly for $i$ a successor ordinal let $i = j+1$, now we have
\medskip
\roster
\item "{$(*)_2$}"  $NF_\lambda(M^a_i,M^a_{i+1},f^b_i(N^b_i),M^b_{i+1})$ 
\newline
[why?  because $NF_{\lambda,\bar \delta^i}(M^a_i,M^a_{i+1},M^c_i,M^b_{i+1})$
by the choice of the \newline
$M^x_\zeta$'s hence by \scite{8.12} we have
$NF_\lambda(M^a_i,M^a_{i+1},M^c_i,M^b_{i+1})$ and as \newline
$M^a_i \le_{\frak K} f^b_i(N^b_i) \le M^b_i,M^c_i$ by \scite{8.13} we 
get $(*)_2$].
\endroster
\medskip

By $(*)_2$ and the uniqueness of smooth amalgamation \scite{8.11} 
there is $f^b_i$ as required. \newline
So without loss of generality $f^b_\alpha$ is the identity, so we have
$N^a_0 = M^a_0,N^a_\alpha = M^a_\alpha,N^b_0 \le_{\frak K} M^b_0,N^b_\alpha
\le_{\frak K} M^b_\alpha$; also as said above $NF_\lambda(M^a_0,M^a_\alpha,
M^b_0,M^b_\alpha)$ holds so by monotonicity i.e. \scite{8.13} we get
$NF_\lambda(N^a_0,N^a_\alpha,N^b_0,N^b_\alpha)$ as required.
\medskip

\noindent
\underbar{Case C}:  General case.

We can find $M^\ell_i$ for $\ell < 3,i \le \alpha$ such that:
\medskip
\roster
\item "{$(a)$}"  $M^\ell_i \in K_\lambda$
\sn
\item "{$(b)$}"  for each $\ell < 3,M^\ell_i$ is $\le_{\frak K}$-increasing
in $i$
\sn
\item "{$(c)$}"  $M^0_i = N^a_i$
\sn
\item "{$(d)$}"  $M^{\ell +1}_{i+1}$ is $(\lambda,\lambda)$-saturated over
$M^\ell_{i+1} \cup M^{\ell + 1}_i$ for $\ell < 2,i < \alpha$
\sn
\item "{$(e)$}"  $NF_\lambda(M^\ell_i,M^\ell_{i+1},M^{\ell + 1}_i,
M^{\ell + 1}_i)$ for $\ell < 2,i < \alpha$
\sn
\item "{$(f)$}"  $M^{\ell + 1}_0$ is $(\lambda,\lambda)$-saturated over
$M^\ell_0$
\sn
\item "{$(g)$}"  for $\ell < 2$ and $i < \alpha$ limit we have
$$
M^{\ell +1}_i \text{ is } (\lambda,\lambda) \text{-saturated over }
\dsize \bigcup_{j < i} M^{\ell +1}_j \cup M^\ell_i
$$
\item "{$(h)$}"  for $i < \alpha$ limit we have
$$
NF_\lambda 
\left( \dsize \bigcup_{j < i} M'_j,M^1_i,\dsize \bigcup_{j < i} M^2_j,
M^2_i \right).
$$

\noindent
[how?  as in the proof of \scite{8.6A}].
\endroster
\medskip

\noindent
Now note:
\medskip
\roster
\item "{$(*)_4$}"  $M^{\ell +1}_i$ is $(\lambda,\text{cf}(\lambda \times
(1 + i)))$-saturated over $M^\ell_i$ \newline
[why?  If $i=0$ by clause $(f)$, if $i$ a successor ordinal by clause $(d)$
and if $i$ is a limit ordinal then by clause (g)]
\sn
\item "{$(*)_5$}"  for $i < \alpha,NF_\lambda(M^0_i,M^0_{i+1},M^2_i,M^2
_{i+1})$ \newline
[why?  we use Case B for $\alpha = 2$ with $M^0_i,M^0_{i+1},M^1_i,
M^1_{i+1}$, \newline
$M^2_i,M^2_{i+1}$ here standing for $N^a_0,N^b_0,N^a_1,N^b_1,
N^a_2,N^b_2$ there].
\endroster
\medskip

Now we continue as in Case B (using $f^a_i = \text{ id}_{N^a_i}$ and defining
by induction on $i$ a $\le_{\frak K}$-embedding $f^b_i$ of $N^b_i$ into
$M^c_i$). \newline
2) By part (1) for each $i$ the sequences $\langle M_{\beta,i}:\beta 
\le \alpha_1 \rangle,\langle M_{\beta,i+1}:\beta \le \alpha_1 \rangle$ we
get $\nonforkin{M_{\alpha_1}}{M_{0,i+1}}_{M_{0,i}}^{M_{\alpha_1,i+1}}$ hence
by symmetry (i.e. \scite{8.11}) we have $\nonforkin{M_{0,i+1}}{M_{\alpha_1,i}}
_{M_{0,i}}^{M_{\alpha_1,i+1}}$.  Applying part (1) to the sequences
$\langle M_{0,j}:j \le \alpha_2 \rangle,\langle M_{\alpha_1,j}:j \le \alpha_2
\rangle$ we get $\nonforkin{M_{0,\alpha_2}}{M_{\alpha_1,0}}_{M_{0,0}}
^{M_{\alpha_1,\alpha_2}}$ hence by symmetry (i.e. \scite{8.11}) we have
$\nonforkin{M_{\alpha_1,0}}{M_{0,\alpha_2}}_{M_{0,0}}
^{M_{\alpha_1,\alpha_2}}$; by monotonicity, i.e \scite{8.13} (or restriction
of the matrix) we get the desired conclusion.   \hfill$\square_{\scite{8.16}}$
\enddemo
\bigskip

\demo{\stag{8.17} Conclusion}   Assume $\langle N^\ell_i:i \le \alpha \rangle$
is $\le_{\frak K}$-increasing continuous for $\ell = 0,1$ where $N^\ell_i
\in K_\lambda$ and $N^1_{i+1}$ is $(\lambda,\kappa_\ell)$-saturated over
$N^0_{i+1} \cup N^1_i$ and \newline
NF$_\lambda(N^0_i,N^1_i,N^0_{i+1},N^1_{i+1})$.

Then $N^1_\alpha$ is $(\lambda,\text{cf}(\dsize 
\sum_{i < \alpha} \kappa_i))$-saturated over $N^0_\alpha \cup N^1_0$ (if
$\alpha$ is a limit ordinal, ``$N^1_{i+1}$ is universal over $N^0_{i+1} \cup
N^1_i$" suffice).
\enddemo
\bigskip

\demo{Proof}  The case $\alpha$ not limit is trivial so assume $\alpha$ is a
limit ordinal.  We choose by induction on $i \le_x \alpha$, a sequence
$\langle M'_{i,\varepsilon}:\varepsilon \le \varepsilon(i) \rangle$ such
that:
\medskip
\roster
\item "{$(a)$}"  $\langle M_{i,\varepsilon}:\varepsilon \le \varepsilon(i)
\rangle$ is (strictly) $<_{\frak K}$-increasing continuous
\sn
\item "{$(b)$}"  $N^0_i \le_{\frak K} M_{i,\varepsilon} \le_{\frak K}
N^1_i$
\sn
\item "{$(c)$}"  $N^0_i = M_{i,0}$
\sn
\item "{$(d)$}"  $\varepsilon(i)$ is (strictly) increasing continuous in
$i$,
\sn
\item "{$(e)$}"  $j < i \and \varepsilon \le \varepsilon(j) \Rightarrow
M_{i,\varepsilon} \cap N^1_j = M_{j,\varepsilon}$
\sn
\item "{$(f)$}"  $\varepsilon(0) = 1,M_{i,1} = N^1_0$
\sn
\item "{$(g)$}"  for $i > 0,\lambda$ divide $\varepsilon(i)$
\sn
\item "{$(h)$}"  $N^i_\ell \le_{\frak K} M_{i+1,\varepsilon(i)+1}$.
\endroster
\medskip

\noindent
If we succeed, then $\varepsilon(\alpha)$ is divisible by $\lambda$ and
$\langle M_{i,\varepsilon}:\varepsilon \le \varepsilon(\alpha) \rangle$ is
(strictly) $<_{\frak K}$-increasing continuous, $M_{\alpha,0} = N^0_\alpha$,
and $M_{\alpha,\varepsilon(\alpha)} \le_{\frak K} N^1_\alpha$ but it includes
$N^1_i$ for $i < \alpha$ hence (as $\alpha$ is a limit ordinal) it includes
$\dsize \bigcup_{i < \alpha} N^1_i = N^1_\alpha$; and by \scite{7.4A} we
conclude that $N^1_\alpha = M_{\alpha,\varepsilon(\alpha)}$ is
$(\lambda,\text{cf}(\alpha))$-saturated over $M_{\alpha,1}$ hence over
$N^0_\alpha \cup N^1_0$ (both $\prec M_{\alpha,i}$).

For $i=0$, and $i$ limit there is not much to do.  For $i$ successor we use
\scite{8.18} below.
\enddemo
\bigskip

\demo{\stag{8.18} Conclusion}  1) If $NF_\lambda(N_0,N_1,N_2,N_3)$ and
$\langle M_{0,\varepsilon}:\varepsilon \le \varepsilon(*) \rangle$ is
$\le_{\frak K}$-increasing continuous, $N_0 \le_{\frak K} M_{0,\varepsilon}
\le_{\frak K} M_2$ \underbar{then} we can find $\langle M_{1,\varepsilon}:
\varepsilon \le \varepsilon(*) \rangle$ and $N'_3$ such that:
\medskip
\roster
\item "{$(a)$}"  $N_3 \le_{\frak K} N'_3 \in K_\lambda$
\sn
\item "{$(b)$}"  $\langle M_{1,\varepsilon}:\varepsilon \le \varepsilon(*)
\rangle$ is $\le_{\frak K}$-increasing continuous 
\sn
\item "{$(c)$}"  $M_{1,\varepsilon} \cap N_2 = M_{0,\varepsilon}$
\sn
\item "{$(d)$}"  $N_1 \le_{\frak K} M_{1,\varepsilon} \le_{\frak K} N'_3$
\sn
\item "{$(e)$}"  if $M_{0,0} = N_0$ then $M_{1,0} = N$.
\endroster
\medskip

\noindent
2) if $N_3$ is universal over $N_1 \cup N_2$, then without loss of generality
$N_3 = N_2$.
\enddemo
\bigskip

\demo{Proof}  1) Straight by uniqueness. \newline
2) Follows by (1).  \hfill$\square_{\scite{8.17}},\,\,\square_{\scite{8.18}}$
\enddemo
\newpage

\head {\S9 Nice extensions in $K_{\lambda^+}$} \endhead  \resetall
\bigskip

\demo{\stag{9.0} Hypothesis}  Assume the hypothesis \scite{8.0}.
\enddemo
\bigskip

\noindent
So by \S8 we have reasonable control on \underbar{smooth} amalgamation in
$K_\lambda$.  We use this to define ``nice" extensions in $K_{\lambda^+}$.
This is treated again in \S10.
\bigskip

\definition{\stag{9.1} Definition}  1) Let $M_0 \le^*_{\lambda^+} M_1$ mean:
\medskip
\roster
\item "{$(a)$}"  $M_\ell \in K_{\lambda^+}$, for $\ell = 0,1$
\sn
\item "{$(b)$}"  we can find $\bar M^\ell = \langle M^\ell_i:i < \lambda^+
\rangle$, a representation of $M^\ell$, so $M^\ell_i \in K_\lambda$ \newline
(and
$M^\ell_i$ is $\le_{\frak K}$-increasing continuously and $M_\ell = \dsize
\bigcup_{i < \lambda^+} M^\ell_i$) such that: \newline
$NF_\lambda(M^0_i,M^0_{i+1},
M^1_i,M^1_{i+1}$) for $i < \lambda^+$.
\endroster
\medskip

\noindent
2)  Let $M_0 <^+_{\lambda^+,\kappa} M_1$ mean $M_0 \le^*_{\lambda^+} M_1$ by
some witnesses $M^\ell_i$ (for $i < \lambda^+,\ell < 2$) such that
$NF_{\lambda,\langle \kappa,1,\kappa \rangle}(M^0_i,M^0_{i+1},
M^1_i,M^1_{i+1})$.  If $\kappa = \lambda$, we omit it.
\enddefinition
\bigskip

\proclaim{\stag{9.2} Claim}  1) If $M_0 \le^*_{\lambda^+} M_1$ and 
$\bar M^\ell = \langle M^\ell_i:i < \lambda^+ \rangle$ is a representation 
of $M_\ell$ (as in \scite{8.16}) \underbar{then} for some club $E$ of 
$\lambda^+$, for every $\alpha < \beta$ from $E$ we have $NF_\lambda
(M^0_\alpha,M^0_\beta,M^1_\alpha,M^1_\beta)$. \newline
2)  Similarly for $<^+_{\lambda^+,\kappa}$; if $M_0 <^*_{\lambda^+,
\kappa_\ell}M_1,\bar M^\ell = \langle \bar M^\ell_i:i < \lambda^+ \rangle$
a representation of $M_\ell$ for $\ell = 1,2$ \underbar{then} for some
club $E$ of $\lambda^+$ for every $\alpha < \beta$ from $E$ we have \newline
$NF_{\lambda,\langle \text{cf}(\alpha),1,\text{cf}(\alpha) \rangle}
(M^0_\alpha,M^0_\beta,M^1_\alpha,M^0_\beta)$. \newline
3)  The $\kappa$ in Definition \scite{9.1}(2) does not matter.  In fact, if
$\langle M^\ell_i:i < \lambda^+ \rangle$ are as in \scite{9.1}(1) then for 
some club $E$ of $\lambda^+$ we have: $\alpha \in E \Rightarrow M^1_\alpha
\cap M_0 = M^0_\alpha$ and $\alpha < \beta \and \alpha \in E \and \beta \in
E \Rightarrow [M^1_\beta$ is cf$(\beta)$-saturated over $M^0_\beta \cup
M^1_\alpha]$. 
\endproclaim
\bigskip

\demo{Proof}  1) Straight by \scite{8.16}. \newline
2) Easy. \newline
3) By \scite{8.17}.  (Could have used \scite{7.6}). 
\hfill$\square_{\scite{9.2}}$
\enddemo
\bigskip

\proclaim{\stag{9.3} Claim} 1) For every $\kappa = \text{ cf}(\kappa) 
\le \lambda$ and $M_0 \in K_{\lambda^+}$ for some $M_1 \in K_{\lambda^+}$ 
we have $M_0 <^+_{\lambda^+,\kappa}M_1$. \newline
2)  $<^*_{\lambda^+}$ are $<^+_{\lambda^+,\kappa}$ are transitive. \newline
3)  If $M_0 \le_{\frak K} M_1 \le_{\frak K} M_2$ 
and $M_0 \le^*_{\lambda^+} M_2$
then $M_0 \le^*_{\lambda^+} M_1$. \newline
4)  [transitivity]  If $M_0 \le^*_{\lambda^+} M_1 <^+_{\lambda^+,\kappa} 
M_2$ then $M_0 <^*_{\lambda^+,\kappa} M_2$.
\endproclaim
\bigskip

\demo{Proof}  1) Let $\langle M^0_i:i < \lambda^+ \rangle$ be a representation
of $M_0$ such that $M^0_{i+1}$ is \newline
$(\lambda,\kappa)$-saturated over $M^0_i$.  We choose by induction on 
$i,M^1_i \in K_\lambda$ such that
\newline
$\langle M^1_i:i < \lambda^+ \rangle$ is $<_{\frak K}$-increasing 
continuously, $M^0_i \le_{\frak K} M^1_i,M^1_i \cap M_0 = M^0_i$ and
\newline
$NF_{\lambda,\langle \kappa,1,\kappa \rangle}(M^0_i,M^0_{i+1},
M^1_i,M^1_{i+1})$. We can do it by \scite{7.5}(4).  \nl
2) Concerning $<^+_{\lambda^+,\kappa}$ use \scite{9.2} and \scite{8.16} (i.e. 
transitivity for smooth amalgamations).  Now the proof for
$<^*_{\lambda^+}$ is similar. \newline
3) By monotonicity for smooth amalgamations; i.e. \scite{8.13}. \newline
4) Check. \hfill$\square_{\scite{9.3}}$
\enddemo
\bigskip

\proclaim{\stag{9.4} Claim}  1) If $M_0 \le^*_{\lambda^+} M_\ell$ for
$\ell = 1,2,\kappa = \text{ cf}(\kappa) \le \lambda$ and $a \in M_2 \backslash
M_0$ \underbar{then} for some $M_3$ and $f$ we have: $M_1 <^+_{\lambda^+,
\kappa}M_3$
and $f$ is an $\le_{\frak K}$-embedding of $M_2$ into $M_3$ over $M_0$ with
$f(a) \notin M_1$, moreover, $f(M_2) \le^*_{\lambda^+}M_3$. \newline
2) [uniqueness]  Assume $M_0 <^+_{\lambda^+,\kappa}M_\ell$ for $\ell = 1,2$
\underbar{then} there is an isomorphism $f$ from $M_1$ onto $M_2$ over
$M_0$.
\endproclaim
\bigskip

\demo{Proof}  We first prove part (2). \newline
2) By \scite{9.2}(1) + (2) there are representations $\bar M^\ell = \langle
M^\ell_i:i < \lambda^+ \rangle$ of $M_\ell$ for $\ell < 3$ 
such that: $M^\ell_i \cap M_0 = M^\ell_0$
and $NF_{\lambda,\langle \kappa,\text{cf}(\kappa \times i),\kappa \rangle}
(M^0_i,M^0_{i+1},M^\ell_i,M^\ell_{i+1})$.
\medskip

\noindent
Now we choose by induction on $i < \lambda$ an isomorphism $f_i$ from 
$M^1_i$ onto $M^2_i$, increasing with $i$ and being the identity over
$M^0_i$.  For $i=0$ use ``$M^\ell_0$ is $(\lambda,\kappa)$-saturated over
$M^0_0$ for $\ell = 1,2$" which holds by \scite{8.1}.  
For $i$ limit take unions, for $i$ successor ordinal use 
uniqueness Claim \scite{8.7}.
\enddemo
\bigskip

\demo{Proof of part (1)}  Let $\kappa = \aleph_0$, by \scite{9.3}(1) there 
are for $\ell = 1,2$ models
$N^*_\ell \in K_{\lambda^+}$ such that $M_\ell <^+_{\lambda^+,\kappa}
N^*_\ell$.  Now let $\bar M^\ell = \langle M^\ell_i:i < \lambda^+ \rangle$
be a representation of $M_\ell$ for $\ell = 0,1,2$ and let
$\bar N^\ell = \langle N^\ell_i:i < \lambda^+ \rangle$ be a representation of
$N^*_\ell$ for $\ell = 1,2$.  By \scite{9.3}(4) and \scite{9.2}(3) 
without loss of generality
$NF_\lambda(M^0_i,M^0_{i+1},M^\ell_i,M^\ell_{i+1})$ for $\ell = 1,2$ and
$NF_{\lambda,\langle \kappa,1,\kappa \rangle}
(M^\ell_i,M^\ell_{i+1},N^\ell_i,N^\ell_{i+1})$ respectively.  Now clearly
$N^\ell_0$ is $(\lambda,\kappa)$-saturated over $M^\ell_0$ hence over
$M^0_0$ (for $\ell = 1,2$) so there is an isomorphism $f_0$ from
$N^2_0$ onto $N^1_0$ extending id$_{M^0_0}$ and $f(a) \notin M^1_0$.
\medskip

We continue as in the proof of part (2).  In the end
$f = \dsize \bigcup_{i < \lambda^+} f_i$ is an isomorphism of $N_2$ onto
$N_1$ over $M_0$ and as $f^1_0(a)$ is well defined and in $N^1_0 \backslash
M^1_0$ clearly \newline
$f(a) = f_0(a) \notin M_1$, as required. \hfill$\square_{\scite{9.4}}$
\enddemo
\bigskip

\proclaim{\stag{9.5} Claim}  If $\delta$ is a 
limit ordinal $< \lambda^{+2}$ and
$\langle M_i:i < \delta \rangle$ is a $\le^*_{\lambda^+}$-increasing
continuous \underbar{then} $M_i \le^*_{\lambda^+} \dsize \bigcup_{j < \delta}
M_j$ for each $i < \delta$.
\endproclaim
\bigskip

\demo{Proof}  We prove it by induction on $\delta$.  
Now if $C$ is a club of $\delta$ with $i \in C$,
then we can replace $\langle M_j:j < \delta \rangle$ by $\langle M_j:j \in
C \rangle$ so without loss of generality $\delta = \text{ cf}(\delta)$, so
$\delta \le \lambda^+$; clearly it is enough to prove $M_0 \le^*_{\lambda^+}
\dsize \bigcup_{j < \delta} M_j$.  Let $\langle M^i_\zeta:\zeta < \lambda^+
\rangle$ be a representation of $M_i$.
\enddemo
\bigskip

\noindent
\underbar{Case A}:  $\delta < \lambda^+$. \newline
Without loss of generality (see \scite{9.2}(1)) for every $i < j < \delta$ and
$\zeta < \lambda^+$ we have: \newline
$M^j_\zeta \cap M_i = M^i_\zeta$ and $NF_\lambda(M^i_\zeta,M^i_{\zeta +1},
M^j_\zeta,M^j_{\zeta +1})$.  Let
$M^\delta_\zeta = \dsize \bigcup_{i < \delta} M^i_\zeta$, so \newline
$\langle M^\delta_\zeta:\zeta < \lambda^+ \rangle$ 
is $\le_{\frak K}$-increasingly continuous sequence of members of $K_\lambda$
with limit $M_\delta$, and for $i < \delta,M^\delta_\zeta \cap M_i =
M^i_\zeta$.  By symmetry (see \scite{8.14}) we have $NF_\lambda(M^i_\zeta,
M^{i+1}_\zeta,M^i_{\zeta +1},M^{i+1}_{\zeta + 1})$ so as $\langle M^i_\zeta:
i \le \delta \rangle$,$\langle M^i_{\zeta +1}:i \le \delta \rangle$ are
$\le_{\frak K}$-increasingly continuous by \scite{8.16} we know 
$NF_\lambda(M^0_\zeta,M^\delta_\zeta,
M^0_{\zeta +1},M^\delta_{\zeta +1})$ hence by symmetry (\scite{8.14}) we have
$NF_\lambda(M^0_\zeta,M^0_{\zeta +1},M^\delta_\zeta,M^\delta_{\zeta +1})$.
\newline
So $\langle M^0_\zeta:\zeta < \lambda^+ \rangle,\langle M^\delta_\zeta:
\zeta < \lambda^+ \rangle$ are witnesses to $M_0 \le^*_{\lambda^+} M_\delta$.
\bigskip

\noindent
\underbar{Case B}:  $\delta = \lambda^+$.

By \scite{9.2}(1) (using normality of the club filter, restricting 
to a club of $\lambda^+$ and renaming), without loss
of generality for $i < j \le 1 + \zeta < 1 + \xi < \lambda^+$ we have
$M^j_\zeta \cap M_i = M^i_\zeta$, and $NF_\lambda(M^i_\zeta,M^i_\xi,M^j
_\zeta,M^j_\xi)$.  Let us define $M^{\lambda^+}_i = \dsize \bigcup
_{j < 1+i} M^j_i$.  So $\langle M^{\lambda^+}_i:i < \lambda^+ \rangle$ is
a representation of $M^\lambda_{\lambda^+} = M_\delta$ and continue as
before. \hfill$\square_{\scite{9.5}}$
\bigskip

\proclaim{\stag{9.6} Claim}  Assume 
$M_0 <^+_{\lambda^+,\kappa} M_2$ and $a \in M_2 \backslash M_0$, 
and for some $N \le_{\frak K} M_0$ we have: 
$N \in K_\lambda$ and
tp$(a,N,M_2)$ is minimal.  \underbar{Then} we can find $M_1$, \newline
$\bar M^0 = \langle M_{0,i}:
i < \lambda^+ \rangle,\bar M^1 = \langle M_{1,i}:i < \lambda^+ \rangle$
such that:
\medskip
\roster
\item "{$(a)$}"  $\bar M^0$ is a $\le_{\frak K}$-representation of $M_0$
\sn
\item "{$(b)$}"  $\bar M^1$ is a representation of $M_1(\in K_{\lambda^+}),
a \in M_{1,i}$, for all $i$
\sn
\item "{$(c)$}"  $M_0 \le_{\frak K} M_1 \le_{\frak K} M_2$
\sn
\item "{$(d)$}"  for $i < \lambda^+$ we have 
$NF_{\lambda,\langle \lambda,1,1 \rangle}
(M_{0,i},M_{0,i+1},M_{1,i},M_{1,i+1})$ \newline
(hence $M_{\ell,i} = M_{\ell + 1,i} \cap M_\ell)$
\sn
\item "{$(e)$}"  $(M_{0,i},M_{1,i},a) \in K^3_\lambda$ is minimal and
reduced.
\endroster
\endproclaim
\bigskip

\demo{Proof}  
Let $\langle M_{0,i}:i < \lambda^+ \rangle,\langle M_{2,i}:i < \lambda^+
\rangle$ be representations of $M_0,M_2$ respectively, as required in 
\scite{9.1}(2) and without loss of generality $N \le_{\frak K} M_{0,0}$ 
and $a \in M_{2,0}$.  
We now choose by induction on $\zeta < \lambda^+$, an ordinal
$i(\zeta)$ and models $M_{1,i(\zeta)},M_{3,i(\zeta)}$ such that:
\medskip
\roster
\item "{$(A)$}"  $i(\zeta) < \lambda^+$ is increasing continuous in $\zeta$
and \newline
$a \in M_{2,i(0)} \backslash M_{0,i(0)},N \le_{\frak K} M_{0,i(0)}$
\sn
\item "{$(B)$}"  $M_{0,i(\zeta)} \le_{\frak K} M_{1,i(\zeta)} \le_{\frak K}
M_{3,i(\zeta)}$ and $M_{2,i(\zeta)} \le_{\frak K} M_{3,i(\zeta)}$
\sn
\item "{$(C)$}"  $a \in M_{1,i(0)}$ and $(M_{0,i(\zeta)},M_{1,i(\zeta)},a)$
is minimal and reduced 
\sn
\item "{$(D)$}"  for $\xi < \zeta$ and $(\ell,m) \in \{(0,1),(0,2),(1,3),
(2,3)\}$ we have \newline
$NF_\lambda(M_{\ell,i(\xi)},M_{\ell,i(\zeta)},M_{m,i(\xi)},
M_{m,i(\zeta)})$
\sn
\item "{$(E)$}"  $M_{1,i(\zeta)},M_{3,i(\zeta)}$ is $\le_{\frak K}$-increasing
continuous in $\zeta$.
\endroster
\medskip

\noindent
For $\zeta = 0$ note that for $i(0) < \lambda^+,a \in M_{2,i(0)}$ and
$M_{2,i(0)}$ is universal over $M_{0,i(0)}$.
\medskip

\noindent
For $\zeta$ limit let $i(\zeta) = \dsize \bigcup_{\xi < \zeta} i(\xi)$ and
$M_{1,i(\zeta)} = \dsize \bigcup_{\xi < \zeta} M_{1,i(\zeta)}$.
\medskip

\noindent
For $\zeta = \xi + 1$, there is $i(\zeta) \in (i(\xi),\lambda^+)$ and
a model $N_\zeta$ such that $M_{1,i(\xi)} \le_{\frak K} N_\zeta 
\in K_\lambda$ and $\le_{\frak K}$-embedding $f$ of $M_{0,i(\zeta)}$ 
into $N_\zeta,f \restriction
M_{0,i(\zeta)}$ the identity and $(f(M_{0,i(\xi)}),N_\zeta,a)$ is minimal
and reduced.  By uniqueness (i.e. Claim \scite{8.1}) we can find such $N$
satisfying $(\exists M)(N \le_{\frak K} M \in K_\lambda \and 
M_{1,i(\zeta)} \le_{\frak K} M)$.  So we can carry the induction.
\medskip

Lastly, by uniqueness of $<^+_{\lambda,\kappa}$ we can make 
$M_3 = \dsize \bigcup_{\zeta < \lambda^+} M_{3,i(\zeta)}$ to be 
$\le_{\frak K} M_2$ as required. \hfill$\square_{\scite{9.6}}$
\enddemo
\bigskip

\definition{\stag{9.7} Definition}  If $(M_0,M_1,a)$ are as in \scite{9.6}
(a)-(e) we say $(M_0,M_1,a)$ is $\lambda^+$-locally reduced nice and minimal
($\lambda^+$-l.r.n.m.).  We omit ``nice" if we omit clause $(d)$.
\enddefinition
\bigskip

\proclaim{\stag{9.8} Claim}  If 
$(M_0,M_1,a)$ is $\lambda^+$-l.r.n.m. \underbar{then}
$(M_0,M_1,a) \in K^3_{\lambda^+}$ is reduced.
\endproclaim
\bigskip

\demo{Proof}  Check.

Using also \scite{7.4A}:
\enddemo
\bigskip

\proclaim{\stag{9.9} Claim}   $M_0 <^+_{\lambda^+,\kappa} M_1$ 
\underbar{if and only
if} we can find $\langle M^*_j,a_j:j < \lambda^+ \times \kappa \rangle$ 
such that:
\medskip
\roster
\item "{$(a)$}"  $M^*_j$ is $\le_{\frak K}$-increasing 
continuous (in $K_{\lambda^+}$)
\sn
\item "{$(b)$}"  $(M^*_j,M^*_{j+1},a_j)$ is $\lambda^+$-l.r.n.m.
\sn
\item "{$(c)$}"  $M^*_0 = M_0,M^*_j \cup M^*_j = M_1$
\sn
\item "{$(d)$}"  for some $N \le M_0,N \in K_\lambda$ and minimal reduced
$p \in S(N)$, for every $j$.
\endroster
\endproclaim
\bigskip

\demo{Proof}  We can find $\langle M^*_j:j \le \lambda^+ \times \kappa)$ 
satisfying clauses $(a),(b)$ and $(d)$.  By \scite{5.4} easily 
$M^*_0 <^+_\lambda M^*_{\lambda^+}$.
Now by the uniqueness (= \scite{9.4}(2)) + categoricity of $K$ in 
$\lambda^+$, we are done. \hfill$\square_{\scite{9.9}}$
\enddemo
\bigskip

\proclaim{\stag{9.10} Claim}  In 
$(K_{\lambda^+},<^*_{\lambda^+})$ we have disjoint amalgamation.
\endproclaim
\bigskip

\demo{Proof}  First redo \scite{9.4} assuming 
$(M_0,M_\ell,a_\ell)$ for $\ell =1,2$
is $\lambda^+$-l.r.n.m., and getting \newline
$a_1 \notin f(M_2),f(a_2) \notin M_1$ (and
we can start with this).  By \scite{9.8} we get \newline
$M_1 \cap f(M_2) = M_3$, so we have
disjoint amalgamation.  By \scite{9.9} and 
chasing arrows we get it in general. \hfill$\square_{\scite{9.10}}$
\enddemo
\bigskip

\remark{Remark}  This is like the proof of disjoint amalgamations in
\scite{5.7}.
\endremark
\newpage

\head {\S10 Non-structure for $\le^*_{\lambda^+}$} \endhead  \resetall
\bigskip

\demo{\stag{10.0} Hypothesis}  Assume the hypothesis \scite{8.0}
and the further model theoretic 
properties deduced since then (including \scite{6.5} but not \scite{6.6}.) 
\enddemo
\bigskip

It would have been nice to prove all disjoint amalgamations in $K_\lambda$ are
$NF_\lambda$, but this is, at this point, not clear.  But as we look upward
(i.e. want to prove statement on $K_{> \lambda^+}$) and
$<^*_{\lambda^+}$ is very nice, it will be essentially just as well if for
$M,N \in K_{\lambda^+}$ we have $M \le_{\frak K} N \Rightarrow M 
\le^*_{\lambda^+} N$.
Our intention is to assume $M^* \le_{\frak K} N^*$ is a counterexample of
this statement and we would like to say that in a sense this implies the 
existence of many
types over $M^*$ so that we can construct many models in $\lambda^{+2}$.
Note: building models in $K_{\lambda^+},K_{\lambda^{++}}$ 
by approximations in $K_\lambda$ is nice if we use the smooth
amalgamation but we do not have it for non-smooth ones.  So we shall use
$M^* \in K_{\lambda^+}$ being saturated so it has many automorphisms.
\bigskip

\proclaim{\stag{10.1} Claim}  
1) Assume $M_1 \le_{\frak K} M_2$ are in $K_{\lambda^+}$.  Then we can find
$M_0 \in K_{\lambda^+}$ such that $M_0 <^+_{\lambda,\kappa} M_1$ and $M_0 
\le^*_{\lambda^+} M_2$. \newline
2) Also we can find $\langle M_{0,i}:i < \lambda^+ \rangle$, an
$\le_{\frak K}$-increasingly continuous sequence of members of $K_{\lambda^+}$
such that $M_{0,i} <^*_{\lambda^+} M_{0,i+1}$ and 
$\dsize \bigcup_{i < \lambda^+} M_{0,i} = M_1$ and \newline
$i < \lambda^+ \Rightarrow  M_{0,i} \le^*_{\lambda^+} M_2$.
\endproclaim
\bigskip

\demo{Proof}  Let $(M^*,N^*) \in K^{3,\text{uq}}_\lambda$ be from 
\scite{8.0}(2).  Let
$\langle M_{\ell,i}:i < \lambda^+ \rangle$ be a representation of $M_\ell$
for $\ell = 1,2$ and without loss of generality $M_{\ell,i+1}$ is $(\lambda,
\lambda)$-saturated over $M_{\ell,i}$ for $\ell =1,2$ and
$M_{2,i} \cap M_1 = M_{1,i}$. \newline
1) Now choose by induction on $i,M_{0,i}$ such that:
\medskip
\roster
\item "{$(a)$}"  $M_{0,i} \le_{\frak K} M_{1,i}$
\sn
\item "{$(b)$}"  $M_{0,i}$ is $\le_{\frak K}$-increasing continuous
\sn
\item "{$(c)$}"  $M_{0,i+1} \cap M_{1,i} = M_{0,i}$
\sn
\item "{$(d)$}"  $M_{1,i+1}$ is $(\lambda,\text{cf}(\lambda))$-saturated
over $M_{1,i} \cup M_{0,i+1}$
\sn
\item "{$(e)$}"  $(M_{0,i},M_{0,i+1}) \cong (M^*,N^*)$.
\endroster
\medskip

\noindent
There is no problem to carry the definition.  Now let
$M_0 = \dsize \bigcup_{i < \lambda^+} M_{0,i}$ so $M_0 <^+_{\lambda^+} M_1$
and $M_0 \le^*_{\lambda^+} M_2$ are checked by their definitions.
\newline
2) We choose by induction on $i < \lambda^+,\langle M^*_{\varepsilon,i}:
\varepsilon \le 1+i \rangle$ such that:
\medskip
\roster
\item "{$(a)$}"  $M^*_{1+i,i} = M_{1,\lambda \times (1 + \varepsilon)
\times i}$
\sn
\item "{$(b)$}"  For each $\varepsilon$ the sequence
$\langle M^*_{\varepsilon,j}:\varepsilon \le j \le i \rangle$
is $\le_{\frak K}$-increasing continuous
\sn
\item "{$(c)$}"  For each $i$, the sequence
$\langle M^*_{\varepsilon,i}:\varepsilon \le 1+i \rangle$
is $\le_{\frak K}$-increasing continuous
\sn
\item "{$(d)$}"  $M^*_{\varepsilon,i} \cap M^*_{\zeta,j} = 
M^*_{\text{min}\{\varepsilon,\zeta\},\text{min}\{i,j\}}$
\sn
\item "{$(e)$}"  $M^*_{\varepsilon +1,i+1}$ is $(\lambda,\text{cf}
(\lambda \times (1 + \varepsilon)))$-saturated over 
$M^*_{\varepsilon,i+1} \cup M^*_{\varepsilon +1,i}$
\sn
\item "{$(f)$}"  $NF_{\lambda,\langle \lambda,1,\lambda \rangle}
(M^*_{\varepsilon,i},M^*_{\varepsilon +1,i},M^*_{\varepsilon,i+1},
M^*_{\varepsilon +1,i+1})$
\sn
\item "{$(g)$}"  for $\varepsilon < 1 + i$ we have
$NF_\lambda(M^*_{\varepsilon,i},M^*_{\varepsilon,i+1},M_{2,\lambda \times
\lambda \times i},M_{2,\lambda \times \lambda \times (i+1)})$.
\endroster
\medskip

\noindent
For $i=0,i$ limit no problem, for $i = j+1$ first choose $N_{i,\varepsilon,
\zeta} \le_{\frak K} M_{1,\lambda \times \lambda \times i + \lambda \times 
\zeta}$ for $\zeta \le \lambda \times (1+i),\le_{\frak K}$-increasing
continuous, $N_{i,\varepsilon,\zeta} \in K_\lambda,N_{i,\varepsilon,0} =
M_{1,j}( = M^*_{1+i,i})$, \newline
$(N_{i,\varepsilon,\zeta},N_{i,\varepsilon,\zeta +1})
\cong (M^*,N^*)$ and $N_{i,\varepsilon,\zeta + 1} \cap M_{2,\lambda \times
\lambda \times i + \lambda \times \zeta} = N_{i,\varepsilon,\zeta}$.

Now by \scite{5.4}, \wilog \, 
$N_{i,\varepsilon,\lambda}$ is $(\lambda,\lambda)$-saturated over 
$M_{1,j}$, and we choose it as $M^*_{1+j,i}$, and we choose
$M^*_{1,i+1}$ as $M_{1+i,i}$ note that clauses (a) and (f) holds.  Now we
can find $N_{\varepsilon,j_1}$ for $j_1 < 1 + j$ as in \scite{8.6A} and use
uniqueness of the $(\lambda,\lambda \times (1+i))$-saturated model over
$M_{1,j}$.   \hfill$\square_{\scite{10.2}}$
\enddemo
\bigskip

\demo{\stag{10.2} Conclusion}  Assume $M \le_{\frak K} N$ are from 
$K_{\lambda^+}$.  If $\langle M_i:i < \lambda^+ \rangle$ is \newline
$\le^*_{\lambda^+}$-increasing
continuous and for each $i$ for some $N_i$ we have \newline
$M_i <^+_{\lambda^+} N_i
\le^*_{\lambda^+} M_{i+1}$ \underbar{then} for some $(M',N')$ we have:

$$
(M',N') \cong (M,N)
$$

$$
M' = \dsize \bigcup_{i < \lambda} M_i
$$

$$
i < \lambda^+ \Rightarrow M_i \le^*_{\lambda^+} N'.
$$
\enddemo
\bigskip

\demo{Proof} By \scite{10.1}(2) and as we know 
$M^a \le^*_{\lambda^+} M^b \le^+_{\lambda^+} M^c \Rightarrow 
M^a <^+_{\lambda^+} M^c$ and the uniqueness
of $M''$ over $M'$ when $M' <^+_{\lambda^+} M''$. 
\hfill$\square_{\scite{10.2}}$
\enddemo
\bigskip

\proclaim{\stag{10.3} Lemma}  If $\le_{\frak K} \restriction 
K_{\lambda^+}$ is not $\le^*_{\lambda^+}$ \underbar{then} 
$I(\lambda^{+2},K) = 2^{\lambda^{+2}}$.
\endproclaim
\bigskip

\demo{Proof}  Let $S \subseteq \{ \delta < \lambda^{+2}:\text{cf}(\delta) =
\lambda^+\}$ be stationary.  We shall below construct a model $M_S \in
K_{\lambda^{+2}}$ such that from the isomorphism type of $M^S_1$ we can
reconstruct $S/{\Cal D}_{\lambda^{+2}}$; this clearly suffices.  Choose
$M \le_{\frak K} N$ from $K_{\lambda^+}$ such that $\neg(M \le^*_{\lambda^+}
N)$ so by \scite{9.3}(3) without loss of generality $|N \backslash M| =
\lambda^+$.
\medskip

We choose by induction on $\alpha < \lambda^{+2}$ a model $M^S_\alpha$
such that:
\medskip
\roster
\item "{$(a)$}"  $M^S_\alpha \in K_{\lambda^+}$ has universe $\lambda \times
(1 + \alpha)$
\sn
\item "{$(b)$}"  for $\beta < \alpha$ we have $M^S_\beta \le_{\frak K}
M^S_\alpha$
\sn
\item "{$(c)$}"  if $\alpha = \beta + 1$, $\beta \notin S$ then $M^S_\beta
<^+_{\lambda^+} M^S_\alpha$
\sn
\item "{$(d)$}"  if $\alpha = \beta + 1,\beta \in S$ then $(M^S_\beta,
M^S_\alpha) \cong (M,N)$
\sn
\item "{$(e)$}"  if $\beta < \alpha,\beta \notin S$ then $M^S_\beta
\le^*_{\lambda^+} M^S_\alpha$.
\endroster
\medskip

We use freely the transitivity (\scite{9.3}(4)) and continuity 
(\scite{9.5}) of $\le^*_\lambda$ and $[M^a \le_{\frak K} M^b \le_{\frak K}
M^c$ in $K_\lambda,\neg(M^a \le^*_{\lambda^+} M^b \Rightarrow \neg
(M^a \le^*_{\lambda^+} M^c)]$ \scite{9.3}(3). \newline
For $\alpha = 0$ no problem.
\medskip

\noindent 
For $\alpha$ limit no problem.
\medskip

\noindent 
For $\alpha = \beta +1,\beta \notin S$ no problem.
\medskip

\noindent 
For $\alpha = \beta + 1,\beta \in S$ so cf$(\beta) = \lambda^+$, let
$\langle \gamma_i:i < \lambda^+ \rangle$ be increasing continuous with limit
$\beta$ and cf$(\gamma_i) \le \lambda^+$, hence $\gamma_i \notin S$.  By
\scite{9.2}(2) without loss of generality $M_{\gamma_i} <^*_{\lambda^+}
M_{\gamma_{i+1}}$.  Now use \scite{10.2} (and the uniqueness 
(\scite{9.4}(2))).   \hfill$\square_{\scite{10.3}}$
\enddemo
\bigskip

\demo{\stag{10.4} Conclusion}  Assume $I(\lambda^{+2},K) < 2^{\lambda^{+2}}$
(in addition to hypothesis \scite{10.0}). \newline
1) $\le^*_{\lambda^+} = \le_{\frak K} \restriction K_{\lambda^+}$.
\newline
2) $(K_{\lambda^+},\le_{\frak K})$, has disjoint amalgamation so no
$M \in K_{\lambda^{+2}}$ is $\le_{\frak K}$-maximal. \newline
3) $K_{\lambda^{+3}} \ne \emptyset$.
\enddemo
\bigskip

\demo{Proof}  1) By \scite{10.3}. \newline
2) By \scite{9.10} (and part (1)). \newline
3) By \scite{10.4}(2) and \scite{2.4}(8) with $\lambda$ there replaced by
$\lambda^+$ here.  \hfill$\square_{\scite{10.4}}$
\enddemo
\bigskip

\centerline {$* \qquad * \qquad *$}
\bigskip

So we have finally proved the main theorem.  Though not directly contributing
to our main theme, we remark some more consequences of $\le_{\frak K} 
\restriction K_{\lambda^+} \ne \le^*_{\lambda^+}$.
\bigskip

\proclaim{\stag{10.5} Claim}  $(*)_0 \Leftrightarrow (*)_1$ where
\medskip
\roster
\item "{$(*)_0$}"  for some $M \le_{\frak K} N$ from $K_{\lambda^+}$, we do
not have $M \le^*_{\lambda^+} N$
\sn
\item "{$(*)_1$}"  for some $M \le_{\frak K} N$ from $K_{\lambda^+}$ we have:
{\roster
\itemitem{ {} }  if $a \in N \backslash M$ \underbar{then}
tp$(a,M,N)$ is not realized in any $M$ \newline
such that $M^* \le^*_{\lambda^+} M \in K_{\lambda^+}$
\endroster}
\endroster
\endproclaim
\bigskip

\definition{\stag{10.6} Definition}  Assume $M_0 <_{\frak K} M_1$ 
are in $K_{\lambda^+}$,and
$\bar M^\ell = \langle M_{\ell,i}:i < \lambda^+ \rangle$ is a \newline
$\le_{\frak K}$-representation of $M_\ell$ for $\ell = 0,1$.  Let

$$
\align
{\Cal S}_0(\bar M^0,\bar M^1) = 
\biggl\{ \delta < \lambda^+:&M_{1,\delta} \cap M_0 =
M_{0,\delta} \tag "{$(a)$}" \\
  &\text{and not } NF_\lambda(M_{0,\delta},M_{0,\delta + 1},M_{1,\delta},
M_1) \biggr\}
\endalign
$$

$$
{\Cal S}_1(M_0,M_1) = {\Cal S}_0(\bar M^0,\bar M^1)/D_{\lambda^+}
\text{ (well defined)}. \tag "{$(b)$}"
$$
\medskip

\noindent
(c) $\qquad \quad J$ is the normal ideal on $\lambda^{++}$ generated by sets 
of the form \newline

$\qquad \quad \,\,{\Cal S}_0(M^0,M^1)$ where $M^0,M^1$ are as above.
\enddefinition
\bigskip

\noindent
\stag{10.7}\underbar{Comment}:  An earlier try for \scite{10.3} was: \newline
1) For every $S \in J$ for some $\bar M^0,M^1$
as in \scite{10.6} we have \newline
$S = S_0(\bar M^0,\bar M^1)$.
\newline
2)  If $S_1 = S_1(M^0,M^1)$ is stationary then for some
$\bar M = \langle M_i:i < \lambda^+ \rangle$ a representation of
$M = \dsize \bigcup_{i < \lambda^+} M_i \in K_{\lambda^+}$ for every $S
\subseteq S_1$ for some $\bar M'$, we have \newline
$\bar M^1 = \langle M^1_i:i <
\lambda^+ \rangle,M_1 <_{\frak K} M^1 = \dsize \bigcup_{i < \lambda^+} M^1_i,
M^1_i \cap M = M_i$ and \newline
$S_0(\bar M,\bar M^1) = S \text{ mod } {\Cal D}_{\lambda^+}$. \newline
3)  If $\le^*_{\lambda^+} \ne \le_{\frak K} \restriction K_{\lambda^+}$ and
$I(\lambda^{+2},K) < 2^{\lambda^{+2}}$ \underbar{then} for
some stationary $S \subseteq \lambda^+$ we have:
\medskip
\roster
\item "{$(a)$}"  ${\Cal D}_{\lambda^+} \restriction S$ is 
$\lambda^{++}$-saturated
\sn
\item "{$(b)$}"  $\bar M^0,\bar M^1$ as in \scite{10.3} implies
$S_1(\bar M^0,\bar M^1) \subseteq S \text{ mod } D_{\lambda^+}$.
\endroster
\medskip

\noindent
4) If $\le^*_{\lambda^+} \ne \le_{\frak K} \restriction K_{\lambda^+}$ and
$2^\lambda < 2^{\lambda^+} < 2^{\lambda^{+2}}$ then $I(\lambda^{+2},K) =
2^{\lambda^{+2}}$.
\bigskip

\demo{Proof}  1) First we prove only ``$S \subseteq S_0(\bar M^0,\bar M^1)$".
Easy as ${\frak K}_{\lambda^+}$ has amalgamation and
\medskip
\roster
\item "{$\bigotimes$}"  if $M_0 \le_{\frak K} M_1 \le_{\frak K} M_2$ are in
$K_{\lambda^+},\bar M^\ell$ \newline
representing the $S_0(\bar M^0,\bar M^1) \subseteq
S_0(\bar M^0,\bar M^2)$.
\endroster
\medskip

\noindent
Now for equality use part (2).  
\newline
2) Similar to the proof of \scite{10.1}. \newline
3) Suppose $S^* = S_1(\bar M^0,\bar M^2)$ is stationary, let
$\bar S = \langle S'_\alpha:\alpha < \lambda^{++} \rangle$ be such that
$S'_\alpha \in J$.  We can build a model $M^{\bar S} \in K_{\lambda^{+2}}$
and a representation $\langle M^{\bar S}_\alpha:\alpha < \lambda^{++} \rangle$
such that

$$
S_1(M_\alpha,M_{\alpha +1}) = S'_\alpha/D_{\lambda^{++}}.
$$
\medskip

\noindent
4) By \S3 (using the proof of part (2). 
\enddemo
\newpage
\shlhetal

    
REFERENCES.  
\bibliographystyle{lit-plain}
\bibliography{lista,listb,listx,listf,liste}

\enddocument

\bye